\magnification = 1200

\input amssym.def
\input amssym.tex
\input epsf

\def \qed {\hfill $\square$}
\def \R {\Bbb R}
\def \Z {\Bbb Z}
\def \C {{\cal C}}
\def \D {{\cal D}}

\def \S {{\cal S}}
\def \ms {\medskip}

\def \vfe {\vfill\eject}
\overfullrule=0pt
\def \diam{\mathop{\rm diam}\nolimits}
\def \dist{\mathop{\rm dist}\nolimits}
\def \inter{\mathop{\rm int}\nolimits}
\def \Min{\mathop{\rm Min}\nolimits}

\centerline{\bf H\"OLDER REGULARITY OF TWO-DIMENSIONAL}
\centerline{ALMOST-MINIMAL SETS IN $\R^n$}
\vskip 0.5cm
\centerline{ Guy David }
\vskip 1cm

\noindent
{\bf R\'{e}sum\'{e}.}
On donne une d\'{e}monstration diff\'{e}rente et sans doute
plus \'{e}l\'{e}mentaire d'une bonne partie du r\'{e}sultat de 
r\'{e}gularit\'{e} de Jean Taylor sur les ensembles presque-minimaux 
d'Almgren. On en profite pour donner des pr\'{e}cisions sur les 
ensembles presque minimaux, g\'{e}n\'{e}raliser une partie du 
th\'{e}or\`{e}me de Taylor aux ensembles presque minimaux de 
dimension~$2$ dans $\R^n$, et donner la caract\'{e}risation attendue
des ensembles ferm\'{e}s $E$ de dimension~$2$ dans $\R^3$ qui sont 
minimaux, au sens o\`{u} $H^2(E\setminus F) \leq H^2(F\setminus E)$
pour tout ferm\'{e} $F$ tel qu'il existe une partie born\'{e}e $B$
telle que $F=E$ hors de $B$ et $F$ s\'{e}pare les points de 
$\R^3 \setminus B$ qui sont s\'{e}par\'{e}s par $E$.

\bigskip \noindent
{\bf Abstract.}
We give a different and probably more elementary proof of a good part of 
Jean Taylor's regularity theorem for Almgren almost-minimal sets 
of dimension $2$ in $\R^3$. We use this opportunity to settle some
details about almost-minimal sets, extend a part of Taylor's result to 
almost-minimal sets of dimension $2$ in $\R^n$, and give the expected 
characterization of the closed sets $E$ of dimension $2$ in $\R^3$ 
that are minimal, in the sense that $H^2(E\setminus F) \leq H^2(F\setminus E)$
for every closed set $F$ such that there is a bounded set $B$ 
so that $F=E$ out of $B$ and $F$ separates points of $\R^3 \setminus B$ 
that $E$  separates.


\ms\ms\ms

\noindent {\bf 1. Introduction}
\medskip

The core of this paper is a slightly different, probably more elementary,
and very detailed proof of a good part of Jean Taylor's regularity theorem [Ta] 
for Almgren almost-minimal sets of dimension $2$ in $\R^3$.

Recall that this result says that these sets are locally 
$C^1$-equivalent to minimal cones (i.e., planes or cones of type
$\Bbb Y$ or $\Bbb T$, as described near Figures 1.1 and 1.2).
We shall only prove the biH\"older equivalence in this text;
the $C^1$-equivalence needs a little bit of extra work, which
is left for a next paper [D3].  

The proof will extend to $2$-dimensional sets in $\R^n$, and will also 
allow us to show that every $MS$-minimal set of dimension $2$ in 
$\R^3$ (where the definition of competitors uses the separation 
constraint in Definition 1.4 below) is a minimal cone.
The issue shows up naturally in the study of the Mumford-Shah
functional (what are the global Mumford-Shah minimizers $(u,K)$ in $\R^3$ 
for which $u$ is locally constant), to the point that the author felt
compelled to announce Theorem 1.9 in [D2]. 

The local biH\"older equivalence to minimal cones can be useful
as a first step to $C^1$ equivalence, as in [D3],  
but could also give enough geometric control on the almost-minimal sets 
to allow one, perhaps, to prove existence results for functionals with a 
dominant area term under suitable topological constraints.

One of the main initial motivations for this text was to help the author 
understand Jean Taylor's paper, and possibly help some readers too. Part of
the reason for this interest was potential applications of the technical
lemmas to the regularity of minimal segmentations for the Mumford-Shah functional 
in dimension 3, but it is not yet certain that the current text will help.

Since this was the triggering reason for this paper, let us define 
$MS$-minimal sets first. We give the definition for $(n-1)$-dimensional 
sets in $\R^n$, but for the theorem we shall restrict to $n=3$.

\medskip
\proclaim Definition 1.1.
Let $E$ be a closed set in $\R^n$. A $MS$-competitor for $E$
is a closed set $F \i R^n$ such that we can find $R>0$ such that
$$
F \setminus B(0,R) = E \setminus B(0,R)
\leqno (1.2)
$$
and
$$
\hbox{$F$ separates $y$ from $z$ whenever } 
y,z \in \R^n \setminus (B(0,R) \cup E)
\hbox{ are separated by } E.
\leqno (1.3)
$$
\medskip

By ``$F$ separates $y$ from $z$", we just mean that $y$ and $z$ lie 
in different connected components of $\R^n\setminus F$. 
Note that if (1.2) and (1.3) hold for $R$, they also hold 
for every $R'>R$. Here $MS$ stands for Mumford-Shah; indeed the 
separation condition (1.3) is the same one as in the definition of 
global minimizers for the Mumford-Shah functional in $\R^n$.

\medskip
\proclaim Definition 1.4.
A $MS$-minimal set in $\R^n$ is a closed set $E \i \R^n$ such that
$$
H^{n-1}(E\setminus F) \leq H^{n-1}(F\setminus E)
\leqno (1.5)
$$
for every $MS$-competitor $F$ for $E$.

\medskip
Here $H^{n-1}$ denotes the $(n-1)$-dimensional Hausdorff measure.
Notice that if $E$ is $MS$-minimal, then
$H^{n-1}(E\cap B(0,R)) < + \infty$ for $R>0$, because
we can compare with $F = [E\setminus B(0,R)] \cup \partial B(0,R)$. 
We shall also use a weaker notion of minimality, based on a more
restrictive notion of competitors. The notion is the same as for 
Almgren's restricted sets [Al].

\medskip
\proclaim Definition 1.6.
Let $E$ be a closed set in $\R^n$ and $d \leq n-1$ be an integer. 
An $Al$-competitor for $E$ is a closed set $F \i R^n$ that can be 
written as $F = \varphi(E)$, where $\varphi : \R^n \to \R^n$ is 
a Lipschitz mapping such that 
$W_\varphi = \{ x\in \R^n \, ; \, \varphi(x) \neq x \}$
is bounded.
\smallskip
An $Al$-minimal set of dimension $d$ in $\R^n$ is a closed set $E \i \R^n$ 
such that $H^{d}(E\cap B(0,R)) < + \infty$ for every $R>0$ and
$$
H^{d}(E\setminus F) \leq H^{d}(F\setminus E)
\leqno (1.7)
$$
for every $Al$-competitor $F$ for $E$.

\medskip
Here it makes sense to define $Al$-minimal sets of all dimensions 
$d <n$; this is not the case for $MS$-minimal sets because of the 
special role played by separation. Also, the reader may be surprised because 
Definition 1.6 seems simpler than if we used the standard definition of
Almgren restricted sets; we shall check in Section 4 
that the two definitions coincide. 

\medskip\noindent{\bf Remark 1.8.}
Let us check that $Al$-competitors for $E$ are automatically 
$MS$-competitors for $E$, and hence $MS$-minimal sets are also 
$Al$-minimal. Let $F$ be an $Al$-competitor for $E$, let 
$\varphi$
and $W_{\varphi}$ be as in Definition 1.6, and choose
$R$ so that $W_{\varphi} \cup \varphi(W_{\varphi}) \i B(0,R)$.
Obviously (1.2) holds. For (1.3), let $y$, 
$z \in \R^n \setminus (B(0,R) \cup E)$ be separated by $E$. Consider the 
continuous deformation given by $\varphi_{t}(x) = (1-t)x + t 
\varphi(x)$ for $x \in \Bbb R^n$ and $0 \leq t \leq 1$. Notice that 
none of the sets $E_{t} = \varphi_{t}(E)$ ever contain $y$ or $z$, 
because all the modifications occur in $B(0,R)$. Then 
$F = \varphi_{1}(E)=E_1$ separates $y$ from $z$, because $E=E_0$ does. 
See for instance 4.3 in Chapter~XVII of [Du], 
on page 360. Actually, that result is only stated when $E$ is a 
compact set of $\Bbb R^n$. To accommodate this minor difficulty, 
we should instead assume that $F$ does not 
separate $y$ from $z$, pick a path $\gamma\i \R^n \setminus F$ from 
$y$ to $z$, choose $S>R$ such that $\gamma \i B(0,S)$, and notice that
$E' = [E \cap B(0,S)]\cup \partial B(0,S)$ is compact and separates 
$y$ from $z$, but its image $\varphi_{1}(E') =  [F \cap B(0,S)]\cup 
\partial B(0,S)$ does not (because of $\gamma$). From the ensuing 
contradiction we would deduce that $F$ is a $MS$-competitor for $E$, as 
needed.

\medskip
Here is a list of $MS$-minimal sets in $\R^3$. First, the empty set,
and planes, are $MS$-minimal. Next let $P$ be a plane, pick a point 
$O \in P$, and let $S$ be the union of three half lines in $P$ that 
start from $O$ and make $120^\circ$ angles with each other at $O$. 
The product of $S$ with a line perpendicular to $P$ is $MS$-minimal. 
We shall call such a set a cone of type $\Bbb Y$; thus cones of type 
$\Bbb Y$ are made up of three half planes with a common boundary $L$, 
and that meet along the line $L$ with $120^\circ$ angles. 
See Figure 1.1.

\hskip 1.7cm  
\epsffile{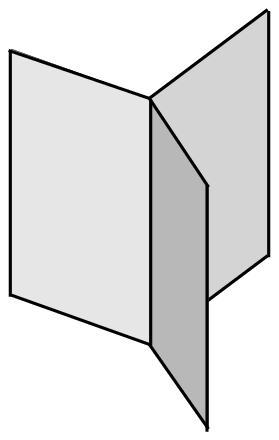} 
\hskip 1.5cm 
\epsffile{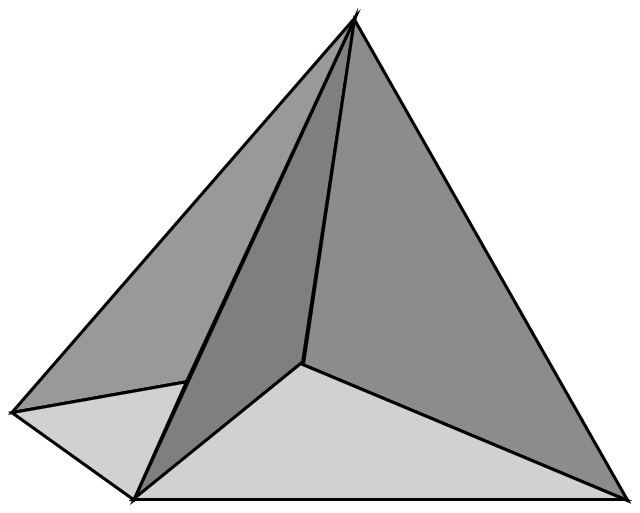} 
\medskip
\noindent\hskip 1.2cm{\bf Figure 1.1.} A cone of type $\Bbb Y$ 
\hskip 2.4cm {\bf Figure 1.2.} A set $T$
\medskip
For the last example, take a regular tetrahedron $R$, call $O$ 
its center, and denote by $T$ the cone centered at $O$ over the 
union of the six vertices of $R$. See Figure 1.2. 
Thus $T$ is composed of six angular faces that meet at $O$, and are 
bounded by four half-lines that leave from $O$ with equal (and 
maximal) angles. We shall call such a set a cone of type $\Bbb T$. 
See K. Brakke's home page http://www.susqu.edu/brakke/
for nicer pictures of $T$ and how $Y$ and $T$ arise as tangent
objects to minimal sets.

The verification that cones of type $\Bbb Y$ and $\Bbb T$ 
are $MS$-minimal can be done 
with a clever integration by parts; see [LM]. 
It turns out that modulo sets of measure $0$, the list is complete.

\bigskip
\proclaim Theorem 1.9.
If $E$ is an $MS$-minimal set in $\R^3$, then there is
a set $E^\ast$ such that $E^\ast \i E$, $H^2(E\setminus E^\ast)=0$,
and $E^\ast$ is the empty set, a plane, or a cone of type $\Bbb Y$ or $\Bbb T$.

\medskip
The fact that we may have to remove a set of vanishing $H^2$-measure
is natural: it is easy to see that the union of a minimal set with a 
closed set of vanishing $H^2$-measure is still minimal.

As we shall see, Theorem 1.9 is a rather simple consequence of the proof
of Jean Taylor's regularity theorem, but to the author's surprise, 
it seems that this was not stated in the literature yet. 
An analogue of this result for size-minimizing currents is stated in 
Section 4 of [Mo1]; the proof sketched there has the same sort of 
ingredients as in the present paper, but does not seem to me to apply 
directly.

\medskip
The situation for $Al$-minimal sets is a little less clear; although
it seem quite probable that Theorem 1.9 extends to $Al$-minimal sets 
in $\R^3$, the author does not know how to prove this without an 
additional assumption, for instance on the separating properties of 
$E$. See Proposition~18.29 and Remark 18.43.

\ms
As was hinted before, Theorem 1.9 is much easier to deduce from
a proof of Jean Taylor's regularity theorem than from its mere 
statement; we shall take this as a first excuse to give, in full details, 
a slightly different proof of that regularity result.

Let us be a little more specific about this. Her theorem concerns
almost-minimal sets of dimension 2 in $\R^3$; we shall give very 
precise definitions of almost-minimal sets, but for the moment
let us just say that they are like the minimal sets of 
Definition 1.6, except that we may add a small term to the right-hand
side of (1.7), and also localize the definition.
Recall that J. Taylor proves that if $E$ is a reduced 
almost-minimal set of dimension $2$ in $\R^3$, with a sufficiently small 
gauge function (a function that controls the small term that we add
to the right-hand side of (1.7)),
then every point of $E$ has a neighborhood where $E$ is 
$C^1$-equivalent to a plane or a cone of type $\Bbb Y$ or $\Bbb T$. 
[See Definition 2.12 for ``reduced sets".]

We shall cut this regularity result in two, and first prove a mere
biH\"older equivalence to a plane, a $\Bbb Y$, or a $\Bbb T$, under 
slightly weaker assumptions on the gauge function. 
See Theorem~16.1 for a precise statement, and observe that in our
definition the biH\"older equivalence comes from a local biH\"older
mapping of the ambient space, so the equivalence also concerns the
way $E$ is embedded in $\R^n$. This is a little simpler to prove 
than the $C^1$ equivalence, and in particular 
does not require epiperimetry estimates. This is also the only part of 
the argument that is needed for Theorem 1.9, and possibly for other 
results as well.

The full $C^1$ equivalence (or even $C^{1+\alpha}$ if the
gauge function is small enough) seems to need a more precise
argument, which uses epiperimetry or some comparison 
with harmonic functions. We intend to give in [D3] 
a special argument for this, which uses a little bit of the 
local biH\"older equivalence through separation properties. 

At this time, only the present paper extends to sets of dimension $2$
in $\R^n$, and it could be that when $n > 3$ the existence of 
epiperimetric inequalities and the local $C^1$ equivalence of 
$E$ to a minimal cone near $x$ depends on the type of tangent cones
that $E$ has at $x$. The author does not even know whether $E$ could
have an interesting one-parameter family of different tangent cones at $x$,
which would of course exclude the $C^1$ equivalence if they are
not isometric to each other.

Theorem 16.1. also applies to the two-dimensional almost-minimal sets 
in $\R^n$, $n > 3$. In this case, we don't know the precise list of minimal 
cones, and by minimal cone of type $\Bbb T$ we shall then mean a 
(reduced) minimal cone $T$ which is neither empty, nor a plane, nor a
minimal cone of type $\Bbb Y$ as above. As we shall see, these are also 
the minimal cones whose density at the origin is larger than $3\pi/2$.
See Section 14 for a rapid description of these cones.
Observe that we do not exclude the possibility that $E$ may 
have two significantly different tangent cones of type $\Bbb T$
at $x$, and then $E$ will be locally biH\"older equivalent to 
either cone.

Our proof of biH\"older equivalence relies on a few constructions 
by hand (for instance, to get the monotonicity of density), a result of 
stability of minimal and almost-minimal sets under Hausdorff limits [D1], 
and an extension of Reifenberg's topological disk theorem
to the case when $E$ is a two-dimensional set in $\R^n$ that stays
close to planes or cones of type $\Bbb Y$ or $\Bbb T$ at all scales 
and locations [DDT]. 
In particular, we never use currents or epiperimetry arguments.

Many of the ingredients in the proofs can be generalized, but the 
author does not know exactly how far the generalizations may lead.
Because of this, we shall try to keep a reasonably large
amount of generality in the intermediate statements, especially 
(but not always) when this is not too costly. For the same sort of reasons, 
the author decided to go ahead and prove in the first sections that various
definitions of almost- and quasi-minimality for sets are equivalent,
or imply others. Hopefully this will help some readers, and the other 
ones will not find it too hard to skip.

\ms 
The plan for the rest of this paper is as follows. 
Part A is mostly concerned in the definitions, regularity, and 
limiting properties of almost-minimal sets.

Section 2 contains the definitions and main regularity properties for a 
slight generalization of Almgren's restricted sets 
(we shall call them quasiminimal sets); 
the small amount of new material will only be used to allow us to simplify
our notions of almost-minimality (in Section 4), and state some results
(like Proposition 5.24 and Theorem 16.1) with a little more generality.

In Section 3 we remind the reader that for sequences of reduced 
quasiminimal sets with uniform constants, limits are also quasiminimal 
with the same constants, and we have lower- and upper semicontinuity
results for the Hausdorff measure. This stays true for the generalized 
quasiminimal sets of Section 2. See Lemmas 3.3 and 3.12 in particular.

Section 4 contains various definitions of almost-minimal sets, and the 
equivalence between the two main ones. Just like Section 2, this 
one is not really needed in view of Theorems 1.9 and 16.1, for instance.

In Part B we prove the almost monotonicity of the density function 
for almost-minimal sets with a small enough gauge function, and results
of approximation by a minimal cone in annuli where the density is almost
constant.

The main result of Section 5 is the fact that if $E$ is a minimal set
of dimension $d$, the density $r^{-d} H^d(E\cap B(x,r))$ is a 
nondecreasing function of $r$ (Proposition 5.16). 
The very standard argument is essentially a comparison of $E\cap B(x,r)$ 
with the cone over $E \cap \partial B(x,r)$. Also see 
Propositions 5.24 and 5.30 for generalizations to almost-minimal sets. 

In Section 6 we show that if $E$ is a reduced minimal set of
dimension 2 in $B(0,b) \i \R^3$ and the density $r^{-d} H^d(E\cap B(0,r))$ 
is constant on $(a,b)$, then $E$ coincides with a reduced minimal cone 
on the annulus $\{ z \, ; \,  a<|z|<b \}$. The proof is surprisingly
painful, in the sense that the author did not manage to avoid
using the minimality of $E$ again. See Proposition 6.2.

Section 7 contains a variant of Proposition 6.2 for almost-minimal 
sets and very small variations of density, obtained from it 
by a compactness argument. See Proposition 7.1 for the initial 
statement in a cone, Proposition 7.24 for the variant in a ball 
which is used most of the time, and Proposition 7.31 for the fact that
if the density of $E$ at $x$ exists, the tangent objects to $E$ at $x$
are minimal cones with the same density.

Part C contains a structural description of the reduced minimal cones
of dimension $2$ in $\R^n$. We show that for such a cone $E$, 
$K = E \cap \partial B(0,1)$ is composed of a finite number of arcs
of great circles, that may only meet by sets of three at their ends, 
and with $120^\circ$ angles. See Proposition 14.1 for
a more precise description (but observe that we shall not give a 
complete list of cones).

We start Section 8 with the simple fact that if the product set $E \times \R^m$
is a reduced $Al$-minimal set of dimension $d+m$ in $\R^{n+m}$, then
$E$ is a $d$-dimensional $Al$-minimal set in $\R^n$. See Proposition 8.3.
We don't know about, nor need the converse (see Remark~8.22).
Proposition 8.3 applies to the blow-up limits of a minimal cone at 
a point of the unit sphere, and we can use this to prove that if $E$ is
a reduced minimal cone of dimension $d$ and $\varepsilon >0$ is small,
then $K = E \cap \partial B(0,1)$ is $\varepsilon$-close (in 
normalized Hausdorff distance) to a minimal set of dimension $d-1$ 
in every small enough ball centered on $K$.

This would lead to some information on $E$, but since we want something
more precise, we start a different approach and try to study the local 
regularity of $K = E \cap \partial B(0,1)$ with our hands. 
In section 9 we show that $K$ is a ``weakly almost-minimal set"
(where ``weakly" comes from the fact that we have to pay some
price when we start from competitors for $K$ and construct a 
competitor for $E$). See Definition 9.1 and Proposition 9.4.

Then we turn more specifically to one-dimensional sets.
In Section 10 we show that the nonempty reduced minimal sets of 
dimension $1$ in $\R^n$ are the lines and the propellers $Y$.
See Theorem 10.1. The proof is also used in Section 11 to show 
that weakly almost-minimal set are often close (in Hausdorff distance)
to lines or propellers (Proposition 11.2). 

In Section 12 we state a simpler, one-dimensional version of the 
Reifenberg-type theorem that will be stated in Section 15. 
This is enough to show that when $K = E \cap \partial B$
as above, $K$ is essentially composed of $C^1$ arcs.
See Proposition 12.6 for a first statement, and Proposition 12.27
for the $C^1$ version.

In Section 13 we compare the cone over an arc of Lipschitz
curve in $\partial B(0,1)$ with the graph of a harmonic function $f$; 
this is a small computation where we use the fact that the main term
in the area functional is the energy $\int |\nabla f|^2$, and whose
conclusion is that the curve cannot yield a minimal cone unless
it is an arc of great circle.

We apply all this in Section 14 and find the description of the 
two-dimensional minimal cones of Proposition 14.1. 

Part D studies the local regularity of the almost-minimal sets.
We start in Section 15 with the description of a result from [DDT],
where a slight generalization of Reifenberg's topological disk theorem 
is proved. This is the result that will be used in Section 16 to 
produce parameterizations of $E$ when we know that is is well 
approximated by two-dimensional minimal cones in small balls.

Section 16 contains the main local regularity results for 
almost-minimal sets of dimension $2$ in $\R^n$. Theorem 16.1 is the main 
statement (the promised generalization of J. Taylor's theorem), 
but some intermediate results may have independent interest, 
such as Lemmas 16.19 and 16.48 (which gives sufficient conditions for $E$ to be 
biH\"older equivalent to a plane in a given ball $B$), 
Proposition 16.24 (which proves the existence of a point of density 
$3 \pi/2$ in $B$), Lemmas 16.25 and 16.51 (sufficient conditions for $E$ 
to be biH\"older equivalent in $B$ to a cone of type $\Bbb Y$),
and Lemma 16.56 and Proposition 18.1 (for the biH\"older equivalence 
to a cone of type $\Bbb T$). Proposition 16.24 is proved in Section 17, 
with a topological degree argument. 

Theorem 1.9 is proved in Section 18, with another topological argument. 
Also see Proposition 18.29 for a version for $Al$-minimal sets
of dimension $2$ in $\R^3$, but with an extra separation assumption.

In Section 19 we describe a set which looks like cone of type 
$\Bbb T$ at large scales, but has no point of type $T$; thus there is no trivial 
way to make the proof of Theorem 1.9 work for $Al$-minimal sets.
In the very short Section 20 we compute the density of the set $T$
of Figure 1.2.

The author wishes to thank Christopher Colin for help with
the topological argument in Section 18, T. De Pauw, J.-C. L\'{e}ger, P. Mattila,
and particularly F. Morgan for interesting discussions about the subject of this 
paper, and gratefully acknowledges partial support from the 
european network HARP.

\bigskip
\noindent {\bf A. ALMOST-MINIMAL SETS}
\medskip
In this first part we shall give a few different definitions of
almost-minimal sets of dimension $d$ in $\R^n$, check that many
definitions are equivalent, discuss some of the main properties
of these sets (including local Ahlfors-regularity and local uniform
rectifiability), and state convergence theorems that will be used
in the rest of this paper.

The equivalence between various definitions of almost-minimality
is not vital for the rest of the paper, but it will be more convenient
to use simple definitions like Definition 1.6, and good to know that they are
equivalent to more classical ones, in particular  those that come more 
naturally from Almgren's definitions in [Al]. 
The equivalence is not particularly subtle or original, but some argument 
is needed, and the author decided that this is a good opportunity
to settle the issue as explicitly as possible.

It appears that for the proof of equivalence, it is useful to introduce
a class of sets which is slightly larger than the class of restricted sets
introduced by F. Almgren, and then check that the main regularity properties 
of restricted sets still hold in this class. This is what  we do in 
Section 2.

In Section 3, we use these properties to show that the Hausdorff measure
has good semicontinuity properties when we restrict to bounded sequences of
generalized Almgren quasiminimal sets. In other words, the Hausdorff
measure of these sets goes to the limit well, 
as in Lemma 3.3 and Lemma 3.12 below, 
and we shall be able to apply this to our various types of 
almost-minimal sets.

The various definitions of almost-minimality, and the equivalence
of the main ones, are given in Section 4. We also use the results of
Section 3 to show that limits of almost-minimal sets are almost-minimal;
see Lemma 4.7.
We try to arrange this text so that the impatient reader can refer 
to this part only occasionally,
for a definition or a believable statement.

\bigskip
\noindent {\bf 2. Generalized Almgren quasiminimal sets}
\medskip

The main results of this paper apply not only to 
minimal sets, but to almost-minimal sets of various types. 
We shall give in Section 4 a few definitions of almost-minimal sets
and clarify the relations between the various classes, but before
we do this we need to talk about Almgren restricted sets (we'll call them 
quasiminimal), and even generalize the notion slightly.

The standard properties described in this section will be used in 
Section 4 to prove the equivalence between different notions of almost-minimality,
and then also in later sections, but only applied to more restricted 
classes of almost-minimal sets.

A good part of this is not needed for the main results of the paper, 
the main point of it is merely to try to give some statements 
in the best possible context. It should be possible for the 
reader to skip Sections 2 and 4 almost entirely if she is just 
interested in minimal sets or happy to use a stronger form of 
almost-minimality.

\smallskip
So we want to discuss minor variations over the notion of 
``restricted set'' of Almgren [Al]. 
For the general setting we are 
given an open set $U \i \R^n$, a quasiminimality constant  $M \geq 
1$, and a diameter $\delta\in(0,+\infty]$. Now $E$ is a closed subset 
of  $U$, and we want to assume that 
$$
H^d(E \cap B) < + \infty
\ \hbox{ for every compact ball } B \i U.
\leqno (2.1)
$$
As in Definition 1.6, we want to compare $E$ with competitors 
$F=\varphi_{1}(E)$, where $\varphi_{t}$, $0 \leq t \leq 1$, is a 
one-parameter family of continuous functions $\varphi_{t}: U \to U$ 
with the following properties:
$$
\varphi_{0}(x) = x  \ \hbox{ for } x\in U,
\leqno (2.2)
$$
$$
\hbox{the function $(t,x) \to \varphi_{t}(x)$, from
$[0,1] \times U$ to $U$, is continuous,}
\leqno (2.3)
$$
$$
\varphi_{1} \hbox{ is Lipschitz,}
\leqno (2.4)
$$
and, if we set 
$$
W_{t} = \{ x\in U ; \varphi_t(x) \neq x \}
\ \hbox{ and then }
\widehat W = \bigcup_{t\in [0,1]} \Big[ W_{t} \cup \varphi_t(W_{t}) 
\big],
\leqno (2.5)
$$
then
$$
\widehat W \hbox{ is relatively compact in $U \ $ and }
\diam(\widehat W) < \delta.
\leqno (2.6)
$$

\medskip
\proclaim Definition 2.7.
We say that the set $E$ is an $(M,\delta)$-quasiminimal set in $U$ 
when $E$ is closed in $U$, (2.1) holds, and
$$
H^d(E \cap W_{1}) \leq M H^d(\varphi_1(E \cap W_{1}))
\leqno (2.8)
$$
for every $\varphi$ such that (2.2)-(2.6) hold.
\medskip

\noindent {\bf Remarks 2.9.} 
\par\noindent {\bf a.}
This definition is almost the same thing as Almgren's definition of a 
restricted set. The only difference is that Almgren does not require 
(2.2) or (2.3). That is, he only considers mappings $\varphi_1 : U 
\to U$ and does not care whether they come from deformations or not. 
Our definition of quasiminimal sets is slightly less restrictive 
(because we add constraints on the mappings $\varphi_{1}$), but this 
is all right, because the regularity results are the same. In fact, 
in most regularity arguments for quasiminimal sets, one uses functions 
$\varphi_1$ such that $W_{1} \cup \varphi_1(W_{1})$ is contained in a 
compact ball $B$ contained in $U$ and with diameter at most $\delta$. 
Then $\varphi_{1}$ automatically comes from a one-parameter family of 
deformations as above, namely 
$\varphi_{t}(x) = (1-t) x + t \varphi_1(x)$, for which $\widehat W$ is
contained in $B$ by convexity.
\par\noindent {\bf b.} 
We require $\varphi_{1}$ to be Lipschitz (in (2.4)) only to make sure 
that (2.8) makes sense, but we do not require any Lipschitz bound for 
$\varphi_{1}$.
\par\noindent {\bf c.} 
It is important that we do not require $\varphi_{1}$ to be injective. 
That is, we want to allow competitors obtained from $E$ by gluing 
pieces together, or making entire regions collapse to one point.
\par\noindent {\bf d.}
We can take $\delta =+\infty$ in Definition 2.7. Our definition of 
$Al$-minimality in Definition~1.6 corresponds to $U=\R^n$, 
$\delta =+\infty$, and $M=1$. In addition, the accounting in (1.7) is 
not the same as in (2.8), but we shall see soon that in this case 
they are  equivalent.
In general $\delta$ and $U$ give more flexibility in the definition; 
they may allow us to localize or account for boundary conditions. For 
instance, catenoids are $(1,\delta)$-quasiminimal sets in $\R^3$, but 
only for $\delta$ small enough.

\medskip
The quasiminimal sets of Definition 2.7 are not yet general enough 
for us, so we introduce an even larger class. This extension will be needed
only because we want to work with different notions of almost-minimal 
sets; the uninterested reader may skip.

\medskip
\proclaim Definition 2.10. Let $U$ and $M$ be as above, and let
$\delta_0 \in (0,+\infty]$ and $h\in[0,1)$ be given. We shall denote by 
$GAQ(M,\delta_0 ,U,h)$ the set of closed subsets $E$ of $U$ such that 
(2.1) holds and
$$
H^d(E \cap W_{1}) \leq M H^d(\varphi_1(E \cap W_{1}))
+ h \delta^d
\leqno (2.11)
$$
for every $\delta \leq \delta_0$ and every
$\varphi$ such that (2.2)-(2.6) hold.
\medskip

Elements of $GAQ(M,\delta_0 ,U,h)$ will be called ``generalized Almgren
quasiminimal sets". The only difference with the quasiminimal sets 
from Definition 2.7 is that we added the error term $h \delta^d$ in 
the right-hand side of (2.11). So $h=0$ corresponds to the previous 
situation. We should always think that $h$ is small, because in all 
the regularity results below, $h$ has to be small enough, depending 
on $n$,  $d$, and $M$.  Also notice that the only situations where we expect a 
real difference between (2.11) and (2.8) (possibly with a larger $M$) 
is when $H^d(E \cap W_{1})$ and $H^d(\varphi_1(E \cap W_{1}))$ are 
very small, i.e., when $\varphi_{1}$ only moves very few points. The 
main point of the discussion below is that these situations can be 
avoided when we prove regularity theorems. 

We like  quasiminimal sets because they have lots of good regularity 
properties. We want to say that these properties still 
hold for generalized quasiminimal sets (always assuming that $h$ is 
small enough), and that the long proofs given elsewhere still work. 
So, for  the rest of this section, we shall try to describe the 
regularity results that we want to use, and at the same time say why 
the proof for quasiminimal sets can be generalized rather easily.

Notice that  if $E$ is a quasiminimal set and $Z$ is any closed set 
such that $H^d(Z)=0$, then $E\cup Z$ is also a quasiminimal set, with 
the same constants. For instance, in the situation of (2.8), we have that
$H^d((E\cup Z)\cap W_{1}) = H^d(E \cap W_{1}) 
\leq M H^d(\varphi_1(E \cap W_{1}))
\leq M H^d(\varphi_1((E \cup Z) \cap W_{1}))$.
This is why we introduce the (hopefully simpler) 
set $E^\ast$ below, and the notion of ``reduced set''. 

\ms\proclaim Definition 2.12.
For each closed subset $E$ of $U$, denote by
$$
E^\ast = \big\{ x\in E \, ; \, H^d(E\cap B(x,r)) > 0
\hbox{ for every } r>0 \big\}
$$
the closed support (in $U$) of the restriction of $H^d$ to $E$. 
We say that $E$ is \underbar{reduced} when $E^\ast=E$.

\ms
It is easy to see that 
$$
H^d(E \setminus E^\ast) = 0;
\leqno (2.13)
$$
indeed, we can cover $E \setminus E^\ast$ by countably many balls 
$B_{j}$ such that $H^d(E \cap B_{j})=0$ (first reduce to
$E \cap K$, where $K$ is compact and does not meet $E^\ast$,
and then use the definition of $E^\ast$ and the compactness of $E\cap K$). 

\medskip
\noindent {\bf Remark 2.14.} If $E$ is quasiminimal, then $E^\ast$ is 
quasiminimal with the same constants. This is very easy, because when 
the $\varphi_{t}$ satisfy (2.2)-(2.6), 
$H^d(\varphi_{1}(E \setminus E^\ast))=H^d(E \setminus E^\ast)=0$ by 
(2.13) and because $\varphi_{1}$ is Lipschitz, so the conditions 
(2.8) and (2.11) are the same for $E^\ast$ as for $E$.
Because of this, it is enough to study reduced quasiminimal 
sets. The same remark will apply to the various types of 
almost-minimal sets introduced in Section 4.

\medskip
We start our list of regularity properties with the basic but 
useful local Ahlfors-regularity of $E^\ast$.

\medskip
\proclaim Lemma 2.15. For each $M \geq 1$ there exist constants $h>0$ 
and $C \geq 1$, that depend only on $n$, $d$, and $M$, such that if 
$E \in GAQ(M,\delta ,U,h)$, $x\in E^\ast$,
and $0<r<\delta$ are also such that $B(x,2r) \i U$, then 
$$
C^{-1} r^d \leq H^d(E \cap B(x,r)) \leq C r^d.
\leqno (2.16)
$$

\medskip
See Proposition 4.1 in [DS] 
for standard quasiminimal sets, but note 
that this is only a reasonably minor modification of a result in 
[Al]. 
Now we should say a few words about why the proof of (2.16) in [DS] 
extends to generalized quasiminimal sets. For the upper bound, there 
is no serious problem; the idea of the proof is to suppose that for 
some cube $Q$, $\diam(Q)^{-d} H^d(E\cap Q)$ is very large, and then 
compare $E$ with a set obtained from $E$ by a Federer-Fleming 
projection on $d$-dimensional skeletons of cubes inside
$Q$.  The numbers at stake in the argument are much larger than 
$\diam(Q)^d$, and the additional small error $h\diam(Q)^d$ does not 
really change the estimates. More specifically, the first place where 
quasiminimality is used (in [DS]) is (4.5), and we should add an 
extra $Ch\diam(Q)^d$ in the right-hand side there. [Let us point out 
to the reader that would wish to check things in [DS], that our 
standard quasiminimal set is called $S$ there, and the quasiminimality 
constant $M$ is called $k$.] The computations can stay the same, up 
to (4.9) and (4.11) where we should also add an extra $Ch\diam(Q)^d$ 
in the right-hand side.
Even with $h=1$, this term is easily eaten up by the other term 
$C 2^{j(n-d)}\diam(Q)^d$ in the right-hand side of (4.9) and (4.11)
(because $j \geq 1$); the rest of the proof is the same.

The proof of the lower bound in (2.16) is probably the place where we 
should be the most careful, since the mass of $E$ nearby could be 
very small and $h\delta^d$ could be the main term in (2.10).

 In the argument of [DS] (which begins a little below (4.21) there), 
we start with a cube $Q$, assume that $\diam(Q)^{-d}H^d(E\cap Q)$ is 
very  small, and then we construct a competitor $F = \varphi(E)$ 
(with the constraints (2.2)-(2.6) with any $\delta > \diam(Q)$). The 
cube $Q$ divides into a thin annular region $A_{j}(Q)$ near the 
boundary  $\partial Q$, and the main part $H_{j}(Q)=Q\setminus 
A_{j}(Q)$. Here $j$ is an integer parameter, and the width of 
$A_{j}(Q)$ is roughly $2^{-j} \diam(Q)$. In the region $H_{j}(Q)$, we 
do a Federer-Fleming projection on skeletons of cubes of size 
$2^{-j} \diam(Q)$, and it turns out that if 
$\diam(Q)^{-d}H^d(E\cap Q)$ is small enough, depending on $j$, the 
projection contains no full face of dimension $d$. This allows us to 
continue the Federer-Fleming projection one more step, up to 
$(d-1)$-dimensional faces. As a result, our final mapping $\varphi$ 
is such that $H^d(\varphi(E\cap H_{j}(Q))) = 0$.  In the intermediate 
region $A_{j}(Q)$, we construct a mapping that matches the desired 
values on the two sides of $\partial A_{j}(Q)$. We do this so that 
$H^d(\varphi(E\cap A_{j}(Q))) \leq C H^d(E\cap A_{j}(Q))$. 
The comparison argument in [DS] yields 
$H^d(E\cap H_{j}(Q)) \leq C M H^d(E\cap A_{j}(Q))$,
where $M$ is the quasiminimality constant. See (4.23)-(4.25) in [DS], 
and the sentence that follows. Here it should be
$$
H^d(E\cap H_{j}(Q)) \leq C M H^d(E\cap A_{j}(Q))
+ C h \diam(Q)^{d},
\leqno (2.17)
$$
because of the extra term in (2.10). It is important that the 
constant $C$ does not depend on $j$, so that we can still take $j$ as 
large as we wish; the only price to pay is that we have to assume
that $\diam(Q)^{-d}H^d(E\cap Q)$ is accordingly small at the 
beginning of the argument, if we want 
our story about no full face in $H_j(Q)$, and then (2.17), to  hold.

The way we can deduce the lower bound in (2.16) from this is slightly 
different from what was done in [DS]. First we observe that it is 
enough to prove this lower bound for $x\in E$ such that
$$
\limsup_{r \to 0} \Big[ r^{-d} H^d(E \cap B(x,r))\Big] \geq c_{0},
\leqno (2.18)
$$ 
where $c_{0}$ is some positive dimensional constant. Indeed, if 
$c_{0}$ is chosen small enough, (2.18) holds for almost-every 
$x\in E$ (see for instance [Ma], Theorem 6.2 on p.89,
or even [D2], p.113). 
Then every point $x\in E^\ast$ is a limit of points $x_k \in E$ such that 
(2.18) holds, and we easily get the desired lower bound 
on $H^d(E \cap B(x,r))$ by using the lower bound that we can get for a 
smaller ball centered at some $x_k$.

Suppose that (2.18) holds, but not the lower bound in (2.16). Then
we can find $r>0$ such that $r < \delta$, $B(x,2r) \i U$, but
$H^d(E \cap B(x,r)) < \alpha r^d$, where $\alpha$ is the inverse of 
$C$ in (2.16), and can be chosen as small as we want. Set 
$\eta = 2^{-d-3}$.
If $\alpha$ is sufficiently small compared to $c_{0}$, we can pick a 
new radius $r_{1} < r$ such that
$$
H^d(E \cap B(x,r_{1})) \leq \alpha r_{1}^d
\leqno (2.19)
$$
but
$$
H^d(E \cap B(x, \eta r_{1})) \geq \alpha (\eta r_{1})^d.
\leqno (2.20)
$$
Indeed, we can try successively all the radii $\eta^k r$, $k \geq 0 \,$; 
if none of them works, we can easily prove by induction that 
$H^d(E \cap B(x,\eta^k r)) < \alpha (\eta^{k}r)^d$
for every $k>0$, but if $\alpha$ is small enough this gives a very 
small lower bound on $t^{-d} H^d(E \cap B(x,t))$ for every $t\leq r$, 
which contradicts (2.18).

Choose $r_{1}$ as in (2.19) and (2.20). Now we can pick a cube $Q$ 
centered at $x$, such that 
$$
B(x,2\eta r_{1}) \i Q \i B(x,r_{1}/2)
\leqno (2.21)
$$
and 
$$
H^d(E\cap A_{j}(Q)) \leq C 2^{-j} H^d(E \cap B(x,r_{1})) 
\leq C 2^{-j} \alpha r_{1}^d.
\leqno (2.22)
$$
Indeed we can easily find $C^{-1} 2^j$ cubes $Q$ such that (2.21) 
holds, and for which the $A_{j}(Q)$ are disjoint, so (2.22) is just 
obtained by picking the one for which $H^d(E\cap A_{j}(Q))$ is 
the smallest. On the other hand,
$$
H^d(E\cap H_{j}(Q)) \geq H^d(B(x,\eta r_{1})) 
\geq \alpha \eta^d r_{1}^d
\leqno (2.23)
$$
by (2.21) and (2.20), so (2.22) says that
$$
H^d(E\cap A_{j}(Q)) \leq C 2^{-j} H^d(E\cap H_{j}(Q)).
\leqno (2.24)
$$

Recall that we can take $j$ as large as we want. If we do so, and 
then take $\alpha$ small enough (to make sure that (2.17) applies to 
$Q$), the term $C M H^d(E\cap A_{j}(Q))$ in (2.17) is smaller than
${1 \over 2} H^d(E\cap H_{j}(Q))$, and so (2.17) says that
$H^d(E\cap H_{j}(Q)) \leq 2Ch \diam(Q)^d$. Recall that $\diam(Q)$ and 
$r_{1}$ are comparable, by (2.21). So if $h$ is now chosen small 
enough, we get the desired contradiction with (2.23). This completes 
our proof of Lemma 2.15. \qed

\medskip
The next important property of quasiminimal sets is the local uniform 
rectifiability of $E^\ast$. Let us not give a precise statement here, 
and refer to Theorem 2.11 in [DS] 
instead. We claim that this extends 
to $E\in GAQ(M,\delta ,U,h)$, provided again that we take $h$ small 
enough, depending on $n$, $d$, and $M$.

The proof in [DS] is rather long and complicated, but fortunately a 
lot of it consists in manipulations on Ahlfors-regular sets, and 
there are only a limited number of places where the quasiminimality 
of $E$ is used.

In Chapter 5, it is used to get a Lipschitz mapping from $E$ to 
$\R^d$, with big image. Quasiminimality is used in (5.7), but not in 
(5.8)-(5.10) (that is, only through Ahlfors-regularity), so the 
effect of  an additional $Ch\diam(Q)^d$ in the right-hand side of 
(5.7) is just an extra $Ch\diam(Q)^d$ in (5.11). If $h$ is small 
enough, this is easily eaten up by the $C^{-1}\diam(Q)^d$ in the 
left-hand side, and we can continue the argument with no further 
modification.

In Chapter 6 we show that if $E \i U \i \R^n$ is a quasiminimal set 
and $g : U \to R^l$ is Lipschitz, the new set 
$\widehat S = \{ (x,y) \, ; \,  x\in S \hbox{ and } y=g(x) \} \i 
\R^{n+l}$
is a quasiminimal set too. Quasiminimality is used in (6.10), and our 
extra $Ch\diam(Q)^d$ gets multiplied by a constant $C$ (that depend 
also on the Lipschitz constant for $g$) in (6.12) and (6.13). Then we 
add things that do not depend on quasiminimality, and eventually we 
get an extra $Ch\diam(Q)^d$ in (6.16). That is, 
$\widehat S \in GAQ(\widehat M,\delta ,\widehat  U, C h)$, where 
$\widehat  U = U \times R^l$, $\widehat M$ is a constant that depends 
on $M$, $n$, $d$ and the Lipschitz constant for $g$, and $C$ also 
depends on these constants.

Chapter 7  only  talks about general Ahlfors-regular sets, and 
Chapter 8 talks about how to deduce the local uniform rectifiability 
of $E^\ast$ from a main lemma, but only uses results from previous 
chapters, with no comparison argument. We can also relax through most 
of Chapter 9, because the first sections concern a deformation lemma 
that has nothing to do with quasiminimal sets. Even Section 9.2, 
whose goal is to apply the deformation lemma to quasiminimal sets, 
starts with a general discussion with dyadic cubes and 
Ahlfors-regularity. Quasiminimality itself appears only on page 75, 
where a deformation $\phi$ is constructed. Thus $\phi$ in [DS] plays 
the role of $\varphi_{1}$ in (2.2)-(2.6), and $W$ in [DS] is the same 
as $W_{1}$ here. Recall also that $S$ and $k$ there play the roles of 
$E$ and $M$ here.

By (9.92), $W \cup \phi(W) \i B'$, where $B'$ is a ball of radius 
$C_{3}N+nN+C$ that is compactly contained in $U$ (by (9.56)). In this 
argument, $N$ is a very large integer, to be chosen near the end, and 
we shall  see soon that for the present discussion we do not really 
need to know who the large constant $C_{3}$ is, because our small 
constant $h$ can be chosen last. We can make sure that our set 
$\widehat W$ from (2.5) is also contained in $B'$, because $B'$ is 
convex (see Remark 2.9.a above).

The quasiminimality of $E$ is used in (9.93), and we should add an 
extra $h \diam(B') \leq h(C+C_{3})N^d$ to the right-hand side, 
because of the additional term in (2.11). 
Then there are a few inclusions and identities, and (9.93) is not 
used until (9.103), where we should also add $ChN^d$ to the 
right-hand side
(let us drop the dependence on $C_{3}$ and other constants). Thus 
(9.103) should be replaced with
$$
H^d(E\cap W) \leq CMN^{d-1} + C hN^d.
\leqno(2.25)
$$
The computations after (9.103) can stay the same; they just use 
previous results to show that $H^d(E\cap W) \geq C^{-1}N^d$ (see 
(9.105) and two lines above it). The desired contradiction is then 
derived by choosing $N$ large enough and, in our case, $h$ small.
This takes care of local uniform rectifiability.

We are also interested in Chapter 10 of [DS], 
where one shows that 
$E^\ast$ also has big projections locally. Quasiminimality shows up 
in (10.22), which should be replaced with
$$
H^d(E\cap W) \leq MH^d(\phi(E\cap W)) + C h r_{1}^d.
\leqno(2.26)
$$
But since we see (independently) in (10.25) that 
$H^d(E\cap W) \geq C^{-1} r_{1}^d$,
the extra term in (2.26) can be eaten up by the term in the left-hand 
side, and then one can proceed as in [DS]. 

This  completes our discussion of the extension  of results from 
[DS], and in particular of Theorem 2.11 there, to generalized 
quasiminimal sets in $GAQ(M,\delta,U,h)$.
But we are not finished yet, because now we want to extend results 
from [D1], 
that concern limits of sequences of quasiminimal sets.

\medskip
\noindent {\bf 3. Limits of quasiminimal sets}
\medskip

We want to show that the Hausdorff measure goes to the limit well along 
bounded sequences of generalized quasiminimal sets; this will be used in 
the next section to prove the equivalence of various definitions, and 
also in later sections, although we shall only apply it to 
almost-minimal sets then.

Let $U \i  \R^n$, $M\geq 1$, $\delta >0$ and $h\geq 0$ be given, and 
let $\{ E_{k} \}$ be a sequence in $GAQ(M,\delta,U,h)$. We shall 
assume that  the $E_{k}$ are reduced, i.e, that $E_{k}^\ast = E_{k}$
(otherwise, we could always replace $E_{k}$ with $E_{k}^\ast$, which 
is a generalized quasiminimal set with the same constants by 
Remark~2.14).

We shall also assume  that $\{ E_{k} \}$ converges to a closed set 
$E$ in $U$. By this we mean that for each compact set $H \i U$, 
$$
\lim_{k \to + \infty} d_{H}(E,E_{k}) = 0,
\leqno (3.1)
$$
where the local variant $d_{H}$ of the Hausdorff distance is defined 
by
$$
d_{H}(E,F) = \sup \{ \dist(x,F) \, ; \, x \in E \cap H \}
+ \sup \{ \dist(x,E) \, ; \, x \in F \cap H \}.
\leqno (3.2)
$$

We use the convention that 
$\sup \{ \dist(x,F) \, ; \, x \in E \cap H \}=0$ when $E \cap H$
is empty.
By the way, if we did not assume the $E_{k}$ to be 
reduced, they could converge to essentially anything, because we 
could not keep track of the $E_{k}\setminus E_{k}^\ast$. Indeed it is 
easy to find sequences of finite sets that converge to any given 
closed set $E$. 

\medskip
\proclaim Lemma 3.3. Let $U, M, \delta, h$, $\{ E_{k} \}$, and $E$ be 
as before, and assume that $h$ is small enough, depending on $n$, 
$d$, and $M$. Then $E$ is reduced,
$$
H^d(E \cap V) \leq \liminf_{k \to +\infty} \, H^d(E_{k} \cap V)
\hbox{ for every open set } V \i U,
\leqno (3.4)
$$
and 
$$
E \in GAQ(M,\delta,U,h)
\hbox{ (i.e., with the same constants).}
\leqno (3.5)
$$
\medskip

When $h=0$, this is the main result of [D1]. 
Here we just have to say 
why the proof of [D1] goes through when $h$ is small. 
As we shall see soon, $E$ is locally Ahlfors-regular, so it is reduced.
For (3.4) (i.e., Theorem 3.4 in [D1]), 
there is nothing to do: the 
proof just consists in putting together results from [DS] 
that we already discussed in Section 2. So  we may now discuss the proof of 
(3.5) (Theorem 4.1 in [D1]).

The verification that $E$ is locally Ahlfors-regular, with big pieces 
of Lipschitz graphs, does not need to be modified (it is again a 
mechanical consequence of the results of Section~2).  The same 
comment applies to the remarks about rectifiability and the existence 
of tangent planes, which leads us near (4.13) in [D1]. 

Then we come to the verification of quasiminimality itself. Let us  
remind the reader of how it goes. We are given a deformation 
$f : U\to U$, which can be written as $f = \varphi_{1}$ for a family 
of mappings $\varphi_{t}$ that satisfy (2.2)-(2.6), where (2.6) holds 
for some $\delta' \leq \delta$, and we want to 
prove (2.11), i.e., that
$$
H^d(E \cap W_{f}) \leq M H^d(f(E \cap W_{f})) + h (\delta')^d,
\leqno (3.6)
$$
where we set $W_{f} = \{ x\in U ; f(x) \neq x \}$. We know that the 
sets $E_{k}$ have the same property, but we cannot take limits 
directly, because $H^d(f(E \cap W_{f}))$ could be much smaller than 
the 
$H^d(f(E_{k} \cap W_{f}))$. A typical unfriendly situation where this 
may happen would be when $f$ is not injective on $E$ and sends two 
different pieces of $E$ to the same set, while it is still injective 
on each of the $E_{k}$.

So we replace $f$ with a different deformation $h$ which does not 
have the same defect. That is, we try to force $h$ to merge points 
systematically near the places where we know that it merges points of $E$.
This is easier to do at places where $E$ is close to a tangent plane 
$P$ and $f$ is close to a linear function $A$, because there we can 
compose $f$ with the orthogonal projection onto $A(P)$. A good part 
of the argument of [D1] consists in doing this in many places, 
gluing the different functions that we construct, and checking that the 
places where we do  the gluing or where we cannot control things do 
not matter too much, because the part of $E$ that lives there has 
small measure.

Fortunately we can keep the construction of $h$ as it is in [D1]; 
we just have to say how the final estimates go. Things start being
dangerous near (4.94), when the function $h$ is finally defined and 
we introduce its variant $h_{1}$. Quasiminimality is used in (4.95), 
with the function $h_{1}$ and the set $E_{k}$. So we add the usual 
$Ch (\delta')^d$ to the right-hand side of (4.95), and get that
$$
H^d(E_{k}\cap W_{h_{1}}) \leq M H^d(h_{1}(E_{k}\cap W_{h_{1}}))
+ h (\delta')^d.
\leqno (3.7)
$$
Then we quietly follow the computations, up to (4.107) which gives an 
upper bound on $H^d(h_{1}(E_{k}\cap W_{h_{1}}))$. We pick the various 
constants as explained after (4.107), so as to get the equivalent but 
simpler (4.108), and then prove that
$$
H^d(h_{1}(E_{k}\cap W_{h_{1}}))
\leq H^d(f(E\cap W_{f})) + 2\eta_{1},
\leqno (3.8)
$$
where  $\eta_{1}$ is as small as we want. In [D1], (3.8) follows from 
(4.108) (our upper bound for $H^d(h_{1}(E_{k}\cap W_{h_{1}}))$) and a 
lower bound (4.109) for $H^d(f(E\cap W_{f}))$. This bound itself 
is the direct consequence of Lemma 4.111.

Unfortunately, we need to say a little more about the proof of 
Lemma 4.111, because it is very badly written. The author claims at 
the beginning of the proof that $E$ is a quasiminimal set, and this 
is a rather stupid thing to say, because this is more or less what 
we try to prove.
So the reader is asked not to pay attention to any of the remarks 
about quasiminimality, and follow quietly the construction of the 
sets $S_{1}$
(just the intersection of $E$ with a ball $B_{1}$), then 
$S_{2} = \{ (x,f(x)) \, ; \, x\in S_{1} \} \i \R^{2n}$, and then 
$S_{3}$, which is the image of $S_{2}$ by some simple affine 
bijection $\varphi$ of $\R^{2n}$, whose goal is just to make some 
tangent plane horizontal.
Then we prove (4.113), which says that $S_{3}$ stays very close to a 
horizontal $d$-plane in some ball. This part only uses the geometry
of the construction, so it is all right.

Now replace $E$ with $E_{k}$; this gives a new set $S_{3,k}$, and 
this set is a quasiminimal set, by the same argument that was falsely 
applied to $E$ in [D1]; this involves constructions from [DS], in 
particular Chapter 5, which also apply to generalized quasiminimal 
sets, by Section 2. In addition, $E_{k}$ also satisfies the distance 
estimate (4.113)
for $k$ large, because the $E_{k}$ converge to $E$.
So we can apply Lemma 2.18 in [D1] (or Chapter 10 of [DS]) 
to get that for $k$ large, the projection of $S_{3,k}$ on 
some $d$-plane contains some fixed ball. Then this is also 
true of $S_{3}$, by compactness, and we can conclude as was 
claimed in [D1]. 

So (3.8) holds, and hence
$$
H^d(E_{k}\cap W_{h_{1}}) \leq M H^d(f(E\cap W_{f})) + 2M\eta_{1}
+ h (\delta')^d
\leqno (3.9)
$$
by (3.7). We also know that
$$
H^d(E \cap W_{h_{1}}) \leq H^d(E_{k}\cap W_{h_{1}}) + \eta_{1}
\leqno (3.10)
$$
for $k$ large, by (3.4). Recall from (4.94) in [D1] 
and the line below it that $W_{h_{1}} \i W_{f}$, 
but we may choose the function $\varphi$ in 
(4.94) so that the difference is as close to $\partial W_{f}$ as we 
want. As a result, we can make $H^d(E \cap W_{f} \setminus W_{h_{1}})$ 
(incidentally, another typo in [D1]) as small as we want. 
When we combine this with (3.10) and (3.9), and then let small constants 
like $\eta_{1}$ tend to $0$, we get that
$$ 
H^d(E \cap W_{f}) \leq M H^d(f(E\cap W_{f})) + h (\delta')^d.
\leqno (3.11)
$$
This is just the same as (2.11), and this completes our proof of 
Lemma 3.3. \qed

\medskip
We shall also need an upper semicontinuity result for $H^d$, which 
unfortunately the author did not care to write down in [D1]. 

\medskip
\proclaim Lemma 3.12. 
Let $U, M, \delta, h$, $\{ E_{k} \}$, and $E$ be 
as in Lemma~3.3. In particular, assume that $h$ is small enough, 
depending on $d$ and $M$. Then
$$
\limsup_{k \to +\infty} \, H^d(E_{k} \cap H)
\leq (1+Ch) M H^d(E \cap H)
\hbox{ for every compact set } H \i U.
\leqno (3.13)
$$
Here $C$ depends only on $d$ and $M$.

\medskip
This will be a consequence of the rectifiability of $E$.
Since $E$ is rectifiable, $E$ has a $d$-dimensional ``approximate 
tangent plane" $P_x$ at $x$ for $H^d$-almost every $x\in E$. 
See for instance [Ma], Theorem  15.19 on page 212. 
Since in addition $E$ is locally Ahlfors-regular (i.e., (2.16) holds), 
$P_x$ is a true tangent plane at $x$. See for instance [D2],  
Exercise 41.21 on page~277. In addition,
$$
\lim_{r \to 0} r^{-d} H^d(E\cap B(x,r)) = c_d
\leqno (3.14)
$$
for $H^d$-almost every $x\in E$, where 
$c_d$ is the $H^d$-measure of the unit ball in $\R^d$.
See again [Ma], Theorem 16.2 on page 222. 
Call $E'$ the set of points $x\in E$ such that $E$ has a tangent plane 
$P_x$ at $x$ and (3.14) holds; thus almost-every point of $E$ lies in $E'$.

Let $\varepsilon > 0$ be small, and denote by $\cal B$ the set of balls
$B= B(x,r)$ such that $x\in E' \cap H$, $2r < \delta$,
$$
B(x,2r) \i \widehat H,
\leqno (3.15)
$$
where $\widehat H$ is some compact subset of $U$
that we choose so that $H \i \inter(\widehat H)$,
$$
H^d(E \cap \partial B) = 0,
\leqno (3.16)
$$
$$
\dist (y,P_x) \leq \varepsilon r
\ \hbox{ for } y \in E \cap B(x,2r),
\leqno (3.17)
$$
and 
$$
|r^{-d} H^d(E\cap B(x,r)) - c_d| \leq \varepsilon.
\leqno (3.18)
$$

Notice that for $x\in E' \cap H$, $B(x,r) \in \cal B$ for almost-every
small enough $r \,$: (3.17) and (3.18) come from the definition of $E'$, and the 
set of small $r$ such that (3.16) fails is at most countable, because $H^d(E)$ 
is locally finite. Thus the balls $B, B \in \cal B$, form what is called a 
Vitali covering of $E'\cap H$. Let us check that there is an at most 
countable collection $\{ B_i \}_{i\in I}$ in $\cal B$ such that the 
$B_i \, , \, i\in I$, are disjoint and almost cover $E'\cap H$ in the sense that
$$
H^d \big([E'\cap H] \setminus  \bigcup_{i \in I} B_i \big) = 0.
\leqno  (3.19)
$$
Let us define successive finite collections $\{ B_i \}_{i\in I_k}$.
Start with $I_0 = \emptyset$. If $I_k$ was already defined, 
set $H_k = [E' \cap H] \setminus \bigcup_{i\in I_k} \overline B_i$. 
For each $x\in H_k$, there are lots of balls $B \in \cal B$ 
centered on $x$, and such that $5B$ does not meet the $B_i$, $i\in I_k$.
This  gives a covering of $H_k$ by balls of $\cal B$; by the standard 
$5$-covering lemma (see the first pages  of [St]), 
we can find an at most countable subcollection $B_i$, $i\in J_k$,
such that the $B_i$, $i\in J_k$, are disjoint, and the $5B_i$, $i\in J_k$,
cover $H_k$. Then
$$
\sum_{i\in J_k} H^d(E\cap B_i)
\geq C^{-1} \sum_{i\in J_k} H^d(E\cap 5B_i)
\geq C^{-1} H^d(H_k),
\leqno  (3.20)
$$
by (2.16) or (3.18) and because the $5B_i$ cover $H_k$.
Notice that 
$$
\sum_{i\in J_k} H^d(E\cap B_i) \leq H^d(E\cap \widehat H) < +\infty
\leqno  (3.21)
$$
because the $B_i$ are disjoint and contained in $\widehat H$ (by (3.15)), 
$\widehat H$ can be covered by a finite collection of compact balls in $U$, 
and by (2.1). So we can pick a finite collection $J'_k \i  J_k$, so that
$$
\sum_{i\in J'_k} H^d(E\cap B_i)
\geq {1 \over 2} \sum_{i\in J_k} H^d(E\cap B_i)
\geq (2C)^{-1} H^d(H_k),
\leqno  (3.22)
$$
and then set $I_{k+1} = I_k \cup J'_k$. This gives a definition
of $I_k$ by induction, and then we take $I = \cup_{k} I_k$.
The $B_i$, $i\in I$,  are disjoint by construction. Then
$$
\sum_{k}\sum_{i\in J'_k} H^d(E \cap B_i)
\leq H^d(E \cap \widehat H) < + \infty,
\leqno  (3.23)
$$
by the proof (3.21). In particular, the series in $k$ converges, and
(3.22) says that $H^d(H_k)$ tends to $0$. But
$$
H^d \big([E'\cap H] \setminus  \bigcup_{i \in I} B_i \big)
\leq H^d \big([E'\cap H] \setminus  \bigcup_{i \in I_k} B_i \big)
= H^d \big([E'\cap H] \setminus  \bigcup_{i \in I_k} \overline B_i \big)
= H^d(H_k),
\leqno  (3.24)
$$
by (3.16). This proves (3.19).

Let us cover $[E \cap H] \setminus \bigcup_{i \in I} B_i$ by a collection
of open balls $B_j$, $j\in J$, such that $\diam(B_j) < \delta$, $2B_j \i U$
for each $j$, and
$$
\sum_{j\in J} \ r_j^d \leq \varepsilon,
\leqno  (3.25)
$$
where $r_j$ denotes the radius of $B_j$. Now the balls $B_i$, 
$i\in I \cup J$, cover $E \cap H$, and by compactness we can find 
finite sets $I' \i I$ and $J' \i J$ such that 
$E \cap H \i \bigcup_{i \in I' \cup J'} B_i$. Since 
the $E_k$ converges to $E$, we also get that
$$
E_k \cap H \i \bigcup_{i \in I' \cup J'} B_i
\leqno  (3.26)
$$
for $k$ large enough, hence
$$
H^d(E_{k} \cap H) \leq \sum_{i \in I' \cup J'} H^d(E_{k} \cap B_i)
\leq \sum_{i \in I'} H^d(E_{k} \cap B_i) + C \varepsilon,
\leqno  (3.27)
$$
by (2.16) and (3.25). Next pick $i \in I'$, and write $B_i = B(x,r)$.
Notice that for $k$ large,
$$
\dist (y,P_x) \leq 2\varepsilon r
\ \hbox{ for } y \in E_k \cap B(x,3r/2), 
\leqno (3.28)
$$
by (3.17) and because the $E_k$ converge to $E$. Also recall that
we only have a finite set $I'$ of indices $i$, so for $k$ large,
(3.28) holds for all of them at the same time. In addition, (3.28) also
holds if we replace $P_x$ with any plane $P$ parallel to $P_x$, and
close enough to it.

We want to use the quasiminimality of $E_k$, so we define 
a Lipschitz deformation $\varphi$. Take 
$$
\varphi(y)=y \ \hbox{ for } \ y\in 
[\R^n \setminus B(x,r)] \cup \big\{ y \in \R^n \, ; \, 
\dist(y,P) \geq 3 \varepsilon r \big\},
\leqno (3.29)
$$
then call $\pi$ the orthogonal projection onto $P$, and set
$$
\varphi(y) = \pi(y) \ \hbox{ on the inside region } \ 
D= \big\{ y \in B(x,(1-\varepsilon) r) \, ; \, 
0 \leq \dist(y,P) \leq 2\varepsilon r \big\}.
\leqno (3.30)
$$
On the remaining annular region 
$A = \big\{ y \in B(x,r) \, ; \, \dist(y,P) < 3 \varepsilon r \big\}
\setminus D$ we choose $\varphi$ so that $\varphi(A) \i \overline B(x,r)$ 
and $\varphi$ is $10$-Lipschitz. For instance, take
$\varphi(x) = \theta(y) y + (1-\theta(y))\pi(y)$ for $y \in A$, 
with $\theta(y) = \Min(1,(2\varepsilon r)^{-1}\dist(y,D))$.
[We added a factor 2 to make sure that $\varphi(A) \i \overline B(x,r)$
if $P$ is close enough to $P_x$, even if it does not exactly go 
through $x$.] 

Set $\varphi_t(y)=t\varphi(y)+(1-t)y$ for $0\leq t \leq 1$;
this family obviously satisfies (2.2)-(2.6), with 
$\widehat W \i \overline B_i$.
Recall that $\diam(B_i) = 2r < \delta$  by definition of $\cal B$.
So we may apply Definition 2.10, and (2.11) holds (with any $\delta' 
\in (2r,\delta)$). That is,
$$
H^d(E_k \cap W_{1}) \leq M H^d(\varphi_1(E_k \cap W_{1}))
+ (2r)^d h. 
\leqno (3.31)
$$
Set $B'=B(x,(1-\varepsilon)r)$ and notice that
if $y \in E_k \cap B'$, then 
$\dist(y,P) \leq 2 \varepsilon r$ by (3.28), so
$\varphi(y) = \pi(y)$ by (3.30). Thus $\varphi(y) = y$
if and only if $y\in P$, and 
$$
E_k \cap W_1 \cap B' = \big\{ y \in E_k \cap B' \, ; \, 
0 < \dist(y,P) \leq 2 \varepsilon r \big\}
 = E_k \cap B' \setminus P.
\leqno (3.32)
$$
There are at most countably many planes $P$ parallel to $P_x$
and for which $H^d(E_k \cap B(x,r) \cap P) > 0$, so we may choose $P$ such that
this is not the case. Then
$$\eqalign{
H^d(E_k \cap B')&= H^d(E_k \cap W_{1} \cap B') 
\leq M H^d(\varphi(E_k \cap W_{1})) + (2r)^d h
\cr& 
\leq M H^d(\varphi(E_k \cap B')) 
+ M H^d(\varphi(E_k \cap B(x,r) \setminus B')) + (2r)^d h
\cr&
\leq M H^d(\varphi(E_k \cap B')) + 10^d M H^d(E_k \cap B(x,r) \setminus B')
+ (2r)^d h
}\leqno (3.33)
$$
by (3.32) and (3.31), because $\varphi_1 = \varphi$, 
and because $E_k \cap W_{1} \i E_k \cap B(x,r) 
\i (E_k \cap B') \cup (E_k \cap B(x,r) \setminus B')$
and $\varphi$ is $10$-Lipschitz. Next let us check that
$$
H^d(E_k \cap B(x,r) \setminus B') \leq C \varepsilon r^d,
\leqno (3.34)
$$
where $C$ may depend on the quasiminimality constants, but not on $\varepsilon$.
We can cover $P_x \cap B(x,r) \setminus B'$ by $C \varepsilon^{-d+1}$ balls 
$B(z_j,\varepsilon r)$ centered on $P_x \cap B(x,r) \setminus B'$. 
By (3.28), the $B(z_j,3\varepsilon r)$ cover $E_k \cap B(x,r) \setminus B'$.
Notice also that the $B(z_j,10\varepsilon r)$ are contained in $U$,
by (3.15). Then $H^d(E_k \cap B(z_j,3\varepsilon r)) \leq C \varepsilon^d r^d$
by (2.16), because the $E_k$ all lie in the same $GAQ(M,\delta,U,h)$, and if 
$h$ is small enough; (3.34) follows.
Now
$$\eqalign{
H^d(E_k \cap B(x,r)) &= H^d(E_k \cap B') + H^d(E_k \cap B(x,r) \setminus B')
\cr&\leq M H^d(\varphi(E_k \cap B')) + C H^d(E_k \cap B(x,r) \setminus B') 
+ (2r)^d h
\cr&\leq M H^d(\varphi(E_k \cap B')) + C \varepsilon r^d + (2r)^d h
}\leqno (3.35)
$$
by (3.33) and (3.34).
Recall from (3.28) and (3.30) that $\varphi(y) = \pi(y)$ for 
$y\in E_k \cap B'$, so $\varphi(E_k \cap B') \i P \cap B(x,r)$ and 
$$
H^d(E_k \cap B(x,r)) \leq M H^d(P \cap B(x,r)) 
+ C \varepsilon r^d + (2r)^d h
\leq M c_d \, r^d + C \varepsilon r^d + (2r)^d h,
\leqno (3.36)
$$
where $c_d$ still denotes the $H^d$-measure of the unit ball in $\R^d$.
On the other hand, (3.18) says that
$$
H^d(E \cap B(x,r)) \geq c_d r^d - \varepsilon r^d.
\leqno (3.37)
$$
Incidentally, we could also have obtained something like (3.37)
directly from (3.17) and the facts that $x \in E$ and $E$ is a reduced 
quasiminimal sets as in (3.5), by the same proof as in 
Lemma~10.10 in [DS]. 
But (3.18) is not that hard to obtain anyway.

Recall that $B_i = B(x,r)$; when we put (3.36) and (3.37) together,
we get that 
$$
H^d(E_k \cap B_i) 
\leq M H^d(E \cap B_i) + C \varepsilon r^d + (2r)^d h
\leq (1+C \varepsilon + Ch) M H^d(E \cap B_i).
\leqno (3.38)
$$
Recall that the $B_i$, $i\in I'$, are disjoint, so 
$$\eqalign{
H^d(E_{k} \cap H) 
&\leq \sum_{i \in I'\cup J'} H^d(E_{k} \cap B_i) 
\leq \sum_{i \in I'} H^d(E_{k} \cap B_i) + C \varepsilon
\cr&
\leq (1+C \varepsilon + Ch) M \sum_{i \in I'} H^d(E \cap B_i)
+ C\varepsilon
\cr&
\leq (1+C \varepsilon + Ch) M H^d\big(E \cap (\cup_{i \in I'}B_i)\big)
+ C\varepsilon
}\leqno (3.39)
$$
by (3.26), (2.16) and (3.25), and then (3.38). 
[We can apply (2.16) with uniform bounds, by (3.14).]

Let $V$ be any open neighborhood of $E\cap H$; we can choose $\cal B$
above so that all the $B_i$ are contained in $V$ (just add an extra 
smallness condition on $r$ near (3.14), or equivalently make 
$\widehat H$ smaller). Then 
$$
H^d(E_{k} \cap H) \leq (1+C \varepsilon+Ch) M H^d(E \cap V) + C\varepsilon
\leqno (3.40)
$$
for $k$ large, by (3.39). Since $\varepsilon$ is as small as we want,
$$
\limsup_{k \to +\infty} H^d(E_{k} \cap H) \leq (1+Ch) M H^d(E \cap V).
\leqno (3.41)
$$
Finally, $H^d(E \cap V)$ is as close to $H^d(E \cap H)$ as we want, 
because $V$ is any open neighborhood of $E\cap H$ and 
$H^d(E \cap V)<+\infty$ for some $V$. So 
$$
\limsup_{k \to +\infty} H^d(E_{k} \cap H) \leq (1+Ch) M H^d(E \cap H),
\leqno (3.42)
$$
and this completes our proof of Lemma~3.12.
\qed

\vfe  
\noindent {\bf 4. Almost-minimal sets}
\medskip

For us, almost-minimal sets will be sets that look more and more like 
minimal sets when you look at them at smaller scales. 
There are a few different possible definitions, even in the same context 
as in Sections 2 and 3 where we only consider Almgren competitors. We shall give 
three of them. The first one (with $A_+$) is the closest to our definition
of quasiminimal sets (Definition 2.7), and it has the advantage that we 
essentially do not need this section or many details of the previous two to use it. 
So the reader that would not be interested in subtleties about classes of 
almost-minimal sets is encouraged to have a look at Definition~4.1 and skip 
the rest.

The other two notions (with $A$ and $A'$) are weaker, and turn out to be 
equivalent to each other (see Proposition~4.10). 
For $A$, which is very close to our notion of 
generalized quasiminimal sets, we have a nice convergence result
(Lemma 4.7). The definition for $A'$ is a little bit simpler, and 
this is the definition that we used for $Al$-minimal sets (in 
Definition~1.6), so the equivalence between the two will make some of
our statements a little simpler.

For all these definitions (with $A_+$, $A$, and $A'$),
we shall measure closeness to minimality with a nondecreasing function 
$h : (0,+\infty) \to [0,+\infty]$, such that 
$\lim_{\delta \to 0} h(\delta) = 0$. We shall call such a function a 
\underbar {gauge function}.

We start with a first notion which is reasonably close to Definition~2.7. 

\medskip
\proclaim Definition 4.1. Let $E$ be a closed subset of $U$. We say 
that $E$ is an $A_+$-almost-minimal set (of dimension $d$) in $U$, 
with gauge function $h$, if (2.1) holds and, for each $\delta >0$ 
and each family $\{\varphi_t \}_{0 \leq t \leq 1}$ of
deformations such that (2.2)-(2.6) hold,
$$
H^d(E \cap W_{1}) \leq (1+h(\delta)) \,  H^d(\varphi_1(E \cap W_{1})),
\leqno (4.2)
$$
where $W_{1} = \{ x\in U ; \varphi_1(x) \neq x \}$, as in (2.5).

\medskip
In other words, we require that for each $\delta > 0$, $E$ be an 
$(1+h(\delta),\delta)$-quasiminimal set in $U$ (as in Definition 2.7).
The main difference is that $M(\delta) = 1+h(\delta)$ tends to $1$
when $\delta$ tends to $0$. Observe that we allowed the possibility 
that $h(\delta)=+\infty$ for $\delta \geq \delta_0$, to make it plain 
that in some cases we do not have any information coming from large scales. 
This is the same as saying that we only consider almost-minimality
at scales $\delta \leq \delta_0$.

Definition 4.1 looks nice, but we shall prefer to use slightly 
weaker versions where the error term is allowed not to depend on 
the set $W_1$ where $\varphi_1(x) \neq x$, but only on $\delta$. The 
next definition is slightly simpler, but still uses the same sort 
of accounting method as in Definitions 2.7 and 4.1.

\medskip
\proclaim Definition 4.3.
We say that the closed set $E$ in $U$ is an $A$-almost-minimal set 
in $U$, with gauge function $h$, if (2.1) holds and, for each $\delta > 0$ 
and each family $\{\varphi_t \}_{0 \leq t \leq 1}$ such that (2.2)-(2.6) hold,
$$
H^d(E \cap W_{1}) \leq H^d(\varphi_1(E \cap W_{1}))
+ h(\delta) \delta^d.
\leqno (4.4)
$$
Again, $W_{1} = \{ x\in U ; \varphi_1(x) \neq x \}$ as in (2.5).

\medskip
\noindent {\bf Remark 4.5.}
So we just replaced the error term $h(\delta) H^d(\varphi_1(E \cap W_{1}))$ 
in (4.2) with the simpler $h(\delta) \delta^d$. 
Let us check that $A_+$-almost-minimal set are 
automatically $A$-almost-minimal sets, at least in any smaller open set
$U_\tau = \{ x \in U \, ; \, \overline B(x,\tau) \i U \}$. 

Choose $\delta < \tau$ such that $h(\delta) < +\infty$, and let 
$\{\varphi_t \}$ satisfy (2.2)-(2.6). Observe that 
$\diam(W_{1}) \leq \delta$, by (2.6), and that 
$\dist(W_{1},\R^n \setminus U) \geq \delta$, because $W_{1} \i U_\tau$
(by (2.6)) and we took $\delta < \tau$. Then 
$H^d(E \cap W_{1}) \leq C \delta^d$, by the local Ahlfors-regularity of
$E^\ast$ (see Lemma~2.15, but we are still in the standard 
quasiminimal case). Here $C$ depends on $M = 1+h(\delta)$, but this is 
all right.

If $H^d(E \cap W_{1}) \leq H^d(\varphi_1(E \cap W_{1}))$, then 
(4.4) holds stupidly. Otherwise, 
$$\eqalign{
H^d(E \cap W_{1}) 
&\leq  (1+h(\delta)) H^d(\varphi_1(E \cap W_{1}))
\leq H^d(\varphi_1(E \cap W_{1})) + h(\delta) H^d(E \cap W_{1})
\cr&\leq H^d(\varphi_1(E \cap W_{1})) + C h(\delta) \delta^d,
}\leqno (4.6)
$$
by (4.2). This is (4.4), with $\widetilde h(\delta) = C h(\delta)$.

The author does not know for sure that Definition 4.1 is strictly 
stronger than Definition 4.3, but it also seems difficult to prove 
that they are equivalent. This is why we spent some time dealing with 
generalized quasiminimal sets in the previous sections. Notice that 
$E$ is an $A$-almost-minimal set in $U$, with gauge function $h$,
if and only if it lies in the class $GAQ(1,\delta,U,h(\delta))$ 
of Definition 2.10 for every $\delta >0$ . If $\delta$ is small,
$h(\delta)$ is as small as we want, and then we can apply all the
results of Sections 2 and 3. In particular, we have the following 
consequence of Lemma 3.3.

\medskip
\proclaim Lemma 4.7. 
Let $\{ E_k \}$ be a sequence of $A$-almost-minimal sets in $U$, with 
the same gauge function $h$. 
Assume that each $E_k$ is reduced (i.e., that $E_k^\ast = E_k$), and that
the sequence $\{ E_k \}$ converges to the closed set $E$ in $U$
(i.e., (3.1) holds for each compact $H \i U$). Then $E$ is a
reduced $A$-almost-minimal set in $U$, with the same gauge function $h$.

\medskip
To be fair, since Lemma 3.3 needs $h$ to be small enough, it seems 
that we should replace $h(\delta)$ with $+\infty$ for $\delta$ large,
but in fact we do not even need to do this. We can apply Lemma 3.3 with 
$\delta$ small, and this gives (3.4) and all the good local properties 
of $E$ and the $E_k$, which we can then use as in the proof of Lemma 3.3
to give the desired result even when $h(\delta)$ is large. But anyway it 
makes no real difference, because we shall not be able to use the 
conclusion when $h(\delta)$ is large. \qed

\medskip
Let us now give a variant of Definition 4.3 where the accounting is
done more simply, as in Definitions 1.4 and 1.6.

\medskip
\proclaim Definition 4.8.
Let $E$ be a closed subset of $U$ such that (2.1) holds. 
We say that $E$ is an $A'$-almost-minimal set in $U$, with 
gauge function $h$, if for each $\delta > 0$ and each family 
$\{\varphi_t \}_{0 \leq t \leq 1}$ such that (2.2)-(2.6) hold,
$$
H^d(E \setminus F) \leq H^d(F \setminus E)
+ h(\delta) \delta^d,
\leqno (4.9)
$$
were we set $F = \varphi_1(E)$.

\medskip
When $E$ and $F$ coincide out of some compact ball $B \i U$,
(4.9) simply says that 
$H^d(E\cap B) \leq H^d(F\cap B)+ h(\delta) \delta^d$;
we wrote (4.9) the way we did just to avoid mentioning the existence
of such a $B$.
Notice that the $Al$-minimal sets of Definition 1.6 are the same as 
$A'$-almost-minimal sets in $\R^n$, with gauge function $h=0$.
Definition 4.8 looks simpler, but is in fact equivalent to 
Definition 4.3.

\medskip
\proclaim Proposition 4.10.
Let $E$ be a closed subset of $U$ such that (2.1) holds. Then
$E$ is an $A$-almost-minimal set in $U$ with gauge function $h$
if and only if $E$ is an $A'$-almost-minimal set in $U$ with 
(the same) gauge function $h$.
\medskip

First suppose that $E$ is an $A'$-almost-minimal set, and let 
$\{ \varphi_t \}$ satisfy (2.2)-(2.6). We know (4.9) and we want to 
check (4.4).

Set $E_{1}=E \cap W_{1}$, $E_{2}=E\setminus E_{1}$, and  
$F=\varphi_{1}(E)$. Then
$$
H^d(E \cap W_{1}) = H^d(E_{1}) 
= H^d(E_{1}\cap F) + H^d(E_{1}\setminus F)
\leqno (4.11)
$$
trivially. Now $E_{2}\i F$ because $\varphi_{1}(x) = x$ on $E_{2}$,
so $E_{1}\setminus F = E\setminus F$. Then 
$$
H^d(E \cap W_{1}) = H^d(E_{1}\cap F) + H^d(E\setminus F)
\leq H^d(E_{1}\cap F) + H^d(F\setminus E) + 
h(\delta)\delta^d,
\leqno (4.12)
$$
by (4.11) and (4.9). Observe that 
$F = \varphi_{1}(E) = \varphi_{1}(E_{1}) \cup \varphi_{1}(E_{2})$ 
and $\varphi_{1}(E_{2}) = E_{2} \i E$ because $\varphi_{1}(x)=x$ on 
$E_{2}$. Then $E_{1}\cap F = E_{1}\cap [\varphi_{1}(E_{1}) \cup 
\varphi_{1}(E_{2})] = E_{1}\cap \varphi_{1}(E_{1})$ because
$E_{1}\cap \varphi_{1}(E_{2}) = E_{1}\cap E_{2} = \emptyset$.
Similarly, $F\setminus E = \varphi_{1}(E_{1}) \setminus E$ because 
$\varphi_{1}(E_{2}) \i E$. Thus (4.12) says that
$$\eqalign{
H^d(E \cap W_{1}) &\leq H^d(E_{1}\cap \varphi_{1}(E_{1}))
+H^d(\varphi_{1}(E_{1}) \setminus E)
+ h(\delta)\delta^d
\cr &
\leq H^d(E_{1}\cap \varphi_{1}(E_{1}))
+H^d(\varphi_{1}(E_{1}) \setminus E_{1})
+ h(\delta)\delta^d
\cr &
= H^d(\varphi_{1}(E_{1})) + h(\delta)\delta^d
= H^d(\varphi_{1}(E\cap W_{1})) + h(\delta)\delta^d,
}\leqno (4.13)
$$
as needed for (4.4). This proves our first inclusion.

Now let us assume that $E$ is an $A$-almost-minimal set,
and prove that it is an $A'$-almost-minimal set with the same
gauge function. Again let $\{ \varphi_t \}$ satisfy (2.2)-(2.6); 
we would like to show that (4.9) holds. We cannot deduce this so
easily from (4.4), because if we are too  unlucky, 
$\varphi_{1}(E\cap W_{1})$ may have a big part in $E\setminus W_{1}$, 
which is accounted for in (4.4) (and helps it hold) but not in (4.9). 
This may happen if $W_{1}$ is somewhat smaller than its image; 
if this is the case, we shall need to modify the $\varphi_{t}$ a bit, 
to make $W_{1}$ larger before we apply (4.4).

Set $\varphi = \varphi_{1}$ to save notation. We 
want to replace $\varphi$ with functions $\varphi_{k}$
defined by
$$
\varphi_{k}(x) = \varphi(x) + \varepsilon_{k} \psi(\varphi(x)) \, v,
\leqno (4.14)
$$
where $v$ is a given unit vector in $\R^n$, $\{ \varepsilon_{k} \}$ 
is a sequence that tends to $0$, and $\psi : U \to \R$ is a correctly 
chosen smooth function with compact support in $U$. Recall  that 
there is a relatively compact set 
$\widehat W \i U$, defined by (2.5) and that satisfies (2.6). Let $S$ 
and $S'$ be compact subsets of $U$ such that
$$
\widehat W \i S \i \inter(S') \i S' \i U
\ \hbox{ and } \ 
\diam(S') < \delta;
\leqno (4.15)
$$
other conditions may show up later. We choose $\psi$ so that
$$
0 \leq \psi(x) \leq 1 \hbox{ everywhere,}
\leqno (4.16)
$$
$$
\psi(x) = 1 \hbox{ for } x \in S, \ \hbox{ and } \ 
\psi(x) = 0 \hbox{ for } x \in U \setminus S'.
\leqno (4.17)
$$
Since we want to apply (4.4) to $\varphi_{k}$, we need to check  
that it
is the endpoint of an acceptable family $\varphi_{k,t}$, $0 \leq t 
\leq 1$. Set
$$
\varphi_{k,t}(x) = \varphi_{t}(x) + t \varepsilon_{k} 
\psi(\varphi_{t}(x)) v. 
\leqno (4.18)
$$

The conditions (2.2)-(2.4) clearly hold for this family, so we just 
have to check (2.6). Set $W_{k,t} = \{ x\in U \, ; \, 
\varphi_{k,t}(x) \neq x \}$. Then $W_{k,t} \i S'$, because if
$x\not\in S'$, (4.15) says that $\varphi_{t}(x)=x$, and then
(4.17) and (4.18) say that $\varphi_{k,t}(x)=x$.
Let us also check that 
$$
\varphi_{k,t}(W_{k,t}) \i S'_{k}=:
\big\{ y \, ; \, \dist(y, S') \leq \varepsilon_{k} \big\}.
\leqno (4.19)
$$
Notice that 
$$
|\varphi_{k,t}-\varphi_{t}(x)| \leq \varepsilon_{k}
\hbox{ for } x\in U,
\leqno (4.20)
$$
by (4.18) and (4.16).
If $x \in W_{t}$, then $\varphi_{t}(x) \in \varphi_{t}(W_{t}) \i S$ 
by (4.15), and so $\varphi_{t,k}(x) \in S'_{k}$ by (4.20). 
If $x \in W_{k,t} \setminus W_{t}$, then 
$\varphi_{t}(x) = x \in W_{k,t} \i S'$, and then 
$\varphi_{k,t}(x) \in S'_{k}$, by (4.20) again. This proves (4.19). 

We are now ready to prove (2.6) for the $\varphi_{k,t}$. Set
$\widehat W_{k} = \bigcup_{t}[W_{k,t} \cup\varphi_{k,t}(W_{k,t})]$.
We just proved that $\widehat W_{k} \i S'_{k}$, and it is obvious 
from (4.15) and the compactness of $S'$ that for $k$ large, $S'_{k}$ 
is a compact subset of $U$ with $\diam(S'_{k}) < \delta$. This 
proves (2.6). So we can apply (4.4), and we get that
$$
H^d(E \cap W_{k,1}) \leq H^d(\varphi_{k}(E \cap W_{k,1})) 
+ h(\delta) \delta^d.
\leqno (4.21)
$$

Now we shall cut our sets into pieces to control  the 
overlaps and eventually get something that looks like (4.9). What  
comes out of the small annular region $A = S' \setminus S$ will be 
estimated brutally, and later thrown out to the trash, so we shall 
concentrate on $E_{1}=E \cap S$ and $E_{2} = E \setminus S'$.  

First observe that $\varphi_{k}(x) = x$ for $x\in E_{2}$, by 
(4.15) and (4.17). Also, $\varphi(E_{1}) \i S$ because if 
$\varphi(x) \neq x$ then $x\in W_{1}$ and $\varphi(x) \in S$ by 
(4.15).
Consequently, $\varphi_{k}(E_{1}) \i S'$ for $k$ large enough, by 
(4.20).
In particular, $\varphi_{k}(E_{1})$ does not meet 
$\varphi_{k}(E_{2})=E_{2}$. 
Observe also that $W_{k,1}$ does not meet $E_{2}=E \setminus S'$. Thus
$$\eqalign{
H^d(\varphi_{k}(E \cap W_{k,1})) 
&= H^d(\varphi_{k}(E \cap W_{k,1} \cap S'))
\cr&
\leq H^d(\varphi_{k}(E_{1} \cap W_{k,1})) + H^d(\varphi_{k}(E \cap A))
\cr&
= H^d(\varphi(E_{1} \cap W_{k,1})) + H^d(\varphi_{k}(E \cap A)),
}\leqno (4.22)
$$
where the last equality comes from the fact that $\varphi_{k}-\varphi$ 
is constant on $E_{1}$, by (4.17) and because we just said that 
$\varphi(E_{1}) \i S$.

It is not necessarily true that $E_{1}\i W_{k,1}$, but the difference 
is small: we claim that
$$
E_{1}\setminus W_{k,1} \i Z_{k}=:
\big\{ x\in E_{1} \, ; \, 0 < |\varphi(x)-x| \leq \varepsilon_{k} 
\big\}.
\leqno (4.23)
$$
Indeed let $x\in E_{1}$ be given, and assume that it does not
lie in the set on the right-hand side. 
If $\varphi(x)=x$, then $\varphi_{k}(x) \neq x$ because 
$\varphi(x) \in S$ and (4.17) says that $\psi(\varphi(x)) = 1$. And 
if $|\varphi(x)-x| > \varepsilon_{k}$, (4.20) says that 
$\varphi_{k}(x) \neq x$ too. So $x\in W_{k,1}$ in both cases;
(4.23) follows. Now
$$\eqalign{
H^d(E_{1}) 
&=  H^d(E_{1}\cap W_{k,1}) + H^d(E_{1}\setminus W_{k,1})
\cr &\leq H^d(\varphi_{k}(E \cap W_{k,1}))
+ h(\delta) \delta^d +H^d(Z_{k})
\cr & \leq
H^d(\varphi(E_{1} \cap W_{k,1})) + H^d(\varphi_{k}(E \cap A))
+ h(\delta) \delta^d +H^d(Z_{k})
\cr & \leq
H^d(\varphi(E_{1})) + H^d(\varphi_{k}(E \cap A))
+ h(\delta) \delta^d +H^d(Z_{k})
}\leqno (4.24)
$$
by (4.21) and (4.23), and then (4.22).
To get closer to (4.9), observe that 
$$
H^d(E\setminus \varphi(E)) = H^d(E\cap S'\setminus\varphi(E))
= H^d(E\cap S') - H^d(E\cap S'\cap\varphi(E))
\leqno (4.25)
$$
because $E_{2}=E\setminus S'$ is contained in $\varphi(E)$. 
Similarly,
$$
H^d(\varphi(E)\setminus E) 
\geq H^d(\varphi(E)\cap S'\setminus E) 
= H^d(\varphi(E)\cap S') - H^d(E\cap S'\cap\varphi(E)).
\leqno (4.26)
$$
Set $\Delta = H^d(E\setminus F) - H^d(F\setminus E) 
=H^d(E\setminus \varphi(E)) - H^d(\varphi(E)\setminus E)$. 
This is the quantity that we want to majorize, and we just 
proved that
$$
\Delta \leq H^d(E\cap S') - H^d(\varphi(E)\cap S').
\leqno (4.27)
$$ 
Recall that $\varphi(E_{1}) \i S$ (see four lines above (4.22)), so 
$\varphi(E)\cap S'$ contains $\varphi(E_{1})$ and (4.27) implies that
$$
\Delta \leq H^d(E\cap S') - H^d(\varphi(E_{1}))
= H^d(E_{1}) + H^d(E\cap A) - H^d(\varphi(E_{1})),
\leqno (4.28)
$$
because $E\cap S' = (E \cap S) \cup (E \cap S'\setminus S)
= E_1 \cup (E \cap A)$ by definitions. Then (4.24) yields
$$
\Delta \leq 
H^d(E\cap A) + H^d(\varphi_{k}(E \cap A)) 
+ h(\delta) \delta^d + H^d(Z_{k}).
\leqno (4.29)
$$

Notice that $\varphi_{k} = h_{k} \circ \varphi$, where $h_{k}$
is Lipschitz, with a Lipschitz constant that stays below $2$ when $k$ 
is large enough. Then for $k$ large,
$H^d(\varphi_{k}(E \cap A)) \leq C_{\varphi} H^d(E \cap A)$, with a 
constant $C_{\varphi}$ that depends on $\varphi$ but not $k$ or 
$\psi$.
On the other hand, when $\varepsilon_{k}$ tends to $0$, the set 
$Z_{k}$ in (4.23) decreases to the empty set, and so 
$\lim_{k \to + \infty}H^d(Z_{k}) = 0$. Thus (4.29) yields
$$
\Delta \leq (1+C_{\varphi}) H^d(E \cap A) + h(\delta) \delta^d.
\leqno (4.30)
$$

Recall that  $A = S' \setminus S$, where we only required that
$S$ and $S'$ be compact and satisfy (4.15). We can do this so that
$H^d(E \cap A)$ is as small as we want, and so 
$\Delta \leq h(\delta) \delta^d$, as needed for (4.9). 
Thus $E$ is an $A'$-almost-minimal set, with the same
gauge function $h$ as for its $A$-almost-minimality.
This completes our proof of Proposition 4.10.
\qed

\bigskip
\noindent{\bf B. MONOTONICITY OF DENSITY AND APPROXIMATION BY CONES}
\bigskip
\noindent {\bf 5. Density is nearly nondecreasing}
\medskip

In this section we show that if $E$ is a minimal set of dimension $d$ 
in $\R^n$, then, for each choice of $x\in \R^n$,
$r^{-d}H^d(E\cap B(x,r))$ is a nondecreasing function of $r$. 
We also give a local version of this, and generalizations to almost-minimal 
sets with sufficiently small gauge functions. 
See Propositions 5.16, 5.24, and 5.30 below. Incidentally, all this 
is supposed to be very classical, and the general idea will be to compare $E$ 
with a cone. We do not do this directly because the cone is not 
exactly a competitor, and to avoid (probably minor)
problems of definitions or regularity. Instead we first prove an integrated 
version.

\medskip \proclaim Lemma 5.1. 
Let $U \i \R^n$ be open, and let $E$ be a reduced $A$- or $A'$-almost-minimal 
set in $U$, with gauge function $h$.
Let $x\in E$ and $0 < \rho_1 < \rho_2$ be such that $B(x,\rho_2) \i U$. Then
$$\eqalign{
H^d(E\cap &B(x,\rho_1)) \cr&\leq h(2\rho_2) (2\rho_2)^d 
+ {1 \over \rho_1 - \rho_2} \int_{E\cap B(x,\rho_2)\setminus B(x,\rho_1)}
{r(y) \cos \theta(y) \over d} \, dH^d(y),
}\leqno (5.2)
$$
where we set $r(y) = |y-x|$, and $\theta(y) \in [0,\pi/2]$ denotes the 
(smallest) angle between the radius $[x,y]$ and the tangent $d$-plane 
to $E$ at $y$.

\medskip
See Definition 4.3 or 4.8 for almost-minimal sets.
Also recall that $E$ is rectifiable, by Theorem 2.11 in [DS] (generalized 
as in Section 2). The rectifiability of $E$ gives an approximate 
tangent plane to $E$  at $y$ for $H^d$-almost every $y\in E$, which is 
enough to define $\theta(y)$. But here $E$ is locally 
Ahlfors-regular too, so the approximate tangent plane is even a true 
tangent $d$-plane (for instance by Exercise 41.21 on page 277 of [D2]). 

To prove the lemma, we may assume that $x=0$. Let $0 < r_1 < r_2$ be 
such that $\overline B(0,r_2) \i U$; we want to compare $E$ with
the competitor $F = \phi(E)$, where $\phi$ is the radial function defined
by $\phi(r \xi) = \varphi(r) \,\xi$ for $r\geq 0$ and $\xi \in \partial 
B(0,1)$, and where $\varphi$ is the piecewise linear function 
such that (as suggested by Figure 5.1)
$$\eqalign{
\varphi(r) = 0 \hbox{ for } 0 \leq &r \leq r_1,
\ \ \varphi(r) = {r_2 (r-r_1) \over r_2-r_1} 
\hbox{ for } r_1 \leq r \leq r_2,
\cr&\hbox{ and } \varphi(r)=r \hbox{ for }  r \geq r_2 \, .
}\leqno (5.3)
$$
Note that the one-parameter family $\phi_t$ defined by 
$\phi_t(y) = (1-t)y + t \phi(y)$ satisfies (2.2)-(2.6) with
any $\delta > 2 r_2$. Here $W_1 = B(0,r_2) \setminus \{ 0 \}$,
and (4.4) says that
$$
H^d(E \cap B(0,r_2)) \leq H^d(\phi(E \cap B(0,r_2))) + h(\delta) \delta^d.
\leqno (5.4)
$$
Notice that (4.9) yields exactly the same thing, after adding
$H^d(E \cap \phi(E) \cap B(0,r_2))$, so we don't need to worry about 
which definition of almost-minimality we take. We shall remember that 
(5.4) holds for all $\delta > 2 r_2$.

Next we need to evaluate $H^d(\phi(E \cap B(0,r_2))$. 
Set $A = B(0,r_2) \setminus B(0,r_1)$. 
We just need to evaluate $H^d(\phi(E \cap A))$, because
$E \cap B(0,r_1)$ is mapped to the origin.

We shall use the rectifiability of $E$.
Let $y \in E$ be given, and assume that $E$ has a tangent $d$-plane at $y$.
Denote by $P$ this $d$-plane and by $P'$ the vector $d$-plane parallel to $P$.
Also set $r(y)=|y|$ and $\widetilde y = y/r(y)$.

We want to compute how the differential $D\phi(y)$ acts on vectors. 
In the radial direction, we just use the fact that for $t > 0$ such 
that $t \widetilde y \in A$,
$\displaystyle\phi(t \widetilde y) = \varphi(t) \, \widetilde y 
= {r_2 (t-r_1) \over r_2-r_1} \, \widetilde y$, by (5.3), so
$D\phi(y)\cdot \widetilde y = [r_2 /(r_2-r_1)] \,\widetilde y$.
For vectors $w$ orthogonal to $\widetilde y$, (5.3) yields
$D\phi(y)\cdot w = [\varphi(r(y))/r(y)] w$.

\medskip 
\hskip 1.6cm 
\epsffile{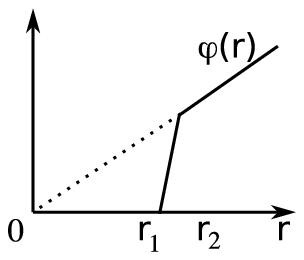}
\hskip 2.9cm
\epsffile{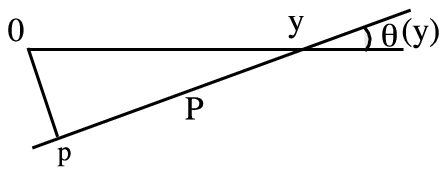}
\smallskip
\noindent\hskip 2.9cm{\bf Figure 5.1.}  
\hskip 4.1cm {\bf Figure 5.2.} 
\medskip 

If $P$ is orthogonal to $[0,y]$ at $y$, $D\phi(y)$ acts on $P'$
by multiplication by $\varphi(r(y))/r(y)$. Otherwise, denote by $p$ 
the orthogonal projection on $P$ of the origin, and set 
$v = {(y-p) \over |y-p|}\in P'$. [See Figure 5.2.] 
By definition, $\theta(y)$ is the angle of $[0,y]$ with $v$, 
so $v = \cos \theta(y) \,\widetilde y + \sin \theta(y) \, w$
for some unit vector $w$ orthogonal to $\widetilde y$. Then
$D\phi(y) \cdot v = [r_2 /(r_2-r_1)] \cos \theta(y) \,\widetilde y + 
[\varphi(r(y))/r(y)]  \sin \theta(y) \, w$ and 
$|D\phi(y) \cdot v|= \beta(y)$, where we set
$$
\beta(y) = \Big\{ [r_2 /(r_2-r_1)]^2\cos^2 \theta(y) 
+ [\varphi(r(y))/r(y)]^2  \sin^2 \theta(y) \Big\}^{1/2}.
\leqno (5.5)
$$
We are ready to apply the area formula to compute $H^d(\phi(E \cap A))$
(in other words, we cover $E \cap A$ with disjoint measurable subsets 
of $C^1$ surfaces, plus a set of vanishing $H^d$ measure, and then 
compute the Hausdorff measure of images on each $C^1$ piece with the 
derivative of $\phi$; we don't need to  worry about overlaps, since $\phi$ 
is bijective on $A$). We get that
$$
H^d(\phi(E \cap B(0,r_2))) =
H^d(\phi(E \cap A)) = \int_{E \cap A} 
[\varphi(r(y))/r(y)]^{d-1} \beta(y) dH^d(y),
\leqno (5.6)
$$
where the extra $[\varphi(r(y))/r(y)]^{d-1}$ comes from the directions
of $P'$ that are orthogonal to $v$. When $P$ is orthogonal to $[0,y]$ at 
$y$, (5.6) also gives the right factor $[\varphi(r(y))/r(y)]^{d}$, 
because $\theta(y)=\pi/2$. We plug this back into (5.4) and get that
$$
H^d(E \cap B(0,r_2)) \leq h(\delta) \delta^d + \int_{E \cap A} 
[\varphi(r(y))/r(y)]^{d-1} \beta(y) dH^d(y).
\leqno (5.7)
$$

It is tempting to let $r_2 - r_1$ tend to $0$ now, but we decided to
take an average before we do this, and this is what will lead to (5.2).
So we fix $0 < \rho_1 < \rho_2$ such that $B(0,\rho_2) \i U$, pick a 
small $\eta > 0$, and for each $t \in (\rho_1,\rho_2)$, we apply (5.7)
with $r_2 = t$, $r_1 = t-\eta$, and $\delta = 2\rho_2$. [Recall that we 
can take any $\delta > 2r_2$ in (5.4) and (5.7).]
Then we average over $t$ and get 
that
$$
H^d(E \cap B(0,\rho_1-\eta)) \leq h(2\rho_2) (2\rho_2)^d 
+ {1 \over\rho_2 - \rho_1} \, I(\eta)
\leqno (5.8)
$$
where
$$
I(\eta) = \int_{t\in (\rho_1,\rho_2)} 
\int_{y\in E \, ; \, t-\eta \leq r(y) \leq t} 
[\varphi_{t}(r(y))/r(y)]^{d-1} \beta_{t}(y) dH^d(y)dt,
\leqno (5.9)
$$
with obvious notations concerning $\varphi_t$ and $\beta_t$. By Fubini,
$$
I(\eta) =  
\int_{y\in E \, ; \,\rho_1-\eta < r(y) < \rho_2} g_\eta(y) \, dH^d(y),
\leqno (5.10)
$$
with 
$$
g_\eta(y) = \int_{t\in (\rho_1,\rho_2)\cap[r(y),r(y)+\eta]}
[\varphi_{t}(r(y))/r(y)]^{d-1} \beta_{t}(y) dt.
\leqno (5.11)
$$
First observe that here $r_2 /(r_2-r_1) = t/\eta$. Then (5.5) yields
$$\eqalign{
\beta_{t}(y) &= \big\{ (t/\eta)^2\cos^2 \theta(y) 
+ [\varphi_{t}(r(y))/r(y)]^2  \sin^2 \theta(y) \big\}^{1/2}
\cr&
\leq \big\{ (t/\eta)^2\cos^2 \theta(y) + 1 \big\}^{1/2}
\leq (t/\eta)\cos\theta(y) + 1 
}\leqno (5.12)
$$
because $\varphi_t(r(y)) \leq r(y)$ and $a^2+b^2 \leq (a+b)^2$. 
This yields
$$
\beta_{t}(y) \leq \eta^{-1} t \cos\theta(y) + 1
\leq \eta^{-1} [(r(y)+\eta) \cos\theta(y) + \eta]
\leqno (5.13)
$$
on the domain of integration in (5.11). Also
$$
{\varphi_{t}(r(y)) \over r(y)} 
= {r_2 \over r(y)} \, {r(y)-r_1 \over r_2 - r_1} 
= {t \over r(y)}\, {r(y)-t+\eta \over \eta}
\leq \eta^{-1}\, {r(y)+\eta \over r(y)}\, 
\big(r(y)-t+\eta\big),
\leqno (5.14)
$$
by (5.3) and because $t \leq r(y)+\eta$ on the domain 
of integration, so
$$\eqalign{
g_\eta(y) &
\leq \eta^{-d}\, [(r(y)+\eta) \cos\theta(y) + \eta] \,
\Big({r(y)+\eta \over r(y)}\Big)^{d-1}
\int_{t\in[r(y),r(y)+\eta]} (r(y)-t+\eta)^{d-1} dt
\cr&
= {1 \over d} \, [(r(y)+\eta) \cos\theta(y) + \eta] \,
\Big({r(y)+\eta \over r(y)}\Big)^{d-1}
\leq {1 \over d} \, r(y) \cos\theta(y) + C \eta,
}
$$
where $C$ is allowed to depend on many things, but not on $\eta$. 
We return to (5.10) and get that $\displaystyle
I(\eta) \leq  {1 \over d} \int_{y\in E \, ; \,\rho_1-\eta < r(y) < \rho_2}
[r(y) \cos\theta(y) + C \eta] \, dH^d(y)$. Then (5.8) becomes
$$\eqalign{
H^d(E \cap &B(0,\rho_1-\eta)) \leq h(2\rho_2) (2\rho_2)^d 
\cr&
+ {1 \over\rho_2 - \rho_1} \, 
{1 \over d} \int_{y\in E \, ; \,\rho_1-\eta < r(y) < \rho_2}
[r(y) \cos\theta(y) + C \eta] \, dH^d(y).
}\leqno (5.15)
$$
Then we let $\eta$ tend to $0$ and get (5.2). This completes our 
proof of Lemma 5.1.
\qed

\medskip\noindent{\bf Remark.} The deformation $\phi$ that we used in 
the proof of Lemma 5.1 is not a bijection, but we could have used 
the homeomorphism $\phi_t(y) = (1-t)y + t \phi(y)$, with $t<1$ very 
close to $1$, and then go to the limit and get (5.2) anyway. So
Lemma 5.1, and Propositions~5.16 and 5.30 below, also hold when 
$E$ lies in the larger class of almost-minimal sets where the 
deformations in (2.2)-(2.6) are also required to be homeomorphisms, 
at least if we assume that $E$ is rectifiable.
I do not think that Proposition 5.24 below extends to this 
context (we shall also use local Ahlfors-regularity there).

\medskip
We want to use Lemma 5.1 to get monotonicity results for 
$r^{-d}H^d(E \cap B(x,r))$. We start with the easier case of minimal 
sets.

\medskip
\proclaim Proposition 5.16.
Let $E$ be a minimal set in $U$. This means that 
$E$ is almost-minimal in $U$, with gauge function $h=0$
(and Definitions 4.1, 4.3, and 4.8 are all equivalent in this case).
Let $x\in U$ and $R>0$ be such that $B(x,R) \i U$. Then
$\theta(x,r) = r^{-d}H^d(E \cap B(x,r))$ is a nondecreasing
function of $r$ on $(0,R]$.

\medskip
Let $B(x,R)$ be as in the statement, and define a locally finite Borel
measure $\mu$ on $[0,R)$ by
$$
\mu(A) = H^d(E \cap \pi^{-1}(A))
\ \hbox{ for every Borel set } A \i [0,R),
\leqno (5.17)
$$
where $\pi$ is the radial projection defined by $\pi(y)=|x-y|$.
Then decompose $\mu$ into its absolutely continuous part $\mu_a$ and 
its singular part $\mu_s$, and denote by $a(r)$ the density of 
$\mu_a$ with respect to the Lebesgue measure. The Lebesgue density
theorem allows us to compute
$a(r)$ at almost-every point $r\in [0,R)$ as a limit of densities, i.e.,
$$
a(r) = \lim_{t \to 0^+} {1 \over t} \, \mu([r-t,r))
= \lim_{t \to 0^+} {1 \over t} \, 
H^d\big(E \cap B(x,r) \setminus B(x,r-t)\big).
\leqno (5.18)
$$
Similarly $a(r) = \lim_{t \to 0^+} {1 \over t} \, 
H^d(E \cap B(x,r+t) \setminus B(x,r))$ for
almost-every $r$. Consequently, if we set
$$
l(r) = H^d(E \cap B(x,r))
\hbox{ for } 0 \leq r < R,
\leqno (5.19)
$$
then $l(r)$ is differentiable almost-everywhere on $(0,R)$, and
$l'(r) = a(r)$ almost-everywhere.
Let us only assume for the moment that $E$ is almost-minimal,
and apply Lemma 5.1 with $\rho_2 = r\in (0,R)$ and $\rho_1 = r-t$. 
We can even drop $\,\cos\theta(y) \,$ and get that 
$$
l(r-t) \leq (2r)^d h(2r) 
+ {r \over td} \, H^d(E \cap B(x,r) \setminus B(x,r-t)).
\leqno (5.20)
$$
Then we let $t$ tend to $0$, use (5.18), and get that
$$
l(r) \leq (2r)^d h(2r) + {r \over d} \, a(r) 
= (2r)^d h(2r) + {r \over d} \, l'(r) 
\ \hbox{ for almost-every } r \in(0,R).
\leqno (5.21)
$$
Here $E$ is even a minimal set, so (5.21) says that 
$l(r) \leq {r \over d} \, l'(r)$.
Recall that $l(r) = \int_{[0,r)} d\mu$; we claim that for $0 < a < b <R$,
$$
\theta(x,b) - \theta(x,a)
= \int_{[a,b)} t^{-d} d\mu(t) - d \int_{[a,b)}  t^{-d-1}l(t) dt.
\leqno (5.22)
$$
Set $g(t)=t^{-d}$; thus $\theta(x,r) = g(r) l(r)$, and we just want 
to show that the derivative of the product 
$gl$ is $g d\mu + g'l$. Let $0 < a < b <R$ be given, and apply Fubini's 
theorem to $I=\int_{ a \leq s \leq t < b} \, g'(s) ds d\mu(t)$. 
One computation gives
$I=\int_{[a,b)} \, [g(t)-g(a)]d\mu(t)$, and the other one yields 
$I=\int_{[a,b)} g'(s)[l(b)-l(s)] ds$. Then
$$\leqalignno{
\int_{[a,b)} g(t) d\mu(t) + \int_{[a,b)} &g'(s)l(s)ds
= \big[I + g(a) \int_{[a,b)} d\mu(t)\big] 
+ \big[-I + l(b) \int_{[a,b)} g'(s) ds \big]
\cr&= g(a) (l(b)-l(a)) + l(b)(g(b)-g(a)) = g(b)l(b)-g(a)l(a).
&(5.23)
}$$
This is just (5.22); note that (5.23) would work with any smooth 
function $g$.

We may now conclude. Recall that $t^{-d}d\mu(t) \geq t^{-d} l'(t) dt 
\geq d \, t^{-d-1} l(t) dt$ because $d\mu(t) \geq d\mu_a(t) = l'(t) dt$ 
and by (5.21) with $h(2r)=0$, so (5.22) says that $\theta(x,\cdot)$ 
is nondecreasing, and Proposition 5.16 follows. \qed

\medskip
\proclaim Proposition 5.24.
There is a constant $\lambda > 1$, that depends only on $n$ and $d$, 
such that the following holds.
Suppose that $E$ is an $A$- or $A'$-almost-minimal set in $U$, with gauge 
function $h$. Let $x\in E^\ast\cap U$ and $R>0$ be such that $B(x,R) \i U$,
and also assume that $h(R)$ is small enough (depending only on $n$ 
and $d$) and that $\int_0^R h(2r) {dr \over r} < +\infty$. Then
$$
\theta(x,r) \, e^{\lambda A(r)}
\hbox{ is a nondecreasing function of } r\in (0,R),
\leqno (5.25)
$$
where we set 
$\displaystyle A(r) = \int_0^r h(2t) {dt \over r}$ for $0 < r < R$.

\medskip
Notice that here we had to take $x\in E^\ast\cap U$ to get the 
conclusion (Recall from Definition~2.12 that $E^\ast$ is the closed
support of $H^d_{|E}$). 
Since $\lim_{r \to 0} A(r) = 0$, $e^{\lambda A(r)}$
tends to $1$ and can be seen as a mild penalty that measures a 
potential defect in monotonicity.

We start the proof as in Proposition 5.16, but we multiply $l(t)$ with
$g(t) = t^{-d}e^{\lambda A(t)}$. By (5.23), the distribution derivative
of $g(t) l(t)$ is $g(t) d\mu(t) + g'(t) l(t) dt$. Also,
$d\mu(t) \geq d\mu_a(t) = a(t) dt = l'(t) dt$, so
$$
\theta(x,b) e^{\lambda A(b)} - \theta(x,a) e^{\lambda A(a)}
\geq \int_a^b [g(t) l'(t) + g'(t) l(t)] dt,
\leqno (5.26)
$$
as in (5.22). So it is enough to check that 
$g(t) l'(t) + g'(t) l(t) \geq 0$ almost-everywhere.
But $g'(t) = - {d \over t} \, g(t) + \lambda A'(t) g(t)
= - {d \over t} \, g(t) + \lambda {h(2t) \over t} g(t)$, so
$$
g(t) l'(t) + g'(t) l(t) = g(t) 
\big[l'(t) - {d \over t} \, l(t) + \lambda {h(2t) \over t} l(t)\big].
\leqno (5.27)
$$
Now (5.21) says that $l'(t) \geq {d \over t} \, [l(t) - (2t)^d h(2t)]$,
so it is enough to check that $(2t)^d \leq \lambda l(t)$.
This comes from Lemma 2.15; we just have to take $\lambda \geq 2^d C$, 
where $C$ is as in (2.16). This is the reason why we have to take 
$x\in E^\ast$ and assume that $h(R)$ is small enough.
\qed

\medskip
The case of $A_+$-almost-minimal sets is a little easier, because we
may replace the error term $h(\delta) \delta^d$ in (5.4) with 
$h(\delta) H^d(\phi(E \cap B(0,r_2))$. [We use (4.2) instead of (4.4).] 
When we follow the computations, we see that we can replace (5.2) with 
$$
H^d(E\cap B(x,\rho_1)) 
\leq {1 + h(2\rho_2) \over \rho_1 - \rho_2} 
\int_{E\cap B(x,\rho_2)\setminus B(x,\rho_1)}
{r(y) \cos \theta(y) \over d} \, dH^d(y).
\leqno (5.28)
$$
Then (5.21) may be replaced with
$$
l(r) \leq  [1+h(2r)] \, {r \over d} \, l'(r) 
\ \hbox{ for almost-every } r \in(0,R),
\leqno (5.29)
$$
and (5.27) allows us to show that $\theta(x,r) e^{\lambda A(r)}$
is  nondecreasing as soon as
$[1+h(2t)]^{-1} \, {d \over t} - {d \over t} + \lambda {h(2t) \over t}
\geq 0$, i.e., when $\lambda \geq d \, \{ 1-[1+h(2t)]^{-1} \}$.
This proves the following intermediate statement for 
$A_+$-almost-minimal sets, where the interest is that we no longer 
need to take a center in $E^\ast$.

\medskip
\proclaim Proposition 5.30. For each $\lambda > 0$ we can find
$\tau >0$, that depends only on $n$, $d$, and $\lambda$, with the 
following property.
Let $E$ be an $A_+$-almost-minimal set in $U$, with gauge function $h$. 
Let $x\in U$ and $R>0$ be such that $B(x,R) \i U$, and suppose that 
$h(2R)\leq \tau$ and $\int_0^R h(2r) {dr \over r} < +\infty$.
Set $A(r) = \int_0^r h(2r) {dt \over t}$ for $0 < r < R$. Then
$$
\theta(x,r) e^{\lambda A(r)}
\hbox{ is a nondecreasing function of } r\in (0,R).
\leqno (5.31)
$$

\medskip\noindent 
{\bf Remark 5.32.} It is tempting to try to prove Propositions 5.16, 
5.24, and 5.30 without using the rectifiability of $E$. It would be
enough to prove that
$$
H^d(E\cap B(x,\rho_1)) \leq h(2\rho_2) (2\rho_2)^d 
+ {1 \over \rho_1 - \rho_2} \int_{E\cap B(x,\rho_2)\setminus B(x,\rho_1)}
{r(y) \over d} \, dH^d(y)
\leqno (5.33)
$$
(or the corresponding variant in the case of Proposition 5.30)
instead of (5.2) (i.e., without the cosine), and then continue as 
before. To do so, we would need to change the argument near (5.5), 
when we compute the effect of $D\phi$ on the tangent plane. 
We observed that $D\phi(y).v = [r_2 /(r_2-r_1)] v$ in the radial direction 
and $D\phi(y).v = [\varphi(r(y))/r(y)] v$ in the other direction, so we 
expect that $\phi$ will at most multiply the Hausdorff measure of $E$
by $[r_2 /(r_2-r_1)][\varphi(r(y))/r(y)]^{d-1}$ at $y$, with the effect of 
replacing $\beta(y)$ in (5.5) with $[r_2 /(r_2-r_1)]$ and eventually 
leading to (5.33).
But the author 
is not aware of a generalization to  unrectifiable sets of the obvious estimate on 
Hausdorff measure of affine images (of rectifiable sets) that would allow 
to conclude here.
We could also try to replace the Hausdorff measure with a more 
friendly variant. Anyway, this may not be so useful, because the 
cosine that shows up in Lemma 5.1 will really be used in the next section.

\bigskip 
\noindent {\bf 6. Minimal sets with constant density}
\medskip

In addition to the monotonicity of density proved in the last section, 
it will be useful to know that when $E$ is a minimal set and 
$$
\theta(r)=r^{-2}H^d(B(0,r))
\leqno (6.1)
$$ 
is constant on some interval $(a,b)$, then $E^\ast$ coincides on 
$B(0,b)\setminus \overline B(0,a)$ with a minimal cone 
(i.e., a minimal set which is a cone) centered at the origin. 

\medskip
\proclaim Theorem 6.2.
Let $E$ be a reduced minimal set of dimension $d$ in $B(0,b) \i \R^n$.
This means that $E$ is a reduced almost-minimal in $B(0,b)$, as in
Definition 4.3 or 4.8 (also see Definition 2.12), with the gauge 
function $h=0$. Assume that we can find $a \in [0,b)$ and a constant 
$\theta_0 \geq 0$ such that 
$$
\theta(r) = \theta_0 \hbox{ for } a < r < b.
\leqno (6.3)
$$
Then there is a reduced minimal cone ${\cal C}$ centered at the origin
such that 
$$
E \cap [B(0,b)\setminus \overline B(0,a)] 
= {\cal C}\cap [B(0,b)\setminus \overline B(0,a)].
\leqno (6.4)
$$

\medskip
Recall that our two  definitions of minimality are equivalent.
We shall not need to know this, because it will be just as simple to
prove the theorem for both notions. Of course Theorem 6.2 also
applies if $E$ is minimal in a domain that contains $B(0,b)$.

Observe also that if we did not assume that $E$ is reduced, we
would still get that $E^\ast \cap [B(0,b)\setminus \overline B(0,a)] 
= {\cal C}\cap [B(0,b)\setminus \overline B(0,a)]$
(instead of (6.4)); this is easy, because $E^\ast$ is a reduced
minimal set, with the same density, by Remark 2.14 and (2.13), so we 
can apply Theorem 6.2 to it.

Let $E$ be as in the theorem. We first look at Proposition 5.16 to see 
if we can derive extra information because (6.3) holds. Set 
$A = B(0,b) \setminus \overline B(0,a) = \pi^{-1}((a,b))$, where
$\pi$ denotes the radial projection defined by $\pi(y) = |y|$.
We want to prove that
$$
\cos\theta(y) = 1 
\hbox{ for $H^d$-almost every } y\in E \cap A,
\leqno (6.5)
$$
where $\cos\theta(y)$ is as in Lemma 5.1.

Let $\nu$ be the finite Borel measure on $(a,b)$ defined by
$$
\nu(S) = \int_{E \cap \pi^{-1}(S)} \cos\theta(y) dH^2(y)
\leqno (6.6)
$$
for $S \i (a,b)$.
Notice that $\nu \leq \mu$, where $\mu$ is as in (5.17) and 
corresponds to replacing $\cos\theta(y)$ with $1$. By (6.3),
$\mu([\alpha,\beta)) = \theta_0 (\beta-\alpha)$ for 
$a<\alpha<\beta<b$, so $d\mu = \theta_0 dx$ on $(a,b)$. Then
$\nu$ is absolutely continuous too. Let $f$ denote its density (with 
respect to the Lebesgue measure). Notice that 
$f(r) = \lim_{t \to 0^+} {1 \over t} \nu([r-t,r))$ almost-everywhere, 
as in (5.18). We apply Lemma 5.1 with $\rho_2=r \in (a,b)$ and 
$\rho_1 = r-t$, as we did for (5.20), and we get that 
$$
l(r-t) \leq {r \over dt} \, H^d(E \cap B(x,r) \setminus B(x,r-t))
\leqno (6.7)
$$
instead of (5.20) (just because $h=0$ here). Then we
let $t$ tend to $0$ and get that $l(r) \leq {r \over d} \, f(r)$ 
(instead of (5.21)). We also have that $f(r) \leq l'(r)$, 
because $l'(r)$ is the density of $\mu$ and $\nu \leq \mu$. 
In addition, $l(r) = r^d\theta(r) = r^d\theta_0$ by (6.1) and (6.3), 
so $l'(r)=dr^{d-1}\theta_0 = dl(r)/r$. Altogether 
$$
l(r) \leq {r \over d} \, f(r) \leq {r \over d} \, l'(r) = l(r)
\ \hbox{ almost-everywhere on $(a,b)$},
\leqno (6.8)
$$
and in particular $f(r)=l'(r)$ almost-everywhere on $(a,b)$. Recall
that $f$ is the density of the absolutely continuous measure $\nu$, so
$$
\nu((a,b)) = \int_a^b f(r)dr = \int_a^b l'(r)dr = \mu((a,b)),
\leqno (6.9)
$$
and (6.5) follows by comparison.

In geometrical terms, (6.5) just says that for $H^d$-almost-every
point $y\in E \cap A$, the tangent plane to $E$ at $y$ 
goes through the origin. This looks like pretty strong evidence that
$E$ coincides with a cone in $A$, but apparently we need some nontrivial 
amount of additional work to get the result. In particular, we shall
use the minimality of $E$ again.

\medskip
Denote by $\pi_S$ the radial projection onto the 
unit sphere. Thus $\pi_S(y) = y /|y|$. We claim that 
$$
H^d(\pi_S(E\cap A)) = 0.
\leqno (6.10)
$$
This is an easy consequence of (6.5). Indeed recall from Section 2 
that $E$ is both rectifiable and locally Ahlfors-regular, 
so by Exercise 41.21 on page 277 of [D2]) 
it has a true tangent $d$-plane at almost every point.
[We don't really need to know this here, but this will
make the discussion simpler.]

By rectifiability, we can cover $E\cap A$, up to a negligible set $N$ 
such that $H^d(N)=0$, by a countable collection of $d$-dimensional 
surfaces $\Gamma_i$ of class $C^1$. 
For almost-every point $y$ of $E \cap A \cap \Gamma_i$, the 
tangent $d$-plane to $E$ at $y$ is the same as the tangent $d$-plane 
to $\Gamma_i$ at $y$, because we may assume that $y$ is a  point of 
positive upper density of $E \cap \Gamma_i$), 
and then (6.5) says that this $d$-plane goes through the origin. 
So $y$ is a critical point of $\pi_S$ on $\Gamma_i$, and 
Sard's theorem (applied to $\pi_S$ on $\Gamma_i$) says that 
$H^d(\pi_S(E \cap A \cap \Gamma_i))=0$. Since $H^d(\pi_S(N))=0$ 
because $\pi_S$ is Lipschitz on $A$, we just need to sum over $i$ 
to get (6.10).

\medskip
\proclaim Proposition 6.11. Let $y\in E \cap A$ be such that $E$
has a tangent $d$-plane $P$ at $y$ which goes through $0$.
Then $E$ contains the radial line segment $L= A \cap \pi_S^{-1}(\pi_S(y))$
through $y$.

\medskip
Let us first check that Proposition 6.11 implies that
$$
E \hbox{ coincides on $A$ with a cone;}
\leqno (6.12)
$$
we shall only see later that this cone is a reduced minimal set. 

Call $E_1$ the set of points $y\in E \cap A$ that satisfy the 
conditions of Proposition 6.11; by the discussion above
(the proof of (6.10)), almost-every point of $E\cap A$ lies in $E_1$. 
In particular, $E_1$ is dense in $E\cap A$. Call $\cal C$ the cone over $E_1$; 
then $E \cap A \i \overline E_1 \i \overline {\cal C}$.
Conversely, $\overline {\cal C} \cap A \i E$, by Proposition~6.11 and 
because $E$ is closed. 
So (6.12) follows from Proposition~6.11.

\medskip
The proof of Proposition 6.11 will take some time. Let $y\in E^\ast \cap A$,
$P$, and $L$ be as in the proposition, and suppose in addition that $L$ is 
not contained in $E$. We want to construct a deformation of 
$E$ and eventually get a contradiction with its minimality. 

We need lots of \underbar{notation} before we start.
Since $L$ is not contained in $E$, we can pick $x\in L \setminus E$.
Let $\varepsilon_0$ be small, to  be chosen near the end. 
Then pick a small radius $r_y$ such that
$$
\dist(z,P) \leq \varepsilon_0 |z-y| 
\ \hbox{ for } z\in E \cap 2B_y,
\leqno (6.13) 
$$
where we set $B_y=B(y,r_y)$ and $2B_y = B(y,2r_y)$. Notice that (6.13) 
holds for $r_y$ small enough, because $P$ is tangent to $E$ at $y$.
We choose $r_y < b/100$ and so small that $2B_y$ is contained in $A$
and does not contain $x$.

Now let $r_x$ be a small radius, smaller than $r_y$ and
such that if $B_x = B(x,r_x)$, then 
$$
2B_x \i A \setminus [E \cup 2B_y].
\leqno (6.14)
$$
We shall also need to take $r_x$ small enough, depending on 
$r_y$, for the final estimate to work.

Without loss of generality, we may assume that coordinates were chosen 
so that $L$ is contained in the positive first axis, and $P$ is the
horizontal $d$-plane where $x_{d+1}= \cdots = x_n = 0$. We know 
that $L \i P$, by (6.5), so this  is coherent. We shall identify 
$\R^n$ with $L' \times P' \times Q$, where $L'$ is the line that 
contains $L$ and the first axis, $P'$ is the orthogonal complement of 
$L'$ in $P$ (so $P'$ corresponds to $x_1 = x_{d+1} = \cdots = x_n = 0$), and
$Q$ is the orthogonal of $P$ in $\R^n$. Thus $Q$ is $(n-d)$-dimensional.

We shall often use this identification to denote points of $\R^n$ by
$z =(z_1,z_2,z_3)$, with $z_1\in L'$, $z_2 \in P'$, and $z_3 \in Q$.
Also, we identify $L'$ with $\R$, so there is an order on $L'$, and we 
assume that $0 < x < y$ in $L'$; the other case would be treaded similarly.

So we want to construct a deformation $f$ and use it to reach a 
contradiction. We start with a first deformation $f_1$ that tries to 
send to $P$ points that lie near $L$ and roughly between $x$ and $y$.
The effect of this map will be to free from $E \setminus P$ 
a small tunnel near $L$. This map $f_1$ will be composed with
other maps $f_j$, $j\geq 2$, that will be defined later; the main goal
is to make most of $E\cap B_y$ disappear, by retracting it on a 
boundary; the tunnel will be used as a clean path that leads from
$x$ to $y$.

\medskip
\proclaim Lemma 6.15.
We can find $\xi \in P' \cap B(x,r_x/100)$ and $\rho \in (r_x/2,2r_x/3)$ 
such that $E \cap A$ does not meet the cone over 
$S_{\xi,\rho} = x + \xi + [Q \cap \partial B(0,\rho)]$.

\medskip
Thus $S_{\xi,\rho}$ is a $(n-d-1)$-dimensional sphere centered
at $x+\xi$. When $d=n-1$, $S_{\xi,\rho}$ is just composed of two points
that lie on different sides of $P$. See Figure 6.1 for two special 
cases when $n=3$.

To prove Lemma 6.15, pick any small constant $\tau > 0$ and set
$Z_\tau = \big\{ (z_1,z_2,z_3) \in \pi_S(A\cap E) \, ; \, 
z_1 \geq \tau \big\}$, where $\pi_S$ still denotes the radial 
projection on the unit sphere. We know from (6.10) that $H^d(Z_\tau) = 0$.
Then define $\psi_\tau : Z_\tau \to P' \times \R$, by 
$\psi_\tau(z_1,z_2,z_3) = (xz_2 /z_1, x|z_3|/z_1)$. Notice that
$\psi_\tau$ is Lipschitz (with a bound that depends on 
$\tau$, but this does not matter) so $H^d(\psi_\tau(Z_\tau))=0$.

Recall that $P' \times \R$ is $d$-dimensional, so we can pick 
$\xi \in P' \cap B(x,r_x/100)$ and $\rho \in (r_x/2, 2r_x/3)$
such that $(\xi,\rho)$ lies in no $\psi_\tau(Z_\tau)$.
Let us check that the cone over $S_{\xi,\rho}$ does not meet
$E \cap A$. Otherwise, we can find $z=(z_1,z_2,z_3) \in E\cap A$
and $\lambda \geq 0$ such that $\lambda z\in S_{\xi,\rho}$.
This last forces $\lambda z_1 = x > 0$ (we identify $x$ with its first 
component), so $z_1 >0$ and $\pi_S(z)$ lies in $Z_\tau$
for $\tau$ small enough. Next 
$\psi_\tau(\pi_S(z)) = \psi_\tau(z/|z|)
= (xz_2 /z_1, x|z_3|/z_1) = (\lambda z_2,\lambda |z_3|)$
because $\lambda z_1 = x$. But $\lambda z\in S_{\xi,\rho}$,
so $\lambda z_2 = \xi$ and $\lambda |z_3|=\rho$, hence 
$\psi_\tau(\pi_S(z)) = (\xi,\rho)$, which is impossible by 
choice of $\xi$ and $\rho$. So the cone over $S_{\xi,\rho}$ does not meet
$E \cap A$; this proves Lemma 6.15. \qed

\medskip 
\hskip -0.3cm 
\epsffile{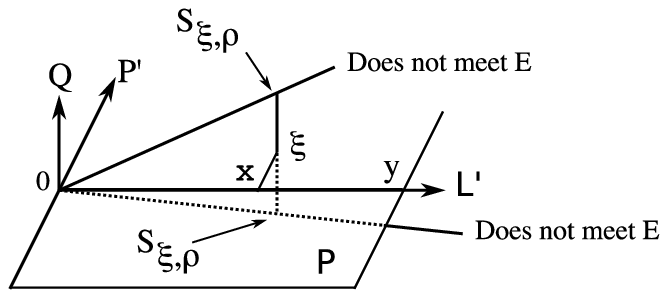}
 \hskip 1.4cm
 \epsffile{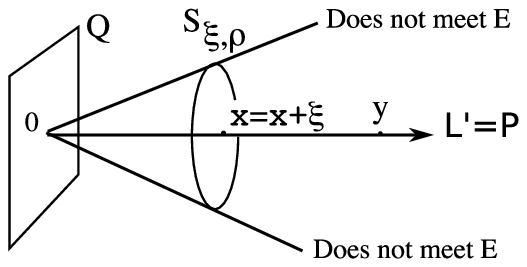}
\smallskip
\hskip 0.9cm{\bf Figure 6.1.a.} $\, n=3$, $d=2$  
\hskip 2.3cm {\bf Figure 6.1.b.} $\, n=3$, $d=1$  
\medskip 

\medskip
Choose a pair $(\xi,\rho)$ as in Lemma 6.15.
For $a \leq s \leq t \leq b$, set
$$
V(s,t) = \big\{ z \in \R^n \, ; \, s \leq z_1 \leq t \big\}.
\leqno (6.16)
$$
With this notation, $S_{\xi,\rho}$ is contained in the vertical
hyperplane $V(x,x)$. Pick $a_1$ and $b_1$ on the first axis, with 
$a < a_1 < b_1 < b$, so that 
$$
\overline B_x \cup \overline B_y \i V(a_1,b_1).
\leqno (6.17)
$$
Set $S_{\xi,\rho}^\varepsilon = \big\{ z\in V(x,x) \, ; \,
\dist(z,S_{\xi,\rho} \leq 2\varepsilon \big\}$ for $\varepsilon>0$.
Also call $\C_\varepsilon$ the cone over $S_{\xi,\rho}^\varepsilon$.
If $r_x$ and $\varepsilon$ are small enough,  
$\C_\varepsilon \cap V(a_1,b_1)$ is a compact subset of $A$.
Then, by compactness and because the cone over $S_{\xi,\rho}$
does not meet $E \cap A$, we can choose $\varepsilon$  so small
that
$$
\C_\varepsilon \cap V(a_1,b_1) \ \hbox{ does not meet } E.
\leqno (6.18)
$$
Of course $\varepsilon$ may be very small, so we shall need
to be careful and avoid letting our final estimates depend on 
$\varepsilon$. Also choose $\varepsilon$ smaller than $r_x/100$ to avoid 
possible geometric complications.

\underbar{Our first map $f_1$} will act inside vertical hyperplanes $V(s,s)$,
and we first want to describe how it works in $V(x,x)$.
The reader may refer to Figure 6.2, which is a picture in $V(x,x)$. 
The dimensions of the picture correspond to
$n=3$ and $d=2$; in general, $V(x,x)$ is $n-1$ dimensional, $P'$ is 
$(d-1)$-dimensional,
and $Q$ is $(n-d)$-dimensional, but the picture is simple in $P'$,
and things will be invariant under rotations in $x+Q$, so Figure~6.2 
still gives a fair idea of what happens.

Set $J=P' \cap \overline B(x+\xi,r_x/10)$ 
and $S'= x + \xi + Q\cap \partial B(0,\rho+\varepsilon)$ (a sphere 
concentric with $S_{\xi,\rho}$, but just a bit larger),
and then denote by $T$ the convex hull of $J \cup S'$. 
Thus 
$$
T = \big\{ (x,z_2,z_3) \, ; \,  |z_2-\xi| \leq 
(\rho + \varepsilon -|z_3|) \, r_x/(10\rho+10\varepsilon)\big\} .
\leqno (6.19)
$$

Set $\widetilde f(z)=z$ in $V(x,x)\setminus T$. 
For $z\in T \setminus \C_\varepsilon$,
write $z=(x,z_2,z_3)$, call $z^\ast = (x,\xi,\rho z_3/|z_3|)$
the closest point of $S_{\xi,\rho}\,$, and then let $\widetilde f(z)$
be the radial projection of $z$ onto $P' \cup \partial T$, with
origin $z^\ast$. That is, $\widetilde f(z)$ lies on 
$P' \cup \partial T$, $z$ lies on the segment $[z^\ast,\widetilde f(z)]$,
and the open segment does not meet $P' \cup \partial T$. 
Obviously $\widetilde f$ is continuous across $\partial T$.

The reader may be worried about the Lipschitz-continuity
of $\widetilde f$ inside $T \setminus \C_\varepsilon$, so let
us say a few words about that. Let us first check that 
$$
|z_3| \leq \rho-\varepsilon
\ \hbox{ when } \  (x,z_2,z_3) \in T\setminus \C_\varepsilon,
\leqno (6.20)
$$
which is already almost obvious from Figure 6.2. 
Otherwise, $|z_2-\xi| 
\leq (\rho + \varepsilon -|z_3|)r_x/(10\rho+10\varepsilon)
\leq 2 \varepsilon r_x/10\rho \leq 4 \varepsilon/10$, by (6.19)
and because $\rho \geq r_x/2$. Recall that 
$z^\ast = (x,\xi,\rho z_3/|z_3|)$, so 
$|z-z^\ast|^2 = |z_2-\xi|^2 + (\rho-|z_3|)^2
\leq 16\varepsilon^2/100 + (\rho-|z_3|)^2
\leq 16\varepsilon^2/100 + \varepsilon^2 < 4\varepsilon^2$.
This is impossible because $z \notin \C_\varepsilon$; (6.20)
follows.

\medskip 
\hskip 3.1cm 
\epsffile{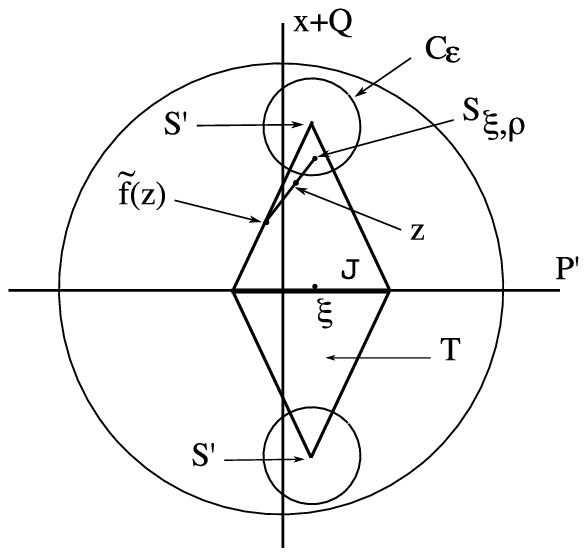}
\smallskip
\centerline{{\bf Figure 6.2.} A picture in $V(x,x)$} 
\medskip 

We may now compute the radial projection of $z$ onto $P'$, 
again with origin $z^\ast$. This point should be of the form
$z'=z^\ast + \lambda (z-z^\ast)$, with $\lambda \geq 1$. 
Since the last coordinate $z^\ast_3 + \lambda (z_3-z^\ast_3)$
vanishes and $z^\ast_3 = \rho z_3/|z_3|$, we get that 
$\lambda = \rho/(\rho-|z_3|)$, which is indeed defined by (6.20)
and larger than $1$. In addition, $\lambda$ and $z'$ are Lipschitz functions
of $z \in T \setminus \C_\varepsilon$. Now the radial projection 
of $z$ onto $\partial T$, again with origin $z^\ast$, is also 
Lipschitz on $T \setminus \C_\varepsilon$. And so is $\widetilde f(z)$,
because the infimum of two Lipschitz functions $\lambda$ is Lipschitz.
Since $\widetilde f$ is continuous across $\partial T$, $\widetilde f$
is Lipschitz on $V(x,x) \setminus \C_\varepsilon$.

Extend $\widetilde f$ to the whole $V(x,x)$, so that it is Lipschitz
and $\widetilde f(T) \i T$. [If this last condition is not satisfied, 
we can always compose the restriction of $\widetilde f$ to $T$ with a
Lipschitz retraction onto $T$.] Finally extend $\widetilde f$ to 
$V(a_1,b_1)$, by homogeneity: set
$$
\widetilde f(z_1,z_2,z_3) = {z_1 \over x} 
\, \widetilde f(x,z_2 x/z_1,z_3 x/z_1).
\leqno (6.21)
$$
Notice that $\widetilde f$ moves the points inside vertical
planes $V(s,s)$.

We are now ready to define $f_1$. Set 
$$
x_0 = x-{r_x \over 5} \ \hbox{ and }\  y_0 = y + {r_y \over 5} \, , 
\leqno (6.22)
$$
and then let $\psi$ be the continuous cut-off function such that 
$$\eqalign{
\psi(t) = 0 \hbox{ out of } &(x_0,y_0), \ 
\psi(t)=1 \hbox{ on } [x,y], \ 
\hbox{ and }
\cr & \psi \hbox{ is affine on the two remaining intervals.} 
}\leqno (6.23)
$$
Finally set
$$
f_1(z) = \psi(z_1) \widetilde f(z) + (1-\psi(z_1)) \, z
\ \hbox{ for } z \in \R^n,
\leqno (6.24)
$$
where $z_1$ still denotes the first coordinate of $z$.
Notice that $z \in V(a_1,b_1)$ when $\psi(z_1) \neq 0$,
by (6.17).

\medskip
Let us record a few properties of $f_1$ before we compose with a 
second deformation. 
As promised, $f_1$ only moves the points inside vertical planes 
$V(s,s)$.
Notice also that
$$
f_1(z) = z \ \hbox{ out of }\ \widehat T \cap V(x_0,y_0),
\leqno (6.25)
$$
where $\widehat T$ denotes the cone over $T$, and that 
$$
f_1\big(\widehat T\cap V(x_0,y_0)\big) \i \widehat T \cap V(x_0,y_0).
\leqno (6.26)
$$
We also want to see where the set $E_1 = f_1(E)$ lies.
Let us check that
$$
E_1 \cap \widehat T \cap V(x,y)
\i \big[\partial \widehat T \cup\widehat{J}\big]\cap V(x,y),
\leqno (6.27)
$$
where $\widehat{J}$ is the cone over $J$
(also notice that $\partial \widehat T$ is the cone 
over the boundary of $T$ in $V(x,x)$). Indeed, let 
$w\in E_1 \cap \widehat T \cap V(x,y)$ be given, and 
let $z\in E$ be such that $f(z)=w$. By (6.25),
$z \in \widehat T$. Also, $z\in V(x,y)$ because $f_1$ moves points 
vertically. By (6.23) and (6.24), $w = \widetilde f(z)$.
By (6.21), $\widetilde f(z) = {z_1 \over x} \, \widetilde f(xz/z_1)$.
By (6.18), $z$ and $xz/z_1$ lie out of $C_\varepsilon$, and then 
$\widetilde f(xz/z_1)$ is the projection on $\partial T \cup J$
of some point $z^\ast$; (6.27) follows.

Near $x$, things are also easy. Indeed $\widehat T \cap V(x_0,x) \i B_x$
by construction, and $B_x$ does not meet $E$ by (6.14), so 
$$
E_1 \cap \widehat T \cap V(x_0,x) =\emptyset
\leqno (6.28)
$$
(again recall that $f_1$ only moves points on planes $V(s,s)$, and
inside $\widehat T$).

Near $y$, we want to use the fact that 
$$
\dist(z,P) \leq \varepsilon_0 |z-y|
\leqno (6.29)
$$
for $z\in E \cap 2B_y$ (by (6.13)) to show that
$$
\dist(w,P) \leq 2\varepsilon_0 |w-y|
\ \hbox{ for } w\in E_1 \cap \overline B_y
\leqno (6.30)
$$

Let $ w\in E_1 \cap \overline B_y$ be given, and let 
$z\in E$ be such that $w=f_1(z)$. Notice that
$z \in V(y-r_y,y+r_y)$ because $f_1$ moves points vertically.
If $f_1(z)=z$, then $\dist(w,P) = \dist(z,P)\leq \varepsilon_0 |w-y|$
by (6.29). Otherwise, $z\in \widehat T$ by (6.25). If $r_x$
is small enough (depending on $x$, $y$, and $r_y$, but this is all 
right), $\widehat T\cap V(y-r_y,y+r_y) \i 2B_y$, so $z\in 2B_y$
and (6.29) holds.

Set $z'=xz/z_1$, and let us first worry about $\widetilde f(z')$.
Recall from (6.18) that $z$ and $z'$ lie out of $C_\varepsilon$;
then $z'\in T\setminus C_\varepsilon$, and (6.20) says that
$|z'_3| \leq \rho-\varepsilon$. The point $z^\ast$ used to 
compute $\widetilde f(z')$ lies ``above" $z'$, and so
$$
\dist(\widetilde f(z'),P) = |\widetilde f(z')_3| \leq |z'_3| = \dist(z',P).
\leqno (6.31)
$$
Let us also check that
$$
|\widetilde f(z')-z'| \leq 6 \dist(z',P).
\leqno (6.32)
$$
If $\dist(z',P) \geq \rho/2$, then (6.32) holds because the diameter
of $T$ is less than $3\rho$ (recall that $\rho \geq r_x/2
\geq 5 \diam(J)/2$). Otherwise, set 
$z^\ast = (x,\xi,\rho z'_3/|z'_3|)$,  as in the definition of
$\widetilde f(z')$, and observe that
$\widetilde f(z')-z^\ast= \lambda (z'-z^\ast)$, for some $\lambda \geq 1$.
In addition, $\lambda \leq \rho/(\rho-|z'_3|)$, which corresponds to 
the radial projection of $z'$ onto $P$; see the computation above (6.21).
Then 
$$
|\widetilde f(z')-z'| \leq (\lambda-1)|z^\ast - z'|
\leq {|z'_3| \over \rho-|z'_3|} \, |z^\ast - z'|
\leq 6 |z'_3| = 6 \dist(z',P)
\leqno (6.33)
$$ 
because $\rho-|z'_3| \geq \rho/2$ and 
$|z^\ast - z'| \leq \diam(T) \leq 3 \rho$.
So (6.32) holds in all cases.

By homogeneity, (6.32) yields $|\widetilde f(z)-z| \leq 6 \dist(z,P)$.
Then 
$$
|f_1(z)-z| \leq 6 \dist(z,P) \leq 6 \varepsilon_0 |z-y|,
\leqno (6.34)
$$
because $f_1(z)$ lies on the segment $[z,\widetilde f(z)]$,
and then because (6.29) holds for $z$.
Similarly, $\dist(\widetilde f(z),P) \leq \dist(z,P)$ by (6.31)
and homogeneity, and then 
$$
\dist(f_1(z),P) \leq \dist(z,P)
\leqno (6.35)
$$
because $f_1(z)$ lies on $[z,\widetilde f(z)]$.

Next $|f_1(z)-y| \geq |z-y|-|f_1(z)-z|\geq (1-6 \varepsilon_0) |z-y|$,
by (6.34). Therefore
$$
\dist(f_1(z),P) \leq \dist(z,P) \leq \varepsilon_0 |z-y|
\leq 2 \varepsilon_0 |f_1(z)-y|
\leqno (6.36)
$$
by (6.35) and (6.29). Recall that $w=f_1(z)$; hence (6.36) proves 
(6.30).

To summarize, $E_1 \cap \widehat T \cap V(x,y)$ is contained in the 
union of the flat floor $\widehat J$ and the triangular tunnel 
$\partial \widehat T$, by (6.27). 
In addition, there is no piece of $E_1$ at the end of $\widehat T$
near $x$ (by (6.28)), and near $y$ we know that $E_1$ stays very close 
to $P$ (by (6.30)), so the end of the tunnel near $y$ is not completely
blocked by $E_1$ either.

\medskip
We are now ready to construct a \underbar{second deformation $f_2$}, 
to be composed with $f_1$. Its goal is to project $E_1 \cap B_y$ onto
$P$. We take
$$
f_2(z)=z \ \hbox{ for } z\in \R^n \setminus B_y,
\leqno (6.37)
$$
and inside $B_y$ we want $f_2$ to look like a radial projection onto 
$P \cup \partial B_y$. It will act a little like $f_1$ above, but
we won't need to move points in vertical planes. On the other hand,
we want to keep the invariance under rotations in $Q$.

Set $S_y = (y/x) S_{\xi,\rho} 
= y + (y/x)\xi + Q \cap \partial B(0,\rho_1)$, where we set
$\rho_1 = y \rho/x$. The precise choice of sphere will be useful later. 
Also set $S_{y,+} = \big\{ z\in \R^n \, ; \, \dist(z,S_y) \leq 
\rho_1/100\big\}$. Notice that $S_y$ and $S_{y,+}$ both lie well inside
$B_y$, because $\rho \leq r_x$ and we take $r_x$ much smaller than 
$r_y$. See Figure 6.3. 

\medskip 
\hskip 2.8cm 
\epsffile{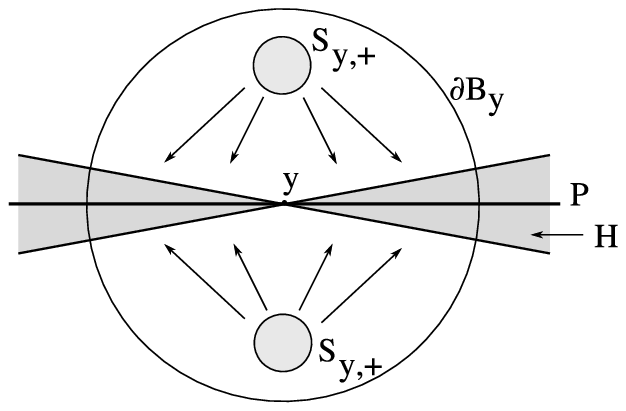}
\smallskip
 {\centerline {\bf Figure 6.3 }} 
 \medskip

For $z\in B_y \setminus S_{y,+}$, denote by $z^\sharp$ the closest point of 
$S_y$. Thus $z^\sharp = (y,(y/x)\xi,\rho_1 z_3/|z_3|)$. Then let
$f_2(z)$ be the radial projection of $z$ onto $\partial B_y \cup P$,
centered at $z^\sharp$. Let us check that $f_2$ is Lipschitz on
$B_y \setminus S_{y,+}$. Let $\tau > 0$ be very small. In the region 
where $|z_3| \leq \rho_1-\tau$, the radial projection of $z$ onto 
$P$ is well defined, and Lipschitz (we even checked below (6.20) that
it is equal to $z^\sharp + \lambda (z-z^\sharp)$, with 
$\lambda = \rho_1(\rho_1-|z_3|)^{-1}$. Then $f_2$ is Lipschitz
on this region. In the region where $|z_3| > \rho_1-\tau$, and if 
$\tau$ is small enough, the projection of $z$ onto $P$ is either not 
defined, or else lies out of $B_y$. Then $f_2(z)$ is the projection
onto $\partial B_y$, and it is Lipschitz on the region. Notice also 
that $f_2(z)=z$ on $\partial B_y$, so it is continuous across 
$\partial B_y$ and Lipschitz on $\R^n \setminus S_{y,+}$.

We now extend $f_2$ to the whole $\R^n$, so that it is Lipschitz
(probably with a very bad constant) and maps $\overline B_y$ to itself
(otherwise, compose its restriction to $\overline B_y$ with a 
Lipschitz retraction onto $\overline B_y$).

We need to know where the set $E_2 = f_2(E_1)$ lies. 
Notice that 
$$
E_2 \setminus \overline B_y = E_1 \setminus \overline B_y
\  \hbox{ and } \ 
E_2 \cap \overline B_y = f_2(E_1 \cap \overline B_y)
\leqno (6.38)
$$
by (6.37) and because $f_2(\overline B_y) \i \overline B_y$. Set
$$
H = \big\{ z \in \R^n \, ; \, 
\dist(z,P) \leq 2 \varepsilon_0 |z-y| \big\};
\leqno (6.39)
$$
thus (6.30) says that $E_1 \cap \overline B_y$ is contained in $H$.
It will be good to know that
$$
E_1 \cap \overline B_y \i \overline B_y \setminus S_{y,+},
\leqno (6.40)
$$
and for this it is enough to check that $S_{y,+}$ does not meet $H$.
Let $z=(z_1,z_2,z_3)\in S_{y,+}$ be given;
then $||z_3|-\rho_1| \leq \dist(z,S_y) \leq \rho_1/100$,
so $\dist(z,P) = |z_3| \geq 99 \rho_1 /100$. On the other hand,
$|z-y| \leq \dist(z,S_y) + \rho_1 + |(y/x)\xi -y|
\leq \rho_1/100 + \rho_1 + y r_x/100x 
\leq \rho_1/100 + \rho_1 + y \rho/50x 
\leq 2 \rho_1$
(see the statement of Lemma 6.15 and recall that $\rho_{1}=y \rho/x$). 
So $z$ lies out of $H$; (6.40) follows. 

Next we check that
$$
E_2 \cap \overline B_y 
\i P \cup [\partial B_y \cap H].
\leqno (6.41)
$$
Let $w\in E_2 \cap \overline B_y$ be given, and let $z\in E_1$
be such that $w=f_2(z)$. By (6.37), $z\in \overline B_y$.
By (6.40), $z$ lies out of $S_{y,+}$, so $w$ is the radial projection
of $z$ onto $P \cap \partial B_y$, with origin $z^\sharp$ as above.
In particular, $w \in P \cup \partial B_y$, and we just have to check
that $w \in H$ if $w\notin P$.

First assume that $z, z^\sharp$, and $y$ are not collinear,
and let $V$ denote the plane that contains them. Also assume that
$P\cap V$ is a line; then $H\cap V$ is a thin double cone centered at 
$y$ that contains $P\cap V$. Now $z^\sharp$ lies out of $H$
(because $z^\sharp \in S_y \i S_{y,+}$), and $z\in H$ by (6.30).
So the half line $[z^\sharp,z)$ starts out of $H$, then gets in $H$
before it hits $z$ (or at the same time), and it does not leave
$H$ before it meets $P$. This last may not happen, but if it happens,
it happens after $[z^\sharp,z)$ hits $z$ (the open line segment
$(z^\sharp,z)$ does not meet $P$ by definition  of $z^\sharp$).
In both cases, $w=f_2(z)$ lies in $H$.

If $P\cap V$ is not a line, it is reduced to $\{ y \}$
(because $z^\sharp$ lies out ou $P$). Then the line through $z$ and 
$z^\sharp$ does not meet $P$ (recall that $z, z^\sharp$, 
and $y$ are not collinear). By definition of 
$z^\sharp = (y,(y/x)\xi,\rho_1 z_3/|z_3|)$, this can only happen if
$|z_3| = \rho_1$. In this case, $\dist(w,P) = \dist(z,P)
\leq 2\varepsilon_0 |z-y| \leq 2\varepsilon_0 r_y$, because 
$z\in H \cap \overline B_y$. In addition, $w\in \partial B_y$
(because it does not lie in $P$), hence $|w-y|=r_y$ and $w \in H$.

Finally suppose that $z, z^\sharp$, and $y$ are collinear.
In this case $z$ lies between $z^\sharp$ and $y$ (because
otherwise $z$ would not lie in $H$; already $z^\sharp$ does not),
and $w=y$ (because $y\in P$ and $[z^\sharp,y]$ does not
meet $\partial B_y$). This completes our proof of (6.41).

Let us also check that $E_2$ does not block our access to
the tube $\widehat T$, i.e., that
$$
E_2 \cap \overline B_y \cap V(x,y) 
\ \hbox{ does not meet } \inter(\widehat T) \setminus P.
\leqno (6.42)
$$
Let $w \in E_2 \cap \overline B_y \cap V(x,y)$ be given,
and suppose it lies in $\inter(\widehat T) \setminus P$.
Let $z\in E_1$ be such that $f_2(z)=w$. Then $z\in \overline B_y$,
because otherwise $w=z$ by (6.37), and (6.27) says that this is 
impossible. In particular, $w$ is obtained from $z$ by the 
usual radial projection process with $z^\sharp$. Then $z\notin P$,
because otherwise $w=z$, and we assumed that $w \notin P$.

Since $w \in V(x,y)$ and $z^\sharp \in V(y,y)$, $z$ lies in $V(x,y)$.
Then (6.27) says that $z$ lies out of $\inter(\widehat T)$.
Call $w',z'$, and $z^\ast$ the radial projections of $w,z,$, and 
$z^\sharp$ on $V(x,x)$. Thus $w' \in \inter(T)$, $z'\notin \inter(T)$,
and $z^\ast = (x,\xi,\rho z'_3/|z'_3|)$ lies in $\inter(T)$.
This is impossible, because $T$ is convex (see near (6.19)) and 
$z'$ lies on $[z^\ast,w]$. So (6.42) holds.

Let us summarize the situation now. Set 
$$
Z= [\widehat T \cap V(x_0,y)] \cup \overline B_y.
\leqno (6.43)
$$
Recall from (6.25), the fact that $\widehat T \cap V(y,y_0) \i B_y$
if $r_x$ is small enough (see (6.22)), and (6.37) that
$$
f_2(z)=f_1(z)=z \ \hbox{ for } z \in \R^n \setminus Z. 
\leqno (6.44)
$$
Let us also  check that 
$$
f_2 \circ f_1(Z) \i Z.
\leqno (6.45)
$$
First take $z\in \widehat T \cap V(x_0,y)$, $f_1(z) \in \widehat T \cap 
V(x_0,y)$, by (6.26) and because $f_1$ moves points vertically.
If $f_1(z)$ lies in $\overline B_y$, 
$f_2 \circ f_1(z) \in \overline B_y \i Z$ because $f_2(\overline B_y)
\i \overline B_y$. Otherwise, 
$f_2 \circ f_1(z) = f_1(z) \in \widehat T \cap V(x_0,y) \i Z$,
by (6.37). Next suppose that $z\in \widehat T \cap V(y,y_0)$; then
$f_1(z) \in \widehat T \cap V(y,y_0) \i B_y$ by (6.26), 
because $f_1$ moves points vertically, and because $r_x$ is very small.
Then $f_2 \circ f_1(z) \in \overline B_y \i Z$, as before.
Finally, if $z\in \overline B_y \setminus [\widehat T \cap V(x_0,y_0)]$,
$f_1(z) = z$ by (6.25), and 
$f_2 \circ f_1(z) = f_2(z) \in \overline B_y \i Z$; (6.45) follows.

We now verify that
$$
E_2 \cap Z \i P \cup [\partial \widehat T \cap V(x,y) \setminus B_y]
\cup [\partial B_y \cap H \setminus (\inter(\widehat T) \cap V(x,y))]. 
\leqno (6.46)
$$
Indeed there is no part in $V(x_0,x)$ because of (6.28) and (6.38). 
The part of $E_2 \cap Z \setminus P$ that lies in 
$V(x,y) \setminus \overline B_y$ lies in $\partial \widehat T$, 
by (6.27) and (6.38). The rest lies in $\overline B_y$, 
hence in $\partial B_y \cap H$ by (6.41), and out of 
$\inter(\widehat T) \cap V(x,y)$ by (6.42).

This is somewhat better than before; we now have a nice floor 
$P \cap Z$, and we managed to surround its interior with a tunnel
that leads from $x$ to $y$ and a big cavity $B_y \setminus P$.

Our \underbar{third mapping $f_3$} will contract most of the floor 
$P \cap Z$ onto its boundary; this will be good, because it will
make a big part of the floor vanish. 
More precisely, set $F= P \cap Z \cap V(x,b)$. 
See Figure 6.4. 
Thus $F$ is the union of  $\widehat J \cap V(x,y)$ and 
$P \cap \overline B_y$, so it looks 
like a long $d$-dimensional cylinder attached to a much larger 
$d$-ball. [Recall that $\widehat J$ is the cone over
the $(d-1)$-ball $J=P' \cap \overline B(x+\xi,r_x/10)$; 
see above (6.19).]
There is a bilipschitz mapping $\psi$ from $F$ to the
cylinder $C = [x,y] \times J \times \{ 0 \}$, and which fixes 
the points of the end 
$F\cap V(x,x) = \{ x \} \times J \times \{ 0 \}$.

\medskip 
 \hskip 3.2cm 
\epsffile{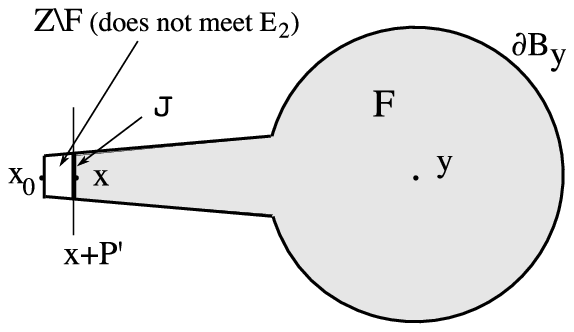}
\smallskip
 \centerline {{\bf Figure 6.4.} A picture in $P$} 
 \medskip

To contract $C$ to a piece of its boundary, we use the
mapping $\varphi$ which to $z\in C$ associates the radial projection
of $z$ onto the box $\big[ [x,y] \times \partial J \times \{ 0 \} \big]
\cup \big[\{ y \} \times J \times \{ 0 \}\big]$ (open at one side)
with center $(x_0,\xi,0)$. 
[We don't care about the precise value; the point is that $x_0 < x$, 
so $\varphi$ is Lipschitz. The contraction that we use on $F$ is
$\psi^{-1} \circ \varphi \circ \psi$. It is Lipschitz on $F$,
even though with a very large constant. We want to set
$$
f_3(z) = \psi^{-1} \circ \varphi \circ \psi(z)
\ \hbox{ for } z\in F
\leqno (6.47)
$$
and 
$$
f_3(z) = z \ \hbox{ for } z\in \overline{\R^n \setminus Z},
\leqno (6.48)
$$
so let us check that this defines a Lipschitz function
$f_3$ on $F \cup (\overline{\R^n \setminus Z})$. The main point is that 
$f_3(z) = z$ on the part of $F$ that touches $\overline{\R^n \setminus Z}$;
the rest of the verification is straightforward. Indeed, we want
to check that $|f_3(z_2)-f_3(z_1)| \leq C |z_2-z_1|$ when $z_1$ and
$z_2$ lie in $F \cup (\overline{\R^n \setminus Z})$. This is trivial if
$z_1$ and $z_2$ both lie in $F$ or in $\overline{\R^n \setminus Z}$,
so we may assume that $z_1\in F$ and $z_2 \in \overline{\R^n \setminus Z}$.
But, due to the simple geometry of $Z$, we can find 
$z_3 \in F \cap \partial Z$ such that $|z_3-z_1| \leq C |z_2-z_1|$, 
and then $|f_3(z_2)-f_3(z_1)| \leq |f_3(z_2)-f_3(z_3)|+ |f_3(z_3)-f_3(z_1)|
\leq |z_3-z_2|+C|z_3-z_1| \leq |z_3-z_1|+|z_1-z_2|+C|z_3-z_1|
\leq C'|z_2-z_1|$, as needed.

Extend $f_3$ to the rest of $Z$, so that 
$$
f_3 \hbox{ is Lipschitz on $\R^n$ and $f_3(Z)\i Z$} 
\leqno (6.49)
$$
(the existence of a Lipschitz extension is a well-known
general fact, and we can arrange that $f_3(Z)\i Z$
by composing the restriction of $f_3$ to $Z$ with a Lipschitz 
retraction onto $Z$ (such a retraction exists because there is a 
bilipschitz mapping of $\R^n$ that maps $Z$ to a cylinder).

Set $E_3 = f_3(E_2)$, and let us check that
$$
E_3\cap Z \i Z_1 \cup Z_2,
\leqno (6.50)
$$
where 
$$
Z_1 = V(x,y) \cap \partial \widehat T \setminus  B_y
\ \hbox{ and } \  Z_2 = \partial B_y \cap H.
\leqno (6.51)
$$
Let $w\in E_3\cap  Z$ and let $z\in E_2$ be such that
$w=f_3(z)$. Notice that $z\in Z$ because $f(z)=z$ for
$z\notin Z$. 
Next, $z$ lies out of $V(x_0,x)$, because $E_2 \cap V(x_0,x)
= E_1 \cap V(x_0,x)$ by (6.38), and $E_1 \cap Z \cap V(x_0,x)
\i E_1 \cap \widehat T \cap V(x_0,x) =\emptyset$ by (6.43) and (6.28).

First suppose that $z\in P$. Then $z\in F = P \cap Z \cap V(x,b)$,
and (6.47) shows that $w = f_3(z) \in \psi^{-1}\big(\big[ [x,y] \times \partial J \times \{ 0 \} \big]
\cup \big[\{ y \} \times J \times \{ 0 \}\big]\big)$.
This last set is the closure of $\partial F \setminus V(x,x)$
(where $\partial F$ denotes the boundary  of $F$ in $P$),
and it is contained in $Z_1 \cup Z_2$.

Now suppose that $z\notin P$; then (6.46) says that $z\in \partial Z$,
so $w=f_3(z)=z$ by (6.48), and $w$ lies in the set in $Z_1 \cup Z_2$ 
by (6.46) again. This proves (6.50). 

Let us also record the fact that
$$
f_3 \circ f_2 \circ f_1 (z) = z \hbox{ for $z \notin Z$, and }
f_3 \circ f_2 \circ f_1 (Z) \i  Z.
\leqno (6.52)
$$

\medskip
When $d=n-1$, this is enough to conclude because 
$H^{n-1}(E_3 \cap Z) \leq H^{n-1}(Z_1 \cup Z_2) 
\leq C \varepsilon_0 r_y^d + C_{x,y} r_x^{d-1}$
(maybe see the proof of (6.80) if you are worried about
the second term), and for this case the reader may
set $f_4(z) = f_5(z) =z$ for all $z$, $Z^\ast =Z$, and
go directly to the line below (6.89) for the final accounting
that leads to a contradiction.

For the general case, we still need to compose 
$f_3 \circ f_2 \circ f_1$ with \underbar{two last mappings}
$f_4$ and $f_5$, whose effect will be to project the potentially
large sets $E_3 \cap Z_1$ and $E_3 \cap Z_2$ onto $d$-dimensional 
skeletons of cubes. In both cases, we shall use a standard 
Federer-Fleming argument.  
 
Some amount of notation concerning dyadic cubes will be needed.
Here dyadic cubes are closed cubes obtained from $[0,2^k]^n$
for some $k\in Z$ by a translation in $2^k \Bbb Z^n$.
For such a dyadic cube $Q$,  we shall denote by $\S_d(Q)$ the 
$d$-dimensional skeleton of $Q$, i.e., the union of the $d$-dimensional 
faces of $Q$.

Also let $j \in \Bbb N$ be given; in this paper $j$ will
be chosen once and for all, sufficiently large for some
geometric conditions explained below to hold, but for all
practical purposes it is a geometric constant that depends 
only on $n$. We denote by $\Delta_j(Q)$ the set of dyadic cubes
of sidelength $2^{k-j}$ that are contained in $Q$, and then by
$\S_{d,j}(Q)$ the union of all the 
$d$-dimensional faces of dyadic cubes $R\in \Delta_j(Q)$. 
In other words, $\S_{d,j}(Q) = \bigcup_{R \in \Delta_j(Q)} \S_d(R)$.

We shall use the following lemma, which is part of  
Proposition 3.1 on page 13 of [DS]. 
We also refer to [DS] for additional background on dyadic cubes
and skeletons.

\ms\proclaim Lemma 6.53.
Let $Q$ be a dyadic cube in $\R^n$,
and let $F \i Q$ be closed, with $H^d(F) < + \infty$.
Then there is a Lipschitz mapping $\varphi_Q : \R^n \to \R^n$
such that
$$
\varphi_Q(z) = z \ \hbox{ for } z\in \overline{\Bbb R^n \setminus Q}.
\leqno (6.54)
$$
$$
\varphi_Q(F \cap Q) \i \S_{d,j}(Q) \cup \partial Q,
\leqno (6.55)
$$
$$
\varphi_Q(R) \i R \ \hbox{ for } R\in \Delta_j(Q),
\leqno (6.56)
$$
and 
$$
H^d(\varphi_Q(F \cap R)) \leq C H^d(F\cap R) 
\ \hbox{ for } R\in \Delta_j(Q).
\leqno (6.57)
$$
The constant $C$ depends only on $n$  and $d$, not on $Q$, $j$,
or $F$.

\ms
We want to start with the function $f_4$ that will take care
of $E_3 \cap Z_1$. Let $k \in \Z$ be such that 
$$
r_x/32 \leq 2^k \sqrt n \leq r_x/16.
\leqno (6.58)
$$

Let $\D_k$ denote the set of dyadic (closed) cubes $Q$ of $\R^n$
with sidelength $2^{k}$ such that 
$\dist(Q,Z_1) \leq 2^{k+1} \sqrt n$.
For each $Q\in \D_k$, we apply Lemma 6.53 to $Q$ and the set
$F = E_3\cap Q$, and get a mapping $\varphi_Q$.

We compose all the mappings $\varphi_Q$, $Q\in \D_k$, 
(in some randomly decided order) and get a mapping $\phi_0$. 
Notice that for a given $z\in \R^n$, at most one $\varphi_Q$
may move $z$, the one (if it exists) for which $z$ lies in the 
interior of $Q$. Even then, $\varphi_Q(z) \in Q$, and 
by (6.54), the other $\varphi_{Q'}$, $Q' \neq Q$, do not move 
$z$ or $\varphi_Q(z)$. Because of this, $\phi_0$ does not change 
if we compose the $\varphi_Q$, $Q\in \D_k$, in any other order. 

Set $E_{3,0} = \phi_0(E_3)$. It should be a little better
than $E_3$ near $Z_1$ because we were able to send
lots of points on the $\S_{d,j}(Q)$, but there might still be some
part of $\phi_0(E_3)$ left on the boundaries $\partial Q$.
To fix this, we shall compose with variants of $\phi_0$ for which the
boundaries $\partial Q$ lie in different places.
Denote by $v_0, v_1, \cdots v_m$ the vectors of $\R^n$ whose coordinates 
are either $0$ at $2^{k-1}$. Thus $m=2^n-1$. To make the notation coherent,
do this so that $v_0 = 0$. 

We now construct $\phi_1$. For each $Q\in \D_k$, we apply
Lemma 6.53 to $Q+v_1$ and the set $F=E_{3,0} \cap [Q+v_1]$.
The fact that $Q$ is not a dyadic cube itself, but just a
translation of such a cube by $v_1$, does not matter, but
we shall use the fact that the cubes $R\in \Delta_j(Q+v_1)$
are still dyadic cubes of sidelength $2^{k-j}$. 

Call $\varphi_{Q,1}$ the map given by Lemma 6.53. 
Then compose the new $\varphi_{Q,1}$ in any order, and get
a new mapping $\phi_1$. Finally set $E_{3,1}=\phi_1(E_{3,0})$.

We iterate this construction in the obvious way: once
we have $\phi_{l-1}$ and $E_{3,l-1}$ for some $l \leq m$, we
apply Lemma 6.53 to the $Q+v_{l}$, $Q\in \D_k$,
with $F=E_{3,l-1} \cap [Q+v_{l}]$. This gives a mapping
$\varphi_{Q,l}$. We compose all these mappings to
get $\phi_{l}$, and set $E_{3,l}=\phi_l(E_{3,l-1})$.

We end with $\phi_m$, and set 
$$
f_4 = \phi_m \circ \cdots \phi_0 \ \hbox{ and }
E_4 = f_4(E_3) = E_{3,m}.
\leqno (6.59)
$$
Our next goal is to control $E_4$ away from $Z_2$, as follows.

\medskip
\proclaim Lemma 6.60.
Let $R$ be a dyadic cube of sidelength $2^{k-j}$
such that 
$$
\dist(R,Z_1 \cup Z_2) \leq r_x/2
\ \hbox{ and } \  \dist(R,Z_2)\geq 2^{k}.
\leqno (6.61)
$$
Then 
$$
H^d(R \cap E_4) \leq C 2^{kd}.
\leqno (6.62)
$$

\ms
Let us first see how $H^d(E_{3,l}\cap R)$ evolves with time.
Set $E_{3,-1}=E_3$ and $\varphi_{Q,0}=\varphi_Q$ to make the notation
uniform. Recall from (6.57) with $F=E_{3,l-1} \cap [Q+v_{l}]$ that
$$
H^d(\varphi_{Q,l}(R\cap E_{3,l-1})) \leq C H^d(R\cap E_{3,l-1})
\leqno (6.63)
$$
for $0 \leq l \leq m$ and $Q\in \D_k$. This is only useful when
$Q$ is the cube such that $Q+v_{l}$ contains $R$, because for the other
cubes $Q'$, (6.54) says that $\varphi_{Q',l}(z)=z$ on $R$. Thus
$$
H^d(\phi_{l}(R\cap E_{3,l-1})) \leq C H^d(R\cap E_{3,l-1}).
\leqno (6.64)
$$

Observe also that $\phi_{l}(R') \i R'$ for every 
dyadic cube $R'$ of sidelength $2^{k-j}$, because (6.54) and (6.56)
say that this is true for every $\varphi_{Q,l}$. 
Thus, if $z\in R'$ and $\phi_{l}(z)$ lies in $R$, $R'$ meets $R$.
That is,
$$
R \cap E_{3,l} = R \cap \phi_{l}(E_{3,l-1})
\i \bigcup_{R' \in {\cal N}(R)} \phi_{l}(R' \cap E_{3,l-1}),
\leqno (6.65)
$$
where we denote by ${\cal N}(R)$ the set of dyadic cubes $R'$ of 
sidelength $2^{k-j}$ that touch $R$. Consequently,
$$
H^d(R \cap E_{3,l})  
\leq \sum_{R' \in {\cal N}(R)} H^d(\phi_{l}(R' \cap E_{3,l-1}))
\leq C \sum_{R' \in {\cal N}(R)} H^d(R'\cap E_{3,l-1}),
\leqno (6.66)
$$
by (6.64). Iterations of this yield
$$
H^d(R \cap E_{4}) = H^d(R \cap E_{3,m}) 
\leq C \sum_{R' \in {\cal N}_{m+1}(R)} H^d(R'\cap E_{3}),
\leqno (6.67)
$$
where the set ${\cal N}_{l}(R)$ of neighbors of order $l$
of $R$ is just defined by induction by ${\cal N}_{1}(R)={\cal N}(R)$
and ${\cal N}_{l+1}(R) = \bigcup_{R' \in {\cal N}_{l}(R)} {\cal N}(R')$.
An obvious induction yields 
$$
\dist(R',R) \leq (l-1) 2^{k-j} \sqrt n
\ \hbox{ for } R' \in {\cal N}_{l}(R),
\leqno (6.68)
$$
so, if $z_0$ is any point of $R$,
$$
\bigcup_{R' \in {\cal N}_{m+1}(R)} R'
\i B\big(z_0, (m+2) 2^{k-j} \sqrt n\big)
= B\big(z_0, (2^{n}+1) \, 2^{k-j} \sqrt n\big)
\i B(z_0, 2^{k-2})
\leqno (6.69)
$$
because $m=2^{n}-1$, and if $j$ is chosen large enough.

When $R$ lies far enough from $Z_1 \cap Z_2$, all this will lead to
(6.62) rather directly. Let us first check that for a dyadic cube $S$ 
of sidelength $2^{k-j}$,
$$
H^d(E_3 \cap S) \leq C 2^{kd}
\ \hbox{ when $\dist(S,Z_1 \cup Z_2) \leq r_x$ and 
$S$ does not meet $Z_1 \cup Z_2$.}
\leqno (6.70)
$$
Indeed $S$ does not meet $E_3 \cap Z$ by (6.50), and then 
$S \cap E_3 \i S \cap E$ by (6.52). So (6.70) follows from
the local Ahlfors-regularity of $E$ (Lemma 2.15 above, but 
only the easier upper bound and in the easier special case 
of locally minimal sets); this is why we assumed
that $\dist(S,Z_1 \cup Z_2) \leq r_x$, to make sure that $S \i B$
for some small ball $B$ such that
$2B \i A \i B(0,b)$. See in particular (6.14), the definition of 
$A$ above (6.5), and the statement of Theorem 6.2.

Now we can already check (6.62) when
$$
2^{k} \leq \dist(R,Z_1 \cup Z_2) \leq r_x/2.
\leqno (6.71)
$$
Indeed, by (6.69) and (6.58), all the cubes $S \in {\cal N}_{m+1}(R)$
satisfy the condition in (6.70), and (6.62) follows from (6.67).
We are thus left with the case of cubes $R$ such that
$\dist(R,Z_1 \cup Z_2) < 2^{k}$ and hence, by (6.61)
$$
\dist(R,Z_1) \leq 2^{k}.
\leqno (6.72)
$$

Pick $z_0 \in R$, denote its coordinates by $z_{0,i} \,$, 
$1 \leq i \leq n$, and then choose $\varepsilon_i \in \{ 0, 2^{k-1} \}$ 
so that $z_{0,i}-\varepsilon_i$ lies at distance at least $2^{k-2}$ from
$2^{k}\Bbb Z$. Recall that there is a $v_l$ such that
$v_l = (\varepsilon_1, \cdots, \varepsilon_n)$, and 
then every coordinate of $z_0 -v_l$ lies at distance at least 
$2^{k-2}$ from $2^{k}\Bbb Z$. Now 
$$
\dist(z_0,v_l + \partial Q) \geq 2^{k-2}
\ \hbox{ for } Q\in \D_k,
\leqno (6.73)
$$
because the points in the faces that compose $\partial Q$
have at least one coordinate in $2^{k}\Bbb Z$.

Now let $S$ denote any cube of ${\cal N}_{m+1}(R)$, and let $Q$
be such that $Q+v_l$ contains $S$. Observe that 
$$\eqalign{
\dist(Q,Z_1) &\leq |v_l| + \dist(S,Z_1) 
\leq |v_l| + 2^{k-2} + \dist(z_0,Z_1) 
\cr& \leq |v_l| + 2^{k-2} + \diam(R) +\dist(R,Z_1)
\cr& \leq 2^{k-1} \sqrt n + 2^{k-2} + 2^{k-j} \sqrt n +\dist(R,Z_1)
\leq 2^{k+1} \sqrt n
}\leqno (6.74)
$$
by (6.69) and (6.72) in particular. So $Q \in \D_k$ (see the 
definition below (6.58)), and (6.55) says that 
$$
\varphi_{Q,l}(E_{3,l-1} \cap S) \i
\varphi_{Q,l}(E_{3,l-1} \cap [v_l +Q]) \i v_l +[\S_{d,j}(Q) \cup \partial Q].
\leqno (6.75)
$$
Now $S\in {\cal N}_{m+1}(R)$, so $S \i B(z_0, 2^{k-2})$ by (6.69),
and (6.73) says that $S$ does not meet $v_l + \partial Q$.
Recall from (6.56) that $\varphi_{Q,l}(S) \i S \,$; then
$\varphi_{Q,l}(E_{3,l-1} \cap S) \i v_l +\S_{d,j}(Q)$, and
$H^d(\varphi_{Q,l}(E_{3,l-1} \cap S)) \leq C 2^{kd}$. Now
$$
H^d(\phi_l(E_{3,l-1} \cap S)) \leq C 2^{kd}
\leqno (6.76)
$$
too, because by (6.54) all the other $\varphi_{Q',l}$ that compose $\phi_l$ 
leave the points of $S$ alone. We shall remember that (6.76) holds for all
$S\in {\cal N}_{m+1}(R)$.

We may now conclude as before, with repeated uses of 
(6.64). That is,
$$\eqalign{
H^d(R \cap E_{4}) &= H^d(R \cap E_{3,m})
= H^d(R \cap \phi_m(E_{3,m-1}))
\cr&\leq \sum_{S \in {\cal N}_{1}(R)} H^d(\phi_m(S \cap E_{3,m-1}))
\leq C \sum_{S \in {\cal N}_{1}(R)} H^d(S \cap E_{3,m-1})
\cr&=C \sum_{S \in {\cal N}_{1}(R)} H^d(S \cap \phi_{m-1}(E_{3,m-2}))
\cr&\leq C \sum_{S \in {\cal N}_{2}(R)} H^d(\phi_{m-1}(S \cap E_{3,m-2})) 
\leq C^2 \sum_{S \in {\cal N}_{2}(R)} H^d(S \cap E_{3,m-2})
\cr \cdots &  
\leq C \sum_{S \in {\cal N}_{m-l+1}(R)} H^d(\phi_{l}(S \cap E_{3,l-1}))
\leq C 2^{kd}
}\leqno (6.77)
$$
by (6.76). This completes our proof of Lemma 6.60.
\qed

\ms
Let us record a few facts about $f_4$ before we move to $f_5$.
Recall that if $R$ is a dyadic cube of sidelength $2^{k-j}$,
then all the mappings that compose $f_4$ send $R$ to $R$, so
$f_4(z) \in R$ for $z\in R$. Consequently,
$$
|f_4(z)-z| \leq 2^{k-j} \sqrt n
\ \hbox{ for } z\in \R^n.
\leqno (6.78)
$$
Also set $Z_1^\ast = \big\{ z\in \R^n \, ; \, 
\dist(z,Z_1) \leq 2^{k+2} \sqrt n \big\}$. Notice that
$Z_1^\ast$ contains all the cubes $Q+v_l$, $Q \in \D_k$.
[See below (6.58) for the definition of $\D_k$.]
Then
$$
f_4(z) = z \hbox{ for } z\notin Z_1^\ast
\ \hbox{ and } \ 
f_4(Z_1^\ast) \i Z_1^\ast
\leqno (6.79)
$$
because the $\phi_l$ only act nontrivially on the cubes $Q+v_l$,
and preserve them. Let us also check that
$$
H^d\big(\big\{ z\in E_4 \cap Z_1^\ast \, ; \, 
\dist(z,Z_2) \geq 2^{k+1} \big\}\big) \leq C_{x,y} \, r_x^{d-1},
\leqno (6.80)
$$
where $C_{x,y}$ depends on $x$, $y$, but not on $r_x$ and $k$.

Call $Y$ the set on the left-hand side of (6.80). 
Cover $Y$ with dyadic cubes $R$ of sidelength $2^{k-j}$. 
Notice that (6.61) holds when $T$ meets $Y$, so 
$H^d(E_4 \cap R) \leq C 2^{kd}$ by Lemma 6.60. Recall from (6.51) 
that $Z_1 = V(x,y) \cap \partial \widehat T \setminus  B_y$
is contained in a thin tube of width $C r_x$ around
the segment $[x,y]$, where $C$ may depend on $x$ and $y$
(typically, through the ratio $|y|/|x|$). Then 
$Z_1$, and also $Y$, can be covered by less than
$C |x-y|/r_x$ cubes $R$ as above (recall from (6.58) that $r_x$ and
$2^k$ are comparable); each cube contributes at most
$C 2^{kd} \leq C r_x^d$ by (6.62), so (6.80) follows. The reader should
not be shocked by the crude estimate above; we shall compensate by
taking $r_x$ really small.

\ms
Lemma 6.60 only controls $E_4=f_4(E_3)$ far from $Z_2$,
and indeed $E_4$ could still contain a big piece of 
$Z_2 = \partial B_y \cap H$. So we shall compose with
a \underbar{very last mapping} $f_5$. The strategy is the same
as for $Z_1$, but this time we shall work at the scale $2^t$,
where $t \in \Bbb Z$  is such that
$$
2^{t-1} \leq \varepsilon_0 r_y \leq 2^t,
\leqno (6.81)
$$
instead of $2^k$. Notice that we shall choose $r_x$ very small, 
depending on $r_y$ and $\varepsilon_0$, so we can assume that
$k$ is much smaller than $t$.

Denote by $\D_t$ the set of dyadic 
cubes $Q$ of sidelength $2^{t}$ such that 
$$
\dist(Q,Z_2) \leq 2^{t+1} \sqrt n
\leqno (6.82)
$$

We shall repeat exactly the same construction
as before. Denote by $w_0, \cdots w_m$ the $2^n$
vectors with coordinates in $\{ 0, 2^{t-1} \}$. We inductively
build successive mappings $\phi_l$, $0 \leq l \leq m$,
and each one is obtained by composing Lipschitz functions 
$\varphi_{l,Q}$, $Q\in \D_t$. We start with $E_{4,-1}=E_4$,
and later set $E_{4,l}=\phi_l(E_{4,l-1})$. Each
$\varphi_{l,Q}$ is obtained by applying Lemma 6.53
to $Q+w_l$ and the set $F=E_{4,l-1} \cap [Q+w_l]$.
All this defines the mappings $\varphi_{l,Q}$ and $\phi_l$ by induction.
At the end, we set $f_5 = \phi_m \circ \cdots \phi_0$
and $E_5 = f_5(E_4) = \phi_m(E_{4,m-1})$.
Here is the analogue of Lemma 6.60.

\medskip
\proclaim Lemma 6.83.
Let $R$ be a dyadic cube of sidelength $2^{t-j}$
such that $\dist(R,Z_2) \leq r_y/2$. Then
$$
H^d(R \cap E_5) \leq C 2^{td}.
\leqno (6.84)
$$

\ms
We shall start with the easier cases where we can prove that
$$
H^d(E_4 \cap S) \leq C 2^{td}
\ \hbox{ for } S \in {\cal N}_{m+1}(R),
\leqno (6.85)
$$
where ${\cal N}_{m+1}(R)$ still denotes the set of 
dyadic cubes of sidelength $2^{t-j}$ that are neighbors
of $(n+1)^{st}$ order of $R$.

Indeed, once we have (6.85), we can obtain (6.84) by the same
argument as before: the analogue of (6.67) holds and yields (6.84).

Let us first take care of the contribution of $E$ in (6.85).
Choose $z_0 \in R$ and let $S\in {\cal N}_{m+1}(R)$ be
given; the proof of (6.69) says that $S \i B(z_0,2^{t-2})$.
Since we want to apply the local Ahlfors-regularity 
of $E$, we should check that $B(z_0,2^{t-1}) \i B(0,b)$
(because $E$ is minimal in $B(0,b)$). And indeed
$B(z_0,2^{t-1}) \i 2B_y$ because 
$\dist(R,\overline B_y) \leq \dist(R,Z_2) \leq r_y/2$
and both $2^t$ and $\diam(R)$ are much smaller than $r_y \,$,
and $2B_y \i A \i B(0,b)$ by definition of $r_y$ (below (6.13)).
The (easier) upper bound in Lemma 2.15 says that 
$$
H^d(E\cap S) \leq  H^d(E\cap B(z_0,2^{t-2})) \leq  C 2^{td}; 
\leqno (6.86)
$$
so it will be enough to control $S \cap E_4 \setminus E$.

We start with the case when $\dist(R,Z_1\cup Z_2) > 2^{t+1}$. 
Recall that $S \i B(z_0,2^{t-2})$,
so $\dist(S,Z_1\cup Z_2) \geq 2^{t}$. In particular 
$S$ does not meet the set $Z_1^\ast$ of (6.79), and (6.79) 
says that $E_4 \cap S=E_3 \cap S$.
In addition, $E_3 \cap S$ does not meet $Z_1 \cup Z_2$,
so (6.50) says that it does not meet $Z$. By (6.52), 
$E_3 \cap S = E_3 \cap S \setminus Z \i E$. Altogether,
$E_4 \cap S \i E$, and (6.85) follows from (6.86).

We may now assume that $\dist(R,Z_1\cup Z_2) \leq 2^{t+1}$,
and our second case will be when $\dist(R,Z_2) \geq 2^{t+1}$.
Let $z\in E_4 \cap S$ be given. Then 
$\dist(z,Z_2) \geq \dist(z_0,Z_2) - 2^{t-2}
\geq \dist(R,Z_2) - \diam(R) - 2^{t-2} \geq 2^{t}$
because $S \i B(z_0,2^{t-2})$. If $z\in Z_1^\ast$,
it lies in the set of (6.80); therefore,
$$
E_4 \cap S \cap Z_1^\ast \leq  C_{x,y} \, r_x^{d-1},
\leqno (6.87)
$$
and we can assume that $z\notin Z_1^\ast$. Then $z\in E_3$, by (6.79).
If $z\in Z$, (6.50) says that $z\in Z_1 \cup Z_2$. This is impossible,
because $Z_1 \i Z_1^\ast$ and $\dist(z,Z_2) \geq 2^{t}$.
So $z\notin Z$. Let $w\in E$ be such that $z=f_3\circ f_2 \circ f_1(w)$.
By (6.52), $w\in Z$, and then $z=w$. So $z\in E$, and (6.86) already
takes care of it. This proves (6.85) and (6.84) in our second case.

We are left with the case when $\dist(R,Z_2) \leq 2^{t+1}$.
In this case we cannot prove (6.85), but instead we proceed as in
the last case of Lemma 6.60 (when (6.72) holds). The point is that
$R$ is so close to $Z_2$ that it lies in the middle of the grid
$\D_t$, and there is an $l \leq m$ such that none of the 
$S \in {\cal N}_{m+1}(R)$ meets a face $w_l +\partial Q$, $Q\in \D_t$.
Then the corresponding images $\phi_l(E_{4,l-1} \cap S)$ are contained
in $d$-dimensional grids $\S_{d,j}(Q)$, whose $H^d$-measures are 
controlled. The proof is the same as before, and completes our proof
of Lemma 6.83.
\qed

\ms
Set $Z_2^\ast = \big\{ z\in \R^n \, ; \, 
\dist(z,Z_2) \leq 2^{t+2} \sqrt n \big\}$. Notice that
$Z_2^\ast$ contains all the cubes $Q+w_l$, $Q \in \D_t \,$;
the argument here and below is just the same as for $Z_1^\ast$, near 
(6.79). Then
$$
f_5(z) = z \hbox{ for } z\notin Z_2^\ast
\ \hbox{ and } \ 
f_5(Z_2^\ast) \i Z_2^\ast
\leqno (6.88)
$$
because the $\phi_l$ only act nontrivially on the cubes $Q+w_l$,
and preserve them. 

Recall from (6.51) that $Z_2= \partial B_y \cap H$, where 
$H$ is a small neighborhood of the $d$-plane $P$, as in (6.39); 
then we can cover $Z_2$ by less than $C \varepsilon_0^{-d+1}$
balls of radius $\varepsilon_0 r_y$. By the same argument, and
because (6.81) says that $2^{t}$ and $\varepsilon_0 r_y$ are 
comparable, we need less than $C \varepsilon_0^{-d+1}$ dyadic
cubes of sidelength $2^{t-j}$ to cover $Z_2^\ast$. All these cubes $R$,
if they touch $Z_2^\ast$, are such that $\dist(R,Z_2) \leq r_y/2$,
so Lemma 6.83 says that $H^d(R \cap E_5) \leq C 2^{td} \leq 
C \varepsilon_0^d r_y^d$. Altogether,
$$
H^d(E_5 \cap Z_2^\ast) \leq C \varepsilon_0 r_y^d \,.
\leqno (6.89)
$$

We are ready for the final comparison.
Set $f=f_5 \circ f_4 \circ f_3 \circ f_2 \circ f_1$, and then
$\varphi_t(z)= tf(z)+(1-t)z$ for $0 \leq t \leq 1$ and $z\in \R^3$.
The maps $\varphi_t$ satisfy the properties (2.2)-(2.4).

Set $Z^\ast = Z \cup Z_1^\ast \cup Z_2^\ast$. Notice that
$$
f(z)=z \hbox{ for } z \notin Z^\ast, \hbox{ and }
f(Z^\ast) \i Z^\ast
\leqno (6.90)
$$
by (6.52), (6.79), and (6.88).
So, if we define the $W_t$ and $\widehat W$ as in (2.5), 
we get that $\widehat W$ is contained in the convex hull of $Z^\ast$, 
and (2.6) holds (with $U = B(0,b)$) because 
$Z^\ast \i B(0,b)$.

We can apply (4.9) or (4.4) (depending on which definition of 
minimality we use), with $h(\delta)=0$. Let us first see
how to derive the final contradiction if we use (4.9). We get that
$H^d(E \setminus E_5) \leq H^d(E_5 \setminus E)$ 
(because here $F = \varphi_1(E) = f(E) = E_5$), or equivalently
$$
H^d(E \cap Z^\ast) \leq H^d(E_5 \cap Z^\ast),
\leqno (6.91)
$$ 
because $E$ and $E_5$ coincide out of $Z^\ast$. Let us check that
$$
H^d(E_5 \cap Z^\ast) 
\leq C \varepsilon_0 r_y^d + C_{x,y} \, r_x^{d-1} \, .
\leqno (6.92)
$$
Let $z \in E_5 \cap Z^\ast$ be given. The case when 
$z\in Z_2^\ast$ is taken care of by (6.89). If
$z\in Z_1^\ast \setminus Z_2^\ast$, (6.88) says that
$z\in E_4$; then $\dist(z,Z_2) \geq 2^{t+2} \sqrt n$
because $z\notin Z_2^\ast$, so $z$ lies in the set of 
(6.80); (6.80) takes care of this case. So we may assume that
$z\in Z^\ast \setminus Z_1^\ast \cup Z_2^\ast$. By 
(6.88) and (6.79), $z \in E_3$. Also, $z\in Z$, so 
(6.50) says that $z\in Z_1 \cup Z_2$. This case is
impossible, because $Z_i \i Z_i^\ast$; (6.92) follows.

On the other hand,
$$
H^2(E \cap Z^\ast) \geq H^2(E \cap B_y) \geq C^{-1} r_y^d,
\leqno (6.93)
$$
by Lemma 2.15, because $y\in E$, $E$ is reduced, and $2B_y \i B(0,b)$.

Incidentally, (6.93) would be easier to prove here because we have 
(6.13). Or we could have used the fact that the upper density of $E$ 
at almost every point is larger than some positive constant, restricted 
to such an $y$ in Proposition 6.11, and then chosen $r_y$ specifically 
so that (6.93) holds. 

We get the desired contradiction with (6.91) or (6.92) by choosing 
$\varepsilon_0$ small, and then $r_x$ very small 
(depending on $x$, $y$, and $r_y$). 

The reader may be worried about the case when $d=1$;
this case is sufficiently simple to be treated directly,
but let us say a few words about it for the sake of completeness.
In this case, $P' = \{ 0 \}$ and lemma 6.15 says that there is 
a $\rho \in (r_x/2, 2r_x/3)$ such that the cone $\widehat S_\rho$ 
over the $(n-2)$-dimensional sphere $S_\rho = x + [Q \cap \partial B(0,\rho)]$ 
in the vertical hyperplane $x+Q$, does not meet $A\cap E$.
At the same time, (6.13) says that the $(n-1)$-dimensional cone $\Sigma$
centered at $y$ and defined by $\dist(z,P) = 2 \varepsilon_0 |z-y|$
does not meet $E \cap B_y \setminus \{ y \}$. We choose $r_x$ so small
that $\widehat S_\rho$ meets $\Sigma$ well inside $B_y$.
Denote by $\Omega$ the region of $V(x,y)$ bounded by 
$\widehat S_\rho \cup \Sigma$; thus $\Omega$ is the intersection of
$V(x,y)$ with a thin double-cone centered at $0$ and $y$.
Notice that $\Omega \cap V(x,x)$ does not meet $E$,
by (6.14). It is easy to contract $\Omega$ onto the line segment $[x,y]$
by radial vertical motion, and since $E$ does not meet $\partial \Omega$ 
except at $y$, we can make sure that the contraction $f$ is Lipschitz 
and does not move the points of $E\setminus \Omega$. 
In addition, $f(E \cap \Omega)$ does not 
meet $B_x$ (again by (6.14)), so we can contract the whole piece 
$f(E \cap \Omega) \i [x,y]$ to the point $y$ along $[x,y]$. 
The final contradiction is now obtained
as in the general case.

Return to the case when $d>1$.
If we use (4.4) (instead of (4.9)), we get that 
$$
H^d(E \cap W) \leq H^d(f(E \cap W)),
\leqno (6.94)
$$
with $W = \{ w \in \R^3 \, ; \, f(w) \neq w \}$. Recall that
$f(z)=z$ out of $Z^\ast$ (by (6.90)), so $W \i Z^\ast$ and 
$$
H^d(f(E \cap W)) \leq H^d(f(E \cap Z^\ast)) \leq H^d(E_5 \cap Z^\ast)
\leqno (6.95)
$$
(by (6.90) again). Notice that 
$$
|f_5(z) - z| \leq 2^{t-j} \sqrt n
\ \hbox{ for } z \in \R^n,
\leqno (6.96)
$$
by the proof of (6.78). Then $E_5 = f_5 \circ f_4(E_3)$
does not meet $B(y,r_y/2)$, by (6.78), (6.96), and because 
$E_3$ does not meet $B_y$ (by (6.50) and the definition (6.43)).
Hence $E \cap B(y,r_y/2) \i E \cap W \setminus E_5$, just because 
$f(z) \neq z$ if $z \in E \setminus E_5$, and
$$\eqalign{
H^d(E \cap B(y,r_y/2)) 
&\leq H^d(E \cap W) \leq H^d(f(E \cap W)) 
\cr&\leq  H^d(E_5 \cap Z^\ast)
\leq C \varepsilon_0 r_y^2 + C_{x,y} \, r_x^{d-1},
}\leqno (6.97)
$$
by (6.94), (6.95) and (6.92). On the other hand, 
$H^d(E \cap B(y,r_y/2)) \geq C^{-1} r_y^d$ by 
Lemma~2.15, so we also get a contradiction when we use (4.4)
in the accounting.

This completes our proof of Proposition 6.11. As was explained before, 
(6.12) follows and so $E$ coincides with a cone inside $A$. \qed

\medskip
We return to the proof of Theorem 6.2. So far we proved that there
is a cone $\cal C$ centered at the origin, that satisfies (6.4). We 
still need to check that $\cal C$ is a reduced minimal set. 
Notice that it is reduced, by (6.4) and because $E$ is reduced.
[See Definition 2.12 for the definition.] Suppose it is not minimal. Then 
we can find a competitor $F= \varphi_1({\cal C})$ where $\{ \varphi_t \}$ 
satisfies (2.2)-(2.6) with $U=\R^3$, and such that
$$
H^2(F\setminus {\cal C}) < H^2({\cal C} \setminus F).
\leqno (6.98)
$$
Observe that we use Definition 4.8 for minimality, but Proposition 4.10
tells us that we could also have used Definition 4.3 (or even Definition 4.1,
since here $h(r)=0$).

Recall from (2.6) that the set $\widehat W$ of (2.5) is bounded. 
Since $\cal C$ is a cone, we may always dilate everything by a small 
constant (this does not change (6.98)), so we may  assume 
that $\widehat W \i B(0,a)$.

We want to use $F$ and the $\varphi_t$ to construct a deformation of
$E$ and contradict its minimality. First observe that 
${\cal C} \cap B(0,b)$ is a deformation of $E$ in $B(0,b)$. Indeed,
pick any $b'\in (a,b)$, and set $g(x) = 0$ for $x\in \overline B(0,a)$,
$g(x) = { b'(|x|-a) \over b'-a}{ x \over |x|}$ for 
$x\in B(0,b') \setminus \overline B(0,a)$, and $g(x) = x$ for 
$x\in B(0,b) \setminus B(0,b')$. 
Then set $f_t(x)=tg(x) + (1-t) \, x \,$
for $0 \leq t \leq 1$. It is easy to see that the $g_t$ satisfy
(2.2)-(2.6) with $U=B(0,b)$, and that 
$g_1(E) = {\cal C} \cap B(0,b)$
(by (6.4)). We now put our two deformations together. That is, we set
$\psi_t(x) = g_{2t}(x)$ for $0 \leq t \leq 1/2$ and 
$\psi_t(x) = \varphi_{2t-1}(g_1(x))$ for $1/2 \leq t \leq 1$, and get 
a new family that satisfies (2.2)-(2.6) with $U=B(0,b)$. 
By minimality of $E$, 
$$
H^2(\psi_1(E) \setminus E) \geq H^2(E \setminus\psi_1(E)).
\leqno (6.99)
$$
But $\psi_1(E) = F \cap B(0,b)$, so (6.99) just says that
$H^2(E) \leq H^2(F \cap B(0,b))$.

Recall from (6.3) that $H^2(E \cap B(0,r)) = \theta_0 r^2$
for $a < r < b$. Also, $H^2({\cal C} \cap A) = H^2(E \cap A)$
by (6.4), so $H^2({\cal C} \cap B(0,r)) = \theta_0 r^2$ 
for $r>0$ and $H^2({\cal C} \cap B(0,b)) = H^2(E)$. Altogether, 
$H^2({\cal C} \cap B(0,b)) = H^2(E) \leq H^2(F \cap B(0,b))$,
which contradicts (6.98). This completes the proof of 
Theorem 6.2. \qed

\bigskip
\noindent {\bf 7. Almost constant density}
\medskip
In this section we use Theorem 6.2 to show that if $E$
is a reduced almost-minimal set of dimension $d$ in $B(0,b) \i \R^n$, 
with sufficiently small gauge function, and if 
$\theta(r) = r^{-d}H^d(E \cap B(0,r))$
varies sufficiently little between $a$ and $b$, then $E$ is very 
close to a reduced minimal cone in $A=B(0,b) \setminus \overline B(0,a)$.

We also have a (simpler) result when $a=0$ (Proposition~7.24),
and we shall see in Proposition 7.31 that every blow-up limit of a 
reduced almost-minimal set of dimension $d$ in $U \i \R^n$ is a 
reduced minimal cone.

\medskip\proclaim Proposition 7.1.
Let $n$, $d$ (the dimensions), $a$, and $b$ be given, with $0 < a < b < +\infty$.
For each $\tau >0$, we can find $\varepsilon > 0$ so that the following
holds. Let $E$ be a reduced almost-minimal set of dimension $d$
in $B(0,b) \i \R^n$, with gauge function $h$, and assume that 
$h(2b) \leq \varepsilon$. Also assume
that $\theta(b) \leq \theta(a) + \varepsilon$. Then there is a reduced 
minimal cone $\C$ centered at the origin such that
$$
\dist (x,\C) \leq \tau 
\ \hbox{ for } x \in E \cap B(0,b-\tau) \setminus B(0,a+\tau),
\leqno (7.2)
$$
$$
\dist (x,E) \leq \tau 
\ \hbox{ for } x \in \C \cap B(0,b-\tau) \setminus B(0,a+\tau),
\leqno (7.3)
$$
$$\eqalign{
|H^d(E\cap B(x,r))-H^d(\C \cap B(x,r))| &\leq \tau
\cr
\ \hbox{ for } x\in \R^n \hbox{ and } r>0
\hbox{ such that } &B(x,r) \i B(0,b-\tau) \setminus B(0,a+\tau),
}\leqno (7.4)
$$
and 
$$
|\theta(r) - r^{-d} H^d(\C \cap B(0,r))| \leq \tau
\ \hbox{ for } a+\tau \leq r \leq b-\tau.
\leqno (7.5)
$$

\medskip
See Definition 2.12 for the notion of reduction, and observe that if
$E$ is not reduced, Remark 2.14 allows us to apply Proposition 7.1
to $E^\ast$. Our statement allows
$\C= \emptyset$, but only if $\theta(r)$ is very small (by (7.5)).
Notice that (7.5) forces $\theta(r)$ to be very close to 
the constant density for $\C$. This is a constraint on $\theta(r)$,
but also can be a way to determine $\C$ in terms of $\theta(r)$.

\ms
We shall prove the proposition by compactness.
Fix $n$, $d$, $a$, and $b$, and consider a sequence $\{ E_n \}$ of reduced
almost-minimal sets on $B(0,b)$, with gauge functions $h_n$ such that
$\varepsilon_n = h_n(2b)$ tends to $0$. It is a standard fact about
the Hausdorff distance that we can extract a subsequence $\{ n_k \}$
so that the sets $E_{n_k}$ converge to a limit $F\i B(0,b)$, as in
the beginning of Section~3 (see for instance [D2], Proposition 34.6 on 
page 214). 

Set $\widetilde h_l = \sup_{n \geq l} h_n$ for $l \geq 0$. Notice that 
$\widetilde h_l$ is nondecreasing (just like the $h_n$), and 
$\lim_{r \to 0} \widetilde h_l(r) =0$, because each $h_n$ tends to $0$
and $h_n(r) \leq \varepsilon_n$ for $r\leq 2b$.
Then Lemma~4.7 (applied to the $E_{n_k}$)
says that for each $l$, $F$ is almost-minimal in $B(0,b)$,
with gauge function $\widetilde h_l$. But since we can take $l$ as large 
as we want and $\widetilde h_l(2b)$ tends to $0$, we get that 
$$
F \hbox{ is a minimal set in } B(0,b). 
\leqno (7.6)
$$
We also get that
$$
H^d(F\cap V) \leq \liminf_{k \to +\infty} H^d(E_{n_k}\cap V)
\leqno (7.7)
$$
for every open set $V \i B(0,b)$, by (3.4) in Lemma 3.3 (applied with
$U=B(0,b)$, $M=1$, $\delta=2b$, and $h=\widetilde h_l(2b)$ for some 
large $l$). 
We also need an inequality in the other direction, so we apply 
Lemma 3.12 to $E_{n_k}$, with the same choice of $U$, 
$M$, $\delta$, and $h$, and then let $l$ tend to 
$+\infty$; this yields
$$\limsup_{k \to +\infty} \, H^d(E_{n_k} \cap H)
\leq H^d(F \cap H)
\hbox{ for every compact set } H \i B(0,b).
\leqno (7.8)
$$
Notice that $H^d(F \cap \partial B(x,r)) = 0$ for almost-every 
$r\in(a,b)$; (7.7) and (7.8) yield
$$\eqalign{
\limsup_{k \to +\infty} \, H^d(E_{n_k} \cap B(x,r))
&\leq H^d(F \cap B(x,r)) = H^d(F \cap \overline B(x,r))
\cr&\leq \liminf_{k \to +\infty} H^d(E_{n_k}\cap B(x,r))
}\leqno (7.9)
$$
for such $r$, so that 
$$
H^d(F \cap B(0,r)) = \lim_{k \to +\infty} H^d(E_{n_k}\cap B(0,r))
\ \hbox{ for almost-every } r\in(a,b).
\leqno (7.10)
$$
Set $\theta_k(r) = r^{-d} H^d(E_{n_k} \cap B(0,r))$ for $r \leq b$,
and suppose in addition that $\theta_k(b) \leq \theta_k(a) + \varepsilon_k$.
Pick $a'$ and $b'$, with $a<a'<b<b'$, such (7.10) holds for $a'$ and $b'$.
Then 
$$\eqalign{
H^d(F \cap B(0,b')) &= \lim_{k \to +\infty} H^d(E_{n_k}\cap B(0,b'))
\leq \liminf_{k \to +\infty} \, H^d(E_{n_k}\cap B(0,b))
\cr& = b^d \liminf_{k \to +\infty} \, \theta_k(b)
\leq b^d \liminf_{k \to +\infty} \, [\theta_k(a) +\varepsilon_k]
\cr& = b^d \liminf_{k \to +\infty} \,\theta_k(a)
= b^d a^{-d} \liminf_{k \to +\infty} \, H^d(E_{n_k}\cap B(0,a))
\cr& \leq b^d a^{-d}\lim_{k \to +\infty} H^d(E_{n_k}\cap B(0,a'))
= b^d a^{-d} H^d(F \cap B(0,a')).
}\leqno (7.11)
$$
We let $b'$ tend to $b$ and get that
$$
\theta(b) =: b^{-d} H^d(F \cap B(0,b)) \leq a^{-d} H^d(F \cap B(0,a'))
= a^{-d} (a')^{d}\theta(a').
\leqno (7.12)
$$
Recall from Proposition 5.16 and (7.6) that $\theta$ is nondecreasing 
on $(a,b)$, so $\theta(a^+) = \lim_{a' \to a^+} \theta(a')$ exists,
and (7.12) says that $\theta(b) \leq \theta(a^+)$. But 
$\theta(b) \geq \theta(a^+)$ by Proposition 5.16, so 
$\theta(b) = \theta(a^+)$ and $\theta$ is constant on $(a,b)$.
By Theorem 6.2, there is a reduced minimal cone $\C$ centered 
at the origin such that
$$
F \cap A = \C \cap A.
\leqno (7.13)
$$
In particular, (7.2) and (7.3) hold for
$E_{n_k}$ and $k$ large enough (and with any $\tau > 0$ given in 
advance). 

It will be good to know that 
$$
H^d(\C \cap \partial B) = 0 \,\hbox{ for every ball $B$,}
\leqno (7.14)
$$
but since this is easy to believe, we shall only
check it in Lemma 7.34, at the end of the section. 
Now let  a ball $B$ be given, with
$B\i B(0,b-\tau/2) \setminus B(0,a+\tau/2)$. 
By (7.14), (7.13), (7.7), (7.8), and the proof of (7.10),
$$
\lim_{k \to + \infty} H^d(E_{n_k} \cap B)
= H^d(F \cap B) = H^d(\C \cap B),
\leqno (7.15)
$$
so the analogue of (7.4) for $E_{n_k}$ and $B$ holds for $k$ large. 

\ms
We are now ready to prove Proposition 7.1. Let $a$, $b$, and $\tau$
be given, and suppose that we cannot find $\varepsilon$ so that the 
proposition holds. Then we can find $E_n \i B(0,b)$ that satisfies the
hypotheses with $\varepsilon_n = 2^{-n}$, but not the conclusion. 
We extract a convergent subsequence, as above, and get a set $F$ such 
that (7.6)-(7.13) hold. Let us show that (7.2)-(7.5) hold for $k$ large;
this will prove the desired contradiction.

We already know that (7.2) and (7.3) hold for $k$ large. 
For (7.4) we only know that it holds for $k$ large when $x$
and $r$ are fixed, so we need a little bit of uniformity.
By (7.14), $H^d(\C \cap B(x,r))$ is a continuous function of
$x$ and $r$, so we can find a finite collection of balls $B_i$, 
$i\in I$, such that $B_i \i B(0,b-\tau/2) \setminus B(0,a+\tau/2)$
for $i\in I$, and such that whenever $B(x,r)$ is a ball such that
$B(x,r) \i B(0,b-\tau) \setminus B(0,a+\tau)$ (as in (7.4)), we can
find $i,j \in I$ such that
$$
B_i \i B(x,r) \i B_j \ \hbox{ and } \ 
H^d(\C \cap B_j \setminus B_i) \leq \tau/5.
\leqno (7.16)
$$

We also know from (7.15) that for $k$ large enough,  
$$
|H^d(\C \cap B_i)-H^d(E_{n_k} \cap B_i)| \leq \tau/5
\ \hbox{ for $i\in I$.}
\leqno (7.17)
$$

Now let $B(x,r)$ be as in (7.4), and let $i,j \in I$ be such that
(7.16) holds. Then 
$$
H^d(\C \cap B_i) \leq H^d(\C \cap B(x,r)) \leq H^d(\C \cap B_j)
\leq H^d(\C \cap B_i) + \tau/5
\leqno (7.18)
$$
and, for $k$ large (depending on $I$, but not on $B(x,r)$),
$$\eqalign{
H^d(\C \cap B_i)-\tau/5 &\leq H^d(E_{n_k} \cap B_i) 
\leq H^d(E_{n_k} \cap B(x,r)) \leq H^d(E_{n_k} \cap B_j) 
\cr&
\leq  H^d(\C \cap B_j) +\tau/5 \leq H^d(\C \cap B_i)+2\tau/5,
}\leqno (7.19)
$$
by (7.17) and (7.16). Thus $|H^d(\C \cap B(x,r)) - H^d(E_{n_k} \cap B(x,r))| 
< \tau$, by (7.18) and (7.19), so (7.4) holds.

We are left with (7.5) to prove. Call $\theta_0$ the constant density
of $\C$; that is, $\theta_0 = r^{-d} H^d(\C \cap B(0,r))$ for $r>0$. 
Notice that $H^d(\C \cap \partial B(0,r)) = 0$ for $a < r < b$, so
$$
\lim_{k \to +\infty} H^d(E_{n_k}\cap B(0,r)) = H^d(F \cap B(0,r))
= H^d(\C\cap B(0,r)) = \theta_0 r^d,
\leqno (7.20)
$$
by (7.7), (7.8), and (7.13). As before, we need some uniformity, 
so we pick a finite (but large) collection $\{ r_i \}_{i\in I}$ 
in $(a,b)$. We know that for $k$ large enough,
$$
|r_i^{-d} H^d(E_{n_k}\cap B(0,r_i)) - \theta_0| \leq \tau/2
\ \hbox{ for $i\in I$.}
\leqno (7.21)
$$
Let $r \in [a+\tau,b-\tau]$ be given. If our collection
$\{ r_i \}$ is dense enough  
we can find $i,j \in I$ such that $r_i \leq r \leq r_j$ and 
$r_j/r_i$ is as close to $1$ as we want. Then for $k$ large
$$
\theta_k(r) = r^{-d} H^d(E_{n_k}\cap B(0,r))
\leq r^{-d} H^d(E_{n_k}\cap B(0,r_j))
\leq r^{-d} r_j^d (\theta_0 + \tau/2) \leq \theta_0 + \tau,
\leqno (7.22)
$$
by (7.21). Similarly, 
$$
\theta_k(r) \geq r^{-d} H^d(E_{n_k}\cap B(0,r_i))
\geq r^{-d} r_i^d (\theta_0 - \tau/2) \leq \theta_0 - \tau,
\leqno (7.23)
$$
and (7.5) holds for $k$ large. This completes the proof of Proposition 7.1 by 
contradiction.
\qed

\ms 
Proposition 7.1 easily extends to the case when $a=0$. 
Let us give a slightly more invariant statement for this to simplify later uses.

\ms\proclaim Proposition 7.24.
For each $\tau >0$, we can find $\varepsilon > 0$ so that the following holds.
Let $E$ be a reduced almost-minimal set in the open set $U \i \R^n$, 
with gauge function $h$. Let $x\in E$ and $r>0$ be such that $B(x,r) \i U$, 
$h(2r) \leq \varepsilon$, and
$$
\theta(x,r) \leq \inf_{0 < t < r/100} \theta(x,t) + \varepsilon,
\leqno (7.25)
$$
where we set $\theta(x,t) = t^{-d} H^d(E\cap B(x,t))$.
Then there is a reduced minimal cone $\C$ centered at $x$, such that
$$
\dist (y,\C) \leq \tau r
\ \hbox{ for } y \in E \cap B(x,(1-\tau)r),
\leqno (7.26)
$$
$$
\dist (y,E) \leq \tau r
\ \hbox{ for } y \in \C \cap B(x,(1-\tau)r),
\leqno (7.27)
$$
and 
$$\eqalign{
|H^d(E\cap B(y,t))-H^d(\C \cap B(y,t))| &\leq \tau r^d
\cr
\ \hbox{ for } y\in \R^n \hbox{ and } t>0
&\hbox{ such that } B(y,t) \i B(x,(1-\tau)r).
}\leqno (7.28)
$$

\ms 
By scale invariance, is enough to prove this when $x=0$, $r=1$, and 
$E$ is almost-minimal in (an open set $U$ that contains) $B(0,1)$.
We can proceed as in the proof of Proposition 7.1 (with $a=0$ and 
$b=1$), up to (7.10) included. Then pick $a_1\in (0,1/100)$ and 
$b' \in (a_1,1)$ such that (7.10) holds for $a_1$ and $b'$. 
By (7.25), $\theta_k(1) \leq \theta_k(a_1)+\varepsilon_k$, so 
$$\eqalign{
H^d(F \cap B(0,b')) &= \lim_{k \to +\infty} H^d(E_{n_k}\cap B(0,b'))
\leq \liminf_{k \to +\infty} \, H^d(E_{n_k}\cap B(0,1))
\cr& 
= \liminf_{k \to +\infty} \, \theta_k(1)
\leq \liminf_{k \to +\infty} \, [\theta_k(a_1) +\varepsilon_k]
=\liminf_{k \to +\infty} \,\theta_k(a_1)
\cr&
= a_1^{-d} \liminf_{k \to +\infty} \, H^d(E_{n_k}\cap B(0,a_1))
= a_1^{-d} H^d(F \cap B(0,a_1)),
}\leqno (7.29)
$$
by (7.10). We let $b'$ tend to $1$ and get that 
$$
H^d(F \cap B(0,1))\leq a_1^{-d} H^d(F \cap B(0,a_1)).
\leqno (7.30)
$$
Now $F$ is minimal in $B(0,1)$, by (7.6), so Proposition 5.16 says that 
$\theta(r) = r^{-d} H^d(F \cap B(0,r))$ is nondecreasing in $(0,1)$.
Then $\theta(1) = \lim_{r \to 1} \theta(r) \geq \theta(a_1)$. But (7.30)
says that $\theta(1) \leq \theta(a_1)$, so $\theta(r)$ is constant on 
$[a_1,1]$, hence also on $(0,1]$ because
we can take $a_1$ as close to the origin as we want.
Theorem 6.2 says that $F$ is a reduced minimal cone, and we can continue the 
argument as in Proposition 7.1.
\qed

\ms
Let us also give a baby version of Proposition~7.24.

\ms\proclaim Proposition 7.31.
Let $E$ be a reduced almost-minimal set in the open set $U$,
and let $x\in E$ be such that
$$
\theta(x) = \lim_{r \to 0} r^{-d} H^d(E\cap B(x,r))
\ \hbox{ exists.}
\leqno (7.32)
$$
Then every blow-up limit of $E$ at $x$ is a reduced minimal cone
$F$ centered at the origin, and $H^{d}(F \cap B(0,1)) = \theta(x)$.

\ms
See Definition 2.12 for ``reduced", and Definition 4.8 or 4.3
for almost-minimality. We do not care about the precise values of
the gauge function $h(r)$ here, we just need it to tend to $0$
when $r$ tends to $0$, but of course (7.32) is easier to obtain
when $h$ is small, by Proposition 5.24. A blow-up limit
of $E$ at $x$ is any closed set in $\R^n$ that can be obtained as
the limit (as in (3.1)) of a sequence $\{ r_{k}^{-1}(E-x) \}$,
with $\lim_{k \to +\infty} r_k = 0$.

To prove the proposition, assume that the sequence 
$\{ r_{k}^{-1}(E-x) \}$ converges to some limit $F$.
Notice that for every small $a>0$, 
$F_k = r_{k}^{-1}(E-x)$ is, for $k$ large enough, a reduced 
almost-minimal set with gauge function $\widetilde h(t) = h(at)$.
Then Lemma 4.7 says that $F$ is a reduced almost-minimal set, 
with gauge function $\widetilde h$. But since we can take 
$a$ as small as we want, and $\lim_{a \to 0} \widetilde h(t)=0$
for each $t$, $F$ is a reduced minimal set in $\R^n$.

We also get that 
$$
H^d(F \cap B(0,r)) = \lim_{k \to +\infty} H^d(F_k\cap B(0,r))
\ \hbox{ for almost every } r > 0,
\leqno (7.33)
$$
by the proof of (7.10). But
$H^d(F_k\cap B(0,r)) = r_{k}^{-d} H^d(E \cap B(x,r_k r))$,
which tends to $r^d \theta(x)$ by (7.32). Thus
$H^d(F \cap B(0,r)) = r^d \theta(x)$ for almost 
every $r > 0$. Theorem 6.2 says that $F$ is a cone, and 
Proposition 7.31 follows.
\qed

\ms
Propositions 7.1 and 7.24 will help us find many places where
$E$ looks a lot like a minimal cone. We shall then use a variant of
Reifenberg's topological disk theorem to get a local parametric 
description of $E$.

\ms
We end this section with the proof of (7.14).

\ms 
\proclaim Lemma 7.34. Let $\C$ be a closed cone in $\R^n$ (centered at the 
origin), and assume that $H^d(\C \cap B(0,1)) < +\infty$.
Then $H^d(\C \cap \partial B) = 0$ for every ball $B$.

\ms
Let $0 < a < b <+\infty$ be given, and set $A=\big\{ x\in \R^n \, ; \, 
a \leq |x| \leq b \big\}$; it is enough to show that
$H^d(\C \cap A \cap\partial B) = 0$. Denote by $S = \partial B(0,1)$
the unit sphere, and set $\varphi(r,\theta)=r\theta$ for
$r\in [a,b]$ and $\theta \in S$. Notice that $\varphi$ is a 
bilipschitz mapping from $[a,b] \times S$ to $A$, and
$\varphi^{-1}(\C \cap A) = [a,b] \times T$, where
$T = S \cap \C$. Then
$$
H^d([a,b] \times T) = H^d(\varphi^{-1}(\C \cap A))
\leq C H^d(\C \cap A) < +\infty.
\leqno (7.35)
$$
Let us first check that 
$$
H^{d-1}(T) \leq C H^d([a,b] \times T) < +\infty,
\leqno (7.36)
$$
where $C$ depends also on $a$ and $b$, but we don't care.

Let $\varepsilon > 0$ be given. By (7.35), there is a covering of 
$[a,b] \times T$ by balls $B_i=B(x_i,r_i)$, $i\in I$,
with $r_i \leq\varepsilon$ for $i\in I$ and  
$$
\sum_{i\in I} r_i^d \leq C H^d([a,b] \times T) + \varepsilon.
\leqno (7.37)
$$
By compactness, we may assume that $I$ is finite. For $t\in [a,b]$,
set $P_t = \{ t \} \times S$. Set $B_{i,t} = B_i \cap P_t$, and denote 
by $I_t$ the set of $i\in I$ such that $B_{i,t} \neq \emptyset$.
Notice that the balls $B_{i,t} = B_i \cap P_t$, $i\in I_t$ form a
covering of $T$ by balls of radius at most $\varepsilon$, so
$$
H^{d-1}_{\varepsilon}(T) \leq c_{d-1} \sum_{i\in I_t} \diam(B_{i,t})^{d-1}
\leq C \sum_{i\in I_t} H^{d-1}(B_{i,t})
= C \sum_{i\in I} H^{d-1}(B_{i,t})
\leqno (7.38)
$$
where $H^{d-1}_{\varepsilon}(T)$ is the usual infimum over coverings
that is used to define $H^{d-1}(T)$. We average this over $[a,b]$ and get 
that
$$\eqalign{
H^{d-1}_{\varepsilon}(T)
&\leq C \int_{[a,b]} \, \sum_{i\in I} \, H^{d-1}(B_{i,t}) dt
= C \sum_{i\in I} \int_{[a,b]} H^{d-1}(B_{i,t}) dt
\cr&\leq C \sum_{i\in I} r_i^d 
\leq C H^d([a,b] \times T) + C\varepsilon
}\leqno (7.39)
$$
by (7.37). Then we let $\varepsilon$ tend to $0$ and we get (7.36).

By (7.36), we can apply Theorem 2.10.45 in [Fe], 
which says that there are positive constants $C_1$ and $C_2$
such that
$$
C_1 H^d(U \times V) \leq H^1(U) H^{d-1}(V) \leq C_2 H^d(U \times V)
\leqno (7.40)
$$
for $U \i [a,b]$ and $V \i T$ measurable. In other words,
the restriction of $H^d$ to $[a,b] \times T$ is comparable
to the product of the Lebesgue measure on $[a,b]$ with
the restriction of $H^{d-1}$ to $F$. 

Now let $B$ be any ball; notice that $\partial B$ never meets a 
radial line more than twice, which means that for $v\in T$,
$\varphi^{-1}(\C \cap \partial B\cap A))$ does not meet
$[a,b] \times \{ v \}$ more than twice. Then
$$\eqalign{
H^d(\C \cap \partial B \cap A)
&\leq C H^d(\varphi^{-1}(\C \cap \partial B \cap A))
\cr&\leq C \int_{v\in T} H^{1}\big([a,b] \times \{ v \} 
\bigcap \varphi^{-1}(\C \cap \partial B\cap A)\big) dH^{d-1}(v)
= 0
}\leqno (7.41)
$$
because $\varphi$ is bilipschitz and by Fubini.
Lemma 7.34 follows. \qed

\bigskip
\noindent {\bf C. MINIMAL CONES OF DIMENSION 2}
\medskip

In this part we study the reduced minimal cones $E$ of dimension $2$ in $\R^n$.
The goal is not to obtain a finite list as when $n=3$, but
a fairly simple description of $K = E\cap \partial B(0,1)$ as a finite
collection of geodesic arcs of circles in $\partial B(0,1)$,
with lengths bounded from below, and that meet only by sets of three 
with $120^\circ$ angles. See Proposition 14.11.

We start with a description of the tangent cones to $E$,
where we show that if $x$ lies in $K = E\cap \partial B(0,1)$, every
blow-up limit of $K$ at $x$ is a one-dimensional reduced minimal set.
See Theorem 8.23 for a more precise statement where $K$ is approximated
by minimal sets in small balls centered at $x$. This 
description is based on the fact that if $E \times \R^m$ is a 
reduced $(m+d)$-dimensional minimal set in $\R^n \times \R^m$,
then $E$ is a $d$-dimensional minimal set in $\R^n$. 
See Proposition 8.3.

With this and the simple version of Reifenberg's topological disk 
theorem which we describe in Section 12, we could easily obtain 
information on the bi-H\"oder structure of $E$, but since it is
preferable to use $C^1$ information, we will show that $K$ satisfies 
a weak form of almost-minimization (see Definition 9.1 and 
Proposition 9.4), and then follow some of the general theory
to prove that these weak almost-minimal sets of dimension $1$ have 
some form of $C^1$ regularity. This is done in Sections 11-13, and
for the local description of $K$ by curves we use a simple version of 
Reifenberg's Topological Disk Theorem (Proposition~12.6).

In the mean time we show in Section 10 that the reduced minimal sets 
of dimension~$1$ are cones. See Theorem 10.1.

We state and prove the description of $K = E \cap \partial B(0,1)$
when $E$ is a reduced minimal cone of dimension $2$ in Section~14. 
The fact that $K$ is composed of $C^1$ curves comes
from the description of the weak almost-minimal sets given in 
Sections 11-13, and then a construction of harmonic competitors for a cone,
which is done in Section 13, allows us to check that all the little 
arcs that compose $K$ are arcs of great circles of $\partial B(0,1)$.

\bigskip
\noindent {\bf 8. Minimality of a product and tangent sets}
\medskip

We start our study of minimal cones with what would be the 
most natural attempt to get a nice description of the $d$-dimensional 
minimal sets by induction on $d$. It will be easy to
control the minimal sets of dimension $1$. Then assume that 
we have a reasonable description of the minimal sets of 
dimension $d$; we can try to use this to get information
on minimal cones of dimension $d+1$, because we shall see in 
this section that blow-up limits of these cones are products 
of a minimal set of dimension $d$ by a line. In turn
a good control on the minimal cones may lead to valuable
information on minimal sets of the same dimension,
by the results of Section 7. We shall only be able to
make this plan work for $2$-dimensional minimal sets, but
we shall write down the first stage in full generality.

In this section, we first show that if the product $E \times \R^m$ 
is a reduced minimal set of dimension $d+m$, then 
$E$ is a reduced dimensional minimal set of dimension $d$.
We shall then use this to show that if the cone over $K$
is a reduced minimal set of dimension $d+1$, then the blow-up
limits of $K$ are minimal sets of dimension $d$, and even that
$K$ is close to minimal sets of dimension $d$ in small balls.

Let us recall what we mean by minimal sets. This is
the same as almost-minimal sets as in Definition 4.3 or 4.8, 
with $h(\delta)=0$ for all $\delta$, but since the definition
simplifies slightly, let us state it again for future reference.

\ms\proclaim Definition 8.1.
The closed set $E \i \R^n$ is a minimal set of dimension $d$
if, whenever $\varphi: \R^n \to \R^n$ is a Lipschitz function
such that $\varphi(x)=x$ out of some compact set of $\R^n$,
then
$$
H^{d}(E \setminus \varphi(E)) \leq H^{d}(\varphi(E) \setminus E).
\leqno (8.2)
$$
We say that $E$ is reduced when $H^d(E\cap B(x,r)) > 0$
for $x\in E$ and $r>0$.

\ms
See Definition 2.12 for the original definition of reduced sets.
Concerning our new definition of minimality, it is easily seen
to be equivalent to Definition 4.3 with the gauge function $h=0$. 
The point is that if $\varphi$ is as in Definition 8.1, we can set
$\varphi_t(x) = t \varphi(x) + (1-t) x$ for $0\leq t \leq 1$;
the family $\{ \varphi_t \}$ satisfies (2.2)-(2.6) for
some $\delta > 0$, and then (8.2) is the same as (4.9) with 
$h(\delta)=0$. Let us also recall that Definitions 4.3 and 4.8
are equivalent, by Proposition 4.10.

The first result of this section is a simple observation 
on the structure of some minimal product sets.

\ms\proclaim Proposition 8.3.
Let $m$, $n$, and $d$ be integers, and let $E \times R^m$ be 
a reduced $(d+m)$-dimensional minimal set in $\R^{n+m}$. Then 
$E$ is a reduced $d$-dimensional minimal set in $\R^{n}$.

\ms
Let $\varphi : \R^n \to \R^n$ be a lipschitz function
such that $\varphi(x)=x$ out of some compact set; we want 
to show that (8.2) holds. Our plan is to use $\varphi$
to construct a Lipschitz function $f : \R^{n+m} \to \R^{n+m}$
and apply the minimality of $E \times R^m$.

Let $R>0$ be large, to be chosen near the end,
and choose a smooth radial cut-off function $\psi$ on
$\R^m$, with 
$$\eqalign{
\psi(y) = 1 \hskip0.5cm &\hbox{ for } \hskip0.3cm |y| \leq R,
\cr
0 \leq \psi(y) \leq 1 \hskip0.5cm &\hbox{ for } 
\hskip0.3cm R < |y| < R+1,
\cr
\psi(y) = 0 \hskip0.5cm &\hbox{ for } \hskip0.3cm |y| \geq R+1,
}\leqno (8.4)
$$
and $|\nabla \psi(y)| \leq 2$ everywhere. 
Then define $g : \R^{n}\times \R^m \to \R^n$ by
$$
g(x,y) = \psi(y) \varphi(x) + (1-\psi(y)) \, x
\ \hbox{ for $x\in \R^n$ and $y\in \R^m$} 
\leqno (8.5)
$$
and set $f(x,y) = (g(x,y),y)$. Notice that
$g(x,y) = x + \psi(y) (\varphi(x)-x)$, and
since $\varphi(x)-x$ is bounded,
$$
\hbox{$g$ and $f$ are Lipschitz, with bounds that do not depend on $R$.}
\leqno (8.6)
$$

Set $W^\sharp = \{ (x,y)\in \R^n \times \R^m \, ; 
\, f(x,y) \neq (x,y) \}$. If $(x,y) \in W^\sharp$,
then $\psi(y) \neq 0$ and $\varphi(x) \neq x$. That is,
$W^\sharp \i W \times B(0,R+1)$, where 
$W = \{ x\in \R^n \, ; \, \varphi(x) \neq x \}$.
In addition, if $B$ is a ball that contains $W \cup \varphi(W)$,
then 
$$
W^\sharp \cup f(W^\sharp) \i \Sigma, \hbox{ where we set }
\Sigma = B \times B(0,R+1).
\leqno (8.7)
$$

By assumption, $E^\sharp = E \times R^m$ is minimal, and
(8.2) says that 
$$
H^{d+m}(E^\sharp\setminus f(E^\sharp))
\leq H^{d+m}(f(E^\sharp)\setminus E^\sharp).
\leqno (8.8)
$$
By (8.7), both sets are contained in $\Sigma$; we add
$H^{d+m}(E^\sharp\cap f(E^\sharp) \cap \Sigma)$ to both sides and get 
that
$$
H^{d+m}(E^\sharp\cap \Sigma)
\leq H^{d+m}(f(E^\sharp\cap \Sigma)).
\leqno (8.9)
$$

We shall need to be a little careful when we compute
Hausdorff measures of product spaces, because the Hausdorff measure
of a product is not always exactly the product of the Hausdorff 
measures. Fortunately, $E^\sharp$ is rectifiable (with dimension $d+m$), 
by Section~2 (or directly [Al]), 
so it will be possible to use the following variant of
the coarea formula, which we borrow from [Fe], 
Theorem~3.2.22. We change the notation a little, and
reduce slightly the generality.

\ms\proclaim Lemma 8.10.
Let $d, m, n$ be nonnegative integers, let $Z\i \R^n$ be a rectifiable
set of dimension $m+d$, with $H^{m+d}(E)<+\infty$, and
let $\pi : Z \to \R^m$ be Lipschitz. Then 
$$
\pi^{-1}(y) \hbox{ is a $d$-dimensional rectifiable 
set for $H^m$-almost every 
$y\in \R^m$,}
\leqno (8.11)
$$
and 
$$\eqalign{
\int_{Z \cap A} apJ_m \pi(z) \, dH^{d+m}(z)
&= \int_{y\in \R^m}\bigg\{\int_{z\in \pi^{-1}(y)} {\bf 1}_{A\cap Z}(z)
\, dH^d(z) \bigg\}\, dH^m(y)
}\leqno (8.12)
$$
for every Borel set $A \i \R^n$, where the approximate Jacobian 
$apJ_m \pi$ of dimension $m$ is defined (almost-everywhere) by 
$apJ_m \pi(z) = ||\wedge_m D\pi(z)||$; 
see 1.3.1 in [Fe]. 

\ms
We apply Lemma 8.10 to $Z = E^\sharp \cap \Sigma$, 
and where $\pi$ is the restriction to $Z$ of the second 
projection $\pi: \R^{n+m} \to \R^m$. Here $\pi$ is even linear,
and since $\pi$ is $1$-Lipschitz and the tangent plane 
to $E^\sharp$ always contains the $m$-plane
$\{ 0 \} \times R^m$ where $D\pi$ is the identity,
the approximate Jacobian $apJ_m \pi$ is equal to $1$.

Also, $\pi^{-1}(y) = E \cap B$
for $y\in B(0,R+1) \i \R^m$, and $\emptyset$ otherwise, 
so (8.11) says that
$$
E \cap B \hbox{ is rectifiable,}
\leqno (8.13)
$$
and (8.12) with $A= \Sigma$ says that
$$
H^{d+m}(E^\sharp\cap \Sigma)
= H^d(E\cap B) \ H^m(\R^m \cap B(0,R+1)).
\leqno (8.14)
$$

Similarly set $\Sigma_1 = B \times B(0,R)$, and observe
that $f(x,y) = (\varphi(x),y)$ for $(x,y) \in \Sigma_1$,
by (8.5) and because $\psi(y)=1$ for $y\in B(0,R)$. Then
$$
f(E^\sharp \cap \Sigma_1)
= f\big[(E \cap B) \times B(0,R)\big]
= \varphi(E \cap B) \times B(0,R)
=[\varphi(E) \cap B] \times B(0,R)
\leqno (8.15)
$$
since $\varphi(E\cap B) = \varphi(E)\cap B$
because by definition of $B$, $\varphi(B) \i B$ and 
$\varphi(x)=x$ out of $B$.

Notice that  $f(E^\sharp)$ is rectifiable, because
$E^\sharp$ is rectifiable and $f$ is Lipschitz.
Apply Lemma 8.10 to the product set 
$f(E^\sharp \cap \Sigma_1) =[\varphi(E) \cap B] \times B(0,R)$ 
(which is rectifiable because it is contained in $f(E^\sharp)$),
and with the same projection $\pi$ as before. The proof
of (8.14) yields
$$\eqalign{
H^{d+m}(f(E^\sharp \cap \Sigma_1))
&=H^{d+m}\big([\varphi(E) \cap B] \times B(0,R)\big)
\cr&= H^d(\varphi(E) \cap B) \ H^m(B(0,R)).
}\leqno (8.16)
$$

Set $\Xi = E^\sharp \cap (\Sigma \setminus \Sigma_1)
= [E \cap B] \times [B(0,R+1) \setminus B(0,R)]$. Now
$$
H^{d+m}(f(\Xi)) \leq C H^{d+m}(\Xi) 
= C H^d(E \cap B) \, H^m(B(0,R+1) \setminus B(0,R))
\leq C R^{m-1}
\leqno (8.17)
$$
because $f$ is Lipschitz and by Lemma 8.10 (applied to 
the rectifiable product set $\Xi \,$; the proof is again the
same as for (8.14)).
Here $C$ depends on the Lipschitz constant for $\varphi$
through $|\nabla f|$, and on $E \cap B$. This is unusual,
but the main point is that $C$ does not depend on $R$, by (8.6). 

Return to the right-hand side of (8.9).
Let $z$ lie in $f(E^\sharp\cap \Sigma)$, and
let $(x,y) \in E^\sharp\cap \Sigma$ be such that
$z=f(x,y)$. If $y\in B(0,R+1) \setminus B(0,R)$,
then $(x,y) \in \Xi$ and $z \in f(\Xi)$.
Otherwise, $\psi(y)=1$, $g(x,y) = \varphi(x)$,
and $z = (\varphi(x),y) \in [\varphi(E) \cap B] \times B(0,R)$. Thus 
$f(E^\sharp\cap \Sigma) \i f(\Xi) \cup [\varphi(E) \cap B] \times B(0,R)$,
and (8.9) yields
$$\eqalign{
H^{d+m}(E^\sharp\cap \Sigma)
&\leq H^{d+m}(f(E^\sharp\cap \Sigma))
\cr&\leq H^{d+m}(f(\Xi)) +  H^{d+m}([\varphi(E) \cap B] \times B(0,R))
\cr&\leq C R^{m-1} +  H^{d+m}(f(E^\sharp \cap \Sigma_1))
\cr&\leq C R^{m-1} + H^d(\varphi(E) \cap B) \,  H^m(B(0,R))
}\leqno (8.18)
$$
by (8.17), (8.15), and (8.16). Now apply (8.14), divide by $H^m(B(0,R+1))$, 
and use (8.18) to get that 
$$
H^d(E\cap B) 
= H^m(B(0,R+1))^{-1} H^{d+m}(E^\sharp\cap \Sigma)
\leq H^d(\varphi(E) \cap B) + C R^{-1}.
\leqno (8.19)
$$
Then let $R$ tend to $+\infty$, remove from both sides 
$H^d(E\cap \varphi(E) \cap B)$, and get (8.2).

\smallskip
We still need to show that $E$ is reduced. So let $x\in E$
and $r>0$ be given. The same computations as above show that
$$
H^{d+m}(E^\sharp \cap [B(x,r) \times B(0,1)]) = 
H^d(E\cap B(x,r)) \, H^m(B(0,1)\cap \R^m). 
\leqno (8.20)
$$
The left-hand side is positive because $E^\sharp$ is reduced, 
so $H^d(E\cap B(x,r)) > 0$, as needed.
This completes our proof of Proposition 8.3.
\qed

\ms\noindent{\bf Remark 8.21.}
The proof of Proposition 8.3 relies on the fact that we work
with minimal sets, so that we were able to use the Lipschitz
norm  of $\varphi$ in the argument above, compensate by
taking $R$ enormous, and get away with it. There will be
an argument of the same type in the next section, 
but we will then be working with almost-minimal sets,
and have to pay by including the Lipschitz properties
of $\varphi$ in the definition of weak almost-minimal sets.

\ms\noindent{\bf Remark 8.22.}
Probably the proof of Proposition 8.3 allows us to show that
if a product $E^\sharp = E \times F$ of a $d_1$-dimensional
closed set $E \i \R^{n_1}$ with a $d_2$-dimensional closed set 
$F \i \R^{n_2}$ is a $(d_1+d_2)$--dimensional minimal set
in $\R^{n_1+n_2}$, and if some mild assumptions on the 
behavior of $H^{d_i}(E_i \cap B(0,R))$ for $R$ large
are satisfied, then each $E_i$ is a minimal set of dimension
$d_i$ in $\R^{n_i}$.

I don't know whether the converse is true, i.e., whether
products of minimal sets are necessarily minimal. There is
an argument in [D2] for the product of a minimal set with
$\R^{n_2}$ (see page 530 and Exercise 76.16 in [D2]), 
but in the different context of Mumford-Shah minimal sets
in codimension 1; the argument does not apply here,
even for $E^\sharp = E \times \R^{n_2}$, because the
sections $\varphi_1(E) \cap [\R^{n_1} \times \{ y \}]$, 
$y\in \R^{n_2}$, are not always competitors for 
$E \times \{ y \}$ in $\R^{n_1} \times \{ y \}$.

\bigskip
\proclaim Theorem 8.23. 
Let the integers $0 < d < n$ be given. For each $\varepsilon >0$,
we can find $r_0 \in (0,1)$ such if $E$ is a reduced minimal
cone of dimension $d+1$ in $\R^n$, $x$ is a point of
$K = E \cap \partial B(0,1)$, $0 < r \leq r_0$,
and $P$ denotes the tangent hyperplane to $\partial B(0,1)$
at $x$, then there is a reduced minimal set $L$ of dimension $d$ 
in $P$ such that 
$$
\dist(y,L) \leq \varepsilon r \hbox{ for } y\in K\cap B(x,r)
\hbox{ and }
\dist(z,K) \leq \varepsilon r \hbox{ for } z\in L\cap B(x,r),
\leqno (8.24)
$$
and 
$$\eqalign{
H^d(L \cap B(y,t-\varepsilon r)) - \varepsilon r^d 
\leq H^d(K \cap B(y,t)) 
\leq H^d(L \cap B(y,t+\varepsilon r)) +\varepsilon r^d 
}\leqno (8.25)
$$
for every ball $B(y,t)$ contained in $B(x,r)$.

\ms
We start with a few comments  before the proof.
The definition of minimal sets was recalled (and slightly
simplified) in Definition~8.1.
A minimal cone is simply a minimal set which happens
to be a cone centered at the origin.
The two possible definitions of reduced minimal sets in $P$ 
(where we see $L$ as a subset of $P$, or directly as a subset 
of $\R^n$) are equivalent. Indeed it is obvious that if
$L$ is minimal in $\R^n$, then it is minimal in $P$
because every deformation in $P$ is a deformation in $\R^n$.
Conversely, suppose that $L$ is not minimal in $\R^n$. Let
$\varphi : \R^n \to \R^n$ be Lipschitz, with $\varphi(x)=x$
out of a compact set, and such that
$H^{d}(E \setminus \varphi(E)) > H^{d}(\varphi(E) \setminus E)$.
Set $\widetilde \varphi = \pi \circ \varphi$, where $\pi$ denotes
the orthogonal projection on $P$. Let $B$ be a large ball centered 
on $P$ such that $\varphi(x)=x$ out of $B$ and $\varphi(B) \i B$.
Thus $E$, $\varphi(E)$, and $\widetilde \varphi(E)$ coincide 
out of $B$. In addition,
$$\eqalign{
H^d(\widetilde \varphi(E) \cap B) 
&= H^d(\pi(\varphi(E)) \cap B) = H^d(\pi(\varphi(E) \cap B)) 
\cr&\leq H^d(\varphi(E) \cap B) < H^d(E\cap B),
}\leqno (8.26)
$$
so $E$ is not  minimal in $P$ either.

Observe also that Theorem 8.23 contains the slightly easier
observation that any blow-up limit of $K=E \cap \partial B(0,1)$ 
at a point $x\in K$ is a reduced minimal set of dimension $d$ in 
the tangent plane to $\partial B(0,1)$ at $x$. For a direct proof
of this observation, simply follow the proof of the theorem up
to (8.30).

Our estimate (8.25) looks a little ugly, compared for instance
to its analogue (7.28); the point is that we don't know that
$H^d(L \cap \partial B(y,t))=0$ (as happened with Lemma 7.34).
Incidentally, if $t \leq \varepsilon r$, simply 
understand $B(y,t -\varepsilon r)$ as the empty set
(and the first half of (8.25) is trivial).
We do not really need (8.25) for the rest of this text
(we use it for (14.2), but we could also use (8.24)
and a compactness argument using Lemma 3.12 there, and the
argument would be a tiny bit simpler),
but we decided to include (8.25) for the sake of completeness
and to avoid potential unpleasant surprises later. 
We shall start with the simpler proof of (8.24) anyway.

\medskip
We now prove the theorem, by contradiction and compactness.
Suppose that we can find $\varepsilon$ such that no $r_0$
works. Then for each integer $k \geq 0$, the statement fails
with $r_0 = 2^{-k}$ and there is a reduced minimal cone $E_k$, 
a point $x_k \in K_k = E_k \cap \partial B(0,1)$,
and a radius $r_k \in (0,2^{-k})$ such that no minimal set
$L_k$ satisfies (8.24) and (8.25). We want to derive a contradiction.

Set $A=(1,0, \cdots ,0)$. 
By rotation invariance, we can assume that $x_k = A$
for every $k$, and then $P$ is always the same vertical plane
through $A$. 
Set $F_k = r_k^{-1} (E_k - A)$, so that now $\{ F_k \}$ is a sequence
of reduced minimal sets through the origin.
We can extract a subsequence, so that after extraction $\{ F_k \}$ 
converges to some closed set $F \i \R^n$. See (3.1) for the 
definition, and recall that the existence of converging subsequences
is a standard consequence of the fact that the Hausdorff distance on 
the set of compact subsets of a compact metric space makes this set
compact too. Then Lemma 4.7 (with the gauge function $h=0$)
says that
$$
F \hbox{ is a reduced minimal set of dimension $d+1$.}
\leqno (8.27)
$$

Next denote by $P_0$ the hyperplane parallel to $P$
through the origin and by $D$ the first axis (the vector
line perpendicular to $P_0$). Let us check that 
$$
F = D \times F^\sharp,
\hbox{ where } F^\sharp = F \cap P_0.
\leqno (8.28)
$$
Let $y\in F$ be given, and let us verify that
$y+D \i F$. Since the $F_k$ tend to $F$, we can find
$y_k \in F_k$ so that $\{ y_k \}$ tends to $y$.
Set $z_k = r_k y_k + A$; notice that $z_k\in E_k$
by definition of $F_k$, and that $z_k$ tends to $A$ because 
$r_k$ tends to $0$. Fix $\lambda \in \R$, and set
$v_k = (1 + \lambda r_k) z_k$. Then $v_k \in E_k$ for
$k$ large, because $E_k$ is a cone, and so
$w_k = r_k^{-1} (v_k - A)$ lies in $F_k$. But
$$\eqalign{
w_k 
&= r_k^{-1} ((1 + \lambda r_k) z_k - A)
= r_k^{-1} ((1 + \lambda r_k)(r_k y_k + A)-A)
\cr&= r_k^{-1} (r_k y_k + \lambda r_k^2 y_k+ \lambda r_k A)
= y_k + \lambda A + \lambda r_k y_k,
}\leqno (8.29)
$$
so $w_k$ tends to $y+\lambda A$  (because $y_k$ tends to $y$
and $r_k$ tends to $0$), and $y+\lambda A \in F$ because
$F$ is the limit of the $F_k$. Thus $y+D \i F$ when $y\in F$;
this proves (8.28).

We deduce from (8.28), (8.27), and Proposition 8.3
(with $m=1$) that 
$$
F^\sharp \hbox{ is a reduced minimal set of dimension $d$ in $P_0$.}
\leqno (8.30)
$$
[This is where a proof of our description of blow-up limits of 
$K$ would end.]

Set $L_k= r_k F^\sharp + A$; this is a reduced minimal set 
of dimension $d$ in $P$. Let us check that for
$k$ large, $L_k$ satisfies (8.24) with respect to
$K_k = E_k \cap \partial B(0,1)$.

First let $y \in L_k \cap B(A,r_k)$ be given,
and set  $z=r_k^{-1} (y - A)$. Thus 
$z\in F^\sharp \cap B(0,1) \i F\cap B(0,1)$.
Let $v\in F_k$ minimize the distance to $z \,$;
thus $|v-z| \leq \varepsilon_k$, where $\varepsilon_k$
measures the distance from $F_k$ to $F$ in $B(0,2)$, say,
and the sequence $\{ \varepsilon_k \}$ converges to $0$.
Set $w=r_k v + A$. Then $w\in E_k$, and $|w-y| = r_k |v-z| 
\leq r_k\varepsilon_k$. Also, $y \in P \cap B(A,r_k)$,
so $|y-1| \leq C r_k^2$ (because $P$ is tangent to $\partial B(0,1)$
at $A$), and $|w-1| \leq |w-y|+|y-1| \leq r_k \varepsilon_k + C r_k^2$.
Now ${w \over |w|} \in K_k$, and $|{w \over |w|}-y| \leq |w-1|+|w-y| 
\leq 2 r_k \varepsilon_k + C r_k^2$, which is less than 
$\varepsilon r_k$ for $k$ large. This proves the second half
of (8.24).

Now consider $y\in K_k \cap B(A,r_k)$, set $z=r_k^{-1} (y - A)$,
and observe that $z\in F_k \cap B(0,1)$ because $y\in E_k$. We can 
find $v\in F$ such that $|v-z| \leq \varepsilon_k$, with 
$\varepsilon_k$ as above. Observe that $\dist(y,P) \leq C r_k^2$
because $P$ is tangent to $\partial B(0,1)$ at $A$, so
$\dist(z,P_0) \leq C r_k$ and $\dist(v,P_0) \leq C r_k+\varepsilon_k$.
Recall from (8.28) that $F = D \times F^\sharp$, so the point
$v'$ at the intersection of $v+D$ and $P_0$ lies in $F^\sharp$,
and of course $|v'-v| = \dist(v,P_0) \leq C r_k+\varepsilon_k$.
Now $w=r_k v'+A$ lies in $L_k$, and 
$|w-y| = r_k |v'-z| \leq r_k(C r_k+2\varepsilon_k) < \varepsilon r_k$
for $k$ large. Thus (8.24) holds for $k$ large. Our proof of 
Theorem 8.23 would end here with a contradiction if we did not 
include (8.25) in the statement.

\medskip
By assumption, there is a ball $B_k=B(y_k,t_k) \i B(A,r_k)$ 
such that (8.25) fails, i.e.,
$$
H^d(L_k \cap B(y_k,t_k-\varepsilon r_k)) 
> H^d(K_k \cap B_k) +\varepsilon r_k^d 
\leqno (8.31)
$$
or
$$
H^d(L_k \cap B(y_k,t_k+\varepsilon r_k)) 
< H^d(K_k \cap B_k) -\varepsilon r_k^d. 
\leqno (8.32)
$$
We want to produce a contradiction. Set 
$$
z_k = r_k^{-1}(y_k-A) \, ,  \hskip0.4cm \rho_k = r_k^{-1}t_k \, ,
\ \hbox{ and } B'_k = B(z_k,\rho_k) = r_k^{-1}(B_k-A).
\leqno (8.33)
$$
Notice that $B'_k \i B(0,1)$ because $B_k \i B(A,r_k)$. 
Since we may replace $\{ E_k \}$ with a subsequence, 
we may assume that $\{ B'_k \}$ converges to some ball
$B(z,\rho) \i \overline B(0,1)$. That is, 
$$
\hbox{$\{ z_k \}$ converges to $z$ and $\{ \rho_k \}$ converges to 
$\rho\geq 0$.}
\leqno (8.34)
$$

We want to use the results of Section 4 to estimate the measure 
of some sets, but we need a small amount of space to move around.
So let $\tau > 0$ be very small (to be chosen later).
We shall start with upper bounds for $H^d(E_k \cap B_k)$.
Identify $D$ with the real line, set $I = [0,1] \i D$,
and apply Lemma 3.12 to the sequence $\{ F_k \}$ of $(d+1)$-minimal
sets, in the compact set 
$$
D_1 = I \times \overline B(z,\rho+\tau).
\leqno (8.35)
$$
We can take $M=1$ and $h=0$ in (3.13), and we get that
$$
\limsup_{k \to +\infty} H^{d+1}(F_k \cap  D_1)
\leq H^{d+1}(F \cap  D_1).
\leqno (8.36)
$$
By (8.27), $F$ is rectifiable, and we can apply
Lemma 8.10 (or Theorem~3.2.22 in [Fe]) 
to the product 
$F \cap  D_1 = [I \times F^\sharp \cap \overline B(z,\rho+\tau)]$,
in the same conditions as above. This yields
$$
H^{d+1}(F \cap  D_1) 
= H^1(I) \, H^d(F^\sharp \cap \overline B(z,\rho+\tau))
= H^d(F^\sharp \cap \overline B(z,\rho+\tau)).
\leqno (8.37)
$$
Let us check that for $k$ large,
$$
H^d(F^\sharp \cap \overline B(z,\rho+\tau))
\leq r_k^{-d} H^d(L_k\cap B(y_k,t_k+\varepsilon r_k)).
\leqno (8.38)
$$
Indeed $\overline B(z,\rho+\tau) \i B(z_k,\rho_k+\varepsilon)$
for $k$ large, by (8.34) and if we choose $\tau < \varepsilon$.
Set $\psi(x) = A+r_k x$ for a moment; we know that 
$L_k = \psi(F^\sharp)$, and $\psi(B(z_k,\rho_k+\varepsilon)) =
B(y_k,t_k+\varepsilon r_k)$ by (8.33). So
$$
\psi(F^\sharp \cap \overline B(z,\rho+\tau)) 
\i \psi(F^\sharp \cap B(z_k,\rho_k+\varepsilon))
= L_k \cap B(y_k,t_k+\varepsilon r_k),
\leqno (8.39)
$$
and (8.38) follows. We deduce from (8.36), (8.37), and (8.38)
that for $k$ large,
$$
H^{d+1}(F_k \cap  D_1) \leq \tau +
r_k^{-d} H^d(L_k\cap B(y_k,t_k+\varepsilon r_k)).
\leqno (8.40)
$$
We also need to compare $H^{d+1}(F_k \cap  D_1)$
with $H^d(K_k \cap B_k)$. Set
$$
E'_k = \big\{ \lambda v \, ; \, 
1+\tau r_k \leq \lambda \leq 1+r_k-\tau r_k
\hbox{ and } v \in K_k \cap B_k \big\}.
\leqno (8.41)
$$
and $E''_k = r_k^{-1}(E'_k - A)$. We claim that for $k$ large,
$$
E''_k \i F_k \cap  D_1.
\leqno (8.42)
$$
Indeed let $\lambda$ and $v$ be as in the definition of $E'_k$;
we want to show that $w=r_k^{-1}(\lambda v -A)$ lies in $F_k \cap  D_1$.
First, $\lambda v \in E_k$ because $v\in K_k \i E_k$ and $E_k$ is
a cone, so $w\in F_k = r_k^{-1}(E_k - A)$. We still need to check
that $w\in D_1$.

Write $v = A + (y_k-A) + (v-y_k)$, and also denote by
$\pi$ the orthogonal projection onto $P_0$. Observe that
$$
|\pi(y_k-A) - (y_k-A)| \leq C r_k |y_k-A| \leq C r_k^2
\leqno (8.43)
$$
because $y_k \in B(A,r_k)$ and $P$ is tangent to 
$\partial B(0,1)$ at $A$, and similarly
$$
|\pi(v-y_k) - (v-y_k)| \leq C r_k |v-y_k| \leq C r_k t_k
\leqno (8.44)
$$
because $v\in B_k=B(y_k,t_k)$ and both points lie in
$K_k \i \partial B(0,1)$. Then
$$\eqalign{
w &= r_k^{-1}(\lambda v -A)
= r_k^{-1} [(\lambda -1)v + (v-A)]
\cr&= r_k^{-1} [(\lambda -1)v + (y_k-A) + (v-y_k)]
\cr&= z_k + r_k^{-1} [(\lambda -1)v + (v-y_k)],
}\leqno (8.45)
$$
because $z_k = r_k^{-1}(y_k-A)$ by (8.33). Then
$$\eqalign{
\big|w-z_k &- r_k^{-1}[(\lambda -1)A-\pi(v-y_k)]\big|
\cr&\leq r_k^{-1}\big[(\lambda -1)|v-A| + |\pi(v-y_k) - (v-y_k)|\big]
\cr&\leq |v-A| + Ct_k \leq r_k + C t_k
}\leqno (8.46)
$$
because $1 \leq \lambda \leq 1+r_k$, by (8.44),
and because $v \in B_k \i B(A,r_k)$ by definition of 
$B_k=B(y_k,t_k)$. Set 
$w'= z + r_k^{-1}(\lambda -1)A + r_k^{-1}\pi(v-y_k)$. Then
$$
|w-w'| \leq |z-z_k| + r_k + C t_k.
\leqno (8.47)
$$

Observe that $z\in P_0$, for instance because it is the
limit of the $z_k = r_k^{-1}(y_k-A)$ and by (8.43).
Thus the projection of $w'$ on $D$ is
$r_k^{-1}(\lambda -1)A \in [\tau,1-\tau]$,
and by (8.47) the first coordinate of $w$ lies in $I$
for $k$ large. 

Similarly, $\pi(w') = z + r_k^{-1}\pi(v-y_k)$, and 
$r_k^{-1} |v-y_k| < r_k^{-1} t_k = \rho_k$ because 
$v \in B_k$ and by (8.33), so $\pi(w') \in B(z,\rho_k)$.
Since $\rho_k$ tends to $\rho$, we deduce from (8.47) that
$\pi(w) \in B(z,\rho+\tau)$ for $k$ large. 
Altogether, $w \in  D_1 = I \times \overline B(z,\rho+\tau)$
(compare with the definition (8.35)),
and we completed our proof of (8.42).

Next we evaluate $H^{d+1}(E''_k)$. Recall that
$E''_k = r_k^{-1}(E'_k - A)$, so 
$$
H^{d+1}(E''_k)=r_k^{-d-1}H^{d+1}(E'_k).
\leqno (8.48)
$$
Notice that $E'_k \i E_k$, which is rectifiable,
so we may apply Lemma~8.10 with $m=1$, $Z = E'_k$ and
to the Lipschitz function $\pi_r(x)=|x|$. 
Here (8.12) can be written as
$$
\int_{E'_k} apJ_1 \pi_r(x) \, dH^{d+1}(x)
= \int_{\lambda}\bigg\{\int_{|x|=\lambda} {\bf 1}_{E'_k}(x)
\, dH^d(x) \bigg\}\, dH^1(\lambda).
\leqno (8.49)
$$
The approximate Jacobian is $apJ_1 \pi_r(x) = 1$, because
$\pi_r$ is $1$-Lipschitz, and the tangent plane to $E_k$ at $x$ 
always contains the radial direction where the differential $D\pi_r$ 
has a norm exactly equal to $1$. Thus (8.49) and (8.41) say that
$$\eqalign{
H^{d+1}&(E'_k)
= \int_{1+\tau r_k \leq \lambda \leq 1+r_k-\tau r_k}
H^d(\lambda(K_k \cap B_k)) \, dH^1(\lambda)
\cr& =  \int_{1+\tau r_k \leq \lambda \leq 1+r_k-\tau r_k}
\lambda^d H^d(K_k \cap B_k) \, dH^1(\lambda)
\geq (r_k-2\tau r_k) H^d(K_k \cap B_k).
}\leqno (8.50)
$$
Altogether,
$$\eqalign{
(r_k-2\tau r_k) H^d(K_k \cap B_k)
&\leq H^{d+1}(E'_k) = r_k^{d+1} H^{d+1}(E''_k)
\leq r_k^{d+1} H^{d+1}(F_k\cap D_1)
\cr&\leq r_k^{d+1} \tau + r_k H^d(L_k\cap B(y_k,t_k+\varepsilon r_k))
}\leqno (8.51)
$$
for $k$ large, by (8.50), (8.48), (8.42), and (8.40).
Thus
$$\eqalign{
H^d(K_k \cap B_k) 
&\leq (1-2\tau)^{-1}[r_k^{d} \tau + 
H^d(L_k\cap B(y_k,t_k+\varepsilon r_k))]
\cr&\leq H^d(L_k \cap B(y_k,t_k+\varepsilon r_k)) 
+ \varepsilon r_k^d,
}\leqno (8.52)
$$
because $H^d(L_k \cap B(y_k,t_k+\varepsilon r_k)) \leq C r^d$
(recall from (8.30) that $L_k = A+r_k F^\sharp$ is a 
$d$-dimensional minimal set, and use Lemma 2.15),
and if $\tau$ is small enough. We just proved that (8.32) fails
for $k$ large.

\ms
Now we want to show that (8.31) fails too. The argument will be
quite similar to the one above. This time we set $I=(0,1) \i D$
and $D_2 = I \times B(z,\rho-\tau)$. [If $\rho < \tau$, then
$\rho_k < \varepsilon$ and hence $t_k < \varepsilon r_k$ for $k$ large, 
and (8.31) fails trivially.]
Apply Lemma 3.3 to the open set $D_2$ and the sequence 
$\{ F_k \}$ of $(d+1)$-minimal sets. We get that
$$
H^{d+1}(F \cap  D_2) \leq 
\liminf_{k \to +\infty} H^{d+1}(F_k \cap  D_2).
\leqno (8.53)
$$
As for (8.37), $F$ is rectifiable and we can apply
Lemma 8.10 to the product 
$F \cap  D_2 = I \times [F^\sharp \cap B(z,\rho-\tau)]$.
We get that
$$
H^{d+1}(F \cap  D_2) 
= H^d(F^\sharp \cap B(z,\rho-\tau)).
\leqno (8.54)
$$
Next we check that for $k$ large,
$$
H^d(F^\sharp \cap  B(z,\rho-\tau))
\geq r_k^{-d} H^d(L_k\cap B(y_k,t_k-\varepsilon r_k)).
\leqno (8.55)
$$
Observe that for $k$ large, $B(z,\rho-\tau)$
contains $B(z_k,\rho_k-\varepsilon)$, by (8.34).
Again set $\psi(x) = A+r_k x$ for $x\in \R^n$; thus 
$L_k = \psi(F^\sharp)$ and $\psi(B(z_k,\rho_k-\varepsilon)) =
B(y_k,t_k-\varepsilon r_k)$ by (8.33). So
$$
\psi(F^\sharp \cap \overline B(z,\rho-\tau)) 
\supset \psi(F^\sharp \cap B(z_k,\rho_k-\varepsilon))
= L_k \cap B(y_k,t_k-\varepsilon r_k),
\leqno (8.56)
$$
and (8.55) follows. Then
$$
H^{d+1}(F_k \cap  D_2) \geq H^{d+1}(F \cap  D_2) -\tau
\geq r_k^{-d} H^d(L_k\cap B(y_k,t_k-\varepsilon r_k)) -\tau
\leqno (8.57)
$$
for $k$ large, by (8.53), (8.54), and (8.55).

We again compare $H^{d+1}(F_k \cap  D_2)$ with 
$H^d(K_k \cap B_k)$. We still set 
$E''_k = r_k^{-1}(E'_k - A)$, but this time with
the larger
$$
E'_k = \big\{ \lambda v \, ; \, 
1-\tau r_k \leq \lambda \leq 1+r_k+\tau r_k
\hbox{ and } v \in K_k \cap B_k \big\}.
\leqno (8.58)
$$
We claim that then
$$
F_k \cap  D_2 \i E''_k \ \ \hbox{ for $k$ large}.
\leqno (8.59)
$$
Let $x\in F_k \cap  D_2$ be given; we need show that
$w=A+r_k x$ lies in $E'_k$. We already know that
$w\in E_k$, so $w/|w| \in K_k$, and it is enough to show 
that 
$$
w/|w| \in B_k \ \hbox{ and } \ 
1-\tau r_k \leq |w| \leq 1+r_k+\tau r_k
\leqno (8.60)
$$
Write $x=\alpha A + z + v$, where $\alpha A$ is the projection
of $x$ on $D$. Thus $0<\alpha<1$, and also
$v\in P_0 \cap B(0,\rho-\tau)$ (because $x\in D_2$ and $z\in P_0$). 
Now $w=A+r_k x = (1+\alpha r_k)A + r_k (z+v)$ so
$$
|w| = 1 + \alpha r_k + O(r_k^2)
\leqno (8.61)
$$
by Pythagoras (recall that $z+v \in P_0$). Then 
$$
w/|w| = (1 - \alpha r_k) \, w + O(r_k^2)
= A + r_k (z+v) + O(r_k^2)
\leqno (8.62)
$$ 
and 
$$
\big|w/|w|-y_k\big| 
\leq r_k|z-z_k| + r_k|v| + O(r_k^2)
\leq r_k|z-z_k| + r_k(\rho-\tau) + O(r_k^2)
\leqno (8.63)
$$
because $y_k = A + r_k z_k$ and $v\in B(0,\rho-\tau)$.
In these estimates, the error terms $O(r_k^2)$ are
less than $C r_k^2$, with $C$ independent of $x$.
If $k$ is large enough, (8.61) implies that
$1-\tau r_k \leq |w| \leq 1+r_k+\tau r_k$ (because $0 \leq \alpha \leq 1$),
and (8.63) says that $\big|w/|w|-y_k\big| < t_k =r_k \rho_k$
(because $z_k$ tends to $z$ and $\rho_k$ tends to $\rho$, by (8.34)).
This proves (8.60) and (8.59).

\smallskip
We still have that $H^{d+1}(E''_k)=r_k^{-d-1}H^{d+1}(E'_k)$
as in (8.48), and Lemma~8.10 now leads to
$$\eqalign{
H^{d+1}(E'_k)
&= \int_{1-\tau r_k \leq \lambda \leq 1+r_k+\tau r_k}
H^d(\lambda(K_k \cap B_k)) \, dH^1(\lambda)
\cr& =  \int_{1-\tau r_k \leq \lambda \leq 1+r_k+\tau r_k}
\lambda^d H^d(K_k \cap B_k) \, dH^1(\lambda)
\cr&\leq (r_k + 2\tau r_k) ( 1+r_k+\tau r_k)^d H^d(K_k \cap B_k).
}\leqno (8.64)
$$
as in (8.49) and (8.50). Thus
$$\eqalign{
H^d(K_k \cap B_k) 
&\geq (1-Cr_k-C\tau) \, r_k^{-1} H^{d+1}(E'_k)
\cr&
= (1-Cr_k-C\tau) \, r_k^{d} \, H^{d+1}(E''_k)
\cr&
\geq (1-Cr_k-C\tau) \, r_k^{d} \, H^{d+1}(F_k\cap D_2)
\cr&
\geq (1-Cr_k-C\tau) H^d(L_k\cap B(y_k,t_k-\varepsilon r_k)) 
-\tau r_k^{d} 
\cr&
\geq H^d(L_k\cap B(y_k,t_k-\varepsilon r_k)) - \varepsilon r_k^{d} 
}\leqno (8.65)
$$
by (8.59) and (8.57), because 
$H^d(L_k\cap B(y_k,t_k-\varepsilon r_k)) \leq C r_k^d$ by
Ahlfors regularity of the minimal set $L_k$ (see the proof of
(8.52)), if $\tau$ is small enough, and for $k$ large.
That is, (8.31) fails for $k$ large. 

We finally managed to show that (8.31) and (8.32) fail for $k$ 
large; this contradiction with the definition of our sequence
$\{ K_k \}$ completes the proof of Theorem~8.13.
\qed

\ms
Theorem~8.13 gives some information on the local behavior of
$K=E \cap \partial B(0,1)$ when $E$ is a reduced minimal cone
of dimension $d+1$. For instance, if we have a good control on 
the minimal sets of dimension $d$, we can try to use a generalization
of Reifenberg's topological disk theorem to get local bi-H\"older 
parameterizations of $K$. This may not be enough to control
$K$ precisely, but it could give some useful topological control.

We shall be interested in the case when $d=1$; then we'll see
in Section 10 that minimal sets of dimension 1 are lines and sets
$Y$, and a simple case of the generalized Reifenberg theorem
shows that $K$ is locally bi-H\"older equivalent to a line or a $Y$.
See Section 15 for the existence of generalized Reifenberg 
parameterizations, Section 12 for the simpler special case 
when $d=1$, and Section 16 for similar arguments in dimension 2. 

But we also want to know more about $K$ when $d=1$. In particular,
we want to check that it is composed of a finite collection of
arcs of great circles, and for this it will  be simpler to have
a Lipschitz control, rather than merely bi-H\"older. 

So we shall start with another approach, where we show that 
$K$ satisfies a weak almost-minimality property, and first prove 
reasonable properties for weak almost-minimal sets. This will
give slightly more precise results than Theorem~8.23, at least 
when $d=1$, and maybe the method can be used in slightly 
different contexts. But to be fair, one could also start from
Theorem~8.23, use Proposition 12.6 to show that $K$ is
locally bi-H\"older equivalent to a line or a $Y$, and then 
find more direct ways to show that $K$ is composed of arcs of
great circles in the regions where it looks a lot like a line.
Some arguments of this type will be used in the next paper [D3], 
to  deduce the $C^1$-regularity of almost-minimal sets of dimension
$2$ in $\R^3$ from a local biH\"older description, and they would be 
simpler for one-dimensional sets.

\bigskip 
\noindent {\bf 9. Weak almost-minimality of 
$\,\bf{K=E \cap \partial B(0,1)}$}
\medskip

Let $E$ be a reduced minimal cone of dimension $d+1$
in $\R^n$, and set $K=E \cap \partial B(0,1)$. We want
to show that $E$ is weakly almost-minimal, as in the next
definition. Later on we shall use this property to
study $E$ when $d=1$.

The following definition is a little awkward, because
it is a compromise between what we can easily prove for
$K$ and what will be needed in later arguments. We
state the definition for subsets of $\partial B(0,1)$,
but we could also work directly in $\R^n$, or in some
affine subspace, with only minor modifications.

\ms\proclaim Definition 9.1.
Let $K$ be a closed subset of $\partial B = \partial B(0,1) \i \R^n$.
We say that $K$ is a weak almost-minimal set of dimension $d$,
with gauge function $h$ if, for each choice
of $x\in \partial B$, $r\in (0,1)$, and a Lipschitz function 
$f : \partial B\to \partial B$ such that 
$$
W_f = \big\{ y\in \partial B \, ; \, f(y) \neq y \big\} \i B(x,r)
\ \hbox{ and $f(W_f) \i \overline B(x,r)$,}
\leqno (9.2)
$$ 
we have that
$$
H^d(K \cap \overline B(x,r))
\leq H^d(f(K) \cap \overline B(x,r))
+ (1+M^{d}) \, h(r) r^d, 
\leqno (9.3)
$$
whenever there is a partition of $K\cap B$ into a finite number 
of pieces where $f$ is $M$-Lipschitz.

\ms
As before, our convention is that a gauge function is a
nonnegative nondecreasing function $h$, with 
$\lim_{r \to 0} h(r) = 0$.
There are two differences with the definitions of
Section~4. The small one is that we work on $\partial B$
instead of an open set in $\R^n$. This is just because
our examples will live on the sphere, but the same definition
would work in $\R^n$ or in some other smooth manifolds.

The main difference is that here we pay an additional price 
that depends on the Lipschitz properties of $f$. 
The story about the partition of $K\cap B$ is a mostly
a way to avoid discussing local Lipschitz constants
or gradients. Note that $M$ may be significantly smaller than 
the Lipschitz norm of $f$, or even of its restriction to $K \cap B$; 
we shall use this advantage in some proofs.
 
We would have a slightly weaker notion of weak almost-minimal set
if we replaced $M^d$ in (9.3) with some larger functions of $M$;
I don't think that this would make much difference in the
properties of weak almost-minimal set that will be proved below.
Also $M^{d+1}$ would be a little easier to get in the proof
of Proposition 9.4. On the other hand $M^d$ looks 
more natural, so we shall keep it in the definition.

\ms\proclaim Proposition 9.4.
Let $E$ be a reduced minimal cone of dimension $d+1$
in $\R^n$, and set $K=E \cap \partial B(0,1)$. Then 
$K$ is a reduced weak almost-minimal set of dimension $d$,
with gauge function $h(r) = Cr$.

\ms
Let us first check that $K$ is reduced. Recall from Definition 2.12 
that this means that 
$$
H^d(K \cap B(x,r)) >0
\leqno (9.5)
$$
for $x\in K$ and $r \in (0,1)$. Let us check this brutally. 
We apply Lemma 8.10 (i.e., Theorem 3.2.22 in [Fe]) 
with $Z = E\cap B(x,r/2)$, $m=1$, $A=Z$, and for $\pi$ the
radial projection defined by $\pi(x) = |x|$.
The rectifiability of $E$ was discussed in Section 2,
in a more general context. So (8.12) yields
$$
\int_{Z} apJ_1\pi(x) dH^{d+1}(x) = 
\int_{\rho} \bigg\{ 
\int_{Z \cap \pi^{-1}(\rho)} dH^d \bigg\} \, d\rho.
\leqno (9.6)
$$
The approximate Jacobian $apJ_1\pi(x)$ is easily computed to be $1$,
because $\pi$ is $1$-Lipschitz and the tangent plane to $E$ always 
contains the radial direction where the Lipschitz constant
for $\pi$ is exactly $1$. Thus (9.6) says that 
$$\eqalign{
H^{d+1}(E\cap B(x,r/2)) &= H^{d+1}(Z) =
\int_{\rho} H^d(Z \cap \pi^{-1}(\rho)) \, d\rho
\cr&= \int_\rho
H^d(E\cap B(x,r/2) \cap \partial B(0,\rho)) \, d\rho.
}\leqno (9.7)
$$
In addition, the only contribution in the right-hand side
comes from $1-r/2 < \rho < 1+r/2$, for which
$E\cap B(x,r/2) \cap \partial B(0,\rho)$ is contained
in $\rho \cdot (K \cap B(x,r))$ because $E$ is a cone
and $K = E \cap \partial B(0,1)$. Thus
$$
H^{d+1}(E\cap B(x,r/2)) \leq 
\int_{1-r/2 < \rho < 1+r/2} H^d(K \cap B(x,r)) \, \rho^d d\rho.
\leqno (9.8)
$$
The left-hand side of (9.8) is positive because $E$ is
reduced, hence $H^d(K \cap B(x,r)) >0$, (9.5) holds,
and $K$ is reduced.

\smallskip
We now prove that $K$ is weakly almost-minimal. Let 
$x\in \partial B$, $r\in (0,1)$, and 
$f : \partial B\to \partial B$ be as in Definition 9.1,
and let us use $f$ to define a competitor for $E$.

Let $\psi : [0,+\infty) \to [0,1]$ be such that 
$$\eqalign{
\psi(t)&=0 \hskip0.4cm \hbox{ for } 0 \leq t \leq {1 \over 2},
\cr
\psi(t)&=1 \hskip0.4cm \hbox{ for } {1 \over 2} + {r \over 10}
\leq t \leq 1 - {r \over 10},
\cr
\psi(t)&=0 \hskip0.4cm \hbox{ for } t \geq 1,
}\leqno (9.9)
$$
$\psi$ is continuous, and $\psi$ is linear in the two 
missing interval.

Write every point of $\R^n$ as $z=\rho y$, with $\rho = |z|$
and $y = y/|y|$ (the choice of $y$ when $z=0$ does not matter),
and set 
$$
\varphi(z) = \psi(\rho) \rho f(y) + (1-\psi(\rho)) \rho y
\ \hbox{ for } z \in \R^n.
\leqno (9.10)
$$
Since $E$ is minimal and $\varphi$ is Lipschitz, 
Definition 8.1 yields
$$
H^{d+1}(E\setminus \varphi(E)) 
\leq H^{d+1}(\varphi(E)\setminus E).
\leqno (9.11)
$$
Let us compute what this means in terms of $K$
and $f(K)$. Observe that
$$
\varphi(z) = z \hbox{ for } z \in \R^n \setminus B(0,1)
\ \ \hbox{ and } \ \ 
\varphi(B(0,1)) \i B(0,1),
\leqno (9.12)
$$
so $\varphi(E)$ coincides with $E$ out of $B(0,1)$
and there will be no contribution to (9.11) coming
from $\R^n \setminus B(0,1)$. The most interesting 
contributions will come from the annulus 
$$
A = \big\{ z \in \R^n \, ; \, {1 \over 2} 
+ {r \over 10} \leq |z| \leq 1 - {r \over 10} \big\}.
\leqno (9.13)
$$
Let us again apply Lemma 8.10 to the mapping 
$\pi : z \to |z|$, from $Z = E \cap A$
to $I=[{1 \over 2}+{r \over 10},1 - {r \over 10}]$.
We get (9.6) as before, again with $apJ_1\pi(x)=1$, so
$$\eqalign{
H^{d+1}(E \cap A) &= H^{d+1}(Z)
= \int_I H^d(Z \cap \pi^{-1}(t)) \, dt
= \int_I H^d(E \cap \partial B(0,t)) \, dt
\cr&= \int_I H^d(t \cdot K) \, dt
= \int_I t^d \, H^d(K) \,  dt
= \alpha H^{d}(K),
}\leqno (9.14)
$$
because $\pi^{-1}(t) = \partial B(0,\rho) \i A$ 
for $t\in I$, and $E$ is a cone. The precise value of 
$\alpha = \int_I t^d dt$ won't matter much.

Observe that $\psi(|z|) = 1$ when $z\in A$, so
$\varphi(z) = \rho f(y)$. Thus
$$
\varphi(E\cap A)\cap \partial B(0,\rho) = \rho f(K)
\ \hbox{for $\rho \in I$,}
\leqno (9.15)
$$ 
and it is empty otherwise. The same computation as for $E \cap A$, 
but with the other rectifiable set $\varphi(E\cap A)$, yields
$$
H^{d+1}(\varphi(E\cap A))
= \int_I H^d(\varphi(E\cap A) \cap \partial B(0,t)) \, dt
= \alpha H^{d}(f(K)).
\leqno (9.16)
$$

We shall also need to estimate the contributions of the
two annuli
$$
A_1 = \big\{ z \in \R^n \, ; \, {1 \over 2} 
\leq |z| \leq  {1 \over 2} + {r \over 10} \big\} 
\ \hbox{ and } \ 
A_2 = \big\{ z \in \R^n \, ; \, 1 - {r \over 10} 
\leq |z| \leq  1 \big\}.
\leqno (9.17)
$$

Set $Z_i = \big\{ z\in E\cap A_i \, ; \, 
\varphi(z) \neq z \big\}$.
We want to estimate the Lipschitz properties of $\varphi$ on $Z_i$,
and these will depend on the Lipschitz properties of $f$
on $K \cap B(x,r)$. Let us assume, as in Definition 9.1,
that we have a partition $K \cap B(x,r) = \cup_j K_j$
of $K \cap B(x,r)$ into a finitely number of sets $K_j$ where
$f$ is $M$-Lipschitz. Note that we may assume that the $K_j$
are Borel sets, because $f$ is automatically $M$-Lipschitz
on the $\overline K_j$, which enables us to first replace $K_j$
with $\overline K_j$, and then build a Borel partition out
of the covering. If $z\in Z_i$, $y=z/|z|$ lies in $K$, because
$z\in E$ and $E$ is a cone, and $f(y) \neq y$ because 
$\varphi(z) \neq z$ (see (9.10)). Then $y \in K \cap B(x,r)$, 
by (9.2). Thus the sets
$$
Z_{i,j} = \big\{ z\in Z_i \, ; \, y = z/|z| \in K_j\big\}
\leqno (9.18)
$$
form a partition of $Z_i$. Let us check that
$$
\varphi \hbox{ is $(25+2M)$-Lipschitz on each $Z_{i,j}$.}
\leqno (9.19)
$$
It will be enough to control 
$g(z) = \varphi(z)-z = \rho \psi(\rho)(f(y)-y)$. 
We write that
$$
g(z_1)-g(z_2) = 
\rho_1\psi(\rho_1)(f(y_1)-y_1)
-\rho_2\psi(\rho_2)(f(y_2)-y_2) = \Delta_1 - \Delta_2 \, ,
\leqno (9.20)
$$
with $\Delta_1 = [\rho_1\psi(\rho_1)-\rho_2\psi(\rho_2)]
(f(y_1)-y_1)$ and $\Delta_2 = \rho_2\psi(\rho_2) 
[ (f(y_2)-y_2)-(f(y_1)-y_1)]$.

Since $|\nabla(\rho \psi(\rho))| \leq \psi(\rho) 
+ \rho |\psi'(\rho)| \leq 1+{10 \over r}$, we get that
$$
|\Delta_1| \leq {11 \over r}|\rho_2-\rho_1||f(y_1)-y_1| 
\leq 22|\rho_2-\rho_1|, 
\leqno (9.21)
$$
because $y_1$ and $f(y_1)$ lie in $B(x,r)$, 
by (9.2). Also, $z_1$ and $z_2$ both lie in $Z_{i,j}$, so 
$y_1$ and $y_2$ lie in $K_j$. Then
$$
|\Delta_2| \leq |(f(y_2)-y_2)-(f(y_1)-y_1)| 
\leq |y_2-y_1| + |f(y_2)-f(y_1)| \leq (1+M)|y_2-y_1|,
\leqno (9.22)
$$
by definition of $M$ and the $K_j$. Altogether 
$$
|g(z_1)-g(z_2)| \leq |\Delta_1| + |\Delta_2|
\leq 22|\rho_2-\rho_1| + (1+M)|y_2-y_1|
\leq (24+2M)|z_2-z_1|
\leqno (9.23)
$$
and (9.19) follows because $\varphi(z) = z + g(z)$.

An immediate consequence of (9.19) is that
$$
H^{d+1}(\varphi(Z_{i,j})) \leq C (1+M)^{d+1} H^{d+1}(Z_{i,j}).
\leqno (9.24)
$$
This would be enough to prove (9.3) with $M^d$ replaced by 
$M^{d+1}$, but since we decided to go for the more restrictive
definition, we shall need to improve (9.24) a little bit.
Observe that $Z_{i,j}$ is rectifiable, because it is contained
in $E$ and by the rectifiability results of Section 2, so
we can use the area formula (see for instance Corollary 3.2.20
in [Fe]). That is, 
$$
H^{d+1}(\varphi(Z_{i,j})) 
\leq \int_{Z_{i,j}} J_\varphi(z) \, dH^{d+1}(z).
\leqno (9.25)
$$
Our improvement on (9.24) will come from a
better control of the $(d+1)$-dimensional Jacobian 
$J_\varphi$ due to a better control on
the derivative of $\varphi$ in the radial direction.

It may comfort the reader to notice that since $E$ is 
Ahlfors-regular and rectifiable, it even has a true tangent 
$(d+1)$-plane at $z$ for $H^{d+1}$-almost every $z$.
This plane $z+P(z)$ contains the radial direction because
$E$ is a cone, and if $z+Q(z)$ denotes the $d$-plane
in $P(z)$ perpendicular to the radial direction, it
is easy to see that $K$ has a tangent plane at
$y=z/|z|$, which is $y+Q(z)$. Note that we already knew
from the first part of Lemma~8.10 (applied as before
to the radial projection $\pi$) that $K$ is rectifiable.

The derivative of $g(z) = \rho \psi(\rho)(f(y)-y)$
in the radial direction is bounded by 
$$
[\psi(\rho)+\rho\psi'(\rho)]|f(y)-y| \leq C r^{-1}|f(y)-y| \leq C,
\leqno (9.26)
$$
by (9.2). The derivative in all the tangential directions
is bounded by $25+2M$, because we restrict our attention to
$Z_{i,j}$ and by (9.19). So $J_\varphi(z) \leq C (1+M)^d$.
Indeed, after some manipulations, we can compute $J_\varphi(z)$ 
as a determinant of $D\varphi$ in some set of bases, and if we start 
the basis of $P(z)$ with a radial vector, the first row of the corresponding
matrix has bounded terms, while the terms in the other rows 
are bounded by $C(1+M)$. Thus (9.25) yields
$$
H^{d+1}(\varphi(Z_{i,j})) \leq C (1+M)^d H^{d+1}(Z_{i,j}).
\leqno (9.27)
$$

Incidentally, we have been trying to avoid talking about
approximate differentials of $f$ and $\varphi$ on $K$
and $E$, but the natural way to deal with (9.25) would use these.
Our discussion leads to $J_\varphi(z) \leq C |D_E\varphi(z)|^d$,
where $|D_E\varphi(z)|$ is the norm of the 
approximate differential of the restriction of $\varphi$ to $E$
(which exits $H^{d+1}$-almost everywhere on $E$).
Then the proof of (9.19) yields $|D_E\varphi(z)| \leq C + |D_Kf(y)|$, 
where $|D_Kf(y)|$ is the norm of the approximate differential of 
the restriction of $f$ to $K$. Thus a more precise version of (9.27)
is that
$$
H^{d+1}(\varphi(Z_i)) 
\leq C \int_{Z_i} (1+|D_K\varphi(y)|^d) \, dH^{d+1}(z).
\leqno (9.28)
$$
Notice that we don't need to cut $K\cap B(x,r)$ or the $Z_i$
into pieces in this context, because $|D_K\varphi(y)|$ is local.

We return to (9.27) and estimate $H^{d+1}(Z_{i,j})$.
Recall from (9.18) that $Z_{i,j}$ is contained in the cone
over $K_j \i K \cap B(x,r)$. It is also contained in 
the thin annulus $A_i$ of (9.17), so $\cup_j Z_{i,j}$ is 
contained in a ball $B_i$ of radius $2r$. Then
$$
\sum_j H^{d+1}(Z_{i,j})
= H^{d+1}\big(\cup_j Z_{i,j}\big)
\leq H^{d+1}(E\cap B_i) \leq C r^{d+1}
\leqno (9.29)
$$
because the $Z_{i,j}$ are disjoint and contained in $E$,
and by the easiest part of Lemma 2.15 (the Ahlfors-regularity 
of the minimal set $E$).

We are almost ready to return to (9.11). Set
$$
S = E\cap B(0,1/2) \cup 
\big[\cup_{i\in \{1,2\}} E \cap A_i \setminus Z_i)\big],
\leqno (9.30)
$$
and check that $\varphi(z)=z$ on $S$. When $z\in B(0,1/2)$,
this comes directly from (9.10) and the fact that $\psi(\rho)=0$
by (9.9). On $E \cap A_i \setminus Z_i$, this is the definition
of $Z_i$ (just below (9.17)). Next
$$\eqalign{
\varphi(E) \cap B(0,1) &= \varphi(E\cap B(0,1)) 
\i \varphi(S) \cup \varphi(E\cap A) 
\cup \big[\cup_i \varphi(Z_i)\big]
\cr& =S \cup \varphi(E\cap A) 
\cup \big[\cup_i \varphi(Z_i)\big]
}\leqno (9.31)
$$
by (9.12), because $E\cap B(0,1) = S \cup (E\cap A) \cup
[\cup_i Z_i]$, and because $\varphi(z)=z$ on $S$. So
$$
H^{d+1}(\varphi(E)\cap B(0,1))
\leq H^{d+1}(S) + H^{d+1}(\varphi(E\cap A))
+ C (1+M)^d r^{d+1}
\leqno (9.32)
$$
by (9.27) and (9.29). Similarly,
$$
H^{d+1}(E\cap B(0,1)) \geq H^{d+1}(S) + H^{d+1}(E\cap A)
\leqno (9.33)
$$
because these sets are disjoint. On the other hand, 
$H^{d+1}(E\setminus \varphi(E)) \leq H^{d+1}(\varphi(E)\setminus E)$
by (9.11), and since $\varphi(E)$ coincides with $E$
out of $B(0,1)$ by (9.12), we get that
$$
H^{d+1}(E\cap B(0,1)) \leq H^{d+1}(\varphi(E)\cap B(0,1))
\leqno (9.34)
$$
by adding $H^{d+1}(E\cap \varphi(E) \cap B(0,1))$ to both sides.
Altogether,
$$\eqalign{
H^{d+1}(E\cap A) &\leq H^{d+1}(E\cap B(0,1)) - H^{d+1}(S) 
\cr&\leq H^{d+1}(\varphi(E)\cap B(0,1)) - H^{d+1}(S) 
\cr&\leq H^{d+1}(\varphi(E\cap A)) + C (1+M)^d r^{d+1}
}\leqno (9.35)
$$
by (9.33), (9.34), and (9.32). By (9.14) and (9.16),
this implies that
$$
\alpha H^{d}(K) \leq \alpha H^{d}(f(K)) + C (1+M)^d r^{d+1},
\leqno (9.36)
$$
with $\displaystyle\alpha = \int_I t^d dt \geq \int_{6/10}^{9/10} t^d dt$
(recall that $I=[{1 \over 2}+{r \over 10},1 - {r \over 10}]$,
see above (9.14)). This is just the same thing as (9.3),
with $h(r)=Cr$. Proposition 9.4 follows.
\qed

\ms\noindent {\bf Remark 9.37.}
We can improve slightly our estimate of $H^{d+1}(\varphi(Z_i))$,
by using (9.28) instead of (9.27). The same computations as above then 
lead to
$$
H^{d}(K) \leq H^{d}(f(K)) 
+ C \int_{Z_1 \cup Z_2} (1+|D_K\varphi(y)|^d) \, dH^{d+1}(z),
\leqno (9.38)
$$
where we still use the convention that $y=z/|z|$.
Recall from just above (9.18) that if $z\in Z_1 \cup Z_2$,
then $f(y) \neq y$. In particular, 
$y\in K \cap B(x,r)$, by (9.2). Also, $z\in A_1 \cup A_2$.
We can apply Lemma 8.10 to $Z=Z_1 \cup Z_2$ (which is 
rectifiable because it is contained in $E$) and the radial
projection $\pi$ from $Z$ to 
$J=[{1 \over 2},{1 \over 2} + {r \over 10}] \cup
[1-{r \over 10},1]$. This time we use sets $A \i Z$ to show that
both sides of (8.12) define equal measures, and we get that
$$\eqalign{
\int_{Z} (1+|D_K\varphi(y)|^d) &\, dH^{d+1}(z)
\cr&= \int_{\rho \in J} \bigg\{\int_{z\in Z\cap \pi^{-1}(\rho)}
(1+|D_K\varphi(y)|^d) \, dH^d(z) \bigg\}\, d\rho,
}\leqno (9.39)
$$
where we omitted the approximate Jacobian from the left-hand side
because it is equal to $1$, as for (9.6) and (9.14).
Recall that $\pi^{-1}(\rho) = \partial B(0,\rho)$
and that $Z$ is contained in the cone over 
$W_f = \{ y \in \partial B(0,1) \, ; \, 
f(y) \neq y \} \i K \cap B(x,r)$. Thus
$$\eqalign{
\int_{z\in Z\cap \pi^{-1}(\rho)} (1+|D_K\varphi(y)|^d) \, dH^d(z)
&\leq \int_{\rho W_f} (1+|D_K\varphi(y)|^d) \, dH^d(z)
\cr&= \rho^d \int_{W_f} (1+|D_K\varphi(y)|^d) \, dH^d(y).
}\leqno (9.40)
$$
We integrate this over $J$ and get that
$$
H^{d}(K) \leq H^{d}(f(K)) 
+ C r \int_{W_f} (1+|D_K\varphi(y)|^d) \, dH^d(y),
\leqno (9.41)
$$
by (9.38) and (9.39). This is slightly more precise and elegant than 
(9.36), but will not be used here.

\bigskip
\noindent {\bf 10. One-dimensional minimal sets in $\R^n$}
\medskip

We want to describe here the reduced minimal sets of 
dimension $1$ in $\R^n$. There are two simple examples,
the lines and the $Y$-sets. In this section, a $Y$-set will be
a union of three half lines contained in an affine plane,
with the same endpoint, and that make angles of $120^\circ$
at that point.

It is fairly easy to see that lines and $Y$-sets are indeed
minimal. See [D2], Lemma~59.1 on page 392, 
but here is a rapid sketch of the argument for the
typical case when $Y$ is a $Y$-set centered at the
origin. Let $Z = \varphi(Y)$ be a deformation of $Y$,
choose $R$ so large that $Y=Z$ out of $B(0,R)$, and denote by 
$a_1, a_2, a_3$ the three points of $Y \cap \partial B(0,R)$.
The images under $\varphi$ of the three segments of 
$Y \cap \overline B(0,R)$ are three
arcs that connect $\varphi(0)$ to the $a_j$. Then (if
$H^{1}(Z\cap B(0,R)) < +\infty$) we can find
a point $z\in Z \cap \overline B(0,R)$ and three simple arcs
in $\overline B(0,R) \cap Z$ that go from $z$ to the $a_j$,
and whose interiors are disjoint. Replacing these by line segments
makes $Z$ shorter, and a variational argument shows that
the optimal case is when $z=0$ and $Z=Y$. 

To the mild surprise of the author, there is no other
one-dimensional reduced minimal set, even in large
ambient dimensions.

\ms\proclaim Theorem 10.1.
Every nonempty reduced minimal set of dimension $1$
in $\R^n$ is a line or a $Y$-set.

\ms
The rest of this section is devoted to the proof.
Even though a direct approach would be possible,
it will be easier for us to apply some general arguments,
and in particular use blow-down limits. 
We start with the simpler case of minimal cones.

\ms\proclaim Lemma 10.2. Every nonempty reduced minimal 
cone of dimension $1$ is a line or a $Y$-set.

\ms
Let $E$ be a reduced minimal cone of dimension $1$ 
(centered at the origin), and set $K = E \cap \partial B(0,1)$. 
Notice that $K$ is finite, because otherwise 
$H^{1}(E\cap B(0,1))$ is infinite. Write 
$K = \{ a_1, \cdots, a_m \}$, and denote by $L_j$ the 
half line through $a_j$. Let us first show that
$$
{\rm Angle}(L_i,L_j) \geq {2\pi \over 3}
\ \hbox{ for } i\neq j.
\leqno (10.3)
$$

Choose $i \neq j$ so that $\theta={\rm Angle}(L_i,L_j)$ 
is minimal, and suppose that $\theta < {2\pi \over 3}$.
We want to get a contradiction by replacing 
$[L_i \cup L_j]\cap \overline B(0,1)$ 
with a shorter fork $F_t$ defined as follows. 
Set $z_t = t(a_i+a_j)$ for $0 < t < 10^{-1}$, say, and 
$F_t = [0,z_t]\cup [z_t,a_i] \cup [z_t,a_j]$.
A simple differentiation of $\ell(t) = H^1(F_t)$ at
$t=0$ shows that 
$H^1(F_t) < H^1([L_i \cup L_j]\cap \overline B(0,1))$
for $t$ small (because $\theta < {2\pi \over 3}$). 
We choose $t$ so that this holds, and set
$$
\widetilde E = \big[ E \setminus \big([L_i \cup L_j]
\cap \overline B(0,1)\big)\big] \cup F_t. 
\leqno (10.4)
$$
Notice that 
$H^{1}(\widetilde E \cap B(0,1)) < H^{1}(E \cap B(0,1))$
and the two sets coincide out of $B(0,1)$, so 
$H^1(\widetilde E \setminus E) < H^1(E \setminus \widetilde E)$.

We want to construct a Lipschitz mapping 
$\varphi : \R^n \to  \R^n$ such that $\varphi(x)=x$ 
out of $B(0,1)$ and $\varphi(E) \i \widetilde E$. 
Since $E$ is minimal, Definition 8.1 will allow us to 
apply (8.2) to $\varphi$, and get that 
$H^1(E \setminus \varphi(E)) \leq H^1(\varphi(E) \setminus E)$,
hence $H^1(E \setminus \widetilde E) 
\leq H^1(\widetilde E \setminus E)$. This is
just the opposite of what we have, so we shall get the
desired contradiction as soon as we have $\varphi$.

Let $L$ denote the half line through $a_i+a_j$, and
set ${\cal C}_\alpha = \big\{ y\in \R^n \, ; \,
\dist(y,L) \leq \alpha |y| \big\}$, with
$\alpha= \dist(a_i,L)$. That is, ${\cal C}_\alpha$
is the smallest revolution cone around $L$ that
contains $L_i$ and $L_j$. Notice that since $|a_i-a_j|$ 
is minimal, $\dist(a_l,L_i) \geq \dist(a_j,L_i)$
for $l \neq i$.
Also, $a_i$ and $a_j$ are furthest possible in 
${\cal C}_\alpha \cap \partial B(0,1)$, so
the only point of ${\cal C}_\alpha \cap \partial B(0,1)$
such that $\dist(y,L_i) \geq \dist(a_j,L_i)$ is $y=a_j$. 
Put in other words, $a_l \notin {\cal C}_\alpha$ unless $l=j$.
Let us choose $\alpha'>\alpha$ so close to $\alpha$
that no $a_l$, $l \neq i,j$,  lies in ${\cal C}_{\alpha'}$
either. 

Set $\varphi(y) = y$ for 
$y\in \R^n\setminus [{\cal C}_{\alpha'}\cap B(0,1)]$.
Thus we still need to define $\varphi$ on 
${\cal C}_{\alpha'} \cap B(0,1)$, and we start 
with $(L_i \cup L_j) \cap B(0,1)$.
Define $\varphi$ on $L_i \cap \overline B(0,1)$ so that
it is Lipschitz and maps $L_i \cap \overline B(0,1)$ into 
$F_t \cap \overline B(0,1)$, with $\varphi(0)=0$ and 
$\varphi(x)=x$ near $a_i$. Proceed similarly with $L_j$.
Because there are small security cones around $L_i$ and $L_j$ 
that are contained in ${\cal C}_{\alpha'}$, our map $\varphi$
is Lipschitz on 
$\big[\R^n\setminus [{\cal C}_{\alpha'}\cap B(0,1)]\big] \cup L_i \cup  L_j$.
We simply use the Whitney extension Theorem
to extend it to the rest of $\R^n$, so that it stays Lipschitz.

Notice that $\varphi$ does not move the half lines $L_l$,
$l \neq i,j$; then $\varphi(E) \i \widetilde E$, where $\widetilde E$
is as in (10.4). This completes our construction of $\varphi$,
and we have noted earlier that the ensuing contradiction proves (10.3).

\smallskip
Next we claim that if $j,k,l$ are distinct, the three
lines $L_j, L_k, L_l$ lie in a plane, and $L_j \cup L_k \cup L_l$
is a $Y$-set. Denote by $P$ the plane that contains $L_j$ and $L_k$,
and identify $P$ with $\Bbb C$. By symmetry, we may assume that
$a_j = e^{i\theta}$ and $a_k = e^{-i\theta}$, with
$0 < \theta \leq \pi/2$. By (10.3), $\theta \geq \pi/3$.
Denote by $\xi$ the orthogonal projection of $a_l$ on $P$; then
(10.3) says that $\langle \xi,a_j \rangle = \langle a_l,a_j \rangle
\leq \cos(2\pi/3) = -1/2$, and similarly $\langle \xi,a_j \rangle
\leq -1/2$. Since $\xi$ lies in the closed unit ball, this forces the
half line through $\xi$ to make angles at least $2\pi/3$ with
$L_i$ and $L_j$. The only possibility is that $\theta = \pi/3$, 
$\xi =-1$ and $a_l=-1$. Our claim follows.

\smallskip
We are now ready to prove the lemma.
If $E \cap \partial B(0,1)$ is empty, then $E$ is empty.
Next $E \cap \partial B(0,1)$ cannot be reduced to one point, 
because a half line is not
minimal (we can contract it along its end to make it shorter).
If it has exactly two points, it is easy to check that 
$E$ is a line (otherwise, deform $E\cap B(0,1)$ into
a segment with the same endpoints, and you'll save some length).
If it has more than two points, choose two of them
(say, $a_i$ and $a_j$); the argument above shows that all the 
other points $a_k$ lie just opposite to $a_i+a_j$, so there is 
only one of them, and $E$ is a $Y$-set. This proves Lemma 10.2.
\qed

\ms
Now we want to  deal with the general case of Theorem 10.1.
We shall use results of the previous parts to save time, but
more direct proofs would be possible (and also the theorems
that we use are easier for one-dimensional sets).

So let $E$ be a nonempty reduced minimal set of dimension $1$.
Let us assume for simplicity that the origin lies in $E$, and
set $\theta(r) = r^{-1} H^{1}(E \cap B(0,r))$. By Lemma 2.15,
$\theta$ is bounded. By Proposition 5.16, $\theta$ is nondecreasing,
so $\theta_\infty = \lim_{r \to +\infty} \theta(r)$ exists and
is finite. 

Let $\tau > 0$ be small, and let $\varepsilon > 0$ be as
in Proposition 7.1, with $a=1$ and $b=2$. 
Then let $R \geq 1$ be such that
$\theta(r) \geq \theta_\infty -\varepsilon$ for $r \geq R$.
For $r \geq R$, apply Proposition~7.1 to $r^{-1} E$; the density
assumption is satisfied because $\theta_\infty -\varepsilon
\leq \theta(r) \leq \theta(2r) \leq \theta_\infty$. The conclusion
is that there is a minimal cone $\C$ centered at the origin,
that is close to $E'=r^{-1} E$ in various respects, but in particular
such that
$$
|\rho^{-1} H^{1}(E' \cap B(0,\rho)) - \rho^{-1} H^{1}(\C\cap B(0,\rho))| 
\leq \tau
\ \hbox{ for } 1+\tau \leq \rho \leq 2-\tau,
\leqno (10.5)
$$
as in (7.5). Take $\rho = 3/2$; we get that
$$
\theta(3r/2) = \rho^{-1} H^{1}(E' \cap B(0,\rho))
\leq \rho^{-1} H^{1}(\C\cap B(0,\rho)) + \tau
\leq 3+\tau
\leqno (10.6)
$$ 
by Lemma 10.2. Denote by $N(t)$ the number of points in
$E\cap \partial B(0,t)$.
By [D3], Lemma~26.1 on page 160, 
or [Fe], Theorem 3.2.22, 
$$
\int_{0}^{3r/2} N(x) dx \leq H^{1}(E\cap B(0,3r/2)) 
= (3r/2) \theta(3r/2) \leq 3(3 +\tau) r/2
\leqno (10.7)
$$
so we can find $\rho\in (r/10,3r/2)$ such that 
$N(\rho) \leq {10 \over 14 r} \, {3(3+\tau) r \over 2} < 4$
(if $\tau$ is small enough). That is, 
$$
E\cap \partial B(0,\rho) \hbox{ has at most three points.}
\leqno (10.8)
$$

\ms\proclaim
Lemma 10.9. If $E$ is a reduced minimal set that 
contains the origin, and $\rho > 0$ is such that (10.8)
holds, then $E$ coincides with a line or a $Y$-set in $B(0,\rho)$.

\ms
Theorem 10.1 will follow as soon as we prove this lemma. 
Indeed, we found arbitrarily large radii $\rho$ such that
(10.8) holds, hence such that $E$ coincides with a line or 
a $Y$-set in $B(0,\rho)$. It is easy to see that $E$ is then
a line or a $Y$-set.

Our proof of Lemma 10.9 will take some time, in particular
because we want to be able to use the same construction
with weak almost-minimal sets. 

Before we really start, let us construct retractions $h_F$ 
onto some simple sets $F$.
Our first case is when $F = Y \cap B$, where $Y$ is a $Y$-set
and $B$ is a closed ball that contains the center of $Y$.
Let $P$ denote the plane that contains $Y$. We first
construct a retraction $h_1 : P \to Y$. Denote by 
$L_1, L_2, L_3$ the three half lines that compose $Y$,
and by $V_i$ the connected component of $P\setminus Y$
opposite to $L_i$. For $y \in V_i$, $h_1(y)$ is the projection 
of $y$ onto $L_j \cup L_k$ along the direction of $L_i$,
where $L_j$ and $L_k$ are the two other half lines.
In other words, $h_1(y)$ is the intersection of $L_j \cup L_k$
with the line through $y$ parallel to $L_i$. It is easy
to see that $h_1$ is $2$-Lipschitz on $V_i$. For instance,
observe that we can split $V_i$ into two equal sectors where
$h_1$ is the oblique projection onto the line $L_j$ or $L_k$,
along a direction that makes an angle $\pi/3$ with that line.
Then $h_1$ is also $2$-Lipschitz on $P$,
because the various definitions coincide with the 
identity on $Y$.

We compose $h_1$ with the retraction $h_2$
from $Y$ to $F$, which maps $y\in F$ to itself,
and maps $y\in Y \setminus F$ to the point of $Y\cap \partial B$
that lies on the same branch of $Y$ as $y$. Obviously
$h_2$ is $1$-Lipschitz.
Finally set 
$$
h_F = h_2 \circ h_1 \circ \pi_P
\leqno (10.10)
$$
on $\R^n$, where $\pi_P$ denotes the orthogonal projection 
onto $P$. Then 
$$
\hbox{$h_F$ is a $2$-Lipschitz retraction of $\R^n$ onto $F$.}
\leqno (10.11)
$$

We also need to construct $h_F$ when 
$F = B \cap (L_1 \cup L_2)$, where $L_1$ and $L_2$ are 
half lines with the same endpoint and that make an angle 
at least $2\pi/3$ at that point, and $B$ is a closed
ball that contains that point. Let us first define a
retraction $h_1$ on the plane $P$ that contains $Z$
(if $L_1\cup L_2$ is a line $L$, choose any plane that contains 
$L$; anyway $h_1 \circ \pi_P$ will end up being the orthogonal 
projection onto $L$). Call $D$ the line in $P$ that 
goes through the common endpoint of $L_1$ and $L_2$
and lies at equal distance from $L_1$ and $L_2$ and, 
for $y\in P$, let $h_1(y)$ be the intersection of 
$L_1\cup L_2$ with the line through $y$ parallel to $D$.
Now $h_1$ is still $2$-Lipschitz on $P$, because
$D$ makes an angle of $\pi/3$ or more with the lines
that contain $L_1$ and $L_2$. We then define $h_F$
as before, by (10.10), and (10.11) holds for the same 
reasons.

\smallskip
Let $E$ be as in Lemma 10.9. We want to treat a 
few cases at the same time, so some notation will
be useful. Set $B=\overline B(0,\rho)$, and denote by $a_j$, 
$1 \leq j \leq N$, the points of $E \cap \partial B$. 
Let us still keep the number $N$ of points in $E \cap \partial B$
open, except that (10.8) says that $N\leq3$.

Then denote by $E_j$ the connected component of 
$E\cap B$ that contains $a_j$. Call $M$ the number of components 
that we get, and renumber the $a_j$ so that the $E_j$, $1 \leq j\leq M$, 
are disjoint, and each $a_k$, $k > M$, lies in some $E_j$ such
that $j\leq M$. Notice that the $E_j$ are compact, so
$$
\dist(E_i,E_j) > 0
\ \hbox{ for } 1 \leq i,j \leq M, i \neq j.
\leqno (10.12)
$$
When $M >1$, we shall also need to dispatch the rest
of $E\cap B$ among the $E_j$, as follows.

\ms \proclaim Lemma 10.13.
We can find a partition of $E\cap B$ into sets $E'_j$, $1 \leq j \leq M$, 
such that $E_j \i E'_j$ for $j \leq M$ and
$$
\dist(E'_i,E'_j) > 0
\ \hbox{ for } i,j \leq M, i \neq j.
\leqno (10.14)
$$

\ms
We leave the proof of Lemma 10.13 for the end of the
section, because it is a pure fact about connected components,
and it does not interfere with the rest of the proof.

\ms
Let us return to the proof of Lemma 10.9.
For $1 \leq j \leq M$, denote by $F_j$ the shortest connected set
that contains all the points of $E_j\cap \partial B$.
We claim that if $E_j\cap \partial B$ has three points, i.e., 
if $E \cap \partial B$ has three points which lie in the same
component of $E \cap B$, then $F_1$ is either the intersection of
$B$ with a $Y$-set whose center lies in $B$, or else we can rename
the three points $a_j$ so that $F = [a_1,a_2] \cup [a_2,a_3]$,
and these segments make an angle at least $2\pi/3$ at $a_2$.

To prove the claim, we consider any connected set $G$ with
finite length that contains the three $a_j$. 
By Proposition 10.37 in [D2] (for instance), 
we can find a simple arc $\gamma$ in $G$ from $a_1$ to $a_2$,
and a simple arc $\gamma_1$ in $G$ from $a_3$ to $a_1$.
We denote by $z$ the first point of $\gamma_1 \cap \gamma$
when we start from $a_3$ and run along $\gamma_1$. Possibly
$z$ is one of the $a_i$, but this is all right. Call
$g_1$ the arc of $\gamma$ from $a_1$ to $z$,
$g_2$ the arc of $\gamma$ from $z$ to $a_2$, and
$g_3$ the arc of $\gamma_1$ from $a_3$ to $z$.
These three arcs are disjoint, except for $z$, and 
they are contained in $G$. Then 
$H^1(G) \geq H^1(\cup_j g_j) = \sum_j H^1(g_j)
\geq H^1(\cup_i [a_i,z])$,
because the $g_i$ are disjoint. So
it is enough to minimize length among sets
of the form $\cup_i [a_i,z]$.
By differentiating in $z$, it is easy to see that 
$H^1(\cup_j g_j)$ is minimal only when $\cup_j g_j$ is
contained in a $Y$-set centered at $z$, or else $z$ is
one of the $a_j$. In this last case, if for instance
$z=a_2$, it is easy to show directly that the segments
$[a_1,z]$ and $[z,a_3]$ make an angle at least $2\pi/3$
at $z$ (otherwise, a $Y$-set would be shorter). 
This proves our claim.
Notice that in both cases, we constructed a Lipschitz
retraction $h_{F_{j}}$ onto $F_j$.

When $E_j$ contains exactly two points $a_i$,
$F_j$ is the line segment between these points, and
we also constructed $h_{F_{j}}$ above (this is the case
when $L_2$ lies opposite to $L_1$). Finally, when 
$E_j \cap \partial B$ is reduced to $a_j$, $F_j = \{ a_j \}$
and we simply take $h_{F_{j}}$ constant equal to $a_j$.
\ms
We are now ready to define a competitor $\varphi(E)$.
First set
$$
\varphi(y)=y \ \hbox{ for } y\in \R^n\setminus B(0,\rho).
\leqno (10.15)
$$
We would have liked to set $\varphi(y) = h_{F_{j}}$
on $E'_j$, but then $\varphi$ would not necessarily
be Lipschitz, because $E \cap B(0,\rho)$ may approach the $a_j$
tangentially, and the definitions of $\varphi$ from
inside and outside $B$ would not match well.
So we shall need a small cut-off function $\psi$.
For this reason we introduce a small positive number
$\varepsilon$, set $d(y) = \dist(y,E\cap \partial B)$, 
and then
$$\eqalign{
\psi(y) &=0  \hskip2.1cm  \hbox{ for } 0 \leq d(y) \leq \varepsilon,
\cr 
\psi(y) &= \varepsilon^{-1}d(y) -1 
\hskip0.5cm \hbox{ for } \varepsilon \leq d(y) \leq 2 \varepsilon,
\cr 
\psi(y) &= 1 \hskip2.1cm  \hbox{ for } d(y) \geq 2\varepsilon.
}\leqno (10.16)
$$
Next set
$$
\varphi_j(y) = y + \psi(y) [h_{F_{j}}(y)-y]
\ \hbox{ for } y\in \R^n,
\leqno (10.17)
$$
and 
$$
\varphi(y)=\varphi_j(y) \ \hbox{ for } y\in E'_j.
\leqno (10.18)
$$
Notice that $\varphi_j(y) = y$ on $F_j$, which contains
all the points of $E'_j \cap \partial B$, because
$F_j$ contains all the points of $E_j \cap \partial B$,
and $E'_j \setminus E_j$  does no meet $\partial B$
(every point of $E\cap \partial B$ lies in some $E_j$).
Thus the two definitions of $\varphi(y)$ coincide on 
$E\cap \partial B$. Let us check that for each $j$,
$$
\varphi_j \hbox{ is $10$-Lipschitz.} 
\leqno (10.19)
$$
Fix $j$, and start the verification in the small balls 
$B_k = \overline B(a_k,2\varepsilon)$ centered at
the points of $E\cap \partial B$.
Notice that $|h_{F_{j}}(y)-y| \leq 6 \varepsilon$ for $y\in B_k$,
because $h_{F_{j}}$ is $2$-Lipschitz and $h_{F_{j}}(a_k)=a_k$
(see (10.11)). Then
$$\eqalign{
|\psi(y) [h_{F_{j}}&(y)-y]-\psi(y') [h_{F_{j}}(y')-y']| 
\cr&\leq \big|[\psi(y)-\psi(y')][h_{F_{j}}(y)-y]\big|
+ \psi(y') \big| [h_{F_{j}}(y)-y]-[h_{F_{j}}(y)-y] \big|
\cr& \leq \varepsilon^{-1} |y-y'| \big| h_{F_{j}}(y)-y\big|
+ 3 \psi(y') |y-y'| \leq 9 |y-y'|
}\leqno (10.20)
$$
for $y$, $y'\in B_k$, again because $h_{F_{j}}$ is $2$-Lipschitz. 
Thus $\psi(y) [h_{F_{j}}(y)-y]$ is $9$-Lipschitz on $B_k$, and 
$\varphi_j$ is $10$-Lipschitz there.
Out of the $B_k$, $\psi(y) = 1$ by (10.16), so 
$\varphi_j(y) = h_{F_{j}}(y)$, which is $2$-Lipschitz. 
This proves (10.19).

Next we check that $\varphi$ is Lipschitz on 
$[\R^n \setminus B(0,\rho)] \cup [E \cap B(0,\rho)]$,
i.e., that 
$$
|\varphi(x)-\varphi(y)| \leq C|x-y|
\leqno (10.21)
$$
for some (very large) $C$. If $x$ and $y$ lie out of 
$B(0,\rho)$, this is clear because of (10.15); so let
us assume that $y\in B(0,\rho)$. Notice that
$y \in E'_j$ for some $j \leq M$, by Lemma~10.13,
so $\varphi(y) = \varphi_j(y)$ by (10.18), hence
$\varphi(y)$ lies in the segment $[y,h_{F_{j}}(y)]$
by (10.17). In particular, we just checked that 
$$
\varphi(y) \in B = \overline B(0,\rho)
\ \hbox{ for } y\in E \cap B(0,\rho)
\leqno (10.22)
$$
because $F_j \i B$.

If $x$ lies out of $B(0,2\rho)$, (10.21) holds
because $|\varphi(x)-\varphi(y)| = |x-\varphi(y)| 
\leq |x-y| + |y-\varphi(y)| \leq |x-y|+2\rho \leq 3|x-y|$.
Next assume that $x\in B(0,2\rho) \setminus B(0,\rho)$. 
If $y$ lies in one of the $B(a_k,\varepsilon)$, then
$\psi(y)=0$, $\varphi(y) = y$, and (10.21) is trivial. 
Otherwise, $y$ lies in the compact set 
$K_\varepsilon = E\cap B \setminus \cup_k B(a_k,\varepsilon)$.
Notice that $K_\varepsilon$ does not meet
$\R^n\setminus B(0,\rho)$, because $E \cap \partial B$
reduces to the $a_k$, so $|y-x| \geq d_\varepsilon = 
\dist(K_\varepsilon,\R^n\setminus B(0,\rho)) >0$, and 
$|\varphi(x)-\varphi(y)| \leq 4 \rho 
\leq 4\rho d_\varepsilon^{-1} |y-x|$ because $\varphi(x)$, 
$\varphi(y) \in B(0,2\rho)$, so (10.21) holds.

We are left with the case when $x,y \in E \cap B(0,\rho)$.
Then $x$ also lies in some $E'_i$. 
If $i=j$, we just use (10.18) and (10.19).
Otherwise, $|\varphi(x)-\varphi(y)| \leq 2 \rho$, and
$|y-x| \geq \dist(E'_j,E'_i) > 0$ by Lemma 10.13, so (10.21) 
holds in this last case as well.

Let us now use the Whitney extension theorem
and extend $\varphi$ to the rest of $\R^n$, so that it
stays Lipschitz. For Lemma 10.9 and Theorem 10.1, we don't 
really need that 
$$
\varphi(B(0,\rho)) \i B
\leqno (10.23)
$$
because there will be no gauge function to pay for,
but let us require this anyway, for later uses of 
the argument. This is easy to arrange; simply compose 
the restriction of $\varphi$ to $B(0,\rho)$ with 
the radial contraction on $B$ defined by
$f(z)=z$ for $z\in B$, and $f(z)=z/|z|$ otherwise.
This does not change the values of $\varphi(y)$ when
$y\in E\cap B(0,\rho)$, by (10.22), and $\varphi$ stays
Lipschitz on $\R^n$ because $\varphi(y)=y=f(y)$ on $\partial B$.

Our argument so far would work with minor modifications
for weak almost-minimal sets. We shall return to this soon,
but for the moment let us do the accounting in the situation
of Lemma 10.9 and Theorem 10.1, where $E$ is a reduced minimal
set. Since $\varphi$ is Lipschitz and $\varphi(x)=x$ out
of $B$, we can apply (8.2), which says that
$H^1(E\setminus \varphi(E)) \leq H^1(\varphi(E)\setminus E)$.
Since $E$ and $\varphi(E)$ coincide out of $B$ 
(by (10.15) and (10.22)), we get that
$$
H^1(E \cap B) \leq H^1(\varphi(E) \cap B) =
H^1(\varphi(E \cap B))
\leqno (10.24)
$$  
(by (10.15) and (10.22) again). Let $y\in E \cap B$
be given. Thus $y$ lies in some $E'_j$. If
$y$ lies out of the $B(a_k,2\varepsilon)$,
$\psi(y) = 1$, and $\varphi(y) = \varphi_j(y) = 
h_{F_{j}}(y) \in F_j$ (by (10.11)). Hence
$$\eqalign{
H^1(\varphi(E \cap B)) &\leq \sum_{j \leq M} H^1(F_j)
+ \sum_j\sum_k H^1(\varphi_j(E'_j \cap B \cap B(a_k,2\varepsilon)))
\cr& \leq \sum_{j \leq M} H^1(F_j) 
+ 10 \sum_j\sum_k H^1(E'_j \cap B \cap B(a_k,2\varepsilon))
\cr& = \sum_{j \leq M} H^1(F_j) 
+ 10 \sum_k H^1(E \cap B \cap B(a_k,2\varepsilon))
}\leqno (10.25)
$$  
by (10.18) and (10.19). Hence
$$
H^1(E \cap B) \leq \sum_{j \leq M} H^1(F_j) 
+ 10 \sum_k H^1(E \cap B \cap B(a_k,2\varepsilon))
\leqno (10.26)
$$
by (10.24). We may now let $\varepsilon$ go to $0$,
and we get that
$$
H^1(E \cap B) \leq \sum_{j \leq M} H^1(F_j)
\leqno (10.27)
$$
because $E\cap \partial B$ is finite.
Now recall that $F_j$ is the connected set
that contains the points of $E_j \cap \partial B$ for which
$H^1(F_j)$ is minimal (the description above gives the uniqueness). 
In particular, $H^1(F_j) \leq H^1(E_j)$ because $E_j$ is connected.
The $E_j$ are distinct and contained in $E\cap B$
by construction, so (10.27) says that $H^1(E_j) = H^1(F_j)$,
whence $E_j = F_j$. In addition, $H^1(E'_j \setminus E_j)=0$.

Suppose for a moment that
$E'_j \setminus E_j$ is not empty. Then we can find
$y\in E'_j \setminus E_j$ and $r>0$ so small that
$B(y,r)$ does not meet the closed sets $E_j = F_j$, nor
$\R^n \setminus B(0,\rho)$ (because $y \in E \cap \partial B$
is impossible). Then $E \cap B(y,r)$ is contained in the
union of the $E'_j \setminus E_j$, so $H^1(E \cap B(y,r)) =0$, 
which contradicts our assumption that $E$ is reduced.
So $E'_j \setminus E_j$ is empty, and $E\cap B$ is the
union of the $F_j$.

If $F_1 \cap \partial B$ contains three points,
then there is no other $F_j$ (by (10.8)), and
$E \cap B =  F_1 \cap B$. In addition, the origin
lies in $F_1$ because we assumed that it lies in 
$E$. If $F_j$ does not come from a $Y$-set, it is of
the form $[a_1,a_2] \cup [a_2,a_3]$. Assume
for instance that $0 \in [a_1,a_2]$. Then 
$[a_1,a_2]$ is a diameter of $B$, and then
its angle with $[a_2,a_3]$ it less than $\pi/2$.
We saw that this is not possible (by minimality of $F_j$).
So $E$ coincides with a $Y$-set in $B$, and we are happy.

If $F_1 \cap \partial B$ contains two points
$a_1$ and $a_2$, $F_2$ either does not exist (if $N$=2)
or is reduced to one point, and $E\cap B(0,\rho) = 
F_1 \cap B(0,\rho) = (a_1,a_2)$.
This case is fine too. The same thing happens if
some other $F_j \cap \partial B$ has exactly two points.

The other cases where no $E_j$ contains more
than one point of $\partial B$ are impossible,
because then the $F_j$ are reduced to one point,
and the origin cannot lie in $E$.

This completes our proof of Lemma 10.9 and Theorem 10.1,
modulo Lemma 10.13 which we prove now.
\qed

\ms\noindent {\bf Proof of Lemma 10.13.}
For each small $\eta >0$, denote by 
$G_{j,\eta}$ the set of points $y \in E\cap B$ that can be 
connected to $E_j$ by an $\eta$-string in $E\cap B$, 
i.e., a finite sequence $\{ y_k \}$ in $E\cap B$, with
$y_0 = y$, $|y_k - y_{k-1}| \leq \eta$, and the
last $y_k$ lies in $E_j$. Let $i \leq M$ be given, with
$i \neq j$, and let us check that for $\eta$ small enough,
$G_{j,\eta}$ does not meet $E_i$.

Suppose not.  We can find a sequence $\{ \eta_l \}$ that
tends to $0$, and for each $l$ an $\eta_l$-string 
$\{ y_k \}$ in $E\cap B$ that goes from some $z_l\in E_i$ 
to some $w_l\in E_j$. Set 
$\Gamma_l = \cup_{k \geq 1} [y_{k-1},y_k]$;
thus $\Gamma_l$ is compact, connected, and contained in $B$.
Replace $\{ \eta_l \}$ with a subsequence for which 
$\{ z_l \}$ converges to some $z$, 
$\{ w_l \}$ converges to some $w$, and 
$\{ \Gamma_l \}$ converges to some compact set 
$\Gamma \i B$. 

Notice that $y\in E_i$ (because $E_i$ is closed),
$z\in E_j$, and also $\Gamma \i E$, because each
point of $\Gamma_l$ lies within $\eta_l$ of $E$,
and $E$ is closed. In addition, we claim  that $\Gamma$ 
is connected. 

Again, suppose this is not the case. Split $\Gamma$ as
$\Gamma = G_1 \cup G_2$, where the $G_i$ are disjoint
nonempty closed subsets. Set $d = \dist(G_1, G_2) >0$, 
and take $l$ so large that every point of $\Gamma_l$
lies within $d/3$ of $\Gamma$, and every point of $\Gamma$
lies within $d/3$ of $\Gamma_l$. 
For $m=1,2$, set $G_{l,m} = \big\{ y\in \Gamma_l \, ; \, 
\dist(y,G_m) \leq d/3 \big\}$. First, $G_{l,m}$ is
closed. It is not empty, because $G_m$ contains some $y$,
then there is a $y'\in \Gamma_l$ such that $|y'-y| \leq d/3$,
and obviously $y'\in G_{l,m}$. Also, $\dist(G_{l,1},G_{l,2})\geq d/3$
by definitions, and finally $\Gamma_l = G_{l,1} \cup G_{l,2}$
because every point of $\Gamma_l$ lies within $d/3$ of
$\Gamma = G_1 \cup G_2$. All this contradicts the connectedness 
of $\Gamma_l$, and our claim that  $\Gamma$ is connected follows. 

We managed to  find a connected set $\Gamma \i E \cap B$ that contains
$z\in E_i$ and $w\in E_j$, but this contradicts our assumption
that $E_i$ and $E_j$ are different connected components of $E\cap B$.
Thus $G_{j,\eta}$ does not meet $E_i$ for $\eta$ small.
Then $\dist(G_{j,\eta},G_{i,\eta}) \geq \eta$, because
otherwise we could put two $\eta$-strings together to
go from $E_i$ to $E_j$.

Now choose $\eta$ so small that all the $G_{j,\eta}$
lie at distance $\geq \eta$ from each other, and 
set $E'_j = G_{j,\eta}$ for $j > 1$ and
$E'_1 = E \cap B \setminus \cup_{j>1} G_{j,\eta}$.
These sets are disjoint because the $G_{j,\eta}$
don't meet. It is obvious that $E'_j \supset E_j$
for $j>1$, and for $j=1$ this is true because
$G_{1,\eta}$ does not meet the other $G_{j,\eta}$.
Their union is $E\cap  B$ by definition of $E'_1$.
We still need to check (10.14). Let us even check that
$\dist(E'_i,E'_j) \geq \eta$ for $i\neq j$.

When $i$ and $j$ are larger than $1$, this comes directly from our 
choice of $\eta$. If $i=1$, we simply observe that if
$y\in E'_1$, then $dist(y,E'_j) = dist(y,G_{j,\eta}) \geq \eta$
because otherwise $y$ would lie in $G_{j,\eta}$.
This completes our proof of Lemma 10.13. 
Lemma 10.9 and Theorem~10.1 follows as well. 
\qed

\bigskip
\noindent {\bf 11. Local regularity for 1-dimensional
weak almost-minimal sets}
\medskip

We still want to describe the two-dimensional minimal 
cones in $\R^n$, and to this effect we shall need information
on the reduced weak almost-minimal sets of dimension 1. 
In this section we let $E \i \partial B(0,1) \i \R^n$ be a closed set,
and we assume that 
$$
\hbox{$E$ is a reduced weak almost-minimal set of dimension $1$,
with gauge function $h$}
\leqno (11.1)
$$
[see Definition 9.1 for weak almost-minimal sets, and 
Definition 2.12 (or (9.5)), for the notion of reduction], 
and we try to derive some amount of regularity from mild
assumptions on $E \cap \partial B(x,\rho)$.
The main result of this section is the following.

\ms
\proclaim Proposition 11.2.
Let $E \i \partial B(0,1) \i \R^n$ satisfy (11.1), and 
suppose that $E\cap \partial B(x,\rho)$ has at most three points.
Also suppose that $r + h(r)$ is small enough (depending only on $n$).
Then for $y\in E \cap B(x,\rho/10)$ and $0 < r \leq \rho/10$, 
we can find a one-dimensional minimal cone $Z$ that contains $y$
(i.e., a line or a $Y$-set through $y$) such that
$$\eqalign{
\dist(z,&Z) \leq  C [r + h(10r)]^{1/2} r
\ \hbox{ for } z\in E\cap B(y,r)
\cr& \hbox{ and }
\dist(z,E) \leq  C [r + h(10r)]^{1/2} r
\ \hbox{ for } z\in Z\cap B(y,r).
}\leqno (11.3)
$$

\ms
Observe that we do not require $Z$ to be centered at $y$ when it is 
a $Y$, and the type of $Z$ may depend on $y$ and $r$ (but in fact not 
wildly). If $E$, $x$, and $\rho$ are as in the proposition and
$h(\rho)$ is small enough, it follows from Proposition 12.6 (the
one-dimensional case of the generalization of Reifenberg' theorem
that will be stated in Section 15)
that $E$ is biH\"older-equivalent to a line or a $Y$-set 
in a ball that contains $B(y,\rho/20)$. 

In the case when $E = K \cup \partial B(0,1)$ for a reduced minimal set 
$K$ of dimension $2$, we could also get the biH\"older equivalence from
Theorem 8.23 and Theorem 10.1, but in this case Proposition 9.4 says that
(11.1) and then (11.3) hold with $h(r) \leq C r$; 
we shall see in Proposition 12.27
that (11.3) then gives a $C^1$ control of $E$ near $y$.
This better control will make it easier to prove that $E$
is composed of arcs of great circles which only meet with
$120^\circ$ angles; see Section 14.

\ms
Let us return to the general case where (11.1) holds, forget
about Proposition 11.2 for the moment, and see 
what the proof of Lemma 10.9 yields in this context.
We let $x\in \partial B(0,1)$ and $\rho\in (0,1/2)$ be given, and
start with the only additional assumption that
$$
E \cap \partial B(x,\rho) \hbox{ is finite.}
\leqno (11.4)
$$

As before, we denote by $a_j$, $j \leq N$, the points of 
$E \cap \partial B(x,\rho)$ (if $E \cap \partial B(x,\rho)
\neq \emptyset$), and by $E_j$ the connected component of $a_j$ 
in $E \cap \overline B(x,\rho)$.
We number the points $a_j$ so that the first
$E_j$, $j \leq M$, are distinct, while each $E_k$,
$k>M$, coincides with some previous $E_j$. And we denote
by $F_j$ the shortest connected set in $\R^n$ that contains
all the points of $E_j \cap \partial B(x,\rho)$.

\ms\proclaim Lemma 11.5. There is a constant $C \geq 1$,
that depends only on $n$, such that 
$$
H^1(E \cap B(x,\rho))
\leq \sum_{j\leq M} H^1F_j) + C \, [h(\rho)+\rho] \, \rho.
\leqno (11.6)
$$

\ms
Let us repeat the construction of the function $\varphi$
as in Lemma 10.9, except that we replace the origin with $x$, 
and we only define $\varphi$ on $\partial B(0,1)$. 
We cannot use $\varphi$ directly, because
of the constraint that competitors lie in $\partial B(0,1)$,
so we shall need to project back onto $\partial B(0,1)$.

Denote by $\xi$ the point of $\R^n$ such that the cone
centered at $\xi$ over $\partial B(x,\rho) \cap \partial B(0,1)$
is tangent to $B$ along $\partial B(x,\rho) \cap \partial B(0,1)$.
That is, $\xi$ lies on the line through $0$ and $x$ and $z-\xi$
is orthogonal to $z-x$ for $z \in \partial B(x,\rho) \cap \partial B(0,1)$.
Notice that $z$ lies on the half line $[0,-x)$ because we may assume 
that $\rho \leq 1$. Set $\eta(z) = z$ for $z\in \partial B(0,1) 
\setminus B(x,\rho)$ and, for $z\in B(z,\rho)$, denote by $\eta(z)$
the radial projection of $z$ onto $\partial B(0,1)\cap B(x,\rho)$ 
centered at $\xi$. That is, $\eta(z)$ is the only point of 
$\partial B(0,1) \cap B(x,\rho)$ that lies on the line through $0$
and $z$. Observe that $\eta$ is $2$-Lipschitz on
$[\partial B(0,1) \setminus B(x,\rho)] \cup B(x,\rho)$; we shall
not need to define it on a larger set.

We set $f(y) = \eta(\varphi(y))$ for $y\in \partial B(0,1)$. 
Observe that we still have that
$$
f(y) = \varphi(y) = y
\ \hbox{ for } y\in \partial B(0,1)\setminus B(x,\rho)
\leqno (11.7)
$$
by (10.15), and $\varphi(y) \in B = \overline B(x,\rho)$
for $y\in B(x,\rho)$, by (10.23), so that in particular
$f(y)$ is well defined,
and
$$
f(B(x,\rho)) \i \overline B(x,\rho) \cap  \partial B(0,1).
\leqno (11.8)
$$

The composition with $\eta$ does not radically change
the Lipschitz properties of $\varphi$, so $f$ is Lipschitz
(by (10.21)), and
$$
\hbox{ the restriction of $f$ to each $E'_j$ is $20$-Lipschitz,}
\leqno (11.9)
$$
by (10.19) and (10.18). 

Let us apply the weak almost-minimality of $E$,
as defined in Definition 9.1. The set
$W_f = \big\{ y\in \partial B(0,1) \, ; \, f(y) \neq y \big\}$ 
is contained in $B(x,\rho)$ by (11.7), and 
$f(W_f) \i \overline B(x,\rho)$ by (11.8). 
In addition, we have a partition of $E\cap \overline B(x,\rho)$
into sets $E'_j$, as in Lemma 10.13, and $f$ is 
$20$-Lipschitz on each $E'_j$ by (11.9). Thus we can take $M=20$ 
in (9.3), which yields
$$
H^1\big(E \cap \overline B(x,\rho)\big)
\leq H^1\big(f(E) \cap \overline B(x,\rho)\big)
+ C \, h(\rho) \rho.
\leqno (11.10)
$$
We want to estimate $H^1\big(f(E) \cap \overline B(x,\rho)\big)$.
We get no contribution from the complement of 
$\overline B(x,\rho)$, by (11.7). 
So it is enough to control the sets $f(E'_j)$
(because $E\cap \overline B(x,\rho) = \cup_j E'_j$).

Recall that $E \cap \partial B(x,\rho)$ is composed
of a finite number of points $a_k$. Observe that
$$\eqalign{
\sum_j \sum_k H^1(f(E'_j \cap B(a_k,2\varepsilon)))
&\leq 20\sum_j \sum_k H^1(E'_j \cap B(a_k,2\varepsilon))
\cr&\leq 20 \sum_k H^1(E \cap B(a_k,2\varepsilon))
}\leqno (11.11)
$$
by (11.9). When $y\in E'_j$ lies out of the $B(a_k,2\varepsilon)$,
we know that $\varphi(y) = h_{F_j}(y) \in F_j$, where
$F_j$ still denotes the shortest connected set
that contains the points of $E_j \cap \partial B(x,\rho)$.
Then $f(y) \in \eta(F_j)$. Altogether,
$$\eqalign{
H^1\big(f(E) \cap \overline B(x,\rho)\big)
& \leq \sum_j H^1(\eta(F_j)) + 
\sum_j \sum_k H^1\big(f(E'_j \cap B(a_k,2\varepsilon)\big)
\cr&
\leq \sum_j H^1(\eta(F_j)) 
+ \leq 20 \sum_k H^1(E \cap B(a_k,2\varepsilon))
}\leqno (11.12)
$$
by (11.11). We replace in (10.10), let $\varepsilon$ tend
to $0$, and get that
$$
H^1(E \cap \overline B(x,\rho))
\leq \sum_j H^1(\eta(F_j)) + C \, h(\rho) \rho
\leq \sum_j H^1F_j) + C \, [h(\rho)+\rho] \, \rho
\leqno (11.13)
$$ 
because $\eta$ barely deforms $F_j$ (recall that the center $\xi$
of the radial projection is far from $x$, and that $F_j$ lies in
the convex hull of the points $a_k$ of $E \cap \partial B(x,\rho)$, 
hence close to the surface $\partial B(0,1)$). This proves Lemma 11.5.
\qed

\ms
\ms\proclaim Lemma 11.14. There is a constant
$C_0\geq 1$, that depends only on $n$, such that
if (11.1) holds, then 
$$
H^1(E \cap B(x_0,r)) \geq r/2
\ \hbox{ for $x_0\in E$ and $r>0$ such that }
r+h(r) \leq C_0^{-1}.
\leqno (11.15)
$$ 

\ms 
Indeed suppose that $x_0\in E$ and $r>0$ are such that
$H^1(E \cap B(x_0,r)) < r/2$. Pick any 
$x\in \partial B(0,1) \cap B(x_0,r/4)$, and set 
$\delta(y) = |y-x|$ for $y \in \R^n$; then 
$H^1(\delta(E \cap B(x_0,r))) < r/2$ because $\delta$ is 
$1$-Lipschitz, so we can find $\rho \in (r/4,3r/4)$ such that 
$\rho \notin \delta(E \cap B(x_0,r))$. 
Since $B(x,\rho) \i B(x_0,r)$, this means that 
$\partial B(x,\rho)$ does not meet $E$.

Thus we can apply Lemma 11.5, with no set $F_j$, and get that
$$
H^1(E \cap B(x,r/4)) \leq H^1(E \cap B(x,\rho)) 
\leq C \, [h(\rho)+\rho] \, \rho 
\leq C \, [h(r)+r] \, r < r/8
\leqno (11.16)
$$
because $h$ is nondecreasing (it is a gauge function)
and if $h(r)+r$ is small enough.

We can iterate this argument, this time keeping the
same point $x$, and we get that
$$
H^1(E \cap B(x,4^{-k}r)) < 4^{-k-1} r
\ \hbox{ for } k \geq 1
\leqno (11.17)
$$
and even, by (11.16) again, 
$$
H^1(E \cap B(x,4^{-k-1}r)) \leq 
C [h(4^{-k}r) + 4^{-k}r] \, 4^{-k} r. 
\leqno (11.18)
$$
Then, if $t$ is small and $k$ is chosen such that
$4^{-k-2} \leq t \leq 4^{-k-1}$, we get that
$$\eqalign{
t^{-1} H^1(E \cap B(x,t)) 
&\leq t^{-1} H^1(E \cap B(x,4^{-k-1})) 
\cr&\leq C t^{-1} [h(4^{-k}r) + 4^{-k}r] \, 4^{-k} r
\leq 4 C [h(16t)+16t].
}\leqno (11.19)
$$
We proved that 
$$
\lim_{t \to 0} t^{-1} H^1(E \cap B(x,t)) =0
\leqno (11.20)
$$
for every $x \in \partial B(0,1) \cap B(x_0,r/4)$. That is,
$B(x_0,r/4)$ contains no point of positive
density for $E$. But it is a simple and standard result 
that as soon as $H^1(E) < +\infty$,
$H^1$-almost every point of $E$ has positive 
upper density; see for instance [Ma], Theorem 6.2 on page 89 
or Lemma 19.1 in [D2]. 
Thus $H^1(E \cap B(x_0,r/4))=0$, and $x\notin E$ because 
$E$ is reduced. Lemma 11.14 follows.
\qed

\ms
We start the local study of $E$ that will lead to Proposition 11.2
with the simpler case when (11.1) holds and
$E \cap  \partial B(x,\rho)$ has two points.
 
\ms
\proclaim Lemma 11.21. There is a constant $C_1 \geq 1$ such that
the following holds. Let $E \i \partial B(0,1)$ satisfy (11.1),
and let $x\in \partial B(0,1)$ and $\rho \in (0,1/2)$ be such that
$E \cap  \partial B(x,\rho) = \{ a_1, a_2 \}$, 
$E \cap B(x,\rho/2) \neq \emptyset$, and $r + h(r) \leq C_1^{-1}$. Then
$$
H^1(E\cap B(x,\rho)) \leq |a_2-a_1| + C [\rho + h(\rho)] \, \rho,
\leqno (11.22)
$$
where $C$ is as in (11.6),
$$
\dist(y,[a_1,a_2]) \leq  C_1 [\rho + h(\rho)]^{1/2} \rho
\ \hbox{ for } y\in E\cap B(x,9\rho/10),
\leqno (11.23)
$$
and 
$$
\dist(z,E\cap B(x,\rho)) \leq C_1 [\rho + h(\rho)]^{1/2} \rho
\ \hbox{ for } z\in [a_1,a_2].
\leqno (11.24)
$$

\ms
Let $E$, $x$, and $\rho$ be as in the lemma.
First observe that if we take $C_1 \geq C_0$, we can
apply Lemma 11.14 to the pair $(x_0,\rho/2)$, where 
$x_0$ is any point of $E \cap B(x,\rho/2)$, and
get that
$$
H^1(E\cap B(x,\rho)) \geq H^1(E\cap B(x_0,\rho/2)) \geq \rho/4.
\leqno (11.25)
$$

Denote by $E_1$ the connected component of $E\cap \overline B(x,\rho)$
that contains $a_1$. First suppose that $E_1$ does not contain $a_2$.
Then Lemma 11.5 applies, with $F_1 = \{a_1\}$ and $F_2 = \{a_2\}$,
and (11.6) says that 
$H^1(E\cap B(x,\rho)) \leq C [\rho + h(\rho)] \, \rho$.
This is impossible (if $\rho + h(\rho)$ is small enough), by (11.25).
So $E_1$ contains $a_2$, we can Lemma 11.5 with $F_1 = [a_1,a_2]$, and
we get (11.22).

Next let $y \in E_1$ be given, set $d=\dist(y,[a_1,a_2])$,
and denote by $Z$ the shortest connected set that contains 
$y$, $a_1$, and $a_2$. Thus 
$$
H^1(Z) \leq H^1(E_1) 
\leq |a_2-a_1| + C [\rho + h(\rho)] \, \rho,
\leqno (11.26)
$$
by (11.22). If $Z = [a_1,a_2] \cup [a_i,y]$ for some $i$, 
we get that $d \leq |y-a_i| \leq  C [\rho + h(\rho)]\, \rho$,
by (11.26). If $Z = [a_1,y] \cup [y,a_2]$, (11.26) and Pythagoras
yield $d \leq C [\rho + h(\rho)]^{1/2} \rho$. The last case is
when $Z = [y,z] \cup [a_1,z] \cup [a_2,z]$ is a piece of $Y$-set,
but then $\dist(z,[a_1,a_2]) \leq C [\rho + h(\rho)]^{1/2} \rho$
because $H^1([a_1,z] \cup [a_2,z]) \leq H^1(Z)$ and by (11.26),
and $|y-z| = H^1(Z) - H^1([a_1,z] \cup [a_2,z])
\leq  H^1(Z) - |a_2-a_1| \leq C [\rho + h(\rho)] \, \rho$.
In this last case too, $d \leq  C [\rho + h(\rho)]^{1/2} \rho$.
So we proved that
$$
\dist(y,[a_1,a_2]) \leq  C_1 [\rho + h(\rho)]^{1/2} \rho
\ \hbox{ for } y\in E_1.
\leqno (11.27)
$$
This proves (11.24), because if $z\in [a_1,a_2]$ is given,
the hyperplane through $z$ which is orthogonal to $[a_1,a_2]$
meets $E_1$ (because $E_1$ connects $a_1$ to $a_2$), and
the intersection lies within $C_1 [\rho + h(\rho)]^{1/2} \rho$
of $z$ by (11.27).

We still need to prove (11.23). Let $y\in E\cap B(x,9\rho/10)$
be given, and set $\delta = \Min(\dist(y,E_1),\rho/10)$. Thus
$B(y,\delta) \i B(x,\rho)\setminus E_1$, and 
$$\eqalign{
H^1(E\cap B(y,\delta)) &\leq H^1(E\cap B(x,\rho)) - H^1(E_1)
\cr&\leq H^1(E\cap B(x,\rho)) - |a_2-a_1| 
\leq C [\rho + h(\rho)] \, \rho,
}\leqno (11.28)
$$
because $E_1$ is connected and contains $a_1$ and $a_2$
and by (11.22). On the other hand, Lemma~11.14 says that
$H^1(E\cap B(y,\delta)) \geq \delta/2$ (if $C_1 \geq C_0$),
so $\delta \leq C [\rho + h(\rho)] \, \rho$. In particular,
$\delta < \rho/10$, so 
$\dist(y,E_1) = \delta \leq C [\rho + h(\rho)] \, \rho$
and then $dist(y,[a_1,a_2]) \leq C [\rho + h(\rho)]^{1/2} \rho$
by (11.27). This proves (11.23) (with a slightly larger $C_1$); 
Lemma 11.21 follows.
\qed

\ms
Lemma 11.21 can be iterated to give more information.
The following lemma gives a proof of Proposition 11.2
in the special case when $E\cap B(x,\rho)$ has at
most two points.

\ms
\proclaim Lemma 11.29.
Let $E \i \partial B(0,1) \i \R^n$ satisfy (11.1), and 
suppose that $E\cap \partial B(x,\rho)$ has at most two points.
Also suppose that $r + h(r)$ is small enough (depending only on $n$).
Then for $y\in E \cap B(x,\rho/10)$ and $0 < r \leq \rho/10$, 
we can find a line $L$ that contains $y$ such that
$$\eqalign{
\dist(z,&L) \leq  C [r + h(10r)]^{1/2} r
\ \hbox{ for } z\in E\cap B(y,r)
\cr& \hbox{ and }
\dist(z,E) \leq  C [r + h(10r)]^{1/2} r
\ \hbox{ for } z\in L\cap B(y,r).
}\leqno (11.30)
$$

\ms
Let $E$, $x$, $\rho$, and $y$ be as in the lemma.
Then we are in the situation of Lemma 11.21 (because
$E \cap B(x,\rho/2)$ contains $y$). By (11.22),
$$
H^1(E \cap B(y,9\rho/10)) \leq 2 \rho + C [\rho + h(\rho)] \, \rho
< 21\rho/10,
\leqno (11.31)
$$
so we can find $r\in (\rho/10,9\rho/10)$ such that
$E \cap \partial B(y,r)$ has at most two points
(otherwise $H^1(E\cap B(y,9\rho/10))\setminus B(y,\rho/10)) 
\geq 24\rho/10$ by the same computation as for (10.7)). 
In fact, $E \cap \partial B(y,r)$ has exactly 
two points, because otherwise Lemma 11.5 says that
$H^1(E\cap B(y,r)) \leq C [r+ h(r)]\, r$, which contradicts
Lemma 11.14. So we can apply Lemma 11.21 again to the pair $(y,r)$.

Let us iterate this with the same $y$. We get a sequence of radii
$r_k$, with $r_0 = r$ and $r_k \in (r_{k-1}/10, 9r_{k-1}/10)$
for $k \geq 1$, such that Lemma 11.21 applies to the pair $(y,r_k)$.
For each $k$, we get a line segment $I_k$ with extremities
on $\partial B(y,r_k)$, such that
$$\eqalign{
\dist(y,&I_k) \leq  C_1 [r_k + h(r_k)]^{1/2} r_k
\ \hbox{ for } y\in E\cap B(y,9r_k/10)
\cr& \hbox{ and }
\dist(z,E\cap B(y,r_k)) \leq  C_1 [r_k + h(r_k)]^{1/2} r_k
\ \hbox{ for } z\in I_k.
}\leqno (11.32)
$$

The sequence $\{ r_{k} \}$ is dense enough; we get that for
$0 < t < 8\rho/10$, we can find a line $L$ such that
$$\eqalign{
\dist(z,L)& \leq  C [t + h(10t)]^{1/2} t
\ \hbox{ for } z\in E\cap B(y,t)
\cr& \hbox{ and }
\dist(z,E) \leq  C [t + h(10t)]^{1/2} t
\ \hbox{ for } z\in L\cap B(y,t).
}\leqno (11.33)
$$
Notice that $\dist(y,L) \leq C [t + h(10t)]^{1/2} t$
because $y \in E$, so the line $L'$ through $y$ and parallel
to $L$ also satisfies (11.33) (with a larger constant).
This proves Lemma 11.29. 
\qed

\ms
In the case of Lemma 11.29, the classical topological disk 
theorem of Reifenberg, in the simple case of sets of dimension 1,
says that $E$ is a H\"older curve near $x$. 
If in addition $h(t) \leq C t^{\alpha}$ for some $\alpha>0$,
say, then (11.30) implies that $E$ has a tangent everywhere
on $B(x,\rho/10)$, and the direction of this tangent is 
H\"older continuous, so that $E$ is in fact a $C^{1,\alpha/2}$ 
curve near $x$.

\ms
Let us now turn to the case when (11.1) holds and
$E \cap  \partial B(x,\rho)$ has three points, and discuss 
the analogue of Lemmas 11.21 and 11.29. Assume that (1.1) holds, that
$$
E \cap  B(x,\rho/10)  \neq \emptyset,
\leqno (11.34)
$$
and that
$$
E \cap  \partial B(x,\rho) = \{ a_{1}, a_{2}, a_{3} \}.
\leqno (11.35)
$$
Recall that $E_{j}$ is the component of $a_{j}$ in
$E \cap \overline B(x,\rho)$. If the three $E_{j}$
are different,  Lemma~11.5 applies with $F_{j}=\{ a_j \}$,
and yields $H^1(E \cap B(x,\rho))
\leq C \, [h(\rho)+\rho] \, \rho$. This is impossible, by
(11.34) and Lemma 11.14. So we may assume that
$E_{1}$ contains $a_{1}$ and $a_{2}$.

If $E_{1}$ does not  contain $a_{3}$, we can apply
the proof of Lemma 11.21, with $F_{1}=[a_{1},a_{2}]$
and $F_2 = \{ a_3 \}$,
and get the conclusions (11.22)-(11.24), and then 
(11.31)-(11.33). In particular, the conclusions of
Lemma 11.29 and Proposition 11.2 are valid.

\ms
So let us assume that $E_{1}$ contains the three
$a_{j}$. Let $F_{1}$ be the shortest connected set
that contains the three $a_{j}$. We know from 
Lemma 11.5 that 
$$
H^1(E \cap B(x,\rho))
\leq H^1(F_1) + C \, [h(\rho)+\rho] \, \rho,
\leqno (11.36)
$$
and we claim that this implies that
$$
\dist(y,F_{1}) \leq C \, [h(\rho)+\rho]^{1/2} \rho
\ \hbox{ for } y\in E_{1}.
\leqno (11.37)
$$
Let $G_{1}$ be a connected set that contains
$y$ and the three $a_{j}$, and whose length is minimal.
Thus (11.36) says that
$$
H^1(G_{1}) \leq H^1(E_{1}) \leq 
H^1F_1) + C \, [h(\rho)+\rho] \, \rho.
\leqno (11.38)
$$

Let us describe $G_{1}$ as we did for $F_{j}$
in Lemma 10.13. As before, we can find a simple arc 
$\gamma_{1} \i G_{1}$ from $a_{1}$ to $a_{2}$, then
a simple arc $\gamma_{2}\i G_{1}$ from $a_{3}$ to
some point $z_{1}\in\gamma_{1}$ (where $z_{1}$
is the first point of $\gamma_{1}\cap \gamma_{2}$
when we run along $\gamma _{2}$ starting from $a_{3}$),
and then a simple arc $\gamma_{3}$ that starts from
$y$ and hits $\gamma_{1}\cup \gamma_{2}$ at some
point $z_{2}$ (and not before). We do not exclude the
cases when $z_{1}$ is one of the $a_{i}$ or $z_{2}=y$.
Set $\gamma=\cup_{1 \leq i \leq 3} \gamma_{i}$.
Since $G_{1}$ has minimal length and $\gamma$
connects $y$  and the three $a_{j}$, $G_{1}=\gamma$.
In addition, the arcs of $\gamma$ between two consecutive
nodes or endpoints are line segments. That is, $G_{1}=\gamma$
is composed of five or less line segments.
The question is how far from $F_{1}$ can it get, knowing (11.38).

Set $G_2 = G_1 \setminus [y,z_2[$. Observe that $G_2$
still connects the three $a_i$, and $H^1(G_2) = H^1(G_1)-|y-z_2|$.
By minimality of $F_1$, $H^1(F_1) \leq H^1(G_2)$, so
$|y-z_2| =  H^1(G_1)- H^1(G_2) \leq H^1(G_1)- H^1(F_1)
\leq [h(\rho)+\rho] \, \rho$, by (11.38). Thus (11.37) will follow as 
soon as we check that 
$$
\dist(z_2,F_{1}) \leq C \, [h(\rho)+\rho]^{1/2} \rho.
\leqno (11.39)
$$
Recall that $z_2$ lies in $\gamma_1 \cup \gamma_2$. First suppose that 
$z_2 \in \gamma_2$. Thus $G_2$ is composed of the four arc segments
$[a_1,z_1]$, $[a_2,z_1]$, $[z_1,z_2]$, and $[z_2,a_3]$, where we do 
not exclude the possibility that some of these segments are reduced
to one point. If we replace $[z_1,z_2] \cup [z_2,a_3]$ with
$[z_1,a_3]$ in $G_2$, we get a shorter connected set $G_3$ that 
contains the three $a_j$, and so 
$H^1(G_2) - H^1(G_3) \leq H^1(G_2) - H^1(F_1) 
\leq [h(\rho)+\rho] \, \rho$, by (11.38). Thus 
$\dist(z_2,[z_1,a_3]) \leq C \, [h(\rho)+\rho]^{1/2} \rho$,
and it will be enough to show that
$$
\dist(z,F_{1}) \leq C \, [h(\rho)+\rho]^{1/2} \rho.
\ \hbox{ for } z\in G_3.
\leqno (11.40)  
$$
In the other case where $z_2 \in \gamma_1$, 
$G_2$ is either composed of the four segments
$[a_1,z_2]$, $[z_2,z_1]$, $[z_1,a_2]$, and $[z_1,a_3]$
(if $z_2$ lies between $a_1$ and $z_1$ on $\gamma_1$),
or $[a_1,z_1]$, $a_3,z_1]$, $[z_1,z_2]$, and $[z_2,a_2]$
(if $z_2$ lies between $z_1$ and $a_2$).
In both cases we set $G_3 = \cup_j [a_j,z_1]$ as before,
and it is still enough to prove (11.40). Notice that
$$
H^1(G_3) \leq H^1(G_1) \leq H^1F_1) + C \, [h(\rho)+\rho] \, \rho
\leqno (11.41)
$$
by (11.38). Also, $F_1 = \cup_j [a_j,z]$, for some $z$ that
may also be equal to one of the $a_j$, and it is enough to show that
$$
|z-z_1| \leq C \, [h(\rho)+\rho]^{1/2} \rho.
\leqno (11.42)
$$
Suppose that $z_1 \neq z$, set $v = (z_1-z)/|z_1-z|$, and
set $f(t) = \sum_j |z+tv-a_j|$ for $t \geq 0$. 
Let us also assume that $z+tv \neq a_j$ for $t\neq 0$, 
so that we can differentiate $f$. 

Let us compute some derivatives. 
Set $r_j(t) = |z+tv-a_j|$ and $u_j(t) = (z+tv-a_j)/r_j(t)$.
Then $r'_j(t) = \langle v,u_j(t) \rangle$ and
$u'_j(t) = v/r_j(t) - (z+tv-a_j) \, r_j(t)^{-2} 
\langle v,u_j(t) \rangle$, so
$$\eqalign{
r''_j(t) &= \langle v,u'_j(t) \rangle = r_j(t)^{-1}
-\langle v, z+tv-a_j\rangle \, r_j(t)^{-2} 
\langle v,u_j(t) \rangle
\cr&= r_j(t)^{-1} \big(1-\langle v,u_j(t) \rangle^2 \big)
}\leqno (11.43)
$$
(where we used again the definition of $u_j(t)$).
Thus $f'(t) = \sum_j \,\langle v,u_j(t) \rangle$
and $f''(t) = \sum_j \, r_j(t)^{-1} 
\big(1-\langle v,u_j(t) \rangle^2 \big) \geq 0$.

Notice that $f$ has a minimum at $t=0$.
Then the half derivative $f'(0_+)$ of $f$ at $0$ is nonnegative.
As long as 
$$
\sum_j \big(1-\langle v,u_j(t) \rangle^2 \big)
> 10^{-2}
\leqno (11.44)
$$
for $0 \leq t < t_0$, we get that $f''(t) \geq (200 \rho)^{-1}$
(because $r_j \leq 2\rho$, since all our points lie in 
$\overline B(x,\rho)$). Then 
$$
f'(t) \geq (200 \rho)^{-1} t
\ \hbox{ for } 0 \leq t \leq t_0.
\leqno (11.45)
$$
If (11.44) fails for some $t_0$, the three 
$\langle v,u_j(t_0) \rangle^2$ are close to $1$, which means 
that the three $u_j(t_0)$ are close to $v$ or $-v$. Recall that
$f'(t) = \sum_j \,\langle v,u_j(t) \rangle$, and $f'(t) \geq 0$
(because $f'(0_+) \geq 0$ and $f'' \geq 0$), so at least two
of the $u_j(t_0)$ are close to $v$. This is even more so for
$t > t_0$ (we move $z+tv$ away from the corresponding two $a_j$), so
$f'(t) \geq 1/2$ for $t\geq t_0$, 
which is even better than (11.45) if $t \leq 2\rho$.

We may now integrate (11.45) or the better estimate between
$0$ and $t_1 = |z_1-z|$ and get that
$f(t_1) - f(0) \geq (400 \rho)^{-1} t_1^2$. But 
$$
f(t_1) - f(0) = \sum_j |z_1-a_j| -\sum_j |z-a_j|
= H^1(G_3) - H^1(F_1) \leq  C \, [h(\rho)+\rho] \, \rho
\leqno (11.46)
$$
by (11.41), so $|z_1-z|^2 = t_1^2 \leq C \, [h(\rho)+\rho] \, \rho^2$,
and (11.42) follows. This completes our proof of (11.37).

\ms
Observe that $H^1(E \cap B(x,\rho)) 
\leq H^1(F_1) + C \, [h(\rho)+\rho] \, \rho 
\leq H^1(E_1) + C \, [h(\rho)+\rho] \, \rho$, 
by (11.36) and because $E_1$ also connects the three $a_j$, so 
$H^1(E \cap B(x,\rho)\setminus E_1) \leq  C \, [h(\rho)+\rho] \, \rho$.
Then
$$
\dist(y,F_{1}) \leq C \, [h(\rho)+\rho]^{1/2} \rho
\ \hbox{ for } y\in E \cap B(x,9\rho/10)
\leqno (11.47)
$$
by (11.37) and Lemma 11.14; the proof is the same as for
the proof of (11.23) at the end of Lemma 11.21.
Next
$$
\dist(y,E\cap B(x,\rho)) \leq \dist(y,E_{1}) 
\leq C \, [h(\rho)+\rho]^{1/2} \rho
\ \hbox{ for } y\in F_{1},
\leqno (11.48)
$$
because $E_1$ connects the three $a_i$ and is contained in a small 
strip around $F_1$, by (11.37). [Cut the strip transversally by a 
segment of length $C \, [h(\rho)+\rho]^{1/2} \rho$, and observe that
if the segment did not meet $E_1$, it would separate some $a_j$
from the other ones.]

Finally notice that since $E$ meets $B(x,\rho/10)$ by (11.34), $F_1$ gets very
close to $B(x,\rho/10)$. If $h(r)$ is small enough, this actually 
prevents $F_1$ from being a simple union of two intervals $[a_i,a_j]$, 
because one of them would be close to a diagonal, and then they would 
make an angle smaller than $2\pi/3$. So $F_1$ is the intersection of 
$\overline B(x,\rho)$ with a $Y$-set. 

\ms
We just finished a verification of the analogue of Lemma 11.21, and we 
are ready to complete the verification of Proposition 11.2 in the 
present remaining case.
Let $y\in E \cap B(x,\rho/10)$ be given; as for (11.31),
$$\eqalign{
H^1(E \cap B(y,9\rho/10)) &\leq H^1(E \cap B(x,\rho))
\leq H^1(F_1) + C \, [h(\rho)+\rho] \, \rho,
\cr&\leq 3 \rho + C [\rho + h(\rho)] \, \rho
< 31\rho/10,
}\leqno (11.49)
$$
by (11.36) and if $h(r)$ is small enough.
Then we can find $r\in (\rho/10,9\rho/10)$ such that
$E \cap \partial B(y,r)$ has at most three points
(otherwise $H^1(E\cap B(y,9\rho/10)\setminus B(y,\rho/10)) 
\geq 32\rho/10$, by the same computation as for (10.7)). 
If $E \cap \partial B(y,r)$ has two points, we
can apply Lemma~11.21 again to the pair $(y,r)$, and conclude as in 
Lemma 11.29. Otherwise, the pair $(y,r)$ satisfies (11.34) and (11.35), 
and we can iterate our argument (with the same $y$).

In all cases we can find a sequence of radii
$r_k$, with $r_0 = r$ and $r_k \in (r_{k-1}/10, 9r_{k-1}/10)$
for $k \geq 1$, so that for each $k$, we can apply 
Lemma 11.21 or the discussion above to the pair $(y,r_k)$.
By (11.23) and (11.24) or (11.47) and (11.48), there is 
minimal cone $Z_k$ (a line or a $Y$-set) such that
$$\eqalign{
\dist(z,&Z_k) \leq  C_1 [r_k + h(r_k)]^{1/2} r_k
\ \hbox{ for } z\in E\cap B(y,9r_k/10)
\cr& \hbox{ and }
\dist(z,E\cap B(y,r_k)) \leq  C_1 [r_k + h(r_k)]^{1/2} r_k
\ \hbox{ for } z\in Z_k\cap B(y,r_k).
}\leqno (11.50)
$$
If $Z_k$ does not go through $y$, we can replace it with a translation 
of $Z_k$, without affecting (11.50).
We also get a similar result (with a constant $C_1$ ten times 
larger) for every $r \in (0,\rho/10)$, by the same 
argument as for (11.33). [Apply (11.50) to $r_k$, where 
$r \leq r_k \leq 10r$.] This completes our proof of Proposition 11.2.
\qed

\ms\noindent{\bf Remark 11.51.}
It is very likely that we can use Lemma 11.5 to prove
a rapid decay estimate on the density of $E$ when 
$E\cap \partial B(x,\rho)$ has more than four points.
Indeed, if $E\cap \partial B(x,\rho)= \{ a_1, \cdots, a_N \}$,
with $N \geq 4$, the length of the shortest connected set that contains
the $a_j$ is at most $\tau N \rho$, where the constant $\tau <1$ 
depends on the dimension, and can be taken very small when $N$
is large.

Once we get such a decay, we can probably find an effective radius
$r_0 > 0$ (that depends on the dimension) such that
for $y \in E$ and $0 < r \leq r_0 \,$, we can find $\rho > r$ so that 
$E\cap \partial B(x,\rho)$ has at most three points. 
Then there is also an effective radius $r_1 > 0$
such that for $\in E$, $E$ is biH\"older-equivalent to a line or a $Y$ 
in a ball that contains $B(y,r_1)$.

\bigskip
\noindent {\bf 12. A baby extension of Reifenberg's theorem}
\medskip

If we want to deduce a good description of $E\cap B(x,\rho/10)$
from the conclusion of Proposition 11.2, the simplest is to use a 
slight improvement of Reifenberg's topological disk theorem, in the
special case of one-dimensional sets, but where we only assume that
$E$ is close enough to lines or $Y$-sets.
The statement below is a special case of the generalization 
of Reifenberg's theorem from [DDT] that will be described in 
Section 15, but for the convenience of the reader, we shall try to give
a short account of how it could be proved.

In this section, $E \i \R^n$ is a compact set.
Denote by $\cal Y$ the set of nonempty one-dimensional minimal 
cones in $\R^n$. That is, $Y \in \cal Y$ if $Y$ is a line or a $Y$-set 
(three half lines in a same plane, with the same endpoint, 
and that make $120^\circ$ angles at that point). We measure 
the closeness of $E$ to sets of $\cal Y$ with
$$
\beta_Y(x,r) = \inf_{Y \in \cal Y} d_{x,r}(E,Y),
\leqno (12.1)
$$
where 
$$\eqalign{
d_{x,r}(E,Y) = 
{1 \over r} \, \sup\big\{ &\dist(y,Y) \, ; \, y\in E \cap B(x,r)\big\}
\cr&
+ {1 \over r} \, \sup\big\{ \dist(y,E) \, ; \, y\in Y \cap B(x,r)\big\}.
}\leqno (12.2)
$$
We shall only use $\beta_Y(x,r)$ and $d_{x,r}(E,Y)$ when $E$ and $Y$
meet $B(x,r)$, so we don't need to worry about the definition of
the supremums in (12.2) when some sets are empty. Notice that
then $\beta_Y(x,r) \leq 2$ (take any set $Y$ through $x$).
We could easily restrict to $Y$-sets (and forget about lines), 
because for each line $L$, there is a $Y$ that coincides with $L$ in
$B(x,4r)$, but it is more natural to keep the lines.
If $\cal Y$ was reduced to the set of lines, the $\beta_Y(x,r)$ would
be the bilateral P. Jones numbers that are often used to study the uniform 
rectifiability of sets. We shall assume that 
the origin lies in $E$ (a normalization), and that
$$
\beta_Y(x,r) \leq \varepsilon
\ \hbox{ for $x\in E \cap B(0,2)$ and $0 < r \leq 2$,}
\leqno (12.3)
$$
where $\varepsilon > 0$ is a small constant that we can choose later.
Let us immediately choose, for $x\in E \cap B(0,2)$ and $0 < r \leq 2$, 
a set $Y(x,r) \in \cal Y$ such that 
$$
d_{x,r}(E,Y(x,r)) \leq \beta_Y(x,r) \leq \varepsilon.
\leqno (12.4)
$$

Our first statement will be a description of $E$ near the origin
as a biH\"older image of a line or a $Y$-set. Later in this section we 
will prove the $C^1$ equivalence under a suitable decay assumption 
on the $\beta_Y(x,r)$. Our initial description will be slightly
different, depending on whether or not
$$
Y(0,2) \hbox{ is a $Y$-set centered in $B(0,12/10)$.}
\leqno (12.5)
$$

\ms\proclaim Proposition 12.6.
For each small constant $\tau > 0$, there is a constant 
$\varepsilon > 0$, that depends only on $n$ and $\tau$,
such that if $E \i \R^n$ is a compact set such that $0 \in E$ and
(12.3) holds, we have the following description. If (12.5)
does not hold, there is a line $L$ through the origin
and a biH\"older mapping $f: L\cap B(0,11/10) \to \R^n$, such that
$$
|f(x)-x| \leq \tau \ \hbox{ for } x\in L\cap B(0,11/10),
\leqno (12.7) 
$$
$$
(1-\tau) |x-y|^{1+\tau} \leq |f(x)-f(y)| 
\leq (1+\tau) |x-y|^{1-\tau}
\ \hbox{ for } x, y \in L\cap B(0,11/10),
\leqno (12.8) 
$$
$$
E \cap B(0,1-\tau) \i f(L\cap B(0,1)) \i E \cap B(0,1+\tau).
\leqno (12.9) 
$$
If (12.5) holds, there is a $Y$-set $Y$ centered at some
$z\in E\cap B(0,13/10)$, and a biH\"older mapping 
$f: Y\cap B(0,14/10) \to \R^n$, such that $f(z)=z$,
$$
|f(x)-x| \leq \tau \ \hbox{ for } x\in Y\cap B(0,14/10),
\leqno (12.10) 
$$
$$
(1-\tau) |x-y|^{1+\tau} \leq |f(x)-f(y)| 
\leq (1+\tau) |x-y|^{1-\tau}
\ \hbox{ for } x, y \in Y\cap B(0,13/10),
\leqno (12.11) 
$$
$$
E \cap B(0,13/10-\tau) \i f(L\cap B(0,13/10)) \i E \cap B(0,13/10+\tau).
\leqno (12.12) 
$$

\ms
Of course the precise values of the explicit constants are not so 
important. The main point is a local description of $E$ in $B(0,1-\tau)$
as a simple curve or a union of three curves that we control well.
We can also extend our mapping $f$ to $\R^n$, so that it is still
a biH\"older homeomorphism of $\R^n$ that satisfies
(12.7) and (12.8) (or (12.10) and (12.11)) for $x,y\in B(0,2)$, 
but we shall not need this. See Section 15 for a more general statement 
and [DDT] for a proof of that result.  
Also see Proposition 12.27 for a $C^1$ description.

\ms
For the proof of Proposition 12.6, we suppose that the reader knows a 
proof of the standard Reifenberg theorem, and we explain how to deduce 
the Proposition 12.6 from that proof.
We start with the observation that the set $Y(x,r) \in \cal Y$
that we chose in (12.4) does not depend too wildly
on $x$ and $r$. Let us check that
$$
d_{y,9t/10}(Y(x,r),Y(y,t)) \leq 100\beta_Y(x,r)+2\beta_Y(y,t)
\leq 102\varepsilon
\leqno (12.13)
$$
when 
$$
x,y \in E \cap B(0,2), \, 0 \leq r \leq 1, \, B(y,t) \i B(x,r), 
\hbox{ and } t \geq r/90. 
\leqno (12.14)
$$

The point will be that $Y(x,r)$ and $Y(y,t)$ are both
very close to $E$ in $B(y,t)$. 
First let $z\in Y(y,t) \cap B(y,9t/10)$ be given. By (12.4),
we can find $z_1 \in E$ such that $|z_1-z| \leq \beta_Y(y,t)t$.
Then $z_1 \in B(y,t) \i B(x,r)$, and we can find $z_2 \in Y(x,r)$
such that $|z_2-z_1| \leq \beta_Y(x,r) r$. Altogether,
$|z_2-z| \leq \beta_Y(y,t)t + \beta_Y(x,r) r \leq \beta_Y(y,t)t +
90\beta_Y(x,r) t$ (by (12.14)).
Similarly, if $z\in Y(x,r) \cap B(y,9t/10)$, we can find
$z_1 \in E$ such that $|z_1-z| \leq \beta_Y(x,r) r \leq 90\beta_Y(x,r) t$, 
so $z_1 \in B(y,t)$ and we can find $z_2 \in Y(y,t)$
so that $|z_2-z_1| \leq \beta_Y(y,t)t$ and 
$|z_2-z| \leq \beta_Y(y,t)t + 90\beta_Y(x,r) t$; (12.13) follows.

Denote by $c(x,r)$ the center of $Y(x,r)$ if $Y(x,r)$ is a $Y$-set,
and set $c(x,r) = \infty$ otherwise. A direct consequence of (12.13) 
is that $c(x,r)$ does not vary too rapidly,
and in particular that
$$
|c(x,r)-c(y,t)| \leq C \big[\beta_Y(x,r)+ \beta_Y(y,t) \big] r 
\leq C\varepsilon r
\leqno (12.15)
$$
when (12.14) holds and $c(x,r)$ or $c(y,t)$ lies in $B(y,t/2)$.

\ms
Let us first prove the proposition when (12.5) fails, i.e., 
when $c(0,2) \notin B(0,12/10)$. 
Then (12.15) says that $c(y,1/10) \notin B(y,1/20)$ for
$y\in E \cap B(0,11/10)$. Repeated applications of (12.15) 
show that $c(y,t) \notin B(y,t/2)$ for $y\in E \cap B(0,11/10)$
and $0 < t \leq 1/10$. In this case we have a good approximation
of $E$ by a line in $B(y,t/4)$ for each $y\in E \cap B(0,11/10)$
and $0 < t \leq 1/10$, and we can conclude with the standard 
Reifenberg theorem.

So we may now assume that (12.5) holds. Let us first look
for the point $z\in E$ where the three curves that compose
$E$ near $B(0,1)$ meet. By definition
of $Y(0,2)$, we can find $z_0 \in E$ such that 
$|z_0-c(0,2)| \leq 2 \beta_Y(0,2) \leq 2 \varepsilon$.
Then 
$$
|c(z_0,1)-c(0,2)| \leq C \big(\beta_Y(0,2)+\beta_Y(z_0,1) \big) 
\leq C\varepsilon, 
\leqno (12.16)
$$
by (12.15). Notice that $|c(z_0,1)-z_0| \leq C \varepsilon$,
so $c(z_0,1)$ lies in the middle of $B(z_0,1)$ and we can
find $z_1 \in E$ such that $|z_{1}-c(z_0,1)| 
\leq \beta_Y(z_0,1) \leq C \varepsilon$. Then we iterate.
By induction, we can find a sequence $\{ z_k \}$ in $E$, so that 
$$
|z_{k+1}-c(z_k,2^{-k})| \leq 2^{-k} \beta_Y(z_k,2^{-k})
\leq C \varepsilon 2^{-k},
\leqno (12.17)
$$
hence
$$
|c(z_{k+1},2^{-k-1})-c(z_k,2^{-k})| \leq C \varepsilon 2^{-k}
\leqno (12.18)
$$
by (12.15), then 
$|c(z_{k+1},2^{-k-1})-z_{k+1}| \leq 2C \varepsilon 2^{-k} < 2^{-k-2}$,
we can find $z_{k+2}$, and so on.
Notice also that
$$\eqalign{
|z_{k+2}-z_{k+1}| 
&\leq |z_{k+2}-c(z_{k+1},2^{-k-1})|+|z_{k+1}-c(z_k,2^{-k})|
\cr& \hskip 3cm+ |c(z_{k+1},2^{-k-1})-c(z_k,2^{-k})|
\leq C\varepsilon 2^{-k}
}\leqno (12.19)
$$
by (12.17) and (12.18), so the $z_k$ converge to some $z\in E\cap B(0,13/10)$.
We shall prove that near $B(0,13/10)$, $E$ is composed of three 
Reifenberg-flat curves that meet at $z$.

Let $k \geq 0$ be given. Observe that 
$|z-c(z_k,2^{-k})| \leq |z-z_{k+1}| + |z_{k+1}-c(z_k,2^{-k})|
\leq C \varepsilon 2^{-k},$ by (12.19) and (12.17). 
Denote by $Y_k$ the translation of $Y(z_k,2^{-k})$ by
$z-c(z_k,2^{-k})$. Thus $Y_k$ is a $Y$-set centered at $z$, and
it also gives a very good approximation of $E$ in $B(z,2^{-k-1})$.
Indeed,
$$
d_{z,2^{-k-1}}(Y_k,E) \leq 2 d_{z_k,2^{-k}}(Y(z_k,2^{-k}),E)
+ 2^{k+1}|z-c(z_k,2^{-k})| \leq C \varepsilon,
\leqno (12.20)
$$
by the same proof as for (12.13). We also
have a similar control out of $B(z,1/2)$.
Indeed
$$
d_{0,2}(Y_0,E) \leq C \varepsilon
\leqno (12.21)
$$
because $d_{0,2}(Y(0,2),E) \leq C \varepsilon$
and $Y_0$ is so close to $Y(0,2)$. This last 
comes from (12.13) which says that 
$d_{z_0,9/10}(Y(0,2),Y(z_0,1)) \leq C\varepsilon$,
and because $|z-z_0| \leq C\varepsilon$.

Next we use (12.20) and (12.21) to get a good description 
of $E$  in the annuli
$$
A_{0} = B(0,5/3)  \setminus B(z,2^{-6})
\ \hbox{ and } \ 
A_k = B(z,2^{-k-3}) \setminus B(z,2^{-k-6}), \, k \geq 1.
\leqno (12.22)
$$
We claim that
$$
c(x,t) \notin B(x,t/2) \ \hbox{ for } 
x\in A_k \hbox{ and } t \leq 2^{-k-7}.
\leqno (12.23)
$$
We start with the case when $2^{-k-8}\leq t \leq 2^{-k-7}$. 
Let us first check that
$$
\dist(w,Y_k) \leq C \varepsilon 2^{-k}
\ \hbox{ for } w \in Y(x,t) \cap B(x,t).
\leqno (12.24)
$$
By (12.4), we can find $w_1 \in E$ so that 
$|w_1-w| \leq \varepsilon t \leq \varepsilon 2^{-k-7}$; 
hence $w_1 \in B(z,2^{-k-2})$ (because $x\in B(z,2^{-k-3})$), 
and (12.20) says that
$\dist(w_1,Y_k) \leq C\varepsilon 2^{-k}$; (12.24) follows.

Recall that $Y_k$ is a $Y$-set centered at $z$ and that
$|z-x| \geq 2^{-k-6}$; then (12.24) prevents $Y(x,t)$ from 
being a $Y$-set with a center in $B(x,t/2)$. That is, 
(12.23) holds for $2^{-k-8}\leq t \leq 2^{-k-7}$.

Finally, if (12.23) holds for some $t \leq 2^{-k-7}$, but
$c(x,t/2) \in B(x,t/4)$, (12.15) says that 
$|c(x,t/2)-c(x,t)| \leq C \varepsilon t$, and 
$c(x,t) \in B(x,t/2)$, a contradiction. This proves (12.23).

Because of (12.23), we can use the standard Reifenberg theorem
to control $E$ in each annulus $A_k$. We get three biH\"older 
parameterizations $f_{i,k}$, defined (for instance) on the three
intervals $J_{i,k}$ that compose $Y_k \cap A_k$, with properties 
similar to (12.7) and (12.8) above, and such that
$E \cap A_k \i \cup_i f_{i,k}(J_{i,k}) \i E$. The argument 
even works on $A_0$, even though $A_0$ is not perfectly symmetric.

Recall that the $A_k$ have large overlaps, so we can easily put together 
our descriptions of $E$ in the various $A_k$. We get three biH\"older curves
$\gamma_i$, that leave from $y$ and end out of $B(0,3/2)$, such that
$E \cap B(0,3/2) \setminus \{ y \} \i \cup_i \gamma_i \i E$. 
These curves are disjoint (they are even far from each other on each $A_k$), 
and they can easily be parameterized (with gluing of the $f_{i,k}$), 
so that they have a single parameterization defined on some $Y$-set 
centered at $z$, and that satisfies the constraints (12.10)-(12.12). 
[Incidentally, we shall not really need (12.11) here.]
Proposition 12.6 follows.
\qed

\ms
Next we want to study the regularity of the branches of $E$
when the numbers $\beta_Y(x,r)$ tend to $0$ uniformly and sufficiently 
fast. Let $E$ satisfy the assumptions of Proposition~12.6, say,
with $\tau = 10^{-2}$, and assume in addition that 
there are numbers $\varepsilon_k > 0$,
$k \in \Bbb N$, such that 
$$
\sum_{k \geq 0} \varepsilon_k < +\infty
\leqno (12.25)
$$
and
$$
\beta_Y(x,2^{-k}) \leq \varepsilon_k
\ \hbox{ for } y\in E\cap B(0,2)
\hbox{ and } k \geq 0.
\leqno (12.26)
$$

Let $L$ or $Y$ and $f$ be as in Proposition~12.6. Set
$\gamma_1 = f(L\cap B(0,11/10))$ when (12.5) does not hold.
When (12.5) holds, denote by $\ell_1$, $\ell_2$, and $\ell_3$,
the three components of $Y\cap B(14/10) \setminus \{ z \}$, 
and set $\gamma_j = f(\ell_j)$ for $1 \leq j \leq 3$.

\ms\proclaim Proposition 12.27.
Let $E$ satisfy the assumptions of Proposition 12.6, 
and also (12.25) and (12.26). Then each $\gamma_j$ is of class $C^1$, 
and the tangent direction $\tau_j(x)$ of $\gamma_j$ at $x\in \gamma_j$ 
is such that
$$
\dist(\tau_j(x),\tau_j(y)) 
\leq C \sum_{k\in \Bbb N \, ; \, 2^{-k}\leq 16|x-y|} \varepsilon_j
\ \hbox{ for } x,y\in \gamma_j
\leqno (12.28)
$$
[see (12.31) for the definition of $\dist(\tau_j(x),\tau_j(y))$]. 
In the case when (12.5) holds, the three $\gamma_j$ 
have tangent half lines at $z$ that are coplanar and 
make $120^\circ$ angles.

\ms
In particular, the curves $\gamma_j$ are of class $C^{1,\alpha}$
if $\varepsilon_k \leq C 2^{-\alpha k}$ for some $\alpha \in (0,1)$.

Let us only prove the proposition when (12.5) holds; the other case is 
merely simpler. Let $x\in \gamma_j$ be given, and let $k \geq 0$ be such
that $x\in A_k$. Also let $t\leq 2^{-k-7}$ be given. Recall from (12.23) 
that $c(x,t) \notin B(x,t/2)$.
Also observe that $Y(x,t)$ comes very near $x$, just by (12.4).
Then $Y(x,t)$ coincides with a line in $B(x,t/2)$, and we denote this line
by $L(x,t)$. Thus
$$
d_{x,t/4}(E,L(x,t)) \leq 4 d_{x,t}(E,Y(x,t)) \leq 4 \beta_Y(x,t) 
\leqno (12.29)
$$
by definition of $Y(x,t)$. It will be enough to restrict to powers 
of $2$. We just proved that 
$$
d_{x,2^{-l-2}}(E,L(x,2^{-l})) \leq 4 \beta_Y(x,2^{-l})
\leq 4 \varepsilon_l
\ \hbox{ for } l \geq k+7,
\leqno (12.30)
$$
by (12.26).

Denote by $\tau_l(x)$ the direction of $L(x,2^{-l})$ for $l \geq k+7$. 
We see $\tau_l(x)$ as a vector in the unit sphere 
$\Bbb S^{n-1}$, determined modulo multiplication by $\pm 1$,
and we use the natural distance 
$$
\dist(\tau,\tau') = {\rm Min}(|\tau-\tau'|,|\tau+\tau'|)
\leqno (12.31)
$$
on the sphere modulo $\pm 1$. 

By (12.30) and its analogue for $l+1$, every point of 
$L(x,2^{-l-1}) \cap B(x,2^{-l-1})$ lies within 
$\varepsilon_{l+1} 2^{-l+1}$ from $E$, and then
within $\varepsilon_{l+1} 2^{-l+1}+ \varepsilon_{l} 2^{-l+2}$ 
from $L(x,2^{-l})$; hence
$$
\dist(\tau_{l+1}(x),\tau_l(x)) 
\leq C \varepsilon_{l} + C \varepsilon_{l+1}.
\leqno (12.32)
$$
The triangle inequality yields $\dist(\tau_l(x),\tau_j(x))
\leq C \sum_{m \geq l} \ \varepsilon_m < + \infty$ when $j\geq l$, 
where the finiteness comes from (12.25). Thus
$\tau(x) = \lim_{l \to +\infty} \tau_l(x)$ exists, and
$$
\dist(\tau(x),\tau_l(x))
\leq C \sum_{m \geq l} \ \varepsilon_m
\ \hbox{ for } l \geq k+7.
\leqno (12.33)
$$

Denote by $L(x)$ the line through $x$ with direction $\tau(x)$,
and let us check that
$$
d_{x,2^{-l-3}}(E,L(x)) \leq C \sum_{m \geq l} \ \varepsilon_m
\ \ \hbox{ for } l \geq k+7.
\leqno (12.34)
$$

First let $w\in L(x) \cap B(x,2^{-l-3})$ be given.
Observe that $\dist(x,L(x,2^{-l})) \leq 2^{-l}\varepsilon_l$ by 
(12.30), so $\dist(w,L(x,2^{-l})) \leq 2^{-l}\varepsilon_l
+ 2^{-l-3} \dist(\tau(x),\tau_l(x)) \leq C 2^{-l}\sum_{m \geq l} \ 
\varepsilon_m$ by (12.33). So we can find $w_1 \in L(x,2^{-l})$
such that $|w_1-w| \leq C 2^{-l}\sum_{m \geq l} \ \varepsilon_m$. Notice 
that $w_1 \in B(x,2^{-l-2})$ (or else some point near $x$ will do even 
better), so by (12.30) we can find $y\in E$
such that $|y-w_1| \leq 2^{-l}\varepsilon_l$. Altogether,
$|y-w| \leq C 2^{-l}\sum_{m \geq l} \ \varepsilon_m$, as needed.

Conversely, let $y\in E\cap B(x,2^{-l-3})$ be given. By (12.30), we can
find $w\in L(x,2^{-l})$ such that $|w-y| \leq 2^{-l} \varepsilon_l$.
Then we can find $w_1 \in L(x)$ such that $|w_1-w| 
\leq C 2^{-l}\sum_{m \geq l} \ \varepsilon_m$. This completes our 
verification of (12.34).

The right-hand side of (12.34) tends to $0$ when $l$ tends to 
$+\infty$, hence $L(x)$ is tangent to $E$ at $x$. 
That is, our curves $\gamma_j$ have tangents everywhere.

\ms
We still need to estimate the oscillation of $\tau(x)$.
Let $y$ be another point of the same $\gamma_j$ as $x$, and assume
for the moment that $y$ lies in the same annulus $A_k$ as $x$, 
and that $|y-x| \leq 2^{-k-11}$. Choose $l$ so that
$2^{-l-4} \leq |y-x| < 2^{-l-3}$; then $l \geq k+7$, and
we can apply (12.30). We also have the analogue of (12.30) for $y$,
i.e., $d_{y,2^{-l-2}}(E,L(y,2^{-l})) \leq 4 \varepsilon_{l}$.
So $E$ is very close to $L(x,2^{-l})$ and $L(y,2^{-l})$ in 
$B(x,2^{-l-3}) \i B(x,2^{-l-2})\cap B(y,2^{-l-2})$, and hence
$$
\dist(\tau_l(x),\tau_l(y)) \leq C \varepsilon_{l} \, ,
\leqno (12.35)
$$
for instance by the proof of (12.32). Then 
$$\eqalign{
\dist(\tau(x),\tau(y)) 
&\leq \dist(\tau(x),\tau_l(x)) + \dist(\tau(y),\tau_l(y))
+ C \varepsilon_{l} 
\cr&\leq C \sum_{m \geq l} \varepsilon_{m}
\leq  C \sum_{2^{-m} \leq 2^{-l}} \varepsilon_{m}
\leq  C \sum_{2^{-m} \leq 16 |z-x|} \varepsilon_{m} \, ,
}\leqno (12.36)
$$
by (12.35) and (12.33), which is enough for (12.28).

Next assume that $y$ lies in the same $A_k$ as $x$, but
that $|y-x| > 2^{-k-11}$. By (12.20) or (12.21) and because $y$ lies
in the same $\gamma_j$ as $x$, there is a chain of points
$x_\ell$, $0 \leq \ell \leq N$, in $E \cap A_k$, with $x_0 = x$, 
$|x_{\ell}-x_{\ell-1}| \leq 2^{-k-11}$ for $\ell \geq 1$, 
$x_N = y$, and $N \leq 2^{15}$. Then 
$\dist(\tau_{k+7}(x_{\ell}),\tau_{k+7}(x_{\ell-1})) 
\leq C \varepsilon_{k+7}$, by the proof of (12.35), and hence
$$
\dist(\tau_{k+7}(x),\tau_{k+7}(y)) \leq C \varepsilon_{k+7} \, .
\leqno (12.37)
$$
Now
$$\eqalign{
\dist(\tau(x),\tau(y)) 
&\leq \dist(\tau_{k+7}(x),\tau_{k+7}(y))
+ \dist(\tau_{k+7}(x),\tau(x))
\cr& \hskip 5.8cm + \dist(\tau_{k+7}(y),\tau(y))
\cr& \leq C \varepsilon_{k+7} \, + C \sum_{m \geq k+7} \varepsilon_{m}
\leq C \sum_{2^{-m} \leq 16|x-y|} \varepsilon_{m}
}\leqno (12.38)
$$
by (12.37), (12.33), and because $|y-x| > 2^{-k-11}$.
Again this is enough for (12.28).

We are left with the case when $y$ lies in an annulus
$A_l \neq A_k$. Without loss of generality, we may assume that 
$l \geq k$. We may also assume that $k$ was chosen as 
large as possible. Thus, by (12.22), $x\in A_k \setminus A_{k+1}
= B(z,2^{-k-3})\setminus B(z,2^{-k-4})$ if $k \geq 1$, and
$x \in B(0,5/2)\setminus B(z,1/16)$ if $k = 0$.

We can use the fairly large overlap of the annuli 
$A_m$ to find a chain of points $x_m$, $k \leq m \leq l$, that all lie in 
the same $\gamma_j$ as $x$ and $y$, and so that each $x_m$ lies in 
$A_m \cap A_{m+1}$. Observe that for $k \leq m \leq l-1$,
$$\eqalign{
\dist(\tau_{m+7}(x_m),\tau_{m+8}(x_{m+1}))
&\leq \dist(\tau_{m+7}(x_m),\tau_{m+7}(x_{m+1}))
\cr& 
\hskip 2.4cm+ \dist(\tau_{m+7}(x_{m+1}),\tau_{m+8}(x_{m+1}))
\cr&\leq C \varepsilon_{m+7} + C \varepsilon_{m+8} 
}\leqno (12.39)
$$
by (12.32) and (12.37) or the proof of (12.35). 
We also have similar estimates on 
\par\noindent 
$\dist(\tau_{k+7}(x),\tau_{k+7}(x_{k}))$
and $\dist(\tau_{l+7}(y),\tau_{l+7}(x_{l}))$,
hence
$$
\dist(\tau_{k+7}(x),\tau_{l+7}(y))
\leq C \sum_{m \geq k+7} \varepsilon_m,
\leqno (12.40)
$$
and then
$$
\dist(\tau(x),\tau(y)) \leq C \sum_{m \geq k+7} \varepsilon_m,
\leqno (12.41)
$$
by (12.33). Recall that $|z-x| \geq 2^{-k-4}$ because 
$x\in A_k \setminus A_{k+1}$. At the same time, $y\in A_l$
for some $l>k$, so $|z-y| \leq 2^{-k-4}$. Since in addition
$z\notin A_k$, we get that $|z-y| \leq 2^{-k-6}$ (compare with (12.22)).
Altogether, $16|x-y| \geq 16|x-z|-16|z-y| \geq 2^{-k-7}$, and 
(12.41) implies (12.28).

We still need to show that the three tangent half lines to
the $\gamma_j$ at $z$ are coplanar and make a $120^\circ$ angle. 
Recall from (12.20) that
$d_{z,2^{-k-1}}(Y_k,E) \leq C \varepsilon$, but the same proof shows 
that 
$$
\lim_{k \to +\infty} d_{z,2^{-k-1}}(Y_k,E) = 0
\leqno (12.42)
$$ 
(that is, (12.25) allows us to replace $\varepsilon > 0$ with any small number 
in the proof of (12.20)). 

Fix a $\gamma_j$ and an annulus $A_k$, and denote by $L_{j,k}$ the
branch of $Y_k$ that gets close to $\gamma_j \cap A_k$. Let
$\tau_{j,k}$ denote the direction of $L_{j,k}$. 
Pick any $x\in A_k \cap \gamma_j$, and apply (12.34) with 
$l=k+7$. We get that
$$\eqalign{
\dist(\tau_{j,k},\tau(x)) &\leq C d_{z,2^{-k-1}}(Y_k,E)
+ C d_{x,2^{-k-10}}(E,L(x))
\cr&\leq C d_{z,2^{-k-1}}(Y_k,E) + C \sum_{m \geq k+7} \ 
\varepsilon_m \, .
}\leqno (12.43)
$$
Now let $k$ tend to $+\infty$. The right-hand side of (12.43) tends to $0$, 
by (12.25) and (12.42). Also, $\tau(x)$ has a limit, by (12.28), so $Y_k$ 
has a limit $Y$. The three branches of $Y$ are tangent half lines for
the corresponding $\gamma_j$ at $z$, by (12.42), and they make
$120^\circ$ angles by definition. In addition, $Y$ lies in a plane
because each $Y_k$ lies in a plane.
This completes our proof of Proposition 12.27.
\qed

\bigskip
\noindent {\bf 13. Harmonic improvement over a small Lipschitz cone}
\medskip

In this section we consider the cone over a small Lipschitz
graph defined on an arc of circle, and show that we almost always reduce its 
surface when we replace it with the graph of a harmonic function with the same
boundary values. This will be useful in the next section,
to show that the reduced two-dimensional minimal cones in $\R^n$ 
are composed of arcs of great circles.

We start with a Lipschitz function $f : [0,T] \to \R^m$, 
and we assume that 
$$
0 < T < \pi \, , \ \ 
f(0) = f(T) = 0,
\ \hbox{ and } \ 
||\nabla f||_{\infty} \leq \tau
\leqno (13.1)
$$
for some small $\tau > 0$. We associate to $f$ a Lipschitz function
$F_0$ on the sector 
$$
D_T = \{ (\rho\cos t,\rho\sin t) \, ; \, 0 \leq \rho \leq 1
\hbox{ and } 0 \leq t \leq T \} \i \R^2
\leqno (13.2)
$$
by
$$
F(\rho\cos t,\rho\sin t) = \rho f(t)
\ \hbox{ for } 0 \leq \rho \leq 1 \hbox{ and } 0 \leq t \leq T.
\leqno (13.3)
$$
Thus the graph of $F$ is a part of the cone over the graph
of the lipschitz function $(\cos t,\sin t) \to f(t)$ that lies above 
$D_T$. We want to compare the surface measure of this graph with what
we get with a harmonic function $G$ with the same values on 
$\partial D_T$.

Since $f(0)=0$, we can extend $f$ to $[-T,T]$ so that it is odd,
and then $f(-T)=f(T)=0$ by (13.1). So
we may write $f$ as a sum of sines, i.e.,
$$
f(t) = \sum_{k \geq 1} \beta_k \sin(\pi k t/T),
\leqno (13.4)
$$
for some coefficients $\beta_k \in \R^m$
such that $\sum_k |k\beta_k|^2 < +\infty$, because
the periodized extension is Lipschitz. Let us check that
$$
{\pi^2 \over 2T} \, \sum_{k \geq 1} k^2 |\beta_k|^2 = ||f'||_2^2
= \int_{[0,T]} |f'(t)|^2 dt.
\leqno (13.5)
$$
Indeed $\displaystyle f'(t) = \,\sum_{k \geq 1} \, 
{\pi k \over T}\, \beta_k \cos(\pi k t/T)$ and the cosine
functions are orthogonal on $[0,T]$ because
$\int_{[0,\pi]} \cos kt \,\cos lt \, dt = {1 \over 2}
\int_{[-\pi,\pi]} \cos kt \,\cos lt \, dt = 0$ for $k \neq l$,
so 
$$
||f'||_2^2 = 
{T \over 2} \sum_{k \geq 1} \Big({\pi k  \over T}\Big)^2 |\beta_k|^2
= {\pi^2 \over 2T} \, \sum_{k \geq 1} k^2 |\beta_k|^2.
\leqno (13.6)
$$

Our new extension $G$ is defined on $D_T$ by
$$
G(\rho \cos t, \rho \sin t) 
= \sum_{k \geq 1} \beta_k \rho^{\pi k/T} \sin(\pi k t/T).
\leqno (13.7)
$$
Notice that it has the same boundary values as $F$,
i.e., that
$$
G(z) = F(z) \ \hbox{ for } z\in \partial D_T.
\leqno (13.8)
$$
Indeed $G(\cos t,\sin t) = f(t)$ by (13.4), and
$G(\rho \cos t, \rho \sin t)  = 0$ when $t=0$ or $t=T$.

We do not really need to know that $G$ is harmonic
on $D_T$, but this is the case, essentially because
(13.7) defines a harmonic function when $T=\pi$ and
$D_T$ is a half disk, while in general $G_T$ is obtained
by composition with the conformal mapping
$z \to z^{\pi/T}$ which sends $D_T$ to a half-disk.
Of course it was reasonable to pick a harmonic function $G$,
because we want its energy to be as small as possible.

\ms\proclaim Lemma 13.9. The functions $F$ and $G$ lie
in $W^{1,2}(D_T)$, with
$$
\int_{D_T} |\nabla G|^2 = {\pi\over 2 }\sum_{k \geq 1} 
\, k |\beta_k|^2 \leq ||f'||_2^2
\leqno (13.10)
$$
and
$$
{1 \over 2} \, ||f'||_2^2
\leq \int_{D_T} |\nabla F|^2 =
\sum_{k \geq 1}  \ {k^2 \pi^2 + T^2 \over 4T} \, |\beta_k|^2
\leq ||f'||_2^2.
\leqno (13.11)
$$

\ms
Let us first compute the energy of $G$. 
We start with the radial derivative 
$$
{\partial G \over \partial \rho}
(\rho \cos t, \rho \sin t) 
= \sum_{k \geq 1} \beta_k \, {\pi k \over T \rho} 
\, \rho^{k\pi /T} \sin(\pi k t/T). 
\leqno (13.12)
$$
For $\rho$ fixed, the sine functions are
orthogonal on $[0,T]$ (again because
$\int_{[0,\pi]} \sin kt \,\sin lt dt = {1 \over 2}
\int_{[-\pi,\pi]} \sin kt \,\sin lt \, dt = 0$ for $k \neq l$ 
by evenness), so 
$$\eqalign{
\int_0^T \Big|{\partial G \over \partial \rho}
(\rho \cos t, \rho \sin t) 
\Big|^2 dt
&= {T \over 2}  \sum_{k \geq 1} |\beta_k|^2 
\Big({\pi k \over T \rho} \Big)^2 \rho^{2\pi k/T}
\cr&= {\pi^2\over 2T} \sum_{k \geq 1} k^2 |\beta_k|^2 
\, \rho^{2\pi k/T} \rho^{-2}
}\leqno (13.13)
$$
and then
$$\eqalign{
\int_{D_T} \Big|{\partial G \over \partial \rho} \Big|^2
&= {\pi^2\over 2T} \ \sum_{k \geq 1} k^2|\beta_k|^2 \,
\int_0^1 \rho^{2\pi k/T} \rho^{-1} d\rho
\cr&=  {\pi^2\over 2T} \ \sum_{k \geq 1} k^2|\beta_k|^2 \, {T \over 2 \pi k}
= {\pi\over 4} \sum_{k \geq 1} k |\beta_k|^2 \, .
}\leqno (13.14) 
$$
Similarly, the derivative in the tangential direction is 
$$
{1 \over \rho} {\partial G \over \partial t} 
(\rho \cos t, \rho \sin t) = 
{1 \over \rho} \sum_{k \geq 1} \beta_k \, \rho^{k\pi /T} \  
{\pi k \over T} \,  \cos(\pi k t/T), 
\leqno (13.15) 
$$
and the same computation (only with sines replaced by cosines) 
leads to
$$
\int_{D_T} \Big| {1 \over \rho} \, {\partial G \over \partial t} \Big|^2
= {\pi\over 4} \sum_{k \geq 1} k |\beta_k|^2.
\leqno (13.16)
$$
We add (13.14) and (13.16) and get the first half of (13.10).
The fact that the result does not depend on $T$ is not too surprising,
because energy integrals are invariant under the conformal mapping 
$z \to z^{\pi/T}$. 
For the second half of (13.10), we simply use (13.5) and the fact that
$k \pi/2 \leq k^2 \pi^2/2T$ for $k\geq 1$, by (13.1).

Let us turn to $F$. Recall from (13.3) et (13.4) that
$$
F(\rho\cos t,\rho\sin t) = \rho f(t)
= \rho \sum_{k \geq 1} \beta_k \sin(\pi k t/T).
\leqno (13.17)
$$
The radial derivative is $\sum_{k \geq 1} \beta_k \sin(\pi k t/T)$
and the integral of its square is 
$\displaystyle{T \over 4} \sum_{k \geq 1} \, |\beta_k|^2$. 
The tangential derivative is
$\displaystyle {1 \over \rho} {\partial F \over \partial t} = 
\sum_{k \geq 1} \, \beta_k {\pi k \over T} \  \cos(\pi k t/T)$, 
and its contribution to the energy is 
$\displaystyle {T \over 4} \sum_{k \geq 1} |\beta_k|^2 \ 
\big( {\pi k \over T} \big)^2$. Altogether
$$
\int_{D_T} |\nabla F|^2 = 
\sum_{k \geq 1}  \ {k^2 \pi^2 + T^2 \over 4T} \, |\beta_k|^2.
\leqno (13.18)
$$
This is the main part of (13.11); for the rest,
just use (13.5) and notice that $0 \leq T^2 \leq \pi^2 k^2$.
Lemma 13.9 follows. 
\qed

\ms
Observe that when $T = \pi$ and $k=1$, the energy terms in (13.10) and 
(13.11) are both equal to ${\pi \over 2} |\beta_1|^2$, which is not 
surprising because for $f(t) = \sin t$, the harmonic extension is 
the same as the homogeneous extension. Here $T < \pi$, and we shall 
see now that all the values of $k$ give an advantage to $G$. 
This is not surprising, since the sine function
in no longer possible, by the constraints in (13.1) and
because a plane that passes through the origin, 
$(1,0,0)$ and $(\cos T,\sin T,0)$, has to be the horizontal plane 
through the origin. 

Set $\lambda = T/\pi < 1$ (by (13.1)), and let us check that
$$
\int_{D_T} |\nabla G|^2 \leq {2\lambda \over 1+\lambda^2} 
\int_{D_T} |\nabla F|^2 
\ \ \hbox{ and } \ \  {2\lambda \over 1+\lambda^2} < 1.
\leqno (13.19)
$$
For the first inequality, we just need to compare the coefficients 
in (13.10) and (13.11), and show that for  $k\geq 1$, 
$$
{k\pi\over 2 } \leq {2\lambda \over 1+\lambda^2} \,
{k^2 \pi^2 + T^2 \over 4T} = {2\lambda \over 1+\lambda^2} \,
{k^2 \pi^2 + \lambda^2 \pi^2 \over 4 \pi \lambda}.
\leqno (13.20)
$$
So we just need to check that
$\displaystyle k \leq {k^2 + \lambda^2 \over 1+\lambda^2}$,
but  $k (1+\lambda^2) \leq k^2 + \lambda^2$ because 
$k^2+\lambda^2 - k (1+\lambda^2) = k(k-1) - \lambda^2(k-1) \geq 0$,
so the first part of (13.19) holds. For the second part, we 
just need to observe that ${2\lambda \over 1+\lambda^2}$ 
is an increasing function of $\lambda \in [0,1]$.

\ms
Let $\Sigma_F$ and $\Sigma_G$ denote the respective graphs 
of $F$ and $G$ (both defined on $D_T$).
We want to use the computations above to compare 
$H^2(\Sigma_F)$ and $H^2(\Sigma_G)$. This will be easy because 
we can assume that $||f'||_{\infty} \leq \tau$, 
which is as small as we want.

\ms \proclaim Lemma 13.21. We have that
$$
H^2(\Sigma_G)-H^2(D_T) \leq {1\over 2}\int_{D_T} |\nabla G|^2
\leq {\lambda \over 1+\lambda^2} \int_{D_T} |\nabla F|^2
\leqno (13.22)
$$
and
$$
H^2(\Sigma_F) - H^2(D_T)
\geq {1\over 2}\, \Big(1-{(1+T)^2\tau^2 \over 2}\Big)
\int_{D_T} |\nabla F|^2.
\leqno (13.23)
$$

\ms
Before we prove the lemma, observe that since we can
take $\tau$ as small as we want (and $T \leq \pi$ anyway),
$1-{(1+T)^2\tau^2 \over 2}$ is as close to $1$ as we want.
Since ${\lambda \over 1+\lambda^2} \leq {1 \over 2}$ by (13.19),
for each choice of $T < \pi$ we can find $\eta \in (0,1)$ such
that (for $\tau$  small enough, depending on $T$ through $\lambda$),
(13.22) and (13.23) imply that
$$
H^2(\Sigma_G)-H^2(D_T) \leq \eta \big[ H^2(\Sigma_F)-H^2(D_T) \big].
\leqno (13.24)
$$
Also,
$$
H^2(\Sigma_F) - H^2(D_T) \geq {1\over 3}\, 
\int_{D_T} |\nabla F|^2 \geq {1\over 6}\, ||f'||_{\infty}^2
\leqno (13.25)
$$
if $\tau$ is small enough, by (13.11).

\ms
Now we prove the lemma.
Notice that $G$ is smooth in the interior of $D_T$; 
the series in (13.12) and (13.15) that give $\nabla G$
converge rather brutally because of the term in $\rho^{k\pi /T}$.
Thus
$$
H^2(\Sigma_G) = \int_{D_T} \Big\{ 1 + |\nabla G|^2 \Big\}^{1/2}.
\leqno (13.26)
$$
Notice that $(1+u)^{1/2} \leq 1+u/2$ for $u \geq 0$, so
$$
H^2(\Sigma_G) \leq \int_{D_T} \Big\{1+{1\over 2}|\nabla G|^2 \Big\}
= H^2(D_T) + {1\over 2}\int_{D_T} |\nabla G|^2,
\leqno (13.27)
$$
which gives the first inequality in (13.22). 
The second one comes from (13.19). 

For (13.23), it will be useful to know that $F$ is 
Lipschitz with a small norm. And indeed
$F(\rho \cos t, \rho \sin t) = \rho f(t)$ by (13.3), so
$$
|\nabla F(\rho \cos t, \rho \sin t)| \leq |f(t)| + |f'(t)| 
\leq (1+T) \, ||f'||_\infty 
\leq (1+T) \, \tau,
\leqno (13.28)
$$
where we also used the fact that $f(0)=0$ (by (13.1)).
Now we use the fact that 
$$
(1+u)^{1/2} \geq 1+{u\over 2}-{u^2\over 4} 
\geq 1 + {u\over 2}\Big(1-{(1+T)^2\tau^2 \over 2}\Big)
\ \hbox { for } 0 \leq u \leq (1+T)^2 \tau^2,
\leqno (13.29)
$$
which we apply with $u = |\nabla F|^2$, to obtain that
$$\eqalign{
H^2(\Sigma_G) - H^2(D_T)
&= \int_{D_T} \Big\{\big\{ 1 + |\nabla F|^2 \big\}^{1/2} - 1 \Big\}
\cr&
\geq {1\over 2}\, \Big(1-{(1+T)^2\tau^2 \over 2}\Big)
\int_{D_T} |\nabla F|^2.
}\leqno (13.30)
$$
This proves (13.23); Lemma~13.21 follows.
\qed

\bigskip
\noindent {\bf 14. Structure of the two-dimensional minimal cones}
\medskip

We are finally ready to give a description of the reduced two-dimensional 
minimal cones in $\R^n$. Recall that by minimal cone, we simply mean
a minimal set (as in Definition~8.1), which is a cone. We say that the 
closed set $E$ is reduced when $E^\ast =E$; see Definition~2.12.
We shall denote by ${\cal Z}$ the set of non empty reduced minimal cones 
of dimension $2$ in $\R^n$, not necessarily centered at the origin, and 
by ${\cal Z}_0$ the subset of elements of ${\cal Z}$ that are centered at 
the origin.

\ms\proclaim Proposition 14.1.
Let $E \in {\cal Z}$ be given, and set $K = E \cap \partial B(0,1)$. 
Then $K$ is a finite union of great circles or arcs of great circles 
$C_j$, $j\in J$. 
The $C_j$ can only meet at their extremities, and each extremity
is a common extremity of exactly three $C_j$, which meet with 
$120^\circ$ angles.
In addition, the length of each $C_j$ is a least $l_0$, where
$l_0 > 0$ depends only on $n$.

\ms
See (14.10) and Lemma 14.12 for additional constraints on $K$.
Thus $K$ is an array of arcs of great circles, with some constraints,
but starting from the ambient dimension $n=4$,
we do not know exactly which arrays give rise to minimal cones. 
When $n=3$, the answer was given by E. Lamarle [La], 
A Heppes [He], and J. Tayor [Ta], 
and indeed some arrays in $\Bbb S^2$ do not give a minimal cone.
See [Ta] and Ken Brakke's home page for descriptions. 
When $n \geq 4$, it is expected that the disjoint union of two disjoint
great circles, contained in almost orthogonal two-planes, gives a 
minimal cone, but the precise condition for this to happen is not known.
There may even be other, combinatorially more complicated, continuous 
families of non mutually isometric minimal cones. 

With the help of Proposition 14.1, we shall be able to give a 
local biH\"older description of the two-dimensional almost-minimal 
sets $E$, but if we want a $C^1$ description near a point $x$, as in 
the original Jean Taylor theorem, it is probable that more information on the 
blow-up limits of $E$ at $x$ will be needed.

\ms
Let us prove the proposition. Take $\varepsilon > 0$ small; 
$\varepsilon = 10^{-2}$ will be more than enough. Apply Theorem 8.23 with $d=1$. 
We get $r_0 > 0$, that depends only on $n$ and $\varepsilon$, such that
if $E \in \cal Z$, then for 
$x\in K=E \cap \partial B(0,1)$ and $0<r<r_0$, we can find a
reduced minimal set $L$ of dimension $1$ such that (8.24) and (8.25) 
hold (that is, $K$ is close to $L$ in $B(x,r)$, both in terms of Hausdorff 
distance and measure). 

Continue with a given choice of $x$ and $r$.
Observe that $L$ is not empty, by (8.24) and because $x\in E$, so Theorem 10.1
says that $L$ is a line or a $Y$-set (not necessarily centered at $x$).
By (8.25) with $y=x$ and $t=r$, we get that
$$
H^1(K\cap B(x,r)) \leq H^1(L \cap B(x,r+\varepsilon r)) +\varepsilon r 
\leq 34r/10.
\leqno (14.2)
$$

Denote by $N(t)$ the number of points in
$K\cap \partial B(x,t)$. By [D3], Lemma 26.1 on page~160, 
or [Fe], Theorem 3.2.22, 
$$
\int_{0}^{r} N(t) dt \leq H^{1}(K\cap B(x,r)) \leq 34r/10,
\leqno (14.3)
$$
so by Chebyshev we can find $\rho\in (r/10,r)$ such that $N(\rho) < 4$.

By Proposition 9.4, $K$ is a reduced weak almost-minimal set,
with gauge function $h(s) \leq Cs$. If $r_0$ was chosen small enough,
Proposition 11.2 applies and yield the good approximation of $K$ by 
one-dimensional minimal cones in every $B(y,t)$, $y\in K \cap B(x,\rho/10)$
and $0 < t < \rho/10$, as in (11.3). This is enough to apply 
Proposition 12.6, with the unit ball replaced with $B(x,\rho/20) 
\supset B(x,r/200)$ (and again, if $r_0$ is small enough). 
We get that there is an open set $V$, with $B(x,r/300) \i V \i 
B(x,r/100)$, where $K$ is either composed 
of a single biH\"older curve, or of three biH\"older curves that 
end at a some point.
In addition, (11.3) says that (12.25) and (12.26) hold, with 
$\varepsilon_k \leq C [2^{-k}+h(2^{-k})]^{1/2}\leq C 2^{-k/2}$. 
So Proposition 11.27 says that each of these curve is of class
$C^{1+\alpha}$. If we prove that every point of these 
curves $C_j$ has a small neighborhood where $C_j$ is a great circle,
Proposition~14.1 will follow, with $l_0 = r_0/1000$, say.

So let $x$ lie in the interior of a curve $C_j$, and choose $r >0$ so 
small that $K\cap B(x,r) = C_j\cap B(x,r)$ is a single $C^{1+\alpha}$ 
curve with two endpoints on $C_j\cap\partial B(x,r)$, and the 
oscillation of the direction of the tangent line to $C_j$ along this 
curve is very small (how small will be decided soon). Without loss of 
generality, we can assume that the two points of $K\cap \partial B(x,r)$ 
are $(1,0,0)$ and $(\cos T, \sin T,0)$ for some small positive $T < \pi$
(and where the last coordinate $0$ lies in $\R^{n-2}$). If the
oscillation of the tangent direction of $K$ in $B(x,r)$ is small enough,
$K\cap B(x,r)$ is a small lipschitz graph over the arc of circle
$J = \{ (\cos t,\sin t,0) \, ; \, 0 \leq t \leq T \}$.

Define the Lipschitz function $f$ on $[0,T]$ by the fact that for
$0 \leq t \leq T$, $(\cos t,\sin t, f(t))$ is the point of $K\cap B(x,r)$ 
above $(\cos t,\sin t, 0) \in J$. We can assume that $||f'||_{\infty}$
is as small as we want, by taking $r$ very small (recall that
$C_j$ is $C^{1+\alpha}$).

Define $F$ on the cone over $J$ by 
$F(\rho\cos t,\rho\sin t) = \rho f(\cos t,\sin t)$ for $\rho \geq 0$
and $0 \leq t \leq T$. Thus $F$ is the homogeneous extension of
the function defined in (13.3). Denote by $C(x,r)$ the cone over
$\partial B(0,1) \cap B(x,r)$ and by $\Gamma$ the graph of $F$.
Thus 
$$
E \cap C(x,r) = \Gamma, 
\leqno (14.4)
$$
because $E$ and $\Gamma$ are cones, and both coincide on the unit 
sphere with $K \cap B(x,r)$. We want to use section 13 to construct 
a competitor for $E$.

Let $G$ be as in (13.7). Also denote by $\pi$ the orthogonal projection
on the horizontal plane, and by $D_T = \{ (\rho\cos t,\rho\sin t,0) \, ; 
\, 0 \leq t \leq T \hbox{ and } 0 \leq \rho \leq 1\}$ the same sector 
as in (13.2). We define a function $\varphi$ on $\Gamma$ by 
$$
\hbox{$\varphi(z)=(\pi(z),G(\pi(z)))$ for
$z \in \Gamma \cap \pi^{-1}(D_T)$
and $\varphi(z)=z$ for $z \in \Gamma \setminus \pi^{-1}(D_T)$.}
\leqno (14.5)
$$

Notice that $\varphi(z)=z$ on $\Gamma\cap\pi^{-1}(\partial D_T)$, 
by (13.8), so $\varphi$ is continuous across 
$\Gamma\cap\pi^{-1}(\partial D_T)$. 
We should also check that $\varphi$ is Lipschitz on $\Gamma$,
or equivalently that $G$ is Lipschitz, but we don't need to have
precise bounds. One could easily get rid of this difficulty by
replacing $G$ with a lipschitz function $G'$, which is so close to $G$ 
in the Sobolev $W^{1,2}$-norm that the strict inequalities of Section 13 
stay true for $G'$. But we can also observe that $G$ is Lipschitz in 
the present situation. Indeed, $f$ is $C^{1+\alpha}$ (by Proposition 11.27), 
and $G$ is obtained from $f$ by a composition with a dilation 
(to get a map defined on $[0,\pi]$), an extension by antisymmetry 
(to get a function on $[-\pi,\pi]$
whose $2\pi$-periodic extension is still $C^{1+\alpha}$), then a 
Poisson extension to the disk (which is Lipschitz because $f$ is 
$C^{1+\alpha}$), which we restrict to the half disk
and compose with the conformal mapping $z \to z^{\pi/T}$ to get a 
function defined on $D_T$, which is still Lipschitz because $T < \pi$. 
Set
$$
\varphi(z)=z \ \hbox{ for } z \in \R^n \setminus C(x,r).
\leqno (14.6)
$$
Recall that $\Gamma \cap \partial C(x,r)$
is the union of the two half lines through 
$(1,0,0)$ and $(\cos T, \sin T,0)$, so (13.8) says that 
$\varphi(z)=z$ there, and $\varphi$ is continuous across 
$\partial C(x,r)$. In fact, due to the fact that $\Gamma$
is transverse to $\partial C(x,r)$, our mapping $\varphi$,
which is now defined on $\Gamma \cup [\R^n \setminus C(x,r)]$,
is Lipschitz there.
In addition, $\varphi(z) = z$ for $|z| \geq 2$. 
So we can extend $\varphi$ to $\R^n$ 
in a Lipschitz way, and so that $\varphi(z) = z$ for $|z| \geq 2$.

Now $E$ is a minimal set, and Definition 8.1 says that
$$
H^2(E\setminus \varphi(E)) \leq H^2(\varphi(E)\setminus E).
\leqno (14.7)
$$

Set $E_1 = E \cap C(x,y) \cap \pi^{-1}(D_T) = \Gamma \cap \pi^{-1}(D_T)$
(by (14.4)). Observe that $\varphi(z)=z$ for $z\in E \setminus E_1$,
by (14.5) and (14.6). Then $E\setminus \varphi(E) 
= E_1\setminus \varphi(E)$ and 
$\varphi(E)\setminus E = \varphi(E_1)\setminus E \i
\varphi(E_1)\setminus E_1$, so (14.7) says that
$$
H^2(E_1\setminus \varphi(E)) =
H^2(E\setminus \varphi(E)) \leq H^2(\varphi(E)\setminus E)
\leq H^2(\varphi(E_1)\setminus E_1),
\leqno (14.8)
$$
and hence
$$
\eqalign{
H^2(E_1) &= H^2(E_1\setminus \varphi(E)) + H^2(E_1 \cap \varphi(E))
\leq H^2(\varphi(E_1)\setminus E_1) + H^2(E_1 \cap \varphi(E))
\cr&
= H^2(\varphi(E_1)\setminus E_1) + H^2(E_1 \cap \varphi(E_1))
= H^2(\varphi(E_1)),
}\leqno (14.9)
$$
because $\varphi(E) \setminus \varphi(E_1) \i \varphi(E\setminus E_1)$ 
does not meet $E_1$. 

On the other hand, $E_1 = \Gamma \cap \pi^{-1}(D_T)$ is the graph
of $F$ over $D_T$, and by (14.5) $\varphi(E_1)$ is the graph of $G$
over $D_T$. With the notation of Section 13, $E_1 = \Sigma_F$
and $\varphi(E_1) = \Sigma_G$ (see above Lemma 13.21).
If we chose $r$ above small enough, (13.1) holds with a 
$\tau$ so small that by Lemma 13.21, (13.24) holds for some 
$\eta \in (0,1)$. That is,
$H^2(\Sigma_G) - H^2(D_T)\leq \eta [H^2(\Sigma_F) - H^2(D_T)]$.
But (14.9) says that $H^2(\Sigma_F) = H^2(E_1) \leq H^2(\varphi(E_1))
= H^2(\Sigma_G)$, so we get that $H^2(\Sigma_F) = H^2(D_T)$.
Then $\int_{D_T} |\nabla F|^2 = 0$, by (13.23) (recall that $\tau$
is small), and $||f'||_2 = 0$, by (13.11). That is, $f = 0$ (by (13.1)).

We just proved that the intersection of the arc $C_j$ with $B(x,r)$
is contained in the horizontal plane. That is, $C_j \cap B(x,r)$ is
an arc of great circle, as needed. Proposition~14.1 follows, as was 
noted before.
\qed

\ms
We complete this section with a little more information on the 
minimal cones in ${\cal Z}_0$. Let $E \i {\cal Z}_0$ be given,
and set $K = E \cap \partial B(0,1)$. Recall from Proposition 14.1
that $K$ is composed of some great circles $C_j$, $j\in J_0$,
and some arcs of great circles $C_j$, $j\in J_1$.

Let $H$ denote the set of extremities of the arcs $C_j$, 
$j\in J_1$. For each $x\in H$, there are exactly three $C_j$ 
that leave from $x \,$; thus there are $3N/2$ arcs, where $N$ is 
the number of points in $H$. In particular, $N$ is even and 
the number of arcs is a multiple of $3$.

Next we claim that there is a constant $\eta_0 > 0$, 
that depends only on $n$, such that 
$$\eqalign{
&\hbox{if $i,j\in J_0 \cup J_1$, $i \neq j$, $x\in C_i$, and
$\dist(x,C_j) \leq \eta_0$, then $i, j \in J_0$,}
\cr& \hskip 1.2cm 
\hbox{ and $C_i$ and $C_j$ have a common endpoint in
$\overline B(x,\dist(x,C_j))$.}
}\leqno (14.10)
$$

We shall take $\eta_0 = r_0/300$, where $r_0$ is as in the 
proof of Proposition 14.1. Let $i$, $j$, and $x\in C_i$
be such that $\dist(x,C_j) \leq \eta_0$. 
As before (see below (14.3)), there is a neighborhood $V$ of 
$B(x,r_0/300)$ where $E$ is composed of one simple curve, 
or of three simple arcs of curves that meet at one point
$a \in B(x,r_0/100)$.
In addition, these arcs are arcs of great circles, and each 
one is contained in a $C_k$, $k\in J_0 \cup J_1$. One of these
arcs is contained in $C_i$ (because $x\in C_i$), and our 
assumption that $\dist(x,C_j) \leq \eta_0$ for some other
$j$ implies that there are three arcs, and another one is 
contained in $C_j$. These three arcs meet at $a$ with 
$120^\circ$ angles, so $|a-x| \leq \dist(x,C_j)$;  
(14.10) follows. 

\ms
The simplest set $E \in {\cal Z}_0$ is a plane, and then $K$ is a
great circle. 

In dimension $n \geq 2p$, $p\geq 2$, we may select
$p$ planes $P_j$ that are orthogonal (or nearly orthogonal) to
each other, and set $E = \cup_j P_j$ (or equivalently, $K$ is
the union of the great circles $P_j \cap \partial B(0,1)$).
It seems that the precise conditions under which $E \in {\cal Z}_0$
is not known. 
More generally, we may always start from a set $E \in {\cal Z}_0$
that lies in $\R^{n-2}$, and add to it the plane $\{ 0 \} \times \R^2$;
this gives a set $E'$ that fits the description above; I don't
know whether $E'$ is necessarily a minimal cone. 

The next simplest example is a set of type $\Bbb Y$, the product of a 
$Y$-set in a plane by a line perpendicular to that plane. It is thus 
contained in a space of dimension $3$, and is composed of three half planes 
with a common boundary $L$, that make $120^\circ$ angles along $L$ (see
Figure 1.1). Notice that it is easy to see that
$$
\hbox{if $E \in {\cal Z}_0$, $K$ is connected, and $H$ has $2$ points, 
then $E$ is a set of type $\Bbb Y$;}
\leqno (14.11)
$$
indeed $E$ is composed of exactly three arcs of great circles
that connect the two points of $H$, these arcs are determined by
their directions at any of the two points of $H$, and we know that
these directions make $120^\circ$ angles.

The simplest set for which $H$ has $4$ points is the cone over
the union of vertices of a regular tetrahedron in $\R^3$ centered at the origin
(see Figure 1.2).

\ms \proclaim Lemma 14.12.
There is a constant $d_T > 3\pi/2$, that depends only on $n$, such that
$$
H^2(E \cap B(0,1)) \geq d_T
\ \hbox{ when $E \in {\cal Z}_0$ is neither a plane nor a set
of type $\Bbb Y$.} 
\leqno (14.13)
$$

\ms
Let us first check that 
$$
H^2(E \cap B(0,1)) > 3\pi/2
\leqno (14.14)
$$
for such sets $E$. If $E$ is composed of more than one full great circles, 
$H^2(E \cap B(0,1)) \geq 2\pi$. Otherwise, $H$ is not empty.
Let $y\in H$ be given, and set $\theta(y,r) = r^{-2} H^2(E\cap B(y,r))$
for $r>0$. Proposition 5.16 says that 
$\theta(y,r) = r^{-2} H^2(E\cap B(y,r))$ is nondecreasing,
so the limits $\theta(y) = \lim_{r \to 0} \theta(y,r)$
and $\theta_\infty(y) = \lim_{r \to +\infty} \theta(y,r)$
exist, and $\theta_\infty(y) \geq \theta(y)$.
In addition, $E$ does not coincide with a cone centered at $y$ in 
big annuli, so Theorem 6.2 says that $\theta_\infty(y) > \theta(y)$.

By Proposition 14.1, the sets $k^{-1}(E-y)$ converge to a minimal cone
of type $\Bbb Y$, and Proposition 7.31 (or directly Lemmas 3.3 and 3.12)
says that $\theta(y) = 3\pi/2$ (the density of a set of type $\Bbb Y$).
Hence $\theta_\infty(y) > 3\pi/2$. Finally, 
$$\eqalign{
H^2(E \cap B(0,1)) &= R^{-2} H^2(E \cap B(0,R))
\geq R^{-2} H^2(E \cap B(y,R-|y|)) 
\cr&= R^{-2}(R-|y|)^2 \, \theta(y,R-|y|)
}\leqno (14.15)
$$
for $R$ large, and the last term tends to 
$\theta_\infty(y) > 3\pi/2$, so (14.14) holds.

Now suppose that (14.13) fails. Then there is a sequence of sets $E_k$ 
in ${\cal Z}_0$, that are neither planes nor sets of type $\Bbb Y$, 
and for which $a_k = H^2(E_k \cap B(0,1))$ tends to $3\pi/2$.
We can replace $\{ E_k \}$ by a subsequence that converges to some 
limit $E$. Obviously, $E$ is a cone, and Lemma 4.7 (with $h=0$) says 
that $E \in {\cal Z}_0$. In addition, $H^2(E \cap B(0,1)) = 3\pi/2$
by Lemmas 3.3 and 3.12, so $E$ is a set of type $\Bbb Y$ by (14.14).

Denote by $H_k$ the set of extremities of the arcs of great circles 
in $E_k$. By Proposition~14.1 or (14.10), the points of $H_k$ lie at 
distance at least $\eta_1 = {\rm Min}(l_0,\eta_0)$ from each other.
If $H_k$ has more than $2$ points and $k$ is large enough,
then some $y_k \in H_k$ lies at distance at least $\eta_1/2$
from the two singular points of $K = E \cap \partial B(0,1)$,
and $K$ is composed of at most an arc of great circle in
$B(y_k,\eta_1/3)$. At the same time, Proposition~14.1 says
that $E_k \cap B(y_k,\eta_1/3)$ is very close, in Hausdorff distance,
to a $Y$-set. This is impossible, because $\{ E_k \}$
converges to $E$.

So $H_k$ has at most two points for $k$ large. If it has exactly
two points, $E_k$ contains a set of type $\Bbb Y$ by the proof of 
(14.11), and since it is not a set of type $\Bbb Y$, 
Proposition~14.1 says that it needs to contain at least another arc 
of circle of length at least $l_0$.
Thus $H^2(E_k \cap B(0,1)) \geq 3\pi/2 + l_0/2$, which does not
tend to $3\pi/2$. If $H_k$ is empty, $E_k$ is composed of at least
two great circles, and $H^2(E_k \cap B(0,1)) \geq 2\pi$. In every
case we reach a contradiction, which proves Lemma~4.12.
\qed

It is believable that we can take 
$d_T = 3 {\rm Argcos}(-1/3) \approx 1.82 \cdot \pi$ in Lemma~4.12, 
which corresponds to the set $T$ based on a tetrahedron in $\Bbb R^3$ 
(see Section 20), but I did not try to check.

\bigskip 
\noindent {\bf D. LOCAL REGULARITY OF $2$-DIMENSIONAL 
ALMOST-MINIMAL SETS}
\medskip

We want to generalize J. Taylor's theorem and give 
a local biH\"older description of $E$ when $E$ is a reduced
almost-minimal set of dimension $2$ in $\R^n$, with sufficiently 
small gauge function.

The proof will rely on a generalization of Reifenberg's Topological
Disk Theorem stated in the next section, together with the 
approximation results of Part B and a small argument
with a topological degree. Most of the local regularity
results are stated and proved in Section 16, except for a 
topological argument which is left for Section 17. 
In Section 18 we prove a regularity result for almost-minimal sets
of dimension $2$  in $\R^3$, which is valid in balls where $E$
is close to a set of type $\Bbb T$, but under an additional separation
assumption. This is enough to establish Theorem 1.9 (i.e., the fact 
that every reduced $MS$-minimal set of dimension $2$ in $\R^3$ is 
a minimal cone). Section 19 contains the description of a set
$E \i \R^3$ that looks like a set of type $\Bbb T$ at large scales,
but does not have the same topology. We rapidly compute the density
of the set of type $\Bbb T$ of Figure 1.2 in Section 20.

\bigskip
\noindent {\bf 15. An extension of Reifenberg's Topological
Disk Theorem}
\medskip

Our statement of the desired extension of Reifenberg's theorem
will use some amount of notation. We shall consider sets of 
dimension $d \geq 2$. The most important case for us is when
$d=2$; we could also take $d=2$, forget about general sets of type
$\Bbb T$, and get a statement as in Section 12, but let us not do this.
We shall allow $d>2$, but we shall keep a collection of cones that
are just obtained as products of a $2$-dimensional cone with
an orthogonal $(d-2)$-dimensional plane.

Let us describe the collection of cones that will be used to 
approximate our set $E$. We start with the collection of 
$d$-dimensional (affine) planes, which we call $\cal P$.

Next we define the collection $\cal Y$ of sets of type $\Bbb Y$.
We start with propellers in a plane, which we called $Y$-sets
in the previous sections, and which are just unions 
of three half lines with a same endpoint and that make $120^\circ$ 
angles at that point. In [DDT],  
we even allow sets $Y$ where the three half lines are only required 
to make angles at least $\pi/2$, but we won't need that generality here.
We obtain a first set of type $\Bbb Y$ as the 
product $\Bbb Y_0 = Y \times V$, where $Y$ is a $Y$-set centered 
at the origin, and $V$ is a $(d-1)$-dimensional vector space that is 
orthogonal to the plane that contains $Y$. We shall call $V$
the spine of $\Bbb Y_0$.
Finally, $\cal Y$ is the collection of sets $\Bbb Y$ of the form
$\Bbb Y = j(\Bbb Y_0)$, where $j$ is an isometry of $\R^n$.
Then the spine of $\Bbb Y$ is the image by $j$ of the 
spine of $\Bbb Y_0$. 

The description of the collection $\cal T$ of sets of type $\Bbb T$
(the sets of type TG3 in [DDT]) 
will take a little longer. For our main application to 
two-dimensional almost-minimal sets, we should take a collection
that includes the minimal cones of dimension $2$, and unfortunately 
this will force us to modify slightly the definitions of [DDT]. 
Our set $\cal T$ will be the collection of sets $T=j(T_0 \times V)$,
where $T_0$ lies in a set ${\cal T}_0$ of $2$-dimensional cones 
in $\R^{n-d+2}$, $V$ is the plane of dimension $d-2$ orthogonal
to $\R^{n-d+2}$ in $\R^n$, and $j$ is an isometry of $\R^n$.

Each $T \in {\cal T}_0$ will be the cone over a set 
$K \i \partial B(0,1)$, with the following properties.
First, $K = \cup_{j\in J} C_j$ is a finite union of great circles, 
or closed arcs of great circles. Denote by $H$ the collection of 
extremities of arcs $C_j$, $j\i J$; each point $x$ of $H$ lies in exactly
three $C_j$, $x$ is an endpoint for each such $C_j$, and the three
$C_j$ make $120^\circ$ angles at $x$. The $C_j$ can only meet at their
endpoints (and hence the full arcs of circles are disjoint from the
rest of $K$). In addition, there is a constant $\eta_0 > 0$ such that
$$
H^1(C_j) \geq \eta_0 \ \hbox{ for } j\in J,
\leqno (15.1)
$$
and 
$$\eqalign{
&\hbox{if $x\in C_i$, and $\dist(x,C_j) \leq \eta_0$ 
for some other $j$, then }
\cr& \hskip 2.4cm \hbox{$C_i$ and $C_j$ have a common 
extremity in $B(x,\dist(x,C_j))$.}
}\leqno (15.2)
$$

In addition, we exclude the case when $T$ is a plane or a
set of type $\Bbb Y$.

When $Y \in {\cal T}_0$, the spine of $E$ is the cone over 
$H$ or, if $K$ is uniquely composed of great circles, the origin.
In general, when $T=j(T_0 \times V)$, the spine of $T$ is
$j(S \times V)$ where $S$ is the spine of $T_0$.

The definition of sets of type $TG3$ given in [DDT] 
was a little different. 
First, we also allowed angles
larger than $\pi/2$ instead of exactly $2\pi/3$, but we don't
need this here. More unfortunately, we forgot the possibility
of adding a collection of full great circles in $K$ (provided
that they stay far from the rest of $K$, as in (15.2)).
We also decided to take the full collection of sets $K$ as above, 
but our statement is a little better if we allow a subcollection
(both in the assumptions and the conclusion).
Nevertheless, the proof given in [DDT] 
goes through without changes, and leads to the statement below.

\smallskip
Before we start, let us record the fact that the class ${\cal T}_0$
of reduced minimal cones of dimension $2$ in $\R^{n-d+2}$ that are 
neither planes nor sets of type $\Bbb Y$ satisfies the requirements above, 
with a constant $\eta_0$ that depends only on $m$. 
This is a consequence of Proposition~4.1, except for (15.2)
which was stated separately in (14.10).

\ms
Let $\cal P$, $\cal Y$, and some choice of $\cal T$ be as above.
Set ${\cal Z} = {\cal P} \cup {\cal Y} \cup {\cal T}$, 
For the statement below, we are given a closed set $E \i \R^n$,
and we assume for simplicity that the origin lies in $E$.
As in Section 12, we measure how close $E$ is to sets of $\cal Z$
with
$$
\beta_Z(x,r) = \inf_{Z \in \cal Z} d_{x,r}(E,Z),
\leqno (15.3)
$$
where $d_{x,r}(E,Z)$ is defined as in (12.2). (Compare with (12.1).)
We assume that 
$$
\beta_Z(x,r) \leq \varepsilon
\ \hbox{ for $x\in E \cap B(0,3)$ and $0 < r \leq 3$,}
\leqno (15.4)
$$
where $\varepsilon \geq 0$ is small, depending on $n, d,$ and
the constant $\eta_0$. This means that for $x\in E \cap B(0,3)$ and 
$0 < r \leq 3$, we can choose $Z(x,r) \in \cal Z$ such that
$d_{x,r}(E,Z(x,r)) \leq \varepsilon$. We want a bilipschitz
parameterization of $E$ by a set of $\cal Z$ near the unit ball.

\ms
\proclaim Theorem 15.5 [DDT].  
Let $\cal Z$ be as above, and let $\tau > 0$ be given.
Then there is a constant $\varepsilon > 0$ that depends only on
$n$, $d$, $\eta_0$, and $\tau$, such that if $E \i \R^n$ is a closed 
set that contains the origin and if (15.4) holds, we can find
$Z \in \cal Z$ and a a biH\"older mapping 
$f: B(0,2) \to \R^n$, such that
$$
|f(x)-x| \leq \tau \ \hbox{ for } x\in B(0,2),
\leqno (15.6) 
$$
$$
(1-\tau) |x-y|^{1+\tau} \leq |f(x)-f(y)| 
\leq (1+\tau) |x-y|^{1-\tau}
\ \hbox{ for } x, y \in B(0,2),
\leqno (15.7) 
$$
$$
B(0,2-\tau) \i f(B(0,2)),
\leqno (15.8) 
$$
and
$$
E \cap B(0,2-\tau) \i f(Z\cap B(0,2)) \i E.
\leqno (15.9) 
$$

\ms
Compared to [DDT], we just changed a little some of the
constants (such as $2$ and $3$), but this does not matter.
\smallskip
It will be useful for us that in addition to giving a parameterization
of $E$ by $Z$ near $0$, $f$ extends to a local homeomorphism
of $\R^n$. In particular, we know that the topology of
$\R^n \setminus E$ near $0$ is the same as for $\R^n \setminus Z$.
\smallskip
We can say a little more on the choice of $Z$.
First observe that $Z$ has to be fairly close to $Z(0,3)$
in Hausdorff distance, because (15.6) and (15.9) say that
it is fairly close to $E$. 

If $Z(0,3)$ is a plane, or coincides with a plane
near $B(0,29/10)$, say, we can show that 
for $x\in E \cap B(0,27/10)$ and $r \leq 1/10$, say,
$Z(x,r)$ is a plane, or at least its spine does not
get close to $x$. Then we are in the standard Reifenberg
situation, and we can take $Z = Z(0,3) \in \cal P$.

Similarly, if $Z(0,3)$ is a set of type $\Bbb Y$ with a spine
through the origin, we can take $Z = Z(0,3) \in \cal Y$. 
Slightly more generally (but we won't need that case in the 
present paper), if $Z(0,3) \in \cal Y$ or coincides with a set of type
$\Bbb Y$ in $B(0,29/10)$, we can show that for $x\in E \cap B(0,27/10)$ 
and $r \leq 1/10$, $Z(x,r)$ is never a set of type $\Bbb T$ centered
near $x$. Then we can take $Z\in \cal Y$, except perhaps 
in the situation above where it is more clever to take a plane.

If $Z(0,3)$ is a set of type $\Bbb T$ centered at the origin, 
then again we can take $Z = Z(0,3)$. More generally, if 
$Z(0,3) \in \cal T$, with a center in $B(0,25/10)$, we can first
show (by checking things at smaller and smaller scales) that
there is a single point $x_0 \in B(26/10)$ such that every $Z(x_0,r)$ 
is a set of type $T$ centered in $B(x_0,r/10)$. Then we use
$Z = Z(0,3)$ (or a translation of it if we want $f$ to map the 
center of $Z$ to $x_0$).

If our class $\cal T$ (modulo identification of isometric sets) is
not discrete, it could be that when $r$ varies, the set $Z(x_0,r)$
changes. This is all right, it just means that some other
sets $Z'\in \cal T$ could be used instead of $Z$ (to parameterize $E$), 
and they are all biH\"older equivalent to each other. 
This is not surprising, since $Z(x_0,r)$ is only 
allowed to vary very slowly with $r$ (because $Z(x_0,r)$ stays so 
close to $E$). The situation would be different if we wanted 
$C^1$ parameterizations. Then we would probably need bounds
on how fast the sets $Z(x_0,r)$ vary and (in the case of minimal
sets), it can be expected that these bounds will be harder to prove
than when $d=2$ and $n=3$.

The remaining case when the center of $Z(0,3)$ lies out of
$B(0,25/10)$, but $Z(0,3)$ does not coincide with a set of type
$\Bbb Y$ in $B(0,29/10)$, could be treated similarly, or just eliminated
by replacing $3$ with a larger constant. We shall not need it here 
anyway.

\ms 
See [DDT] for a proof of Theorem 15.5. 
The proof is a little ugly because there are more cases to
study, but we essentially use Reifenberg's scheme, plus 
an organization of the construction in layers (first find
the center if there is one, then find the image of the
spine, then define $f$ on $Z$, and finally extend to the rest of
$\R^n$).

There is no immediate obstruction to having an additional layer
in the complexity of our sets,
but the situation is already complicated enough, and we were 
lucky that, just because there is at most one center $x_0$,
we did not have to understand what happens when
the isometry class of $Z(x_0,r)$ changes with $r$.

\ms
We shall get our local descriptions of the almost-minimal sets 
by brutal applications of Theorem 15.5 in the case when $Z(0,3)$
is a cone centered at $0$. The following definition will make it
more convenient.

\ms\proclaim Definition 15.10. 
Let $E$ be a closed set in $\R^n$. We way that $B(0,1)$ is a 
biH\"older ball for $E$, with constant $\tau \in (0,1)$, 
if we can find a cone $Z \in \cal Z$ centered at $0$, and 
$f : B(0,2) \to \R^n$, with the properties (15.6)-(15.9). 
We say that $B(0,1)$ is of type $\Bbb P$, $\Bbb Y$, or $\Bbb T$, 
depending on whether $Z$ lies in $\cal P$, $\cal Y$, or $\cal T$.
Finally, we say that $B(x,r)$ is a biH\"older ball for $E$
(with the same parameters) when $B(0,1)$ is a biH\"older ball 
for $r^{-1}(E-x)$.

\ms
With this definition, we have the following consequence of 
Theorem 15.5 and the discussion above. We keep the same assumptions
on $\cal Z$.

\ms
\proclaim Corollary 15.11. 
For each small $\tau > 0$, we can find $\varepsilon > 0$,
that depends only on $n$, $d$, $\eta_0$, and $\tau$, such that if 
$E \i \R^n$ is a closed set, $x\in E$, and $r > 0$ are such that
$$\eqalign{
&\hbox{for $y\in E \cap B(x,3r)$ and $0 < t \leq 3r$, we can find}
\cr& \hskip 3cm
Z(y,t) \in {\cal Z}  \hbox{ such that } 
d_{y,t}(E,Z(y,t)) \leq \varepsilon,
}\leqno (15.12)
$$
and in addition $Z(x,3r)$ is a cone centered at $x$,
then $B(x,r)$ is a biH\"older ball for $E$, with constant $\tau$,
of the same type as $Z(x,3r)$, and we can even take $Z$ isometric
to $Z(x,3r)$ in Definition 15.10.

\bigskip 
\noindent {\bf 16. A local regularity theorem for two-dimensional 
almost-minimal sets in $\R^n$} 
\medskip

We are now ready to state a generalization to higher
dimensions of the biH\"older part of Jean Taylor's 
local regularity theorem [Ta]. 

\ms\proclaim Theorem 16.1. 
Let $U$ be an open set in $\R^n$ and $E$ a closed subset of $U$. 
Suppose that $E$ is a reduced  $A$-almost-minimal set in $U$, 
with a gauge function $h$ such that 
$$
\int_{0}^{r_1} h(r) {dr \over r} < +\infty \ \hbox{ for some } r_1 > 0.
\leqno (16.2)
$$
Then for each $x_0 \in E$ and every  choice of $\tau \in (0,1)$, 
there is an $r_0 > 0$ such that $B(x_0,r_0)$ is a biH\"older 
ball for $E$, with constant $\tau$.

\ms
See Definition 4.3 for $A$-almost-minimal sets. 
Notice that we could equivalently have assumed that 
$E$ is a reduced  $A'$-almost-minimal set with the same gauge 
function $h$, as in Definition 4.8, because Proposition 4.10 says that
the two definitions are equivalent. Recall that we say that $E$
is reduced when $E^\ast = E$, where $E^\ast$ is the closed support 
of the restriction of $H^2$ to $E$; see Definition 2.12.

The notion of biH\"older ball is defined in Definition 15.10, and here
we take for $\cal Z$ the collection of two-dimensional reduced minimal 
cones in $\R^n$, as described in Section 14. Thus, in short, Theorem 16.1
says that every point of $E$ has a small neighborhood where $E$ is 
biH\"older equivalent to a minimal cone of dimension $2$.

\smallskip
Jean Taylor [Ta] also gives the $C^1$ equivalence when $n=3$, 
under conditions on $h$ that are stronger than (16.2) (but it
is not clear that she tried to optimize this).
We shall see in [D3] 
how to deduce this from Theorem 16.1 (or some weaker separation 
result). It is not clear to the author that this more precise result
will generalize almost-minimal sets of dimension $2$ in $\R^n$, 
$n \geq 4$; the answer may depend on which type of blow-up limits 
for $E$ at $x$ exist. 

\smallskip
Curiously, the condition (16.2) is only used to prove control the
density, as in (16.6) at the beginning of the argument; 
I don't know whether it is sharp.

\smallskip
It should be observed that neither here nor in [Ta] 
(I think) do we get any precise lower bound for $r_0$. 
The proof will say that for each $x\in E$, the density 
$\theta(x,r)$ eventually gets very close to its limit at $r=0$, 
and then $B(x,r)$ is a biH\"older ball for $E$, but we do not know in 
advance when this happens. Some of our intermediate results, such as 
Lemma~16.19 and 16.48, Proposition 16.24, Lemmas 16.25, 16.51, and 16.56,
and later on Proposition 18.1, are more quantitative, though. 

\ms
We shall prove Theorem 16.1 now, except for Proposition~16.24, 
which uses a little bit of degree theory and will be treated 
in Section 17. We shall systematically assume in this section that
$$\eqalign{
&\hbox{$E$ is a reduced  $A$-almost-minimal set in $U$,} 
\cr& \hskip 2cm
\hbox{with a gauge function $h$ such that (16.2) holds.}
}\leqno (16.3)
$$
Set
$$
\theta(x,r)=r^{-2} H^2(E\cap B(x,r))
\leqno (16.4)
$$ 
for $x\in E$ and $r>0$. Set 
$$
A(r) = \int_{0}^{r} h(2t) {dt \over t} \, .
\leqno (16.5)
$$
By (16.2), we can apply Proposition 5.24 
and get that 
$$
\theta(x,r) e^{\lambda A(r)}
\ \hbox{ is a nondecreasing function of } r,
\leqno (16.6)
$$
as long as $B(x,r) \i U$ and $h(r)$ stays small enough,
and where $\lambda$ is a positive constant.

Notice that $A(r)$ tends to $0$ when $r$ tends to $0$
(by (16.2)), so (16.6) allows us to define
$$
\theta(x) = \lim_{r \to 0} \theta(x,r)
\leqno (16.7)
$$
for $x\in E$. As we shall see soon, there are restrictions 
on the values that $\theta(x)$ can take. 

\ms
Proposition 7.31 says that every blow-up limit of $E$
at $x$ is a reduced minimal cone, whose density is $\theta(x)$.

If $\theta(x) = \pi$, all the blow up limits are planes, and we shall
call $x$ a $P$-point.
If $\theta(x) = 3\pi/2$, all the blow up limits are minimal cones of type 
$\Bbb Y$, and we shall call $x$ a $Y$-point.
Otherwise, we shall say that $x$ is a $T$-point.
In this case do not know whether the various blow-up limits of 
$E$ at $x$ are all isometric to each other, but they all have 
the same density $\theta(x)$. Lemma 14.12 then says that 
their density is at least $d_T$ (where $d_T > 3\pi/2$ is a
constant that may depend on $n$), so 
$$
\theta(x) \geq d_T
\ \hbox{ when $x$ is a $T$-point.}
\leqno (16.8)
$$

When $n=3$, there is only one type of minimal cone of type $\Bbb T$
(i.e., which is not a plane or a set of type $\Bbb Y$), which is the one 
introduced in Remark 1.8 and Theorem 1.9.
Then $\theta(y) = 3 {\rm Argcos}(-1/3) \approx 1.82 \cdot \pi$
when $x$ is a $T$-point. See the computation in Section 20.  
In higher dimensions, there may be lots
of minimal cones with different densities $\geq d_T$. 
But this will not matter much here, because the $T$-points are isolated.

We start with a simple consequence of (16.6) that will be used
repeatedly. We claim that 
$$
\theta(x) \leq e^{\lambda A(r)} \theta(x,r)
\leqno (16.9)
$$
for $x\in E$ and $r > 0$ such that $B(x,r) \i U$ and
$h(r)$ is small enough. Indeed notice that for 
$0 < s \leq r$, 
$$
\theta(x,s) \leq e^{\lambda A(s)} \, \theta(x,s) 
\leq e^{\lambda A(r)} \, \theta(x,r);
\leqno (16.10)
$$
we let $s$ tend to $0$ and get the result.

The following consequence of Proposition 7.24 will also
be used quite often.

\ms\proclaim Lemma 16.11. 
For each $\tau > 0$, we can find $\varepsilon(\tau) > 0$ such that
if $x\in E$ and $r> 0$ are such that
$$
B(x,r) \i U, \ \ \ 
h(2r) \leq \varepsilon(\tau), \ \ \ 
\int_{0}^{2r} h(r) {dr \over r} < \varepsilon(\tau),
\leqno (16.12)
$$
and
$$
\theta(x,r) \leq \theta(x) + \varepsilon(\tau),
\leqno (16.13)
$$
then for every $\rho \in (0,9r/10]$ there is a reduced minimal cone 
$Z(x,\rho)$ centered at $x$, such that
$$
d_{x,\rho}(E,Z(x,\rho)) \leq \tau
\leqno (16.14)
$$
and
$$\eqalign{
&\big| H^2(E \cap B(y,t)) - H^2(Z(x,\rho)\cap B(y,t)) \big|
\leq \tau \rho^2
\cr&\hskip 2cm
\hbox{ for $y\in \R^n$ and $t>0$ such that $B(y,t) \i B(x,\rho)$.}
}\leqno (16.15)
$$
If $\theta(x) = \pi$, $Z(x,\rho)$ is a plane; if $\theta(x) = 3\pi/2$, 
$Z(x,\rho)$ is a set of type $\Bbb Y$, and if 
$\theta(x) > 3\pi/2$, $Z(x,\rho)$ is a set of type $\Bbb T$.

\ms
Let $\rho \in (0,9r/10]$ be given, and let us try to apply 
Proposition 7.24, with the constant $\tau'=\tau/2$, to the ball 
$B(x,10\rho/9)$. The main assumption is (7.25), which holds because 
$$\eqalign{
\theta(x,10\rho/9) 
&\leq e^{\lambda A(10\rho/9)} \, \theta(x,10\rho/9) 
\leq \theta(x,r) e^{\lambda A(r)}
\leq \theta(x,r) e^{\lambda \varepsilon(\tau)}
\cr&
\leq [\theta(x) + \varepsilon(\tau)] \, 
e^{\lambda \varepsilon(\tau)}
\leq \theta(x) + C\varepsilon(\tau)
\cr& 
\leq e^{\lambda A(t)} \,\theta(x,t) + C \varepsilon(\tau) 
\leq e^{\lambda \varepsilon(\tau)} \,\theta(x,t) + C \varepsilon(\tau) 
\cr&
\leq \theta(x,t) + C\varepsilon(\tau)
\leq \theta(x,t) + \tau'
}\leqno (16.16)
$$
for $t < \rho$, by (16.6), (16.5), (16.12), (16.13), (16.9),
and (16.12), and if $\varepsilon(\tau)$ is small enough.

The other assumptions are satisfied because of 
(16.12), so Proposition 7.24 applies and gives a reduced minimal 
cone $\C$ with the properties (7.26)-(7.28). We take
$Z(x,\rho) = \C$, (16.14) follows from (7.26) and (7.27),
and (16.15) follows from (7.28), and because we may assume that
$(1-\tau') > 9/10$. 

We still need to show that the type of $Z(x,\rho)$
is determined by $\theta(x)$. Denote by 
$D =  H^2(Z(x,\rho) \cap B(x,1))$ the density
of $Z(x,\rho)$. Then
$$
\big| D - \theta(x,\rho) \big|
= \rho^{-2} \big| H^2(\C \cap B(x,\rho)) - H^2(E \cap B(x,\rho)) \big|
\leq \tau
\leqno (16.17)
$$
by (16.15), applied with $B(y,t)=B(\rho,r)$. 

But $\theta(x,\rho) \leq \theta(x) + C\varepsilon(\tau)$ 
by the first two lines of (16.16) 
(with $10\rho/9$ replaced with $\rho$)
and $\theta(x) \leq e^{\lambda A(\rho)} \theta(x,\rho)
\leq\theta(x,\rho) + C \varepsilon(\tau)$
by (16.9) and (16.2), we get that 
$|\theta(x,\rho) - \theta(x)| \leq C \varepsilon(\tau)$, and then
$$
| D - \theta(x) |
\leq | D - \theta(x,\rho) | + |\theta(x,r) - \theta(x)|
\leq \tau + C \varepsilon(\tau).
\leqno (16.18)
$$
So $|D - \theta(x)|$ is as small as we want. But 
$Z(x,\rho)$ is a plane as soon as $D < \pi$; by Lemma~14.12,
it is a set of type $\Bbb Y$ if $\pi < D < d_T$, and it is 
of type $\Bbb T$ if $D > 3\pi/2$; Lemma~16.11 follows.
\qed

\ms
The following lemma will allow us to deal with $P$-points
in Theorem 16.1

\ms\proclaim Lemma 16.19. For each choice of $\tau \in (0,1)$ 
there is a constant $\varepsilon_1 > 0$ such that if $x_0 \in E$ and 
$r_0 > 0$ are such that $B(x_0,9r_0) \i U$,
$h(18r_0) \leq \varepsilon_1$, 
$$
\int_0^{18r_0} h(t) \, {dt \over t} \leq \varepsilon_1 \, ,
\ \hbox{ and } \ \theta(x_0,9r_0) \leq \pi+\varepsilon_1, 
\leqno (16.20)
$$
then $B(x,r)$ is a biH\"older ball of type $P$ for $E$, 
with constant $\tau$, for every choice of $x\in E \cap B(x_0,r_0)$ and 
$0 < r \leq r_0$.

\ms
Most of the time we shall use this with $x=x_0$, but 
we shall need the additional uniformity provided by the present
statement in Section 17.
 
Let $x_0$, $r_0$, $x$, and $r$ be as in the statement.
We want to apply Corollary 15.11 to the pair $(x,r)$.
Notice however that here all the sets $Z(y,t)$ in (15.12) 
will be planes, so Corollary 15.11 comes directly from
Reifenberg's original theorem.

So we need to show that (15.12) holds for the pair $(x,r)$. 
That is, for each $y\in E \cap B(x,3r)$ and $0 < t \leq 3r$, 
we need to find a plane $Z(y,t)$ such that
$d_{y,t}(E,Z(y,t)) \leq \varepsilon$, where $\varepsilon$
comes from Corollary 15.11 and depends on $\tau$.
For the rest of the argument, we shall only need to know that 
$y\in E \cap B(x_0,4r_0)$ and $0 < t \leq 3r_0$ 
(which is true because $x\in E \cap B(x_0,r_0)$ 
and $0 < r \leq r_0$).

\smallskip
Let $\tau_1$ be small (much smaller than $\tau$, and to be chosen 
soon), and apply Lemma 16.11 to $B(x_0,9r_0)$. The assumptions
are satisfied if $\varepsilon_1 \leq \varepsilon(\tau_1)$, and we get 
that for $0 < \rho  \leq 8r_0$, there is a plane $Z(x_0,\rho)$ through
$x_0$, with the properties (16.14) and (16.15). We are only interested in 
$P = Z(x_0,8r_0)$, and (16.15) which says that
$$\eqalign{
&\big| H^2(E \cap B(y,t)) - H^2(P\cap B(y,t)) \big|
\leq 64 \tau_1 r_0^2
\cr&\hskip 2cm
\hbox{ for $y\in \R^n$ and $t>0$ such that $B(y,t) \i B(x_0,8r_0)$.}
}\leqno (16.21)
$$
Take $y\in E \cap B(x_0,4r_0)$ as above, and $t = 4 r_0$.
We get that
$$\eqalign{
\theta(y,4r_0) &= (4r_0)^{-2} H^2(E \cap B(y,4r_0))
\cr&
\leq (4r_0)^{-2} \big[ H^2(P \cap B(y,4r_0)) + 64 \tau_1 r_0^2 \big]
\cr&
 \leq \pi + 4 \tau_1 
}\leqno (16.22)
$$
by (16.21) and because $P$ is  a plane.

Notice that $\theta(y) \geq \pi$ because $y\in E$ and $E$
is reduced (see the discussion below (16.7)). 
On the other hand,
$$
\theta(y) \leq e^{\lambda A(4r_0)} \theta(y,4r_0)
\leq e^{\lambda \varepsilon_1} \, \theta(y,4r_0)
\leq e^{\lambda \varepsilon_1} \, (\pi + 4 \tau_1 )
\leq \pi + 4 \tau_1 + C \varepsilon_1 < 3\pi/2
\leqno (16.23)
$$
by (16.9) and (16.20), (16.22), and if $\tau_1$ and $\varepsilon_1$
are small enough. Hence $\theta(y)=\pi$.

Apply Lemma 16.11 to the ball $B(y,4r_0)$ and with 
$\tau = \varepsilon$, where $\varepsilon$ comes from 
Corollary~15.11 as above. The hypotheses are satisfied, 
by (16.19) and (16.22) (and if $\varepsilon_1$ and $\tau_1$ 
are small enough). We obtain that for $0 < t < 36 r_0/10$,
there is a plane $Z(y,t)$ through $y$ such that (16.15) holds.
This means that $d_{y,t}(E,Z(y,t)) \leq \varepsilon$, and this 
is exactly what was needed to apply Corollary 15.11.
This completes our proof of Lemma 16.19.
\qed

\ms
Lemma 16.19 gives the local regularity of $E$
near every $P$-point (and for this we may even take $x=x_0$).
Next we want to take care of the $Y$-points. The proof above
does not work directly in this case, because we used once or twice
the fact that the current density was close to $\pi$, and $\pi$
is the smallest possible density at a point of $E$. 
Thus it will be good to know that when $E$ is close to a set 
of type $\Bbb Y$, we have $Y$-points around to which we shall 
be able to apply Lemma 16.11 because their density is not too small.
This is the point of the next proposition.

\ms\proclaim Proposition 16.24.
There is constant $\varepsilon_2 > 0$, that depends only on $n$, 
such that if (16.3) holds, $B(x,r) \i U$, $h(2r) \leq \varepsilon_2$, 
$\int_0^{2r} h(t) dt/t \leq \varepsilon_2$, and there is a reduced 
minimal cone $Y$ of type $\Bbb Y$ centered at $x$ 
and such that $d_{x,r}(E,Y) \leq \varepsilon_2$, then 
$E \cap B(x,r/100)$ contains (at least) a $Y$-point.

\ms
The proof of Proposition 16.24 relies on Lemma 16.11
(which we use to get the local regularity of $E$ if there is no
$Y$-point), and a topological argument that counts the points
of intersection of $E$ with some $(n-2)$-dimensional spheres.
We shall do the proof in Section~17, and in the mean time we
proceed with the regularity of $E$ near a $Y$-point.

\ms\proclaim Lemma 16.25. 
For each choice of $\tau \in (0,1)$, there is a constant 
$\varepsilon_3 > 0$ such that if $x_0 \in E$ is a $Y$-point, 
$r_0 > 0$, 
$$
B(x_0,50r_0) \i U, \ \ 
h(100r_0) \leq \varepsilon_3, \ \ 
\int_0^{100r_0} h(t) {dt\over t} \leq \varepsilon_3,
\leqno (16.26)
$$
and 
$$
\theta(x_0,50r_0) \leq 3\pi/2+\varepsilon_3 \, ,
\leqno (16.27)
$$
then $B(x_0,r_0)$ is a biH\"older ball of type $\Bbb Y$
for $E$, with constant $\tau$.

\ms
It is enough to prove the lemma when $x_0=0$ and $r_0=1$, 
because it is invariant under rotations and dilations. 
Let $\tau_2$ be small, to be chosen later, and apply Lemma~16.11
to the ball $B(0,50)$; the assumptions (16.12) and (16.13) are satisfied 
if $\varepsilon_3$ is small enough. We get that for $0 < \rho \leq 45$,
there is a reduced minimal cone $Y(\rho)=Z(0,\rho)$ of type $\Bbb Y$,
centered at $0$, such that (16.14) and (16.15) hold. That is,
$$
d_{0,\rho}(E,Y(\rho)) \leq \tau_2
\leqno (16.28)
$$
and
$$\eqalign{
&\big| H^2(E \cap B(y,t)) - H^2(Y(\rho)\cap B(y,t)) \big|
\leq \tau_2 \rho^2
\cr&\hskip 2cm
\hbox{ for $y\in \R^n$ and $t>0$ such that $B(y,t) \i B(0,\rho)$.}
}\leqno (16.29)
$$
Let us check that 
$$
\theta(y,10) \leq 3\pi/2 + C \tau_2
\ \hbox{ for } y\in E\cap B(0,35).
\leqno (16.30)
$$
Set $Y = Y(45)$, and observe that
$H^2(Y\cap B(y,t)) \leq 3\pi r^2/2$ for all $y$ and $t$
because the density is largest for balls centered on the 
spine of $Y$, for instance by Proposition 5.16.
Then apply (16.29) with $\rho = 45$ and $t=10$
(so that $B(y,t) \i B(0,45)$). We get that
$$\eqalign{
\theta(y,10) &
= 10^{-2} H^2(E \cap B(y,t))
\cr&\leq 10^{-2} H^2(Y \cap B(y,t)) + 10^{-2} \tau_2 \rho^2
\cr&\leq 3\pi/2 + C \tau_2 \, ,
}\leqno (16.31)
$$
as needed for (16.30). In addition, (16.9) says that
$$
\theta(y) \leq e^{\lambda A(10)} \, \theta(y,10) 
\leq e^{\lambda \varepsilon_3} \, \theta(y,10) 
\leq \theta(y,10) + C\varepsilon_3
\leqno (16.32)
$$
by (16.5) and (16.26), so (16.31) says that
$\theta(y) \leq 3\pi/2 + C \tau_2 + C\varepsilon_3 < d_T \ $
if $\tau_2$ and $\varepsilon_3$ are small enough. By (16.8),
$y$ cannot be a $T$-point. So we just checked that
$$
\hbox{there is no $T$-point in $E\cap B(0,35)$.}
\leqno (16.33)
$$

We shall need some control on 
$$
E_Y = \Big\{ x\in E \, ; \, x \hbox{ is a $Y$-point of } E \Big\}
=\Big\{ x\in E \, ; \, \theta(x)= {3\pi \over 2} \Big\}.
\leqno (16.34)
$$

Let $\tau_3$ be small, to be chosen soon, and let us apply
Lemma 16.11 with $\tau = \tau_3$, and to the ball $B(y,10)$,
where we now assume that $y\in  E_Y \cap B(0,35)$.
If $\varepsilon_3$ and $\tau_2$ are small enough, 
the assumption (16.12) follows from (16.26), and (16.13)
follows from (16.31) because $\theta(y) = 3\pi/2$ when
$y\in  E_Y$. We get that for 
$$
y\in  E_Y \cap B(0,35) \ \hbox{ and } \ 0 < \rho \leq 9,
\leqno (16.35)
$$
there is a reduced minimal cone $Y(y,\rho)$ of type $\Bbb Y$, 
centered at $y$, such that
$$
d_{y,\rho}(E,Y(y,\rho)) \leq \tau_3
\leqno (16.36)
$$
and
$$\eqalign{
&\big| H^2(E \cap B(x,t)) - H^2(Y(y,\rho)\cap B(x,t)) \big|
\leq \tau_3 \rho^2
\cr&\hskip 2cm
\hbox{ for $x\in \R^n$ and $t>0$ such that $B(x,t) \i B(y,\rho)$.}
}\leqno (16.37)
$$

\smallskip
Notice that (16.36) gives the sort of control on $E \cap B(y,\rho)$ 
that we need in order to apply Corollary 15.11, but only in balls 
$B(y,\rho)$ that are centered on $E_Y$. 
We also need to control the balls that are centered on $P$-points.

Let $x\in E \cap B(0,3) \setminus E_Y$ be given, and set
$$
d(x) = \dist(x,E_Y).
\leqno (16.38)
$$
Observe that $d(x) \leq 3$, because the origin lies in $E_Y$
by assumption. Also, $d(x) > 0$, because the proof of Lemma 16.19
gives a neighborhood of $x$ which contains no $Y$-point (see below (16.23)).
We do not even need to know this, because if $d(x)=0$ we 
can use the cones $Y(y,\rho)$ to control $E$ in the balls centered at $x$,
so we do not need the construction below. See the small argument above
(16.42).

Pick $y\in E_Y$ such that $|y-x| \leq 11d(x)/10$, and set
$\rho = 2d(x)$. Then $y\in E_Y \cap B(0,10)$ and 
$\rho \leq 6$, so (16.35) holds and we can set $Y = Y(y,\rho)$, 
where $Y(y,\rho)$ is as in (16.36) and (16.37).
Denote by $L$ the spine of $Y$; that is, $L$ is the common 
boundary of the three half planes that compose $Y$. It is a 
line through $y$. Let us check that
$$
\dist(x,L) \geq 3d(x)/4.
\leqno (16.39)
$$

Suppose that (16.39) fails. Choose $z\in L$ so that $|z-x| \leq 3d(x)/4$, 
and apply Proposition~16.24 to the ball $B(z,d(x)/10)$. 
We need to check the assumption that 
$d_{z,d(x)/10}(E,Y) \leq \varepsilon_2$.
Indeed, $B(z,d(x)/10) \i B(y,\rho)$ because
$|z-y| \leq |z-x| + |x-y| \leq 3d(x)/4 + 11d(x)/10$ and
$\rho = 2d(x)$. Then
$d_{z,d(x)/10}(E,Y) \leq 10 d(x)^{-1} \rho d_{y,\rho}(E,Y) 
\leq 20\tau_3 < \varepsilon_2$ if $\tau_3$ is small enough,
as needed.

So Proposition~16.24 applies, and there is a point of $\xi \in E_Y$ inside
$B(z,10^{-3}d(x))$. But then $|\xi - x| \leq 10^{-3}d(x) + |z-x|
\leq [10^{-3}+3/4] \, d(x) < d(x)$. This is incompatible  with the definition
of $d(x)$, and so (16.39) holds.

Notice that $B(x,d(x)/2)$ is contained in $B(y,\rho)$
(because $|y-x| \leq 11d(x)/10$ and $\rho = 2d(x)$), so (16.37) says that
$$\eqalign{ 
\theta(x,d(x)/2) &= 
(d(x)/2)^{-2} H^2(E\cap B(x,d(x)/2) 
\cr&
\leq (d(x)/2)^{-2} H^2(Y \cap B(x,d(x)/2)) + \tau_3 (d(x)/2)^{-2} \rho^2
\cr& \leq \pi  + 4 \tau_3 d(x)^{-2} \rho^2 \leq \pi + 16 \tau_3,
}\leqno (16.40)
$$
because (16.39) says that $Y$ coincides with a plane in $B(x,d(x)/2)$.
Recall that $\theta(x)=\pi$ because we assumed that $x \in E \setminus E_Y$
(and by (16.33)). Hence, if $\tau_3$ is small  enough, (16.40) allows us to 
apply Lemma 16.11 with $\tau = \varepsilon$,
where $\varepsilon$ comes from Corollary 15.11. 
We get that for $0 < r < d(x)/3$, there is a plane $Z(x,r)$ through
$x$ such that $d_{x,r}(E,Z(x,r)) \leq \varepsilon$, as in (16.14).

\ms
We are now ready to deduce Lemma 16.25 from Corollary 15.11.
We want to get the constant $\tau$ in the conclusion of Lemma 16.25, 
so we need to prove (15.12), with a constant $\varepsilon$ that depends 
on $\tau$. For each choice of $x\in E \cap B(0,3)$ and $0 \leq r \leq 3$, 
we need to find a minimal cone $Z(x,r)$ such that 
$$
d_{x,r}(E,Z(x,r)) \leq \varepsilon.
\leqno (16.41)
$$

When $x \in E_Y$, the pair $(x,r)$ satisfies (16.35), so we can use
$Z(x,r) = Y(x,r)$, and (16.41) follows from (16.36) if $\tau_3 < 
\varepsilon$.  Note that $Z(0,3)$ is a set of type $\Bbb Y$ centered
at the origin, as required by the last assumption of Corollary 15.11.

When $x\in E \cap B(0,3) \setminus E_Y$ and $r < d(x)/3$,
we can use the set $Z(x,r)$ that was found a few lines ago,
below (16.40). 

We are left with the case when 
$x\in E \cap B(0,3) \setminus E_Y$ and $d(x)/3 \leq r \leq 3$.
Our argument also applies when $d(x)=0$, if we pretend that we 
don't know that $x\in E_Y$ when $d(x)=0$.
Let $y\in E_Y$ be such that $|y-x| \leq d(x)+r$, and set
$\rho =  2r+d(x)$. Observe that 
$|y| \leq |x| + d(x) + r \leq 9$ because $x \in B(0,3)$, 
$d(x) \leq \vert x \vert$, and $r \leq 3$, and also $\rho  \leq 9$, 
so (16.35) holds, and we can set $Z(x,r) = Y(y,\rho)$.
Now $B(x,r) \i B(y,\rho)$ because $\rho \geq r + |y-x|$,
so
$$
d_{x,r}(E,Z(x,r)) \leq r^{-1} \rho \, d_{y,\rho}(E,Z(x,r))
\leq r^{-1} \rho \, \tau_3 \leq 5 \tau_3
\leq \varepsilon
\leqno (16.42)
$$
by (16.36) and if $\tau_3$ is small enough. This completes 
our proof of Lemma 16.25.
\qed

\ms 
The next lemmas are a little bit of a digression, and are not needed
for the proof of Theorem 16.1.
We shall use a compactness argument
to replace the density assumptions in Lemmas 16.19 and 16.25
with assumptions of good approximation of $E$ by a minimal cones.

\ms\proclaim Lemma 16.43. 
For each small $\delta > 0$, there is a constant  $\varepsilon > 0$
(that depends only on $n$ and $d$) such that if
$E$ and $F$ are reduced almost minimal sets of dimension $d$
in the open set $U \i \R^n$, with gauge function $h$, and if 
$$
B(x,10r/9) \i U, \ \ 
h(20r/9) \leq \varepsilon, 
\ \hbox{ and } \ 
d_{x,10r/9}(E,F) \leq \varepsilon,
\leqno (16.44)
$$
then
$$
H^d(E \cap B(x,r)) \leq H^d(F \cap B(x,(1+\delta)r)) 
+ \delta r^d.
\leqno (16.45)
$$

\ms
Suppose that the lemma fails, and choose for each $k \geq 0$
sets $E_k$ and $F_k$ that provide a counterexample with
$\varepsilon = 2^{-k}$.
By homogeneity, we can assume that $x=0$ and $r=1$
for all $k$. Thus $E_k$ and $F_k$ are reduced almost minimal sets
in $B(0,10/9)$. It is a standard fact about the Hausdorff metric
that we can find a subsequence (which we shall still denote
by $\{ E_k \}$ to save notation) such that $\{ E_k \}$
converges to a limit $E$ locally in $B(0,10/9)$.
See for instance Proposition 34.6 on page 214 of [D2]. 
Since $d_{0,10/9}(E_k,F_k) \leq 2^{-k}$ by (16.44), 
$\{ F_k \}$ also converges to $E$ locally in $B(0,10/9)$
(and we don't even need to extract a second subsequence).

Our assumption is that $E_k$  and $F_k$ are almost minimal
in a domain that contains $B(0,10/9)$, and with a gauge function 
$h_k$ such that $h_k(20/9) \leq 2^{-k}$.
The sequence $\{ F_k \}$ satisfies the assumptions
of Lemma 3.3 (and even of Theorem 3.4 in [D1]; 
we do not need generalized quasiminimal sets here), so
$$
H^d(E\cap B(0,1+\delta)) \leq \liminf_{k \to +\infty} 
H^d(F_k\cap B(0,1+\delta))
\leqno (16.46)
$$
by (3.3), applied with $V=B(0,1)$. Similarly, $\{ F_k \}$ satisfies the 
assumptions of Lemma~3.12, with $M=1$ and the number $h$ as small 
as we want. This is because (the end of) Remark 4.5 says that 
$F_k \in GAQ(M,\delta,U,h)$ with $M=1$, $\delta = 20/9$, $U = B(0,10/9)$,
and $h=h_k(20/9)$. Then (3.13) with $H = \overline B(0,1)$ says that
$$
\limsup_{k \to +\infty} H^d(E_k  \cap\overline B(0,1))
\leq H^d(E \cap \overline B(0,1)).
\leqno (16.47)
$$
Altogether, 
$\limsup_{k \to +\infty} H^d(E_k  \cap  B(0,1))
\leq \liminf_{k \to +\infty} H^d(F_k\cap B(0,1+\delta))$,
and (16.45) holds for $k$ large. This contradiction completes
our proof of Lemma 16.43.
\qed

\ms
Here is a variant of Lemma 16.19.

\ms
\proclaim Lemma 16.48.
For each choice of $\tau \in (0,1)$, there is a constant 
$\varepsilon > 0$ such that if $x_0 \in E$, $r_0 > 0$, 
$$
B(x_0,10r_0) \i U, \ \ 
h(20r_0) \leq \varepsilon, \ \ 
\int_0^{18r_0} h(t) {dt\over t} \leq \varepsilon,
\leqno (16.49)
$$
and there is a plane $P$ such that
$$
d_{x_0,10r_0}(E,P) \leq \varepsilon \, ,
\leqno (16.50)
$$
then for every choice of $x\in B(x_0,r_0)$
and $0 < r \leq r_0$, $B(x,r)$ is a biH\"older ball of type $\Bbb P$
for $E$, with constant $\tau$.

\ms
Indeed, let $\varepsilon_1$ be as in Lemma~16.19.
If $\varepsilon$ is chosen small enough we can apply
Lemma~16.43 to $B(x_0,9r_0)$, $\delta = \varepsilon_1/10$,
and $F= P$, and get that $\theta(x_0,9r_0) \leq \pi + \varepsilon_1$. 
Then the conclusion follows from Lemma 16.19
(if in addition $\varepsilon \leq \varepsilon_1$).
\qed

\ms
Now we give a variant of Lemma 16.25.

\ms\proclaim Lemma 16.51.
For each choice of $\tau \in (0,1)$, there is a constant 
$\varepsilon_4 > 0$ such that if $x_0 \in E$, 
$r_0 > 0$, 
$$
B(x_0,60r_0) \i U, \ \ 
h(100r_0) \leq \varepsilon_4, \ \ 
\int_0^{100r_0} h(t) {dt\over t} \leq \varepsilon_4,
\leqno (16.52)
$$
and there is a reduced minimal cone $Y$ of type $\Bbb Y$,
centered at $x_0$, and such that
$$
d_{x_0,60r_0}(E,Y) \leq \varepsilon_4 \, ,
\leqno (16.53)
$$
then $B(x_0,r_0)$ is a biH\"older ball of type $\Bbb Y$
for $E$, with constant $\tau$.

\ms
Let $x_0$ and $r_0$ be as in the statement. Let 
$\varepsilon_5$ be small, to be chosen soon, and
apply Proposition 16.24 to the ball $B(x_0, \varepsilon_5 r_0)$.
Notice that 
$$
d_{x_0,\varepsilon_5 r_0}(E,Y) 
\leq 60\varepsilon_5^{-1} d_{x_0,60 r_0}(E,Y)
\leq 60\varepsilon_5^{-1} \varepsilon_4 < \varepsilon_2
\leqno (16.54)
$$
if $\varepsilon_4$ is small enough, and the other assumptions
of Proposition 16.24 are satisfied more trivially, so there is a point
$x_1 \in E \cap B(x_0, \varepsilon_5 r_0)$ which is of type $\Bbb Y$.
If $\varepsilon_4$ is small enough, (16.53) and (16.52) allow us to 
apply Lemma 16.43 to $B(x_1,500r_0/9)$ and get that
$$
\theta(x_1,50r_0) \leq 3\pi/2 + \varepsilon_3 \, ,
\leqno (16.55)
$$
where $\varepsilon_3$ is as in Lemma 16.25.

The pair $(x_1,50 r_0)$ satisfies the hypotheses of
Lemma 16.25, so $B(x_1,r_0)$ is a bih\"older ball of
type $\Bbb Y$. Since this is not exactly the conclusion
of Lemma 16.43, let us say a little more. The proof of 
Lemma 16.25 gives that for $x \in E \cap B(x_1,3r_0)$
and $0 < r \leq 3 r_0$, there is a reduced minimal cone
$Z(x,r)$ such that $d_{x,r}(E,Z(x,r)) \leq \varepsilon$
(see (16.41) and below), and then we use Corollary 15.11
to conclude. But if $\varepsilon_5$ is small enough 
(and $\varepsilon_5 = 10^{-2}$ would be more
than enough), we also get $Z(x,r)$ for $x\in E \cap B(x_0,3r_0)
\i B(x_1,3r_0+\varepsilon_5 r_0)$. In addition, (16.44) allows us
to use $Y$ itself (if $\varepsilon_4$ is small enough) as the set
$Z(x_0,r_0)$, so the last condition in Corollary 15.11 is also 
fulfilled. Finally Corollary 15.11 says that $B(x_0,r_0)$ 
also is a bih\"older ball of type $\Bbb Y$, and Lemma 16.51
follows.
\qed

\ms 
Return to the slow proof of Theorem 16.1. Lemmas 16.19
and 16.25 show that $E$ has some biH\"older regularity near
$P$-points and $Y$-points. We deal with remaining case of $T$-points
now.

\ms\proclaim Lemma 16.56. 
For each choice of $\tau \in (0,1)$, there is a constant 
$\varepsilon_6 > 0$ such that if $x_0 \in E$ 
is a $T$-point, $r_0 > 0$, 
$$
B(x_0,10r_0) \i U, \  
h(20r_0) \leq \varepsilon_6, \  
\int_0^{20r_0} h(t) {dt\over t} \leq \varepsilon_6,
\leqno (16.57)
$$
and 
$$
\theta(x_0,10r_0) \leq \theta(x_0) +\varepsilon_6
\leqno (16.58)
$$
then $B(x_0,r_0)$ is a biH\"older ball of type $\Bbb T$
for $E$, with constant $\tau$.

\ms
The proof will go as for Lemma 16.25. By rotation and
dilation invariance, we may assume that $x_0=0$ and $r_0=1$.
Let $\tau_4$ be small, to be chosen later, and apply 
Lemma~16.11 with $\tau = \tau_4$ and to $B(0,10)$. If 
$\varepsilon_6 \leq \varepsilon(\tau_4)$, we get that for
$0 < \rho \leq 9$, there is a reduced minimal cone
$T(\rho)$ of type $\Bbb T$, centered at $0$, and such that
$$
d_{0,\rho}(E,T(\rho)) \leq \tau_4
\leqno (16.59)
$$
and 
$$
\eqalign{
&\big| H^2(E \cap B(y,t)) - H^2(T(\rho)\cap B(y,t)) \big|
\leq \tau_4 \rho^2
\cr&\hskip 2cm
\hbox{ for $y\in \R^n$ and $t>0$ such that $B(y,t) \i B(0,\rho)$.}
}\leqno (16.60)
$$
as in (16.14) and (16.15).

We want to control the density of $T(\rho)$ in small balls.
Set $\eta = 10^{-1}{\rm Min}(\eta_0,l_0)$, where 
$\eta_0$ is as in (14.10) and $l_0$ comes from 
Proposition 14.1. It will be useful to know that for 
$\rho > 0$ and $x\neq 0$,
$$
T(\rho) \hbox{ coincides with a cone of type $\Bbb Y$ in } 
B(x,\eta |x|).
\leqno (16.61)
$$
Indeed suppose that $B(x,\eta |x|)$ meets two of the faces
that compose $T(\rho)$. Call $F_1$ and $F_2$ these faces, and
pick $y_j$ in $F_j \cap B(x,\eta |x|)$. Set $z_j = y_j/|y_j|$;
thus $z_j$ lies in the arc of great circle 
$C_j = F_j \cap \partial B(0,1)$.
Observe that $|z_1-z_2| < 3 \eta < \eta_0$, so (14.10) 
says that the $C_j$ have a common endpoint $a\in B(z_1,3\eta)$. 
If $B(x,\eta |x|)$ meets some other face $F_3$ of $T(\rho)$, 
the corresponding arc $C_3$ also has a common endpoint with $F_1$
in $B(z_1,3\eta)$. This endpoint is $a$, because the length of $C_1$
is at least $l_0$. In other words, all the faces that meet 
$B(x,\eta |x|)$ are contained in the set of type $\Bbb Y$
which coincides with $T(\rho)$ in a small neighborhood of $a$;
(16.61) follows.

Next let $x\in E \cap B(0,8) \setminus \{ 0 \}$ be given,
and set $\rho(x) = 9 |x| / 8$. Thus $T(\rho(x))$ is defined,
and
$$\eqalign{
\theta(x,\eta |x|) &= (\eta |x|)^{-2} H^2(E \cap B(x,\eta |x|))
\cr&
\leq (\eta |x|)^{-2} H^2(T(\rho(x)) \cap B(x,\eta |x|))
+ (\eta |x|)^{-2}\tau_4 \rho(x)^2
\cr& 
\leq 3 \pi/2 + 2 \eta^{-2}\tau_4
}\leqno (16.62)
$$
by (16.60), because $B(x,\eta |x|) \i B(0,\rho(x))$, and by 
(16.61). Then 
$$
\theta(x) \leq e^{\lambda A(\eta |x|)} \theta(x,\eta |x|)
\leq e^{\lambda \varepsilon_6} [3 \pi/2 + 2 \eta^{-2}\tau_4]
\leq 3 \pi/2 + C \varepsilon_6 + 2 \eta^{-2}\tau_4
\leqno (16.63)
$$
by (16.9), (16.57), and (16.62). In particular, $\theta(x) < d_T$
if $\varepsilon_6$ and $\tau_4$ are small enough, and
$$
0 \hbox{ is the only $T$-point of } B(0,8).
\leqno (16.64)
$$

Next assume that our $x\in E \cap B(0,8)$ is a $Y$-point.
Let $\tau_5$ be small, and apply Lemma~16.11 with the constant
$\tau_5$, and to the ball $B(x,\eta |x|)$. The assumption
(16.12) follows from (16.57), and (16.13) follows from (16.62)
because $\theta(x) = 3\pi/2$ (if $\varepsilon_6$ and $\tau_4$ 
are small enough, depending on $\tau_5$). We get that for
$0 < r \leq 9 \eta |x|/10$, there is a reduced minimal cone
$Y(x,r)$ of type $\Bbb Y$, centered at $x$, and such that
$$
d_{x,r}(E,Y(x,r)) \leq \tau_5
\leqno (16.65)
$$
and 
$$
\eqalign{
&\big| H^2(E \cap B(y,t)) - H^2(Y(x,r)\cap B(y,t)) \big|
\leq \tau_5 r^2
\cr&\hskip 2cm
\hbox{ for $y\in \R^n$ and $t>0$ such that $B(y,t) \i B(x,r)$.}
}\leqno (16.66)
$$

We also need to approximate $E$ near the points of type $P$.
So let $z\in E \cap B(0,3)$ be a $P$-point, and set
$d(z) = \dist(z,\{ 0 \} \cup E_Y)$. Let us show that
for $0 < r < 9\eta d(z)/10$,
$$
\hbox{there is a plane $P(z,r)$ through $z$ such that 
$d_{z,r}(E,P(z,r)) \leq \varepsilon$.}
\leqno (16.67)
$$

Let us first assume that
$$
d(z) \geq \eta |z|/2.
\leqno (16.68)
$$

Set $\rho = 2|z| \leq 6$, and $T = T(\rho)$.
We shall use $T$ to compute the density of $E$ near $z$. 
Denote by $L$ the spine of $T$ 
(the cone over the the union of the extremities of the arcs of circles 
that compose $T \cap \partial B(0,1)$). Let us show that 
$$
\dist(z,L) \geq d(z)/2.
\leqno (16.69)
$$
Suppose instead that we can find $\xi \in L$,
with $|\xi - z| = \dist(z,L) < d(z)/2$. 
Set $t=\eta d(z)/4$; we want to apply 
Proposition 16.24 to $B(\xi,t)$, so we need a set of type $\Bbb Y$.
By (16.61), there is a cone $Y$ of type $\Bbb Y$ that coincides 
with $T$ in $B(\xi,\eta|\xi|)$. Observe that $Y$ is centered
at $\xi$, because $\xi$ lies in the spine of $T$.
Next $|\xi| \geq |z|-|\xi - z| \geq |z|-d(z)/2 \geq d(z)/2$
(because $d(z) \leq |z|$), so $t =\eta d(z)/4 \leq \eta |\xi|/2$
and $Y$ coincides with $T$ in $B(\xi,2t)$. Thus
$d_{\xi,t}(E,Y) = d_{\xi,t}(E,T)$.

In addition, 
$|\xi| \leq |z| + |\xi - z| \leq |z| + d(z)/2 \leq 3|z|/2$ 
and $\rho = 2|z|$, so $B(\xi,t) \i B(0,\rho)$. Hence
$$\eqalign{
d_{\xi,t}(E,Y) &= d_{\xi,t}(E,T) 
\leq t^{-1} \rho \, d_{0,\rho}(E,T)
\cr&
\leq t^{-1} \rho \,\tau_4  
= 8 \eta^{-1} d(z)^{-1} |z| \,\tau_4  
\leq 16 \eta^{-2} \tau_4  
}\leqno (16.70)
$$
by (16.59) and (16.68).

If $\tau_4$ is small enough, (16.70) allows us to apply 
Proposition 16.24 to $B(\xi,t)$, and get a $Y$-point $y$
in $B(\xi,t/100)$. Recall that $|\xi - z| < d(z)/2$;
then $|y-z| \leq d(z)/2 + t < d(z)$, which contradicts
the definition of $d(z)$. So (16.69) holds.

Next we check that
$$
T \hbox{ coincides with a plane in } B(z,\eta d(z)).
\leqno (16.71)
$$
Otherwise, $B(z,\eta d(z))$ meets two faces of $T$.
Let $z_1, z_2 \in B(z,\eta d(z))$ lie in different
faces $F_1$ and $F_2$, and set $w_i = z_i/|z_i|$.
First observe that 
$$
|z_i - z| \leq \eta d(z) \leq \eta |z|,
\leqno (16.72)
$$
so $|z_i| \geq |z|/2$, and then
$$
|w_1-w_2| \leq 2|z_1-z_2|/|z| \leq 4\eta d(z)/|z|
\leqno (16.73)
$$
because the radial projection onto $\partial B(0,1)$
is $(2/|z|)$-Lipschitz on $\R^n \setminus B(0,|z|/2)$.

Next, $w_1$ and $w_2$ lie in different arcs of circles 
$C_i = F_i \cap \partial B(0,1)$ in the description of 
$T \cap \partial B(0,1)$, and 
$$
\dist(w_1,C_2) \leq |w_1-w_2| \leq 4 \eta d(z)/|z| \leq 4 \eta < \eta_0
\leqno (16.74)
$$
because $\eta \leq \eta_0/10$ (see below (16.60)).
So (14.10) says that $C_1$ and $C_2$ have a common endpoint 
$a$ such that $|a-w_1| \leq |w_1-w_2| \leq 4 \eta d(z)/|z|$.
Of course $|z_1| a$ lies in the spine $L$, so
$$
\dist(z_1,L) \leq \big|z_1 - |z_1| a \big|
\leq |z_1| |w_1-a| \leq 4 \eta |z_1| d(z)/|z| 
\leq 8 \eta d(z).
\leqno (16.75)
$$
This is impossible, because $z_1 \in B(z,\eta d(z))$
and $\dist(z,L) \geq d(z)/2$ (by (16.69)). So (16.71) holds, and
$H^2(T \cap B(z,\eta d(z))) \leq (\eta d(z))^2 \, \pi$.
Recall that $B(z,\eta d(z)) \i B(0,\rho)$ because 
$\rho = 2|z|$ and $d(z) \leq |z|$, so (16.60) says that
$$\eqalign{
\theta(z,\eta d(z)) 
&\leq (\eta d(z))^{-2} H^2(T \cap B(z,\eta d(z))) 
+  (\eta d(z))^{-2} \tau_4 \rho^2
\cr&
\leq \pi + (\eta d(z))^{-2} \tau_4 \rho^2
\leq \pi + 4\eta^{-4} |z|^{-2} \tau_4 \rho^2
\leq \pi + 16 \eta^{-4} \tau_4
}\leqno (16.76)
$$
by (16.68). 

Recall that $z$ is a $P$-point; 
if $\tau_4$ is small enough, (16.76) allows us to
apply Lemma~16.11 to $B(z,\eta d(z))$, where this time
we take $\tau = \varepsilon$, with $\varepsilon$
coming from Corollary~15.11. We get that for
$0 < r \leq 9\eta d(z)/10$, there is a plane $P(z,r)$ through
$z$ such that $d_{z,r}(E,P(z,r)) \leq \varepsilon$.
In other words, (16.67) holds.

\smallskip
We are left with the case when $z \in E \cap B(0,3)$ 
is a $P$-point as before, but (16.68) fails, i.e.,
when $d(z) \leq \eta |z|/2$.

Let $x\in E_Y$ be such that 
$|x-z| \leq 11 d(z)/10 \leq \eta |z|/20$. Thus $x\in B(0,4)$.
Also set $r=3d(z)/2$. Then $r\leq 3\eta |z|/4 < 9 \eta |x|/10$,
and we have a set $Y = Y(x,r)$ of type $\Bbb Y$ that satisfies 
(16.65) and (16.66). Let $L$ denote the spine of $Y$; then
$\dist(z,L) \geq d(z)/4$, by the proof of (16.39), and so
$Y$ coincides with a plane in $B(z,d(z)/4)$. Now
$$\eqalign{
\theta(z,d(z)/4) 
&\leq (d(z)/4)^{-2} H^2(Y \cap B(z,d(z)/4)) + (d(z)/4)^{-2}\tau_5 r^2
\cr&\leq \pi + (d(z)/4)^{-2}\tau_5 \, (3d(z)/2)^2
= \pi + 36 \tau_5
}\leqno (16.77)
$$
because $B(z,d(z)/4) \i B(x,r) = B(x,3d(z)/2)$, and by (16.66).
By Lemma 16.11, we can find planes $P(z,t)$, $0 < t \leq d(z)/5$, 
such that $d_{z,t}(E,P(z,t)) \leq \varepsilon$. That is, we have
(16.67) for a slightly larger range of radii.

\smallskip
We are now ready to use our various cones to check the assumptions
of Corollary~15.11. We need to define a minimal cone $Z(x,r)$ 
for $x\in E \cap B(0,3)$ and $0 < r \leq 3$, so that
$$
d_{x,r}(E,Z(x,r)) \leq \varepsilon
\leqno (16.78)
$$
as in (15.12). 

When $x=0$, we take $Z(x,t) = T(\rho)$, and (16.78) follows from
(16.59). Notice, by the way, that $Z(0,3)$ is a set of type $\Bbb T$ 
centered at $0$, as required by Corollary 15.11 if we want a 
biH\"older ball of type $\Bbb T$.

Now suppose that $x \in E_Y$. If $0 < r < 9\eta |x|/10$, we can 
use the set $Y(x,r)$ constructed near (16.65), and (16.78) follows 
from (16.65). Otherwise, we take $Z(x,r) = T(|x|+r)$, which
is defined because $|x|+r \leq 6$. 
Notice that $B(x,r) \i B(0,|x|+r)$, so
$$\eqalign{
d_{x,r}(E,Z(x,r)) &\leq r^{-1} (|x|+r) \, d_{0,|x|+r}(E,T(|x|+r))
\cr&\leq r^{-1} (|x|+r) \, \tau_4 \leq 2 \eta^{-1} \tau_4
}\leqno (16.79)
$$
by (16.59) and because $|x| \leq 10 \eta^{-1} r/9$. Again 
(16.78) holds if $\tau_4$ is small enough.

We are left with the case when $y$ is a $P$-point.
If $0 < r < 9 \eta d(x)/10$, we can use the plane 
$P(x,r)$ of (16.67), and (16.78) follows directly
from (16.67). So we may assume that $r \geq 9 \eta d(x)/10$.

If $r+2d(x) \leq \eta |x|/3$, we choose $y\in E_Y$ such that 
$|x-y| < 2 d(x)$, and try to take $Z(x,r) = Y(y,r+2d(x))$.
[If $d(x)=0$, apply the argument below with $d(x)$ replaced with a 
very small number (compared to $r$ and $|x|$).]
Let us check that $Z(x,r)$ is defined. First, 
$|x-y| \leq 2 d(x) \leq \eta |x|/3$, so $y\in B(0,4)$
(because $x\in B(0,3)$), and also $|y| \geq 2 |x|/3$. Then
$r+2d(x) \leq \eta |x|/3 \leq \eta |y|/2$,
which is enough to define $Y(y,r+2d(x))$ (see above (16.65)).
Now
$$\eqalign{
d_{x,r}(E,Z(x,r)) &\leq r^{-1}(r+2d(x)) \, d_{y,r+2d(x)}(E,Z(x,r))
\cr&= r^{-1}(r+2d(x)) \, d_{y,r+2d(x)}(E,Y(y,r+2d(x))) 
\cr&\leq r^{-1}(r+2d(x)) \,\tau_5 \leq 3 \eta^{-1} \tau_5
}\leqno (16.80)
$$
because $B(x,r) \i B(y,r+2d(x))$, by (16.65), and because
$r \geq 9 \eta d(x)/10$. Thus (16.78) holds if $\tau_5$ is small enough.

Finally assume that $r+2d(x) > \eta |x|/3$. Take $Z(x,r) = T(r+|x|)$,
which is defined because $r \leq 3$ and $x\in B(0,3)$. Then
$$\eqalign{
d_{x,r}(E,Z(x,r)) &\leq r^{-1} (r+|x|) \, d_{0,r+|x|}(E, T(r+|x|))
\leq 2 r^{-1} (r+|x|) \, \tau_4
}\leqno (16.81)
$$
by (16.59). Since $|x| \leq 3 \eta^{-1} (r+2d(x))
\leq 3 \eta^{-1} (r + 3 \eta^{-1} r)$ because
$r \geq 9 \eta d(x)/10$, (16.78) follows from (16.81).

We completed the verification of (16.78) in all cases, 
Corollary 15.11 applies, and Lemma 16.56 follows.
\qed

\ms
Let us now see how to deduce Theorem 16.1 from our various lemmas.
Let $U$ and $E$ be as in the statement, and let $x_0\in E$ be given.
If $x_0$ is a $P$-point, we can choose $r_0$ so small that the 
assumptions of Lemma 16.19 are satisfied, and then $B(x_0,r_0)$
is a biH\"older ball of type $\Bbb P$.
If $x_0$ is a $Y$-point, we use Lemma 16.25 instead.
Otherwise, $x_0$ is a $T$-point, and we use Lemma 16.56.
This completes our proof of Theorem~16.1, except for the
fact that the two last lemmas depend on Proposition 16.24, 
which will be proved in the next section. 
\qed

\ms
\noindent {\bf Remark 16.82.} We do not know whether the obvious 
analogue of Lemma 16.43 for $T$-points holds. 
That is, suppose that the locally minimal set $E$ is very close 
to a set of type $\Bbb T$ in a ball; is there necessarily a 
point of type $\Bbb T$ near the center?
However, the answer is yes if $n=3$ and we also assume that
$E$ separates locally $\R^3$ into at least three big connected
components. See Proposition 18.1.

\bigskip
\noindent {\bf 17. A topological argument for the existence of 
$Y$-points}
\medskip

This section is devoted to the proof of Proposition 16.24.
So we are given a reduced almost-minimal set $E$ of dimension $2$
in the domain $U \i \R^n$, with a gauge function $h$ such that
(16.2) holds, and a pair $(x,r)$ such that 
$$
\hbox{$B(x,r) \i U$, $h(2r) \leq \varepsilon_2$, and
$\int_0^{2r} h(t) dt/t \leq \varepsilon_2$.}
\leqno (17.1)
$$
We also assume that there is a reduced minimal cone $Y$ of type $\Bbb Y$, 
centered at $x$ and such that 
$$
d_{x,r}(E,Y) \leq \varepsilon_2,
\leqno (17.2)
$$ 
where $\varepsilon_2$ is s small constant that depends on $n$,
and we want to prove that $E \cap B(x,r/100)$ contains
(at least) a $Y$-point.
By rotation and dilation invariance, we may assume that $x=0$
and $r=200$. Thus
$$
\int_0^{400} h(t) dt/t \leq \varepsilon_2
\ \hbox{ and } \ 
d_{0,200}(E,Y) \leq \varepsilon_2,
\leqno (17.3)
$$ 
by (17.1) and (17.2).

We shall prove the proposition by contradiction. So we assume that 
$B(0,2)$ contains no $Y$-point, and try to prove a contradiction. 

Let us check that for any given $\varepsilon > 0$, we have that
$$
H^2(E\cap B(y,t)) \leq H^2(Y\cap B(y,t)) + \varepsilon
\ \hbox{ for $y\in B(0,2)$ and 
$10^{-1} \leq t \leq 1$}
\leqno (17.4)
$$
if $\varepsilon_2$ is small enough. 

Let $\varepsilon'$ denote the small constant that we get
when we apply Lemma 16.43 with $\delta = \varepsilon/10$,
and let us check that (17.4) holds as soon as 
$\varepsilon_2 \leq \varepsilon'/1800$. Let 
$y\in B(0,2)$ and $10^{-1} \leq t \leq 1$ be given;
observe that 
$$
d_{y,10t/9}(E,Y) 
\leq {9 \cdot 200 \over 10 t} d_{0,200}(E,Y)
\leq {9 \varepsilon \cdot 200} < \varepsilon'
\leqno (17.5)
$$
because $B(10t/9) \i B(0,200)$ and by (17.3), and more trivially
$h(20t/9) \leq h(20/9) \leq \int_0^{400} h(t) dt/t \leq \varepsilon_2
\leq \varepsilon'$ by (17.3), so Lemma 16.43 applies (with $F=Y$)
and we get that
$$
H^2(E \cap B(y,t)) \leq H^2(Y \cap B(y,(1+\delta)t)) + \delta t^2
\leq  H^2(Y \cap B(y,t)) + 10 \delta t^2
\leqno (17.6)
$$
(the measure of $Y \cap B(y,(1+\delta)t)) \setminus B(y,t)$
is easy to control, in particular because $y$ lies very close
to $Y$ by (17.3)). Now (17.4) follows from (17.6) (recall that 
$\delta = \varepsilon/10$).

\ms
Let us use (17.4) to check that
$$
\hbox{every point of $E \cap B(0,2)$ is a $P$-point,}
\leqno (17.7)
$$ 
or equivalently that $B(0,2)$ contains no $T$-point.
Let $x\in B(0,2)$ be given; then
$$\eqalign{
\theta(x) &\leq e^{\lambda A(1)} \,\theta(x,1)
= e^{\lambda \int_0^{2} h(t) dt/t} \,\theta(x,1)
\leq e^{\lambda \varepsilon_2} \,\theta(x,1)
\cr&
\leq e^{\lambda \varepsilon_2} \, [H^2(Y \cap B(x,1)) + \varepsilon]
\leq e^{\lambda \varepsilon_2} \, (3\pi/2 + \varepsilon) < d_T
}\leqno (17.8)
$$
by (16.9), (16.5), (17.3), (17.4), because the density of a 
set of type $\Bbb Y$ never exceeds $3\pi/2$, and if
$\varepsilon$, and then $\varepsilon_2$ are small enough. 
Recall that $d_T$ is the constant from Lemma 14.12 and
that $d_T > 3 \pi/2$. Now (16.8) says that $x$ is not
a $T$-point; (17.7) follows.

\smallskip
Next choose a small constant $\tau > 0$; we claim that 
we can find $\eta \in (0,1)$ such that
$$\eqalign{
B(x,r) &\hbox{ is a biH\"older ball of type $P$ for $E$,}
\hbox{ with constant $\tau$,} 
\cr& \hskip 0.8cm
\hbox{ for every $x\in B(0,3/2)$ and } 0 < r \leq 10\eta.
}\leqno(17.9)
$$

Indeed, set $K = E \cap \overline B(0,3/2)$ and, 
for $x_0 \in K$, choose $r_0 = r(x_0)$ so small that the assumptions of 
Lemma 16.19 are satisfied. This is possible, because $\theta(x_0) = \pi$
by (17.7). Then cover $K$ by a finite number of balls $B(x_0,r(x_0))$.
Finally take $\eta$ smaller than the tenth of the smallest radius $r(x_0)$.
If $x\in K$, we can find $x_0$ as above such that $x \in B(x_0,r(x_0))$,
and then Lemma 16.19 says that $B(x,r)$ is a biH\"older ball of type $P$ 
for $E$, with constant $\tau$, for $0 < r \leq r_0$. This proves the 
claim (17.9).

\medskip
Let us first give a vague idea of the plan that we would like to 
follow to get the desired contradiction.
First, we would find a sphere $S$ of codimension $2$
in $B(0,3/2)$ that meets $E$ exactly three times.
We would deform this sphere through a family $S_t$ of spheres that
eventually collapse to a point. Then we would use (17.9) to show 
that along this deformation, the number of points of $E\cap S_t$ 
can only jump by multiples of two; this is impossible when we start
with $3$ intersections and end with none.

Things will be a little more complicated because we need to put
our objects in general position before we compute the number of 
intersections; we shall also need to relate our numbers
to degrees of mappings on spheres to get some stability.

To compute the number of points in $E \cap S_t$ (maybe after we
modify $S_t$ slightly), one way would be to define a one-parameter
family of functions $f_t$ on $S$, so that $S_t = f_t(S)$, give 
equations of $E$ locally, and compute how many points $x\in S$ 
satisfy the equations. We shall find it more convenient to proceed 
in the dual way, i.e., give equations of the sets $S_t$, and count 
the number of solutions of these equations in $E$.
It seems a little more convenient because it is apparently easier 
to localize in $E$.

I wish to thank Christopher Colin for his help in my initial struggle
with degrees, and in particular the suggestion that I use equations.

Let us fix the setting and define a sphere $S$. 
Write $\R^n = \R^3 \times \R^d$, with $d=n-3$.
By rotation invariance, we may assume that $Y = Y_0 \times \{ 0 \}$,
where $Y_0 \i \R^3$ is a vertical propeller (or $Y$-set). That is, we select
three points $a_1$, $a_2$, and $a_3$ in the unit circle
of $\R^2$, so that if $L_j$ denotes the  half line that starts from $0$
and goes trough $a_j$, the three $L_j$  make $120^\circ$ angles
at the origin, and then we set 
$Y_0 = [L_1 \cup L_2 \cup L_3] \times \R\i \R^3$.

Denote by $P$ the horizontal hyperplane 
$P = \R^2 \times \{ 0 \} \times \R^d$, and set
$S = \partial B(0,1) \cap P$. Observe that 
$S \cap Y = \{ a_1, a_2, a_3 \}$; $S$ is the 
sphere of codimension $2$ that we would like to use
in the argument suggested above, but we need to modify
it slightly near the $a_j$ if we want it to meet $E$
exactly three times. 

Fix $j \leq 3$ and use (17.3) to choose $z_j \in E$ 
such that $|z_j - a_j| \leq 300 \varepsilon_2$. 
Choose $\varepsilon$ in (17.4) smaller than $\varepsilon_1/4$, 
where $\varepsilon_1$ is as in Lemma 16.19 
(and where, as before, the small constant $\tau$
just needs to be reasonably small, depending on the geometry). Then 
$\theta(z_j,1/2) \leq 4 H^2(Y \cap z_j,1/2) + 4 \varepsilon
\leq \pi + 4 \varepsilon$ if $\varepsilon_2$ is small 
enough, by (17.4) and because $B(z_j,1/2)$ only meets one 
face of $Y$. Then we can apply Lemma 16.19 to the pair
$(z_j,1/18)$, and we get that 
$$
B(z_j,r)
\hbox{ is a biH\"older ball of type $P$ for $E$ for } 
0 < r \leq 1/18.
\leqno(17.10)
$$

\ms 
We start with a rough model, which will be used as a
pattern, where the $S_t$ is $(d+1)$-dimensional sphere with 
radius $(1-t)$, and centered at $b_t = t(a_1+a_2)/100$, say. 
We choose a final center $b_1=(a_1+a_2)/100$ purposefully away 
from $Y$, to be sure that $S_t$ does not meet $E$ for $t$ close 
to $1$. A formula for the pair of equations of $S_t$ is $f_t(x)=0$, 
where
$$
f_t(x) = (x_2, |x-b_t|^2 - (1-t)^2).
\leqno(17.11)
$$
Let us record the fact that
$$
|f_1(x)| \geq |x-b_1|^2 \geq \dist(b_1,E)^2 \geq 10^{-5}
\ \hbox{ for } x\in E,
\leqno(17.12)
$$
by (17.3).

We shall not use the one-parameter family $\{ f_t \}$ to go from 
$f_0$ to $f_1$, and instead we shall use a finite number of 
intermediate functions,  so that each time the modifications occur 
on a small biH\"older ball for $E$. So we shall use a partition of unity.
Denote by $\{ x_j \}$, $1 \leq j \leq l$, a maximal collection of
points in $E \cap B(0,3/2)$, with the property that 
$|x_i-x_j| \geq \eta$ when $i \neq j$. Let $\widetilde\varphi_j$
be a smooth bump function supported in $B(x_j,2\eta)$ and such that
$\widetilde\varphi_j(x) = 1$ for $x\in \overline B(x_j,\eta)$
and $0 \leq \widetilde\varphi_j(x) \leq 1$ everywhere. 
Notice that $\sum_j \widetilde\varphi_j(x) \geq 1$ for 
$x\in E \cap B(0,3/2)$, because $x$ lies in some
$\overline B(x_j,\eta)$ by maximality of the family $\{ x_j \}$.
Also let $\widetilde\varphi_0$ be a smooth function on $\R^n$ such that
$\widetilde\varphi_0(x) = 0$ when $|x| \leq 3/2-\eta$, 
$\widetilde\varphi_0(x) = 1$ when $|x| \geq 3/2$,
and $0 \leq \widetilde\varphi_0(x) \leq 1$ otherwise.
Now $\sum_{j=0}^l \widetilde\varphi_j(x) \geq 1$ on $E$,
and we set 
$$
\varphi_j(x) = \widetilde\varphi_j(x) \,\Big\{
\sum_{j=0}^l \widetilde\varphi_j(x)\Big\}^{-1}
\ \hbox{ for $x\in E$ and $0 \leq j \leq l$.}
\leqno(17.13)
$$
Thus the $\varphi_j$ keep the same support properties as the
the $\widetilde\varphi_j$, and in particular
$$
\varphi_j
\hbox{ is supported in $B(x_j,2 \eta)$ for $j \geq 1$,}
\leqno(17.14)
$$
and in addition
$$
\sum_{j=0}^l \varphi_j(x) = 1 
\hbox{ for $x\in E \, , \ \ $ and } \ \ 
\sum_{j=1}^l \varphi_j(x) = 1
\hbox{ for $x\in E \cap B(0,3/2-\eta)$.}
\leqno(17.15)
$$
(because $\varphi_0(x)=0$ there).
Our first approximation of a discrete path from $f_0$ to $f_1$
is given by the maps
$$
g_k = f_0 + \sum_{0 < j \leq k} \varphi_j \, (f_1-f_0),
\leqno (17.16)
$$
with $0 \leq k \leq l$. Thus $g_0 = f_0$ and
$$
g_l(x) = f_1(x) \hbox{ for }x\in E \cap B(0,3/2-\eta),
\leqno (17.17)
$$
by (17.15). Notice that for $k \geq 1$,
$$
g_k(x) - g_{k-1}(x) = \varphi_k(x) (f_1(x)-f_0(x)) 
\ \hbox{ is supported in $B(x_k,2 \eta)$,}
\leqno (17.18)
$$
by (17.14).
Recall that we would like to compute the number of solutions in 
$E\cap B(0,3/2)$ of the systems $g_k(x)=0$.
This will be easier to do with the small modifications
of $f_0$ and the $g_k$ that we shall construct next, so that
they only have a finite number of simple roots. 
We start with a modification of $f_0$

\ms\proclaim Lemma 17.19.
There is a continuous function $h_0$ on $E$, such that
$$
|h_0(x) - f_0(x)| \leq 10^{-6} \ \hbox{ for } x\in E,
\leqno (17.20)
$$
$h_0$ has exactly three zeroes $\xi_1$, $\xi_2$,
and $\xi_3$ in $E$, one in each $B(a_i,1/60)$, 
and in addition
$$
\hbox{each $\xi_j$ is a non-degenerate simple zero of $h_0$.}
\leqno (17.21)
$$

\ms
Let us first define (17.21). We shall say that $\xi \in E$ is a 
non-degenerate simple zero of some continuous function $h$ 
on $E$ if $h(\xi) = 0$ and there is a ball $B(\xi,\rho)$ and a 
bi-H\"older map $\gamma$ from $E \cap  B(\xi,\rho)$ to an 
open set $V$ of a plane, such that $h \circ \gamma^{-1}$ is 
of class $C^1$ on $V$ and the differential of $h \circ \gamma^{-1}$
at $\gamma(\xi)$ is of rank $2$.

Notice that when this holds, there is a small neighborhood $W$
of $\gamma(\xi)$ such that $h \circ \gamma^{-1}$ is a 
homeomorphism from $W$ to its image in $\R^2$, so
if we replace $\rho$ with a sufficiently small radius $t>0$,
we get that
$$\eqalign{
&\hbox{the restriction of $h$ to }
E\cap B(\xi,t)  \hbox{ is a homeomorphism }
\cr & \hskip 1.1cm \hbox{from $E\cap B(\xi,t)$ to the open subset }
h(E\cap B(\xi,t)) \hbox{ of } \R^2.
}\leqno (17.22)
$$

We now prove the lemma. We shall only modify $f_0$
near the three points $a_j$ of $S\cap Y$.
Recall from (17.10) that $B(z_j,1/18)$ is a biH\"older ball of 
type $\Bbb P$. Thus there is a plane $P_j$ centered at $z_j$ and 
a mapping $\psi_j : B(z_j,1/9) \to \R^n$, with properties
like (15.6)-(15.9) (see Definition~15.10). 
In particular, $|\psi_j(x)-x| \leq \tau/18$ for
$x\in B(z_j,1/9)$ and
$$
E \cap B(z_j,10^{-1}) \i \psi_j(P_j \cap B(z_j,1/9)) \i E.
\leqno (17.23)
$$

We keep $h_0 = f_0$ out of the three $B(z_j,1/30)$.
In $B(z_j,1/60)$, we replace $f_0$ with 
$h_0 = f_0 \circ \psi_j^{-1}$. In the intermediate
regions, we interpolate, i.e., we set 
$$
h_0(x)= \alpha(x) f_0(x) + (1-\alpha(x)) f_0 \circ \psi_j^{-1}(x),
\leqno (17.24)
$$
with $\alpha(x) = 60|x-z_j|-1$. Notice that
$|h_0(x) - f_0(x)| \leq |f_0(x) - f_0 \circ \psi_j^{-1}(x)|
\leq C \tau$ for $x \in B(z_j,1/30)$, because 
$|\psi_j(x)-x| \leq \tau/18$ and the derivative of $f_0$ is bounded
there; thus (17.20) holds.

Recall from (17.11) that $f_0(x) = (x_2,|x|^2 - 1)$; 
then $|f_0(x)| \geq 1/500$ for 
$x\in Y \setminus \cup_j B(x_j,10^{-2})$
(essentially, because $S \cap Y = \{ a_1, a_2, a_3 \}$
and $|z_j-a_j| \leq 300\varepsilon_2$ (see above (17.10)). 
This implies that 
$$
|h_0(x)| \geq 10^{-3} \hbox{ for }
x\in E \cap B(0,2) \setminus \cup_j(B(z_j,1/60)),
\leqno (17.25)
$$
because $f_0$ varies slowly, and by (17.3) and (17.20). 
Thus $h_0$ only has zeroes in the $B(z_j,1/60)$.

We still need to check that $h_0$ has exactly one zero in
each $B(z_j,1/60)$, and that (17.21) holds.
Set $\gamma_j(x) = \psi_j^{-1}(x)$ for
$x\in B(z_j,1/60)$. By (17.23) and because $\psi_j$ is biH\"older,
$\gamma_j$ is a homeomorphism from $E \cap B(z_j,1/60)$
to its image, which is an open subset of $P_j$.

Recall that $h_0 = f_0 \circ \psi_j^{-1} = f_0 \circ \gamma_j$
in $B(z_j,1/60)$. Thus, for $\xi \in E \cap B(z_j,1/60)$,
$\xi$ is a zero of $h_0$ if and only if $\gamma_j(\xi)$
is a zero of $f_0(x) = (x_2,|x|^2 - 1)$ in $P_j$. So we study
$f_0$ in $P_j$.

By (17.3), we can find $z'_j \in E$ such that 
$|z'_j - (1+10^{-2}) \, a_j| \leq 300 \varepsilon_3$, and
$z''_j \in E$ such that 
$|z''_j - a_j - 10^{-2} e_2| \leq 300 \varepsilon_3$,
where $e_2 = (0,1,0, \cdots ,0)$. Set
$w=\gamma_j(z_j)$, $w'=\gamma_j(z'_j)$, and $w''=\gamma_j(z''_j)$.
These three points lie in $P_j$, so $P_j$ is the plane
that contains $w$ and the two vectors
$v_1 = 100(w'-w)$ and $v_2 = 100(w''-w)$.

Notice that $|w-z_j| \leq \tau/18$, by (15.6), and similarly
for $w'$ and $w''$. Then 
$|f_0(w)| \leq 10 |w-a_j| \leq C \varepsilon_3 + C\tau$,
because $f_0(a_j)=0$. 

Since $|v_1-a_j| \leq C\varepsilon_3 + C\tau$ and
$|v_2-e_2| \leq C\varepsilon_3 + C\tau$, when we can compute the
differential of $f_0$ in $P_j$, we can replace $v_1$ with $a_j$
and $v_2$ with $e_2$, and make errors smaller than 
$C\varepsilon_3 + C\tau$. But 
$D_{a_j} f_0(x) = (0,2\langle x,a_j \rangle)$ and
$D_{e_2} f_0(x) = (1,2x_2)$. 
These stay within $10^{-1}$ of $(0,2)$ and $(1,0)$ when
$x\in P_j \cap B(a_j,1/50)$, so, 
if $\varepsilon_3$ and $\tau$ are small enough, the differential of 
$f_0$ in $P_j \cap B(a_j,1/50)$ is invertible and almost constant,
and $f_0$ has a unique zero $\zeta_j$ in $P_j \cap B(a_j,1/50)$,
which even lies in $B(a_j,10^{-3})$ (recall that
$|f_0(w)| \leq C \varepsilon_3 + C\tau$).

Thus $h_0 = f_0 \circ \psi_j^{-1}$ has a unique zero in $B(z_j,1/60)$,
namely, $\xi_j = \psi_j(\zeta_j)$. In addition, (17.21) holds, precisely because
$f_0 = h_0 \circ \psi_j = h_0 \circ \gamma_j^{-1}$ is $C^1$ on
$P_j$, with a differential of full rank at $\zeta_j$. This completes 
our proof of Lemma 17.19.
\qed

\ms\proclaim Lemma 17.26. 
We can find continuous functions $\theta_k$,
$1 \leq k \leq l$, such that 
$$
\theta_k \hbox{ is supported in } B(x_k,3\eta),
\leqno (17.27)
$$
$$
||\theta_k||_\infty \leq 2^{-k} 10^{-6},
\leqno (17.28)
$$
and, if we set
$$
h_k = h_{k-1} + \varphi_k \, (f_1-f_0) + \theta_k
\leqno (17.29)
$$
for $1 \leq k \leq l$, then 
$$\eqalign{
&\hbox{each $h_k$ has only a finite number of zeroes in $E$,}
\cr&\hskip 2.4cm
\hbox{ which are all non-degenerate simple zeroes of $h_k$.}
}\leqno (17.30)
$$

\ms
Thus we keep the same formula as for $g_k$ in 
(17.18), but we add small perturbations $\theta_k$ to
obtain (17.30).

We shall construct the $h_k$ by induction. 
Let $k \geq 1$ be given, and suppose that we already constructed 
$h_{k-1}$ so that (17.30) holds. Notice that this  is
the case when $k=1$, by Lemma 17.19.

Notice that $h_{k-1} + \varphi_k \, (f_1-f_0)$ coincides with
$h_{k-1}$ out of $B(x_k,2\eta)$, by (17.14). 
We select a very thin annulus 
$$
A = \overline B(x_k,\rho_2) \setminus B(x_k,\rho_1),
\hbox{ with } 2\eta < \rho_1 < \rho_2 < 3 \eta, 
\leqno (17.31)
$$
that does not meet the finite set of zeroes of $h_{k-1}$. 
Recall from (17.9) that $B(x_k,10\eta)$ is a biH\"older ball
of type $\Bbb P$, so there is a biH\"older mapping 
$\psi_k : B(x_k,20\eta) \to \R^n$ and a plane $P_k$ through
$x_k$ such that $|\psi_k(x) - x| \leq 10\eta \tau$ for 
$x\in B(x_k,20\eta)$, as in (15.6), and 
$$
E \cap B(x_k,19 \eta) \i \psi_k(P_k \cap B(x_k,20\eta)) \i E,
\leqno (17.32)
$$
as in (15.9). 

We shall take $\theta_k$ supported in $B(x_k,\rho_2)$,
so (17.27) will hold, and also the desired control on the 
zeroes of $h_k$ out of $B(x_k,\rho_2)$ will come from the
induction assumption on $h_{k-1}$, since $h_k = h_{k-1}$ out of 
$B(x_k,\rho_2)$.

We shall also take $||\theta_k||_\infty$ smaller than 
$2^{-k} 10^{-6}$, so that (17.28) holds, and smaller than
$\inf_{x\in A} |h_{k-1}(x)| > 0$, so that (17.29) (with the fact that
$0 \leq \varphi_k(x) \leq 1$ everywhere) implies that $h_k$
has no zero in $A$.

We still need to control $h_k$ in $B(x_k,\rho_1)$.
Set $\gamma(x) = \psi_k^{-1}(x)$ for $x\in E \cap B(x_k,\rho_1)$.
By (17.32) and because $\psi_k$ is biH\"older-continuous on
$B(x_k,20\eta)$, $\gamma$ is a biH\"older homeomorphism
from $E \cap B(x_j,\rho_1)$ to some open subset $V$ of the 
plane $P_k$, and its inverse is the restriction of $\psi_k$
to $V$.

By density of the $C^1$ functions (in the space of bounded
continuous functions on $V$, with the ``sup" norm), we can choose 
$\theta_k$ with the constraints above, and such that
$$
h_k \circ \psi_k
\hbox{ is of class $C^1$ on $V$.}
\leqno (17.33)
$$
We can even keep the option of adding a very small constant $w \in \R^2$
to $\theta_k$ in $E \cap B(x_j,\rho_1)$ (and interpolating continuously 
on $A$). Let us first check that for almost every choice of $w$,
$$
h_k \hbox{ has a finite number of zeros in }
E \cap B(x_j,\rho_1)
\leqno (17.34)
$$
for the modified $h_k$. Set
$$
Z_y = \big\{ z\in V \, ; \, h_k \circ \psi_k(z) = y \big\}
\leqno (17.35)
$$
for $y\in \R^2$. By (17.33), we can apply the coarea theorem to 
$h_k \circ \psi_k$ on $V$. We get that
$$
\int_V J(z) \, dH^2(z) = \int_{y \in \R^2} \, H^0(Z_y) \, dH^2(y),
\leqno (17.36)
$$
where $J$ is a bounded Jacobian that we don't need to 
compute (apply (8.12) with $m=2$ and $d=0$). Then 
$Z_y$ is finite for almost every $y \in \R^2$. If we choose
$w$ such that $Z_{w}$ is finite and add $-w$ to $\theta_k$ in 
$E \cap B(x_j,\rho_1)$, then the new $Z_0$ will be finite, i.e., 
(17.34) will hold.

We also care about the rank of the differential. 
Sard's theorem says that the set of critical values 
of $h_k \circ \psi_k$ has vanishing Lebesgue measure in 
$\R^2$; if we choose $w \in \R^2$ such that $w$ is not
a critical value, and add $-w$ to $\theta_k$ in 
$E \cap B(x_j,\rho_1)$, then the differential of the new
$h_k \circ \psi_k$ at every zero of $h_k \circ \psi_k$ has rank $2$.

We select $w$ very small, with the properties above, and add $-w$ t
o $\theta_k$ in $E \cap B(x_j,\rho_1)$; then $h_k$ has a finite number 
of zeros in $E \cap B(x_j,\rho_1)$, and they are all 
non-degenerate simple zeroes of $h_k$ (see the definition
below Lemma 17.19 and recall that $h_k \circ \gamma^{-1}
= h_k \circ \psi_k$). So (17.30) holds, and Lemma 17.26
follows.
\qed

\ms
Denote by $N(k)$ the number of zeroes of $h_k$ in $E\cap B(0,2)$.
By Lemma 17.9, $N(0) = 3$. Let us check that for the last index $l$,
$N(l) = 0$. First observe that
$$
h_l - h_0 = \sum_{1 \leq k \leq l} (h_k - h_{k-1})
= \sum_{1 \leq k \leq l} \varphi_k \, (f_1-f_0)
+ \sum_{1 \leq k \leq l} \theta_k
\leqno (17.37)
$$
by (17.29). If $x\in E \cap B(0,4/3)$, (17.15) says that
$\sum_{1 \leq k \leq l} \varphi_k (x) = 1$, so
$$
h_l(x) = h_0(x) + f_1(x)-f_0(x) + \sum_{1 \leq k \leq l} \theta_k(x)
\leqno (17.38)
$$
and 
$$
|h_l(x)| \geq |f_1(x)| - |h_0(x)-f_0(x)| 
- \sum_{1 \leq k \leq l} |\theta_k(x)|
\geq 10^{-5} - 10^{-6} - \sum_{1 \leq k \leq l}  2^{-k} 10^{-6} > 0
\leqno (17.39)
$$
by (17.12), (17.20), and (17.28).

If $x\in E\cap B(0,2) \setminus B(0,4/3)$, (17.15) says that
$\sum_{1 \leq k \leq l} \varphi_k (x) = 1-\varphi_0(x)$, so
$$
h_l(x) = h_0(x) + [1-\varphi_0(x)][f_1(x)-f_0(x)] 
+ \sum_{1 \leq k \leq l} \theta_k(x),
\leqno (17.40)
$$
and
$$
\vert h_l(x) - f_0(x) - [1-\varphi_0(x)][f_1(x)-f_0(x)] \vert
\leq |h_0(x)-f_0(x)| + \sum_{1 \leq k \leq l} |\theta_k(x)|
\leq 2 \cdot 10^{-6}
\leqno (17.41)
$$
by (17.20) and (17.28). But the second coordinate of 
$f_0(x) + [1-\varphi_0(x)][f_1(x)-f_0(x)]$ is
$$
|x|^2 - 1 + [1-\varphi_0(x)] [|x-b_1|^2 - |x|^2 + 1]
= \varphi_0(x) [|x|^2 - 1] + [1-\varphi_0(x)]  |x-b_1|^2
\geq 1/3
\leqno (17.42)
$$
by (17.11) and because $|x| \geq 4/3$. So $h_l(x) \neq 0$
in this case as well, $h_l$ has no zero in $E\cap B(0,2)$, 
and $N(l) = 0$.

We shall reach the desired contradiction as soon 
as we prove that
$$
N(k) - N(k-1) \hbox{ is even for } 1 \leq k \leq l.
\leqno (17.43)
$$

This is the moment where we shall use some degree theory.
Since it is easier to study the degree for mappings between spheres,
we first replace $h_k$ and $h_{k-1}$ with mappings
between $2$-spheres.

Recall that $h_{k-1}$ does not vanish on $A$
(by definition of $A$), and that we chose $||\theta_k||_\infty$ 
so small that this stays true for $h_k$ (see (17.29) and recall
that $\varphi_k = 0$ on $A$, by (17.31) and (17.14)). Set
$$
m_t(x) = h_{k-1}(x) + t [ h_{k}(x) + h_{k-1}(x)] 
\hbox{ for $x\in E \cap \overline B(x_k,\rho_2)$ and $0 \leq t \leq 1$.}
\leqno (17.44)
$$
Thus $m_0 = h_{k-1}$ and $m_1 = h_k$ on 
$E \cap \overline B(x_k,\rho_2)$.
Notice that $m_t(x) = h_{k-1}(x) + t \theta_k(x)$ 
for $x\in E \cap A$ and $0 \leq t \leq 1$, by (17.29), so
$m_t(x) \neq 0$ if we took $||\theta_k||_\infty$ small enough.

Choose $\beta_k > 0$ such that $|m_t(x)| \geq \beta_k$ 
for $x\in E \cap A$ and $0 \leq t \leq 1$.
Set $S_{\infty} = \R^2 \cup \{ \infty \}$, which we see as a 
$2$-dimensional sphere, and define $\pi : \R^2 \to S_{\infty}$
by
$$
\pi(x) = \infty \ \hbox{ if } |x| \geq \beta_k
\ \hbox{ and } \ \ \pi(x) = {x \over \beta_k - |x|} 
\ \hbox{ otherwise.}
\leqno (17.45)
$$
Then set 
$$
p_t(x) = \pi(m_t(x))  
\ \hbox{ for $x\in E \cap \overline B(x_k,\rho_2)$ and $0 \leq t \leq 1$.}
\leqno (17.46)
$$
Notice that $p_t(x)$ is a continuous function of $x$
and $t$, with values in $S_\infty$. By definition of $\beta_k$,
$$
p_t(x) = \infty
\hbox{ for $x\in E\cap A$ and $0 \leq t \leq 1$}.
\leqno (17.47)
$$

It will also be more convenient to replace the
domain $E \cap \overline B(x_k,\rho_2)$ with a subset of the plane.
Let us use the same mapping $\psi_k$ as above (near (17.32)),
and its inverse $\gamma$, which is a biH\"older homeomorphism
from $E \cap B(x_k,\rho_2)$ to some open subset $V'$ of $P_k$
(for the same reason as before). For $0 \leq t \leq 1$, set
$$
q_t(x) = p_t(\psi_k(x)) \ \hbox{ for } x\in V' \ \hbox{ and } \ 
q_t(x) = \infty 
\hbox{ for } x\in P_k \setminus V'.
\leqno (17.48)
$$
Let us check that $q_t$ is continuous on $P_k \times [0,1]$.
It is continuous on $V' \times [0,1]$, because $p_t$ is continuous
on $E \cap B(x_k,\rho_2)\times [0,1]$, and also on 
$[P_k \setminus \overline V']\times [0,1]$, because it is constant 
there. We are left with $\partial V' \times [0,1]$.
But for $x\in \partial V'$, $\psi_k(x)$ lies in 
$E \cap \partial B(x_k,\rho_2)$ (again by (17.32) and because $\psi_k$
is biH\"older), so some neighborhood of
$\psi_k(x)$ in $\overline B(x_k,\rho_2)$ is contained in $A$, 
hence $p_t(w) = \infty$ there (by (17.47)), and $q_t(y) = \infty$ 
near $x$ (and for all $t$).

Extend $q_t$ to the sphere $S'=P_k \cup \infty$ by setting 
$q_t(\infty) = \infty$; obviously $q_t$ is still continuous
on $S' \times [0,1]$.

Now we have two mappings $q_0$ and $q_1$ from the $2$-sphere
$S'$ to the $2$-sphere $S_\infty$, so we can compute their
degrees. See for instance in [Do],  
to which we shall systematically refer because it contains
the needed information rather explicitly.
First observe that since $q_0$ and $q_1$ are homotopic, 
they have the same degree (see Proposition 4.2 on page 63
of [Do]).  

Fix $t \in \{ 0,1 \}$, and try to compute the degree of $q_t$.
First observe that the degree (which roughly speaking counts 
the total number of inverse images of any point of the target 
space $S_{\infty}$, counted with orientation and multiplicity) 
can be computed locally near any point $Q$ of the target 
space; see Corollary 5.6 on page 67 of [Do]. 
Here we shall use $Q=0$. Let us check that for $x\in S'$,
$$
\hbox{if $q_t(x)=0$, then $\psi_k(x)$ is a zero of $h_{k-1}$ if $t=0$,
and of $h_k$ if $t=1$.}
\leqno (17.49)
$$
Indeed, if $q_t(x)=0$, then $x\in V'$ (otherwise $q_t(x) = \infty$), 
hence $q_t(x) = p_t(\psi_k(x))$ (by (17.48), so $p_t(\psi_k(x))=0$, 
and then $m_t(\psi_k(x))=0$ (by (17.46) and (17.45)).
Now (17.49) follows from (17.44). 

By (17.30), both $h_{k-1}$ and $h_{k}$ have a finite number of
zeroes, so the same thing holds for $q_t$. 
Then the degree of $q_t$ (computed at $0$) can be computed 
locally: Proposition 5.8 of [Do] says that if 
we cover $S'$ by finitely many open sets $Y_l$,
so that each $Y_l$ contains exactly one zero $y_l$ of
$q_t$, then the degree of $q_t$ is the sum of the
degrees of the restrictions of $q_t$ to $Y_l$. 
In addition, Proposition 5.5 of [Do] allows us to
replace each $Y_l$ with any smaller open neighborhood 
$W_l$ of $y_l$. Let us check that if $W_l$ is small enough, 
$$
q_t \hbox{ is a homeomorphism from $W_l$ to $q_t(W_l)$.}
\leqno (17.50)
$$
Indeed, set $\xi_l = \psi_k(y_l)$; by (17.49), $\xi_l$ is a zero of 
$h = h_{k-1}$ or $h_k$, and (17.30) says that it is a non-degenerate
simple zero. This means that there is a small ball
$B(\xi_l,\rho)$ and a bi-H\"older map $\gamma_l$ from 
$E \cap  B(\xi_l,\rho)$ to an open set $V$ of a plane, 
such that $h \circ \gamma_l^{-1}$ is 
of class $C^1$ on $V$ and the differential of 
$h \circ \gamma_l^{-1}$ at $\gamma_l(\xi_l)$ is of rank $2$.
But it was observed near (17.22) that if we take $\rho$ small enough, 
the restriction of $h$ to $E\cap B(\xi,\rho)$ is a homeomorphism
from $E\cap B(\xi,\rho)$ to the open subset $h(E\cap B(\xi,\rho))$ 
of $\R^2$. Recall that $h(\xi)=0$, so
if $\rho$ is small enough, $\pi \circ h$ is a homeomorphism 
from $E\cap B(\xi,\rho))$ to $\pi \circ h(E\cap B(\xi,\rho))$.

Take $W_l = \gamma(E\cap B(x,\rho))
= \psi_k^{-1}(E\cap B(x,\rho))$ above; then 
$\pi \circ h \circ \psi_k$ is a homeomorphism from $W_l$ to 
$\pi \circ h(E\cap B(x,\rho)) \i S_\infty$.

At the same time, $y_l$ is a zero of $q_t$, so by the proof of
(17.49), $q_t(x) = p_t \circ \psi_k(x) = \pi \circ m_t  \circ \psi_k(x)
= \pi \circ h \circ \psi_k(x)$ near $y_l$
(by (17.48), (17.46), and (17.44)); (17.50) follows.

By (17.50) and Example 5.4 in [Do], the degree at $0$ of the  
restriction of $q_t$ to $W_l$ is $\pm 1$. We add up all these
degrees and get that the degree of $q_t$ (on the sphere) is equal, 
modulo $2$, to the number of zeroes of $q_t$. Since $q_0$ and $q_1$
are homotopic, they have the same degree, and $q_0$ and $q_1$ have 
the same number of zeroes modulo $2$. 

We are now ready to prove (17.43). Since $h_{k-1}$ and $h_k$
coincide out of $B(x_k,\rho_2)$ and have no zero in $E \cap A$,
we just need to consider the zeroes in $E \cap B(x_k,\rho_1)$.
If $\xi \in E \cap B(x_k,\rho_1)$ is such that $h_{k-1}(x) = 0$, 
then $p_0(x) = 0$ by (17.44) and (17.46), and 
$q_0(\psi_k^{-1}(x)) = p_0(x) = 0$ by $0$ by (17.48). So every
zero of $h_{k-1}$ gives a zero for $q_0$. The converse is true, 
by (17.49). Similarly, the number of zeroes for $h_k$ in 
$E \cap B(x_k,\rho_1)$ is the same as the number of zeroes for
$q_1$. This proves (17.43) and, as explained near (17.43),
Proposition 16.24 follows.
\qed

\bigskip
\noindent {\bf 18. $MS$-minimal sets in $\R^3$}
\medskip

In this section we return to sets of dimension $2$ in $\R^3$, 
and in particular prove Theorem~1.9, which says that every nonempty 
reduced $MS$-minimal set in $\R^3$ is a plane, a cone of type $\Bbb Y$, 
or a cone of type $\Bbb T$.

We start with a variant of Proposition 16.24 for points of type
$T$, which is valid for (standard) almost-minimal sets, but needs
an extra separation assumption.  
This proposition yields an analogue of Theorem 1.9 for $Al$-minimal sets,
again with an extra separation assumption when $E$ looks like
a minimal cone of type $\Bbb T$ near infinity (see Proposition~18.29).

We can then apply Proposition 18.29 to $MS$-minimal sets,
because the separation assumption is easy to get in this
context, and obtain Theorem 1.9.

The results of this section and the previous ones
also apply to $MS$-almost-minimal sets, but to avoid new 
discussions about definitions, we simply mention this in 
Remark 18.44 at the end of the section.

We start with a variant of Proposition 16.24 that finds
a $T$-point when the reduced almost-minimal set $E$
is close enough to a set of type $\Bbb T$ in a ball.

\ms\proclaim Proposition 18.1.
There is constant $\varepsilon_7 > 0$ with the following property.
Let $E$ be a reduced $A$-almost-minimal set of dimension $2$
in an open set $U \i \R^3$, with the gauge function $h$. 
Suppose that  
$$
B(x,r) \i U, \hskip 0.2cm
h(2r) \leq \varepsilon_7, \hskip 0.2cm
\int_0^{2r} h(t) dt/t \leq \varepsilon_7, 
\leqno (18.2)
$$
and that there is a minimal cone $T$ of type $\Bbb T$, centered at $x$,
such that 
$$
d_{x,r}(E,T) \leq \varepsilon_7.
\leqno (18.3)
$$
Denote by $a_1$, $a_2$, $a_3$, and $a_4$ the points of 
$\partial B(x,r/2)$ whose distance to $T$ is largest, and  
by $V(a_j)$ the connected component of $a_j$ in $B(x,r)\setminus E$.
Suppose in addition that at least three of the $V(a_j)$ are distinct.
Then $E \cap B(x,3r/4)$ contains at least a $T$-point.

\ms
See Definition 4.3 for $A$-almost-minimal sets, and recall that
Definition 4.8 gives an equivalent notion, by Proposition 4.10.
See Definition 2.12 for the ``reduced sets".
Also recall that in $\R^3$, there is only one sort of set of type 
$\Bbb T$: the cones based on a tetrahedron that are described near
Figure 1.2. A $T$-point is therefore a point $x\in E$ such that every 
blow-up limit of $E$ at $x$ is of this type (see the definition 
and discussion near (16.8)).

The author doesn't know whether the separation assumption can be removed. 
At least the most obvious attempt fails: see Section 19. 

\smallskip
By scale invariance, it is enough to prove the proposition when
$x=0$ and $r=2$. Our plan is to proceed by contradiction, so we shall 
assume that there is no point of type $T$ in $E\cap B(0,3/2)$ and try to 
reach a contradiction. This assumption will make it easier for us to 
study $E_Y \cap B(0,1)$, where $E_Y$ still denotes the set of 
$Y$-points of $E$.

Let $\alpha > 0$ be small, to be chosen later.
Let us first prove that if $\varepsilon_7$ is small enough,
$$
\dist(y,L) \leq 2\alpha
\ \hbox{ for } y\in E_Y \cap B(0,3/2),
\leqno (18.4)
$$
where $L$ denotes the spine of $T$ (a union of four half lines).
First observe that if $\varepsilon_7$ is small enough,
$$
H^2(E\cap B(y,\alpha)) \leq H^2(T\cap B(y,\alpha)) + 10^{-1} \alpha^2 
\ \hbox{ for $y\in E \cap B(0,3/2)$,}
\leqno (18.5)
$$
because we can apply Lemma~16.43 to $B(y,\alpha)$ with
$\delta > 0$ as small as we want. The details are the 
same as for (17.4). If $y\in E \cap B(0,3/2)$ and
$\dist(y,L) \geq 2\alpha$, then $H^2(T\cap B(y,\alpha)) \leq \pi \alpha^2$,
and (18.5) says that $\theta(y,\alpha) \leq 11\pi/10$. Then 
$$
\theta(y) \leq e^{\lambda A(\alpha)} \theta(y,\alpha)
= e^{\lambda \int_{0}^{\alpha} h(2t) {dt \over t}} \,\theta(y,\alpha)
\leq e^{\lambda \varepsilon_7} \,\theta(y,\alpha)
< 3 \pi/2
\leqno (18.6)
$$
by (16.9), (16.5), and (18.2), so $y$ is a $P$-point.
So (18.4) holds.

\medskip
Next we want to show that every point of $E_Y \cap B(0,4/3)$ has
a neighborhood where $E_Y$ is a nice curve. So let
$y\in E_Y \cap B(0,4/3)$ be given. Lemma 16.25 says that 
for $r$ small, $B(y,r)$ is a biH\"older ball of type $\Bbb Y$ for $E$ 
(with $\tau =10^{-4}$, say)
Let $r=r_y$ be such a radius, denote by $f_y: B(y,2r) \to \R^3$ 
the corresponding biH\"older mapping, and let $Y_y$ be the minimal
cone of type $Y$ centered at $y$ that goes with it. Here (15.6)-(15.9)
say that
$$
|f_y(x)-x| \leq 2\tau r \ \hbox{ for } x\in B(y,2r),
\leqno (18.7)
$$
$$
(1-\tau) |(x-z)/r|^{1+\tau} \leq |(f(x)-f(z))/r| 
\leq (1+\tau) |(x-z)/r|^{1-\tau}
\ \hbox{ for } x, z \in B(y,2r),
\leqno (18.8) 
$$
$$
B(y,19r/10) \i f(B(y,2r),
\leqno (18.9) 
$$
and
$$
E \cap B(y,19r/10) \i f(Y_y\cap B(y,2r)) \i E.
\leqno (18.10) 
$$

\ms\proclaim Lemma 18.11.
Denote by $L_y$ the spine of $Y_y$. Then
$$
E_Y \cap B(y,r/2) \i f_y(L_y \cap B(y,r)) \i E_Y \cap B(y,3r/2).
\leqno (18.12)
$$

\ms
We start with the first inclusion.  Let $z \in E_Y \cap B(y,r/2)$
be given, and set $w=f_y^{-1}(z)$.
Clearly $w\in B(y,r)$ (by (18.7)); let us assume that $w\not\in L_y$
and find a contradiction. The idea will be that this implies that
$\R^3 \setminus E$ has only two local connected components 
near $w$, while the fact that $z \in E_Y$ implies that there are three.
But we want to check this more carefully. 
Let us first check that 
$$\eqalign{
&\hbox{for each $\rho_2 > 0$, we can 
find $\rho > 0$ such that for $z_1$, $z_2$, $z_3 \in B(z,\rho) 
\setminus E$,}
\cr&\hskip 0.5cm
\hbox{we can find a path 
$\gamma \i \overline B(z,\rho_2) \setminus E$, that 
connects two of the $z_j$.}
}\leqno (18.13)
$$

So let $z_1$, $z_2$, and $z_3$ lie in $B(z,\rho) \setminus E$.
Set $w_j=f_y^{-1}(z_j)$ (these are well defined if $\rho \leq r$,
by (18.9)). By the first half of (18.8),
$(1-\tau) (|w_j-w|/r)^{1+\tau} \leq |z_j-z|/r \leq \rho/r$, 
or equivalently 
$$|w_j-w| \leq \rho_1, \ \hbox{ with }
\rho_1 = r \, [(1-\tau)^{-1}\rho/r]^{1/(1+\tau)}.
\leqno (18.14)
$$

If $\rho$ is small enough, $4\rho_1 \leq \dist(w,L_y)$, and 
$Y_y$ coincides with a plane in $B(w,2\rho_1)$. Notice  that
the $w_j$ lie out of $Y_y$, because they lie in $B(y,2r)$ and
otherwise $z_j =  f_y(w_j)$ would lie in $E$ by the second 
part of (18.10). So one of the segments $[w_j,w_k]$ does not meet $Y_y$
(pick $j$ and $k$, $j\neq k$, so that $w_j$ and $w_k$ lie on the
same side of the plane that coincides with $Y_y$ in $B(w,2\rho_1)$).

Set $\gamma = f_y([w_j,w_k])$; this is (the support of) a path
from $z_j$ to $z_k$. Recall that $z=f_y(w)$.
By (18.14) and the second half of (18.8), 
$\gamma \i \overline B(z,\rho_2)$, where 
$$
\rho_2 = (1+\tau) \, r \, (\rho_1/r)^{1-\tau}
\leq (1+\tau) \, r \, [(1-\tau)^{-1}\rho/r]^{(1-\tau)/(1+\tau)}.
\leqno (18.15)
$$
It is a little unpleasant that $\rho_2$ is not equivalent to $\rho$
(because $f_y$ is merely biH\"older-continuous, and not 
quasisymmetric), but at least it is as small as we want, which will be 
enough for (18.13).

Suppose that $\gamma$ meets $E$ at some point $x$.
Since $\gamma \i B(y,19r/10)$, the first part of (18.10) says
that $x=f_y(\xi)$ for some $\xi\in Y_y\cap B(y,2r)$. But $f_y$ is injective,
so $\xi\in [w_j,w_k]$ and this is impossible because $[w_j,w_k]$ does
not meet $Y_y$. So $\gamma$ does not meet $E$. This proves (18.13).

\medskip
Next we use the fact that  $z\in E_Y$, so Lemma 16.25 says that 
for $t$ small, $B(z,t)$ is a biH\"older ball of type $Y$ for $E$.
Choose $\rho_2$ so small that this is the case for 
$t \leq 2\rho_2$. Denote by $f_{z,t}: B(z,2t) \to \R^3$ 
the corresponding biH\"older mapping and by $Z_{z,t}$ the associated 
cone of type $\Bbb Y$ centered at $z$. 

Let $\rho$ be as in (18.13);  we can safely assume that 
$\rho < \rho_2$, so $f_{z,\rho}$ and $Z_{z,\rho}$ are defined.
Pick three points $\zeta_j \in B(z,\rho/2)$, lying in the three 
different connected components of $\R^3 \setminus Z_{z,\rho}$, 
and set $z_j= f_{z,\rho}(\zeta_j)$. By the analogue of (18.7)
the $z_j$ lie in $B(z,\rho) \setminus E$, and by (8.13) there is a path 
$\gamma \i \overline B(z,\rho_2) \setminus E$ that connects 
two of the $z_j$.

Since $\gamma$ may get way out of $B(z,\rho)$, we cannot
find a contradiction immediately in $B(z,\rho)$, and we need 
to zoom out again. But we want to do this carefully, because
we want to keep the information that $E$ should separate the
three $z_j$ because $Z_{z,\rho}$ separates the $\zeta_j$.

Let us construct three new paths leaving from the $z_j$. 
Denote by $\zeta_{0,j}$ the point of $\partial B(z,\rho)$
that lies in the same connected component of 
$\R^3 \setminus Z_{z,\rho}$ as $\zeta_j$, and as far as possible
from $Z_{z,\rho}$. Set $l_{0,j} = [\zeta_j,\zeta_{0,j}]$. Thus $l_{0,j}$
does not meet $Z_{z,\rho}$ (the components of 
$\R^3 \setminus Z_{z,\rho}$ are convex), so 
$g_{0,j} = f_{z,\rho}(l_{0,j})$ does not meet $E$
(because it is contained in $B(z,3\rho/2)$ by the analogue
of (18.7), and by the first part of the analogue of (18.10)).

Now we iterate, and construct, for $k \geq 1$ such that
$2^k\rho \leq 2\rho_2$, line segments $l_{k,j}$ and arcs $g_{k,j}$,
so that in particular
$$
g_{k,j} \i B(z,3 \cdot 2^{k-1}\rho) \setminus E.
\leqno (18.16)
$$
Let $k \geq 1$ be given, with $2^k\rho \leq 2\rho_2$.
Thus we can use $f=f_{z,2^k\rho}$ and $Z=Z_{z,2^k\rho}$.
Denote by $z_{k-1,j}$ the final endpoint of $g_{k-1,j} \,$. If 
(18.16) holds for $k-1$ (which, by the way, is the case when $k=1$),
then $z_{k-1,j} \in B(z,3 \cdot 2^{k-2}\rho) \setminus E$, and
$\zeta'_{k-1,j} = f^{-1}(z_{k-1,j})$ is defined and lies in
$B(z,2^{k}\rho) \setminus Z$, by the analogues for $f$ of (18.7) 
and (18.10). Choose $l_{k,j}=[\zeta'_{k-1,j},\zeta_{k,j}]$, where
$\zeta_{k,j}$ is the point of $\partial B(z,2^{k}\rho)$ that lies in the
same component of $\R^3 \setminus Z$ as $\zeta'_{k-1,j}$ and as far as 
possible from $Z$. As before (for $l_{0,j}$), $l_{k,j}$ does not meet $Z$ 
and is contained in $\overline B(z,2^{k}\rho)$. 
Finally set $g_{k,j} = f(l_{k,j})$;
then (18.16) holds, with the same proof as for $g_{0,j}$.
Thus we can construct the $g_{k,j}$ by induction.

Next we check that for each $k$, the three $\zeta_{k,j}$
lie in different components of $\R^3 \setminus Z_{z,2^k\rho}$. This
is true for $k=0$, because the three $\zeta_j$ lie in 
different components of $\R^3 \setminus Z_{z,\rho}$, 
and $\zeta_{0,j}$ lies in the same component as $\zeta_j$.
Now suppose that $k \geq 1$ and that this is true for $k-1$. Observe that
$|\zeta'_{k-1,j}-\zeta_{k-1,j}| \leq |\zeta'_{k-1,j}-z_{k-1,j}| 
+|z_{k-1,j}-\zeta_{k-1,j}| \leq 2^{k+1} \tau \rho$, by (18.7)
for $f$ and $f_{z,2^{k-1}\rho}$. In addition, $Z$ is quite close
to $Z_{z,2^{k-1}\rho}$, because they are both close to $E$
in $B(z,2^{k-1}\rho)$  (use the analogue of (18.10)). 
The three $\zeta_{k-1,j}$ lie far from $Z_{z,2^{k-1}\rho}$ (by construction), 
and in different components of $\R^3 \setminus Z_{z,2^{k-1}\rho}$ 
(by induction assumption). Then the $\zeta'_{k-1,j}$ also lie far from $Z$, 
and in different components of $\R^3 \setminus Z$. This is also true
for the $\zeta_{k,j}$, because $\zeta_{k,j}$ lies in the same 
component as $\zeta'_{k-1,j}\,$. This is what we wanted. 

Let $m$ denote the largest exponent $k$ such that 
$2^k \rho \leq 2\rho_2$; thus $2^m \rho > \rho_2$.
We just proved that the three $\zeta_{m,j}$ lie in different 
components of $\R^3 \setminus Z_{z,2^{m}\rho}\,$.
Denote by $g_j$ the path from $z_j$ to $z_{j,m}$ obtained by 
concatenation of all the $g_{j,k}$, $0 \leq k \leq m$.
Then $g_j \i B(z,3 \cdot 2^{m-1}\rho) \setminus E$, by (18.16).
We also have a path $\gamma \i \overline B(z,\rho_2) \setminus E$ 
that connects two of the $z_j$, and when we add to it the two
corresponding $g_j$, we get a path
$\widehat\gamma \i B(z,3\cdot 2^{m-1}\rho) \setminus E$ 
that connects two $z_{j,m}$. Set 
$\Gamma = f_{z,2^{m}\rho}^{-1}(\widehat\gamma)$. This path
does not meet $Z_{z,2^{m}\rho} \,$, by (18.10) and the usual argument.
Yet it connects two different $\zeta_{m,j}$, a contradiction.
This completes our proof of the first inclusion in (18.12).

\ms
Now we consider $w\in L_y \cap B(y,r)$ and check that 
$z= f_y(w) \in E_Y$ (it is clear that $f_y(w) \in B(y,3r/2)$,
by (18.7)). By (18.10), $E$ is locally biH\"older-equivalent to a 
set of type $\Bbb Y$ near $z$, and in particular $\R^3 \setminus E$ 
has three local components near $z$. 
More precisely, choose three points $\zeta_1$, $\zeta_2$
and $\zeta_3$ in $B(y,r) \setminus Y_y$, that lie in different
connected components of $\R^3 \setminus Y_y$. For each small
$\rho$ and each $j \leq 3$, choose a point 
$\zeta_{j,\rho} \in B(w,\rho) \setminus Y_y$, that lies in the
same component of $\R^3 \setminus Y_y$ as $\zeta_j$.
The segment $l_{j,\rho} = [\zeta_{j},\zeta_{j,\rho}]$
does not meet $Y_y$ and is contained in $B(y,r)$, so (18.10) 
says that $\gamma_{j,\rho} = f_y(l_{j,\rho})$ does not meet $E$.
In addition, (18.8) says that the extremities 
$z_{j,\rho} = f_y(\zeta_{j,\rho})$ lie in $\overline B(z,\rho_1)$, 
where $\rho_1 = r (1+\tau) |\rho/r|^{1-\tau}$. Notice that 
$\rho_1$ tends to $0$ when $\rho$ tends to $0$.

If $z$ is a $P$-point, Lemma 16.19 says that for $\rho_1$ small enough,
$B(z,\rho_1)$ is a biH\"older ball of type $P$ for $E$. 
Denote by $f : B(z,2\rho_1) \to \R^3$ the corresponding mapping, 
and $P$ the corresponding plane through $z$. 
The three $f^{-1}(z_{j,\rho})$ lie in $B(z,11\rho_1/10) \setminus P$,
by the analogues of (18.7) and (18.10) and because $z_{j,\rho} \notin E$,
and at least two of them lie on the same side of $P$. Choose $j$ and 
$k \neq j$ such that $\Gamma = [f^{-1}(z_{j,\rho}),f^{-1}(z_{k,\rho})]$
does not meet $P$. Then $f(\Gamma) \i B(z,3\rho_1/2) \setminus E$
(by (18.7) and (18.10)). Then $\gamma = \gamma_{j,\rho} \cup f(\Gamma) 
\cup \gamma_{k,\rho} \i B(y,r) \setminus E$, so its inverse image
$f_y^{-1}(\gamma)$ does not meet $Y_y$. But $\gamma$ contains
$\zeta_{j}$ and $\zeta_k$, which lie in different components of 
$\R^3 \setminus Y_y$. This contradiction proves that $z$ is not a $P$-point.

We assumed that there is no $T$-point in $B(0,3/2)$, and 
$y\in E \cap B(0,4/3)$, so it is a $Y$-point. [We could also
exclude the case when $y$ is a $T$-point, by an argument similar
to the proof of the first inclusion in (18.12).]
This proves the second half of (18.12), and Lemma 18.11 follows.
\qed

\bigskip
We shall also need to know that near $\partial B(0,1)$, $E_Y$
is composed of three nice arcs of curves. Let $T$ be the
minimal cone of (18.3), and denote by $L$ the spine of $T$.
Thus $L$ is composed of four half-line $L_j$, and we choose 
the names so that the point $a_j$ in the statement of 
Proposition 18.1 lies exactly opposite to $L_j$. 
Call $z_j$ the point of $L_j \cap \partial B(0,1)$;
thus $z_j = - a_j$. Observe that $T$ coincides with a set 
of type $\Bbb Y$ in $B(z_j,10^{-1})$. Choose $y_j \in E$ such that
$$
|y_j-z_j| \leq 2 d_{0,2}(E,T) \leq 2 \varepsilon_7
\leqno (18.17)
$$
(by (18.3)). If $\varepsilon_7 < 10^{-4}\varepsilon_4$, 
where $\varepsilon_4$ comes from Lemma~16.51 (applied with 
$\tau = 10^{-4}$, say), we can apply Lemma~16.51 to $B(y_j,10^{-3})$,
where we use a translation by $y_j-z_j$ of the set of type $\Bbb Y$ that
coincides with $T$ near $z_j$, to make sure that the spine goes 
through $y_j$. We get that $B(y_j,10^{-3})$ is a biH\"older ball 
of type $\Bbb Y$ for $E$, with constant $\tau$.

Denote by $f_{j} : B(y_{j},2 \cdot 10^{-3}) \to \R^3$ the 
corresponding biH\"older mapping, by $Y_j$ the associated
minimal cone of type $\Bbb Y$ and centered at $y_{j}$, and by 
$\ell_{j}$ the spine of $Y_j$.
The proof of Lemma 18.11 yields
$$
E_Y \cap B(y_j,10^{-3}/2) \i f_j(\ell_{j} \cap B(y_j,10^{-3})) 
\i E_Y \cap B(y_j,3\cdot 10^{-3}/2).
\leqno (18.18)
$$

Set $g_j = \ell_{j} \cap B(y_j,10^{-3})$ and 
$\gamma_j = f_{j}(g_j)$; these are nice little arcs, 
and (18.18) says that $\gamma_j \i E_Y$. Notice that $\gamma_j$
crosses the annulus $A = B(0,1-10^{-4})\setminus B(0,1+10^{-4})$, 
because $|f_j(x)-x| \leq \tau 10^{-3}$ for 
$x\in g_j$ and the extremities of  $g_j$ lie on both sides of $A$ 
(and reasonably far from it). Let us also check that 
$$
E_{Y} \cap A \i \bigcup_{j=1}^4 \gamma_j \, . 
\leqno (18.19)
$$
Indeed, if $y\in E_{Y} \cap A$, (18.4) (with $\alpha = 10^{-4}/2$)
says that $\dist(y,L) \leq 10^{-4}$, so $y$ lies in one of the 
$B(y_j,10^{-3}/2)$ (because it lies in $A$). 
Hence (18.19) follows from (18.18).

\ms
The situation in $B(0,1)$ will also be easy to control.
First observe that
$$
E_Y\cap \overline B(0,1) \hbox{ is compact.}
\leqno (18.20)
$$
Indeed, let $\{ y_k \}$ be a sequence in $E_Y \cap \overline B(0,1)$,
and assume that it converges to $y$. Then $y\in E$, because $E$
is closed in $U$. If $y$ is a $P$-point, Lemma 16.19
says that every point $x\in E$ near $y$ is a $P$-point again.
[Incidentally, we do not need the whole proof; we could deduce the
same result from Lemma 16.11 (and in particular (16.15) with the
plane $Z(x,\rho)$), plus (16.9) to control $\theta(x)$.]
This is impossible, because $y_k \in E_Y$, so $y$ is not a $P$-point.
It cannot be a $T$-point either, because we assumed that there is
no $T$-point in $B(0,3/2)$. So $y \in E_Y$, as needed for (18.20).

\ms
Select a small biH\"older ball $B(y,r_y)$ for each 
$y\in E_Y\cap \overline B(0,1)$, as we did for (18.7)-(18.10), 
and then cover $E_Y\cap \overline B(0,1)$ by a 
finite collection of $B(y,r_y/10)$, $y\in J$. Thus $J$ is finite and
contained in $E_Y\cap \overline B(0,1)$. Call $\eta$ the minimum of
$10^{-5}$ and the smallest of the $r_y/10$, $y\in J$.

Now pick one of the four $\gamma_j$. Recall that $\gamma_j$ crosses $A$. 
Let us continue the branch of $\gamma_j$ that goes inside $B(0,1)$, 
by little steps of diameter comparable to $\eta$, by a curve in $E_Y$. 
For each extension that we get, and if the endpoint $\xi$ of the extension 
lies in $\overline B(0,1)$, we can find $y\in J$ such that 
$\xi \in B(y,r_y/10)$. Then by Lemma 18.11 there is a neighborhood  of $\xi$,
that contains $B(\xi,\eta)$, where $E_Y$ coincides with the nice simple 
curve $f_y(L_y\cap B(y,r_y/2))$. We continue our extension of $\gamma_j$
by an arc of $f_y(L_y\cap B(y,r_y/2))$ of diameter $\eta$, in the
direction that was not covered (if it exists). That is, 
since the curve comes from the outside boundary of $A$,
it has to cover at least one side of $f_y(L_y\cap B(y,r_y/2))$.
We continue the path in the other direction.

We continue the construction like this, as long as we do not leave 
$B(0,1)$. Observe that when we hit a point $x$, our curve really 
crosses a neighborhood of $x$ in $E_Y$. As a result, our curve is simple:
otherwise, consider the first time where we hit a point $x$ that
was already on the curve; we cannot come from any of the two
branches of $E_Y$ that leave from $x$, because they were covered some
time ago and $x$ is the first point of return, but then there is no 
access left, so we cannot have returned to $x$. In fact, this argument shows 
that as soon as our curve meets $x$, our next extension leaves 
$B(x,\eta)$, and then the curve never returns to $B(x,\eta)$ again.

As a result of this (and because $B(0,1)$ does not contain too many 
disjoint balls of radius $\eta/2$), our extension process has 
to stop after a finite number of steps. This means that our extension
of $\gamma_j$ eventually hits $\partial B(0,1)$ again. By (18.19), 
this happens in one of the $\gamma_k$. Since our curve is simple,
$k \neq j$. So we proved that for each $j$, 
$$
\hbox{there is a path in $E_Y \cap \overline B(0,1)$ that connects 
$\gamma_j$ to some other $\gamma_k$.}
\leqno(18.21)
$$

Next we want to study the connected components of
the complement of $E$ near $E_Y$. For $y\in E_Y \cap B(0,4/3)$,
denote by ${\cal A}(y)$ the collection of connected components $V$
of $B(0,2) \setminus E$ such that $y \in \overline V$.
Recall that each $y\in E_Y \cap B(0,4/3)$ has a small neighborhood 
where $E_Y$ is biH\"older-equivalent to a set of type $\Bbb Y$ 
as in (18.7)-(18.10). Consequently, ${\cal A}(y)$ has at most $3$ 
elements, and is a locally constant function of $y\in E_Y \cap B(0,4/3)$.

We can control ${\cal A}(y)$ near the four $y_j$. For $1 \leq j \leq 4$,
pick $w_j \in \gamma_j \cap \partial B(0,1)$.

\ms \proclaim Lemma 18.22.
For each $j$,
$$
{\cal A}(w_j) = \big\{ V(a_i) \, ; \, i \neq j \big\},
\leqno (18.23)
$$
where the $V(a_i)$, $1 \leq i \leq 4$, are as in 
the statement of Proposition 18.1.

\ms
We do not say here that the three components mentioned 
in (18.23) are different; we shall take care of this later.

Recall that for $1 \leq i \leq 4$, $a_i = - z_i$, where $z_i$ 
is the point of the spine of $T$ that was used to define $y_i$; 
see a little above (18.17). And $V(a_i)$ is the connected component 
of $a_i$ in $B(0,2) \setminus E$.

Fix $j$ and $i \neq j$. Denote by $b_i$ the point
of $[a_i,z_j] \cap \partial B(z_j,10^{-5})$. Notice that
$$
\dist([a_i,b_i],E) \geq \dist([a_i,b_i],T) - 2 \varepsilon_7
\geq 10^{-6}
\leqno (18.24)
$$
by simple geometry and (18.3), so $b_i \in [a_i,b_i] \i V(a_i)$.

To relate $V(b_i) = V(a_i)$ to ${\cal A}(w_j)$, we shall use the 
same biH\"older ball $B(y_j,10^{-3})$ as for (18.18), and the 
corresponding mapping $f_j : B(y_{j},2 \cdot 10^{-3}) \to \R^3$.
Recall that 
$$
|f_j(x)-x| \leq \tau 10^{-3} \ \hbox{ for } x\in B(y_{j},2 \cdot 10^{-3}),
\leqno (18.25) 
$$
$$
B(y_j,(2-\tau)10^{-3}) \i f_j(B(y_j,2 \cdot 10^{-3})),
\leqno (18.26) 
$$
and
$$
E \cap B(y_j,(2-\tau)10^{-3}) 
\i f_j(Y_j\cap B(y_j,2 \cdot 10^{-3}))) \i E.
\leqno (18.27)
$$
by (15.6), (15.8) , and (15.9). 

Recall that $b_i \in \partial B(z_j,10^{-5}) \setminus E$
(by (18.24)), and that $|y_j-z_j| \leq 2\varepsilon_7$
(by (18.17)). Then $c_i = f_j^{-1}(b_i)$ is well defined
(by (18.26)), lies in $B(y_j,10^{-4})$ (by (18.25)), 
and out of  $Y_j$ (by (18.27)). Denote by $W_i$ the connected 
component of $B(y_j,10^{-4})\setminus Y_j$ that contains $c_i$, 
and let us check that
$$
\hbox{ the three $W_i$, $i \neq j$, are different.}
\leqno (18.28) 
$$
Indeed $Y_j$ is quite close to $T$ near $y_j$: every
point $\xi \in Y_j \cap B(y_j,2 \cdot 10^{-3})$ lies 
within $\tau 10^{-3}$ of $f_j(\xi) \in E$, by (18.27),
hence $\dist(\xi,T) \leq \tau 10^{-3} + 2 \varepsilon_7$,
by (18.3). Now the $b_i$ lie reasonably far from 
$T$ and $Y_y$, and on different sides of $Y_y$ (by (18.24)),
and this stays true for the $c_i$ (by (18.25)); this proves
(18.28).

Recall that $w_j \in \partial B(0,1) \cap \gamma_j$,
and $\gamma_j = f_j(g_j)$, where $g_j = \ell_j \cap B(y_j,10^{-3})$
(see below (18.18)). Let $\zeta_j \in g_j$ be such that 
$f_j(\zeta_j)=w_j$. Then $|\zeta_j - 1| \leq \tau 10^{-3}$,
by (18.25). The distance from $\zeta_j$ to the spine of $T$
(and here this means the half line through $z_j$) is less than
$C (\tau+\varepsilon_7)$, because $\zeta_j \in g_j$ and $Y_j$ is quite close to
$T$ near $y_j$. Then $|\zeta_j-z_j| \leq C (\tau+\varepsilon_7)$ 
(because $|z_j| =|-a_j|=1$), so $\zeta_j$ and $w_j$ both lie in 
$B(y_j,10^{-4})$ (recall that $|y_j-z_j| \leq 2 \varepsilon_7$).

For each $i \neq j$, $c_i \in B(y_j,10^{-4}) \setminus Y_j$
and $\zeta_j \in B(y_j,10^{-4})\cap \ell_j \,$; then
$I = [c_i,\zeta_j) \i B(y_j,10^{-4}) \setminus Y_j$ and 
$f_j(I) \i B(0,2) \setminus E$, by (18.27).
Hence $f_j(I) \i V(b_i) = V(a_i)$ (because $f_j(c_i)=b_i$). 
Since $f_j(I)$ gets as close to $w_j = f(\zeta_j)$ as we want, 
$V(a_i)$ is one of the elements of ${\cal A}(w_j)$.

Conversely, let $V$ be one of the elements of 
${\cal A}(w_j)$, and pick $x \in V \cap B(w_j,10^{-4})$.
Then $f_j^{-1}(x)$ is defined by (18.26), lies in
$B(y_j,10^{-3})$ by (18.25) and because $w_j \in B(y_j,10^{-4})$,
and lies out of $Y_j$ by (18.27). Thus $f_j^{-1}(x)$ lies in
some $W_i$, by (18.28). Then $J = [f_j^{-1}(x),c_i]$ lies in $W_i$ by
convexity of $W_i$, and $f(J)$ does not meet $E$ by (18.27). 
Its extremities $x$ and $b_i=f_i(c_i)$ lie in the same
component of $B(0,2)\setminus E$. That is, $b_i$ lies in $V$,
hence $V = V(b_i) = V(a_i)$.

This completes our proof of Lemma 18.22.
\qed

\ms
We are ready to conclude.
By assumption, at least three of the four $V(a_i)$ are different.
Let us assume that these are the $V(a_i)$, $1 \leq i \leq 3$.
By (18.21), we can connect $w_4 \in \gamma_4$ to some other 
$w_k$ by a path in $E_Y \cap \overline B(0,1)$; without loss of 
generality, we can assume that $k=3$. 
Since ${\cal A}(x)$ is locally constant on $E_Y$, 
${\cal A}(w_4) = {\cal A}(w_3)$ or, by Lemma~18.22,
$\{ V(a_1),V(a_2),V(a_3) \} = \{ V(a_1),V(a_2),V(a_4)\}$.
Since the $V(a_j)$, $1 \leq j \leq 3$, are distinct, 
$V(a_4)=V(a_3)$. 

Let us also apply (18.21) with $j=2$;
we get that ${\cal A}(w_2) = {\cal A}(w_l)$ for some
$l \neq 2$. Since $V(a_2) \in {\cal A}(w_l)$, it lies in 
${\cal A}(w_2)$ as well, so it is equal to 
$V(a_1)$ or $V(a_3)$ (we already know that $V(a_4)=V(a_3)$).
This contradicts our assumption that three $V(a_j)$ are
different, and completes the proof of Proposition 18.1.
\qed

\bigskip
Next we prove a variant of Theorem 1.9 for $Al$-minimal sets, 
but with an additional assumption that will allow us to apply
Proposition 18.1. Theorem 1.9 will then follow rather easily.

\ms\proclaim Proposition 18.29.
Let $E$ be a reduced $Al$-minimal set of dimension $2$ in $\R^3$
(as in Definition~1.6). If there is a minimal cone $T$
of type $\Bbb T$, centered at the origin, and a sequence 
$\{ t_k \}$ such that 
$$
\lim_{k \to +\infty} t_k = +\infty \ \hbox{ and } \ 
\lim_{k \to +\infty} d_{0,t_k}(E,T) = 0,
\leqno (18.30)
$$
assume in addition that there is a subsequence $\{ t_{k_l} \}$
for which the following separation assumption holds.
Denote by $a_i$, $1 \leq i \leq 4$, the four points of 
$\partial B(0,t_{k_l}/2)$ that lie the furthest from $T$. 
Thus the $- a_i$ are the intersections 
of $\partial B(0,t_{k_l}/2)$ with the spine of $T$. Denote by 
$V(a_i)$ the connected component of $a_i$ in 
$B(0,t_{k_l})\setminus E_{k_l}$. Then (assume that)
$$
\hbox{at least three of the four components $V(a_i)$ are distinct.}
\leqno (18.31)
$$
Then $E$ is a reduced minimal cone, i.e., $E$ is empty, a plane,
or a cone of type $\Bbb Y$ or $\Bbb T$.

\ms
Since we deal with two dimensional sets in $\Bbb R^3$, the
reduced minimal cones of type $\Bbb T$  are the cones based on the
edges of a regular tetrahedron, as defined near Figure 1.2.

Recall that an $Al$-minimal set in $\R^3$ is also the same as an 
almost-minimal set in $\R^3$, with gauge function $h=0$, as
in Definition~4.8 (or even Definitions 4.1 and 4.3, see Proposition~4.10),
so we may apply to $E$ all the results from the previous sections.

By Proposition 5.16, $\theta(x,r) = r^{-2} H^2(E\cap B(x,r))$ is a 
nondecreasing function of $r$ for each $x\in \R^3$. It is also 
bounded, by the easy part of (2.16), so it has finite limits
$$
\theta(x) = \lim_{r \to 0} \theta(x,r)
\ \hbox{ and }\ 
\theta_{\infty} = \lim_{r \to +\infty} \theta(x,r).
\leqno (18.32)
$$
It is easy to see that $\theta_{\infty}$ does not depend
on $x \,$; indeed $B(0,r-|x|) \i B(x,r) \i B(0,r+|x|)$ for $r$ large, so
$(r-|x|)^2 \, \theta(0,r-|x|) \leq r^2 \theta(x,r) \leq (r+|x|)^2 
\,\theta(0,r+|x|)$, which yields $\lim_{r \to +\infty} \theta(x,r)
=\lim_{r \to +\infty} \theta(0,r)$ after dividing by $r^2$ and taking 
limits.

We know from the discussion below (16.7) that $\theta(x)$ can only take 
the values $0, \pi, 3\pi/2$, and $d_+ = {\rm Argcos}(-1/3)$, corresponding 
to the various types of blow-up limits that $E$ may have at $x$. 
When $x \in E$, $\theta(x)=0$ is excluded by (2.16).

But we also have that $\theta_{\infty} \in \{ 0, \pi, 3\pi/2, d_+ \}$.
Indeed consider a blow-in sequence $\{ E_k \}$, where
$E_k = t_k^{-1} E$ for some sequence $\{ t_k \}$ that tends
to $+\infty$. We can extract a subsequence (which we shall still
denote by $\{ E_k \}$), which converges to a limit $F$ (in Hausdorff 
distance on compact subsets of $\R^3$, as in Section 3).
Recall that each $E_k$ is a minimal set, by homogeneity.
By Lemma 4.7, $F$ is a reduced minimal set in $\R^3$. By (3.4),
$$\eqalign{
H^2(F \cap B(0,r)) &\leq \liminf_{k \to +\infty} H^2(E_k \cap B(0,r)) 
= \liminf_{k \to +\infty} t_k^{-2} H^2(E \cap B(0,t_k r)) 
\cr&= \liminf_{k \to +\infty} r^2\theta(0,t_kr) 
= r^2\theta_{\infty} \, .
}\leqno (18.33)
$$
Similarly, Lemma 3.12 says that 
$$\eqalign{
H^2(F \cap \overline B(0,r)) 
&\geq \limsup_{k \to +\infty} H^2(E_k \cap \overline B(0,r)) 
\cr&= \limsup_{k \to +\infty} t_k^{-2} H^2(E \cap \overline B(0,t_k r)) 
= r^2\theta_{\infty} \, ,
}\leqno (18.34)
$$
so $H^2(F \cap B(0,r)) = r^2\theta_{\infty}$ for $r>0$. By Theorem 6.2,
$F$ is a reduced minimal cone. By [Ta], 
$F$ is the empty set, a plane, or a cone of type $\Bbb Y$ or $\Bbb T$. 
Then $\theta_{\infty} \in \{ 0, \pi, 3\pi/2, d_+ \}$, as announced.

Let us distinguish between cases, depending on the value of $\theta_{\infty}$.
If $\theta_{\infty}=0$, $\theta(0,r) = 0$ for all $r$, $E = \emptyset$,
and we don't need to worry. If $\theta_{\infty}=\pi$, $\theta(0,r) > 0$ for 
some $r$, and $E$ is not empty. Pick $x\in E$. We know that 
$\theta(x) > 0$, by (2.16), so $\theta(x) \geq \pi = \theta_{\infty}$.
Then $\theta(x,\cdot)$ is constant and  Proposition 6.2 says that 
$E$ is a cone (in fact, a plane).

Now suppose that $\theta_{\infty}=3\pi/2$. As before, we can find a 
sequence $\{ t_k \}$ that tends to $+\infty$, so that the sets
$E_k = t_k^{-1}E$ converge to a limit $F$, and then $F$ is a 
reduced minimal cone. In fact, (18.33) and (18.34) tell us that the density
of $F$ is $3\pi/2$, so $F$ is a cone of type $\Bbb Y$ centered at the origin.

Each $t_k^{-1}E$ is a reduced minimal set, and the sequence
converges to $F$. For $k$ large, it satisfies the assumptions of
Proposition 16.24 (in the unit ball), so there is a $Y$-point in
$t_k^{-1}E$, hence also in $E$. For this point $x$,
$\theta(x) = 3\pi/2 = \theta_\infty$. Again, $\theta(x,\cdot)$ is constant
and Proposition~6.2 says that $E$ is a minimal cone (here a set
of type $\Bbb Y$).

\ms 
We are left with the case when $\theta_{\infty}=d_+$. The same 
argument as before gives a blow-in sequence 
$\{ E_k \} = \{ t_k^{-1} E \}$ that converges to a reduced 
minimal cone $T$ of type $\Bbb T$ centered at the origin. This time, 
we do not have a good enough analogue of Proposition 16.24 that
would allow us to conclude, and this is why we added an extra hypothesis.

The sequence $\{ t_k \}$ satisfies the assumption (18.30),
so we have (18.31). This is exactly what we need to apply 
Proposition 18.1 to the ball $B(0,t_{k_l})$ (for $l$ large enough, 
so that (18.3) holds). We get the existence of a point of type $T$
in $E \cap B(0,t_{k_l})$. Then $\theta(x) = d_+ = \theta_\infty$
at this point, $\theta(x, \cdot)$ is constant, and $E$ is a minimal 
cone (of type $\Bbb T$), as before.
This completes the proof of Proposition 18.29.
\qed

\bigskip
\noindent {\bf Proof of Theorem 1.9.}
\ms
Let us prove Theorem 1.9. Let $E$ be a $MS$-minimal set
in $\Bbb R^3$ (see Definition 1.4). By Remark 1.8, 
$E$ is $Al$-minimal. By Remark 2.14, $E^\ast$
(the closed support of the restriction to $E$ of $H^2$, as in 
Definition 2.12) is a reduced  $Al$-minimal set. We just need to show that
$E^\ast$ is a reduced minimal cone (see (2.13)).

We want to apply Proposition 18.29 to $E^\ast$, so we
need to check the separation assumption. So we assume that 
$\{ t_k \}$ is a sequence such that (18.30) holds for $E^\ast$, 
and we want to prove (18.31) (for the same sequence, but $k$ large; 
we won't need to take a real subsequence).

Let us show the dependence in $k$ in our notation.
Call $a_{j,k}$, $1 \leq j \leq 4$, the four points of
$\partial B(0,t_k/2)$ that lie furthest from $T$,
and denote by $V_{j,k}$ the connected component of $a_{j,k}$ in 
$B(0,t_{k})\setminus E^\ast$. We want to show that for $k$ large, 
the $V_{j,k}$ are all different.

Suppose instead that two of the $V_{j,k}$ are the same,
for instance $V_{1,k}$ and $V_{2,k}$.
We want to show that this is impossible, and for this we intend to
construct a better competitor $F$, obtained by removing a big piece
of the wall near $a_{k} = {1 \over 2}(a_{1,k} + a_{2,k})$.
The comparison will be easier once we know a little more about $E$ 
near $a_{k}$.

Set $E_k = t_k^{-1} E^\ast$, and $a = t_k^{-1}a_k$
(this is a fixed point). By (18.30), the $E_k$ converge to $T$
in $B(0,1)$, so Lemma 3.3 says that
$$
\limsup_{k \to +\infty} H^2(E_k \cap B(a,10^{-1})) 
\leq H^2(T \cap B(a,10^{-1})) \leq 10^{-2} \pi,
\leqno (18.35)
$$
where the last inequality holds because $T$ coincides with a plane 
in $B(a,10^{-1})$.

Pick $x_k \in E_k$ as close to $a$ as possible,
and apply Lemma 16.19 to $E_k$ in $B(x_k, r_k)$,
where we choose $r_k$  such that $9r_k = 10^{-1}-|x_k-a|$
and hence $B(x_k,9r_k) \i B(a,10^{-1})$. 
The assumption (16.20) holds for $k$ large, because
$$
\limsup_{k \to +\infty} H^2(E_k \cap B(x_k,9r_k)) 
\leq \limsup_{k \to +\infty} H^2(E_k \cap B(a,10^{-1})) 
\leq 10^{-2} \pi
\leqno (18.36)
$$
by (18.35), and $9r_k$ tends to $10^{-1}$. So Lemma 16.19 applies, and
says that $B(x_k,10^{-2})$ is a biH\"older ball of type $\Bbb P$ for $E_k$.
Let $f_k : B(x_k,2/100) \to \R^3$ be the corresponding mapping,
and $P_k$ the associated plane. Let us use $f_k$ to show that 
$$
B(x_k,10^{-2})\setminus E_k \i 2^{-k} V_{1,k} \, ;
\leqno (18.37)
$$
this is believable because we assumed that the components $V_{1,k}$ 
and $V_{2,k}$ coincide, but let us check the details.
Let $\zeta_1$ and $\zeta_2$ be the two points of
$\partial B(x_k,10^{-2})\setminus P_k$ that lie at distance 
$10^{-2}$ from $P_k$ (so they are on opposite sides of $P_k$).
If $x\in B(x_k,10^{-2})\setminus E_k$, 
$f_k^{-1}(x)$ is defined (by (15.8)), lies in 
$B(x_k,{11 \over 1000})$ (by (15.6)) and away from $P_k$
(by (15.9)), so one of the two line segments 
$I = [f_k^{-1}(x),\zeta_i]$ does not meet $P_k$.
Then $f_k(I)$ is an arc in $B(x_k,{12 \over 1000}) \setminus E$
(by (15.6) and (15.9)) that goes from $x$ to
$z_i = f(\zeta_i)$. Thus (18.37) will follow if we show that
both $z_i$ lie in $2^{-k} V_{1,k}$.
Now $\dist(\zeta_i,P_k) = 10^{-2}$, so 
$$\eqalign{
\dist(z_i,E_k \cap &B(x_k,(2-\tau)10^{-2}))
\geq \dist(z_i,f(P_k \cap B(x_k,2 \cdot 10^{-2})))
\cr&\geq \dist(\zeta_i, P_k \cap B(x_k,2 \cdot 10^{-2})) - 2 \tau 10^{-2}
\geq (2-2\tau) 10^{-2}
}\leqno (18.38)
$$
by (15.9), (15.6), and the definition of $\zeta_i$. Since
$$
\dist(z_i,E_k \setminus B(x_k,(2-\tau)10^{-2}))
\geq (2-\tau)10^{-2} - |z_i - x_k| \geq (1-2 \tau) 10^{-2}
\leqno (18.39)
$$
by (15.6) and because $\zeta_i \in \partial B(x_k,10^{-2})$,
we get that 
$$
\dist(z_i,T) \geq \dist(z_i,E_k) - d_{0,1}(E_k,T)
\geq 10^{-2}/2
\leqno (18.40)
$$
for $k$ large, because the $E_k$ converge to $T$
in $B(0,1)$ (by (18.30)). Then $z_i$ lies in 
the same component as $a_1$ or $a_2$ in $B(0,1) \setminus E_k$, 
or equivalently $2^k z_i \in V_{1,k} \cup V_{2,k}$.
But we assumed that $V_{1,k} = V_{2,k} \,$, so 
$z_i \in 2^{-k} V_{1,k} \,$, as needed for (18.37).

Set $F_k = E\setminus B(r_k x_k,10^{-3}r_k)$. This set may be
a little larger than $r_k (E_k \setminus B(x_k,10^{-3})$ if
$E$ is not reduced, and we need to pay some attention to
the difference because did not verify that $E^\ast$ is
a $MS$-minimal set. Let us check that 
$$
F_k \hbox{ is a $MS$-competitor for $E$},
\leqno (18.41)
$$
as in Definition 1.1. Let us take $R = r_k$ in that definition;
(1.2) is clear, and we just need to prove (1.3). 
Let $z_1, z_2$ lie in $\R^3 \setminus [E \cup B(0,r_k)]$.
We need to show that if $E$ separates $z_1$ from $z_2$, $F_k$ separates 
them too. Suppose it does not, and pick a path $\gamma$ from $z_1$ to $z_2$, 
that does not meet $F_k$. Since $\gamma$ meets $E$ and the only
difference is $E \cap B(r_k x_k,10^{-3}r_k)$, $\gamma$ has to go through 
$\partial B$, where we set $B=B(r_k x_k,2 \cdot 10^{-3}r_k)$.

Call $z'_1$ the first point of $\gamma \cap \partial B$ when 
we go from $z_1$ to $z_2$ along $\gamma$, and $z'_2$ the last point of 
$\gamma \cap \partial B$. The part of $\gamma$ between
$z_1$ and $z'_1$ does not meet $E$, because it does not meet $F_k$ 
(and $F_k$ coincides with $E$ out of $B$).
Similarly, the part of $\gamma$ between $z_2$ and $z'_2$ does not meet 
$E$ either. So we just need to show that $E$ does not separate 
$z'_1$ from $z'_2$ and we shall reach the desired 
contradiction with our assumption that $E$ separates $z_1$ from $z_2$.

Notice that the $z'_{i}$ lie on $\partial B \setminus E$, hence also
on $\partial B \setminus E^\ast$. Then the $r_k^{-1} z'_i$ lie
in $\partial B(x_k,2 \cdot 10^{-3}) \setminus E_k$ (recall that
$E_k = r_k^{-1} E^\ast$). By (18.37), the $z'_i$ lie in $V_{1,k}$,
which means that there is a path $\gamma_1 \i B(0,r_k) \setminus E^\ast$ 
that goes from $z'_{1}$ to $z'_2$.

We can even make $\gamma_1$ piecewise linear.
Now $\gamma_1$ may meet $\Delta = E \setminus E^\ast$, but 
observe that $H^2(\Delta) =0$ by (2.13).
Let us modify $\gamma_1$ slightly, so that it does not meet 
$E$. First we move the endpoints of the segments
that compose $\gamma_1$ very slightly (to make sure that the new
path will not meet $E^\ast$), so that the endpoints lie out of $\Delta$
too. This is very easy to arrange, because the complement of $\Delta$
is dense. Note that we can leave $z'_1$ and $z'_2$ alone, since they 
already lie out of $E$.

For each of the segments that compose our new path, the orthogonal 
projection of $\Delta$ on a plane perpendicular to that segment has 
vanishing Lebesgue measure because $H^2(\Delta) =0$, so we may translate the
segment a tiny bit so that it no longer meets $\Delta$. We do this
for each segment, and then connect the new endpoints together by very
short line segments to make a new path $\gamma_2$ from $z'_1$ to $z'_2$.
Now $\gamma_2$ does not meet $E$. So $z'_1$ and $z'_2$ lie in the same
component of $\R^3 \setminus E$, we reached the desired contradiction, and
(18.41) holds.

Since $E$ is $MS$-minimal, Definition 1.4 says that
$H^2(E\setminus F_k) \leq H^2(F_k\setminus E) = 0$.
But this is not true: 
$$
H^2(E\setminus F_k) = 
H^2(E \cap B(r_k x_k,10^{-3}r_k) 
= r_k^2 H^2(E_k \cap B(x_k,10^{-3}) > 0
\leqno (18.42)
$$
by (2.13) and (2.16), for instance.

This last contradiction completes our proof of (18.31); 
Proposition 18.29 applies, and Theorem 1.9 follows.
\qed

\ms\noindent {\bf Remark 18.43.}
So we got a complete description of $MS$-minimal sets in $\R^3$,
but unfortunately not of $Al$-minimal sets. It is tempting to try
to deduce the separation assumption in Proposition 18.29 from the
closeness of $E$ to a cone of type $\Bbb T$ and its local 
description of $E$ in Theorem 16.1. But we shall see in Section 19 
that it is not too hard to construct sets $E \i \R^3$ that are very 
close to a cone of type $\Bbb T$ in the unit ball, 
satisfy the conclusions of Theorem 16.1, but do not have 
any $T$-point and only separate $\R^3$ locally into two connected 
components only. So a proof of the analogue of Theorem~1.9
for $Al$-minimal sets would need to use $Al$-minimality one more time.

\ms\noindent {\bf Remark 18.44.}
It is possible to define a notion of $MS$-almost-minimal sets in $U\i\R^3$,
as we did for $A'$-almost-minimal sets in Definition 4.8. That is,
if $E$ is a closed subset of $U \i \R^3$ such that 
$H^2(E \cap \overline B) < + \infty$ for every compact ball $\overline B \i U$, 
we say that $E$ is a $MS$-almost-minimal set in $U$, with gauge function $h$, 
if 
$$
H^2(E \setminus F) \leq H^2(F \setminus E) + h(r) r^2
\leqno (18.45)
$$
whenever $\overline B(x,r) \i U$ and $F$ is another closed set in $U$ 
such that $F \setminus \overline B(x,r)=E\setminus \overline B(x,r)$,
and 
$$
\hbox{$F$ separates $y$ from $z$ whenever } 
y,z \in U\setminus (\overline B(x,r) \cup E)
\hbox{ are separated by } E.
\leqno (18.46)
$$
When $h = 0$ and $U = \R^3$, we get the minimal sets of Definition 1.4.

We claim (but we shall not check) that the following things are true.
First, $E^\ast$ is a $MS$-almost-minimal set in $U\i\R^3$, with the
same gauge function as $E$. See ``Coral pairs" on pages 58-60 of [D2], 
but the  main point is our proof of (18.41). 

Next, all the results from the previous sections apply to 
$MS$-almost-minimal sets. We cannot say
directly that this is because all $MS$-almost-minimal sets are 
$A'$-almost-minimal, even though we can localize the proof
of Remark 1.8. Indeed the $Al$-competitors $F=\varphi_1(E)$, where
the family $\{\varphi_t \}$ satisfies (2.2)-(2.6), are not 
necessarily $MS$-competitors as above, because even though 
$\widehat W$ in (2.6) is compact and contained in $U$, it is not
sure that it is contained in a compact ball $\overline B(x,r) \i U$.
But we claim that this is not a problem, that all sets $\widehat W$
used in and for this paper are contained in compact balls $\overline B(x,r) \i U$
with a diameter that is not too large.

Finally, we can even remove the separation assumption in Proposition 18.29
when we deal with $MS$-almost-minimal sets, by the final argument in our
proof of Theorem 1.9.

\bigskip 
\noindent {\bf 19. A set $E \i \R^3$ that looks like a cone of type $\Bbb T$
but has no $T$-point}
\medskip

In this section we try to describe a closed set $E \i \R^3$ that coincides 
outside of the unit ball with a set of type $\Bbb T$ centered at the origin, 
but is locally $C^1$-equivalent to a plane or a set of type $\Bbb Y$ 
near each of its points. 
This set is probably never minimal, so it is only a counterexample to what would
be a wrong way to prove an analogue of Theorem 1.9 for $Al$-minimal
sets. It is quite possible that this set is the same as the one depicted in
Figure 1.1.1 in [LM] and Figure 3 of [Mo2], 
following a suggestion of R. Hardt.

Even though we shall try to give a reasonable description, it is 
probable that the best way for the reader to convince herself that nothing 
goes wrong with the topology is to construct her own paper model with 
scissors and glue.

So we start from the minimal set $T$ centered at the origin, and call
$L$ its spine. Thus $L$ is composed of four half lines $L_{j}$,
$1 \leq j \leq 4$. Call $x_j$ the point of $L_j \cap \partial B(0,1)$.
We first construct two faces of $E$ and
the set $E_{Y}$ of $Y$-points of $E$. 

Call $P_1$ the plane that contains $L_{1} \cup L_{2}$,
and then pick a nice smooth curve $\Gamma_{1}\i P_1$ that has
$L_{1} \setminus B(0,1/2)$ and $L_{2} \setminus B(0,1/2)$ as its ends.
The simplest is to add to $L_{1} \setminus B(0,1/2)$ and 
$L_{2} \setminus B(0,1/2)$ the arc of circle that goes from
$L_{1} \cap \partial B(0,1/2)$ to $L_{2} \cap \partial B(0,1/2)$
and is tangent to each $L_i$ at the point of $L_{i} \cap \partial B(0,1/2)$.
But the reader could take a smoother, equivalent curve.
Call $F_1$ the convex component of $P_1 \setminus \Gamma_1$;
this will be one of the faces of $E$; it is just a little smaller
than the corresponding face of $T$.

Similarly, call $P_4$ the plane that contains $L_{3} \cup L_{4}$,
and pick a nice curve $\Gamma_4 \i P_4$ that contains
$L_{3} \setminus B(0,1/2)$ and $L_{4} \setminus B(0,1/2)$.
Also call $F_4$ the convex component of $P_4 \setminus \Gamma_4$.

At this point, we have two faces of $E$ ($F_1$ and $F_4$), and we 
also decide that $E_Y = \Gamma_1 \cup \Gamma_4$. Observe that here
$\Gamma_1$ and $\Gamma_4$ are neither knotted nor braided; this will 
make our discussion and the construction of paper models simplest, but
the case when the $\Gamma_i$ are knotted and linked 
could be more interesting , because then it could be harder to prove
that $E$ is not minimal.

Before we continue with the construction, we should observe that when
we do all the construction below, we shall be able to preserve the 
two symmetries, with respect to the planes $P_1$ and $P_4$ 
respectively.

There is a little more of $E$ that we could construct without thinking.
We know that a neighborhood of each $\Gamma_i$ will be composed of 
three thin stripes that touch $\Gamma_i$. We already know one of these
stripes, which comes from $F_i \,$, and we could brutally construct the two
other ones: we would just need to draw two little curves vaguely parallel
to $\Gamma_i$, and make sure that the two new stripes meet $F_i$ neatly
along $\Gamma_i$, with $120^\circ$ angles. This is what is suggested by 
Figure 19.1. We shall not really need this, 
but maybe knowing this will help the reader follow the 
description of $F$ below.

The remaining set $F=E\setminus [\overline F_1 \cup \overline F_4]$ 
will be a little harder to describe.
The contact with $F_1$ and $F_4$ is along the 
$\Gamma_i$, and through the four stripes alluded to above. 
Now $F$ is actually a single connected face (we shall see that 
it is even oriented). Let us use the symmetry of $E$ and $F$ 
with respect to $P_1$ to simplify our description.

\medskip 
 \hskip 0.4cm
\epsffile{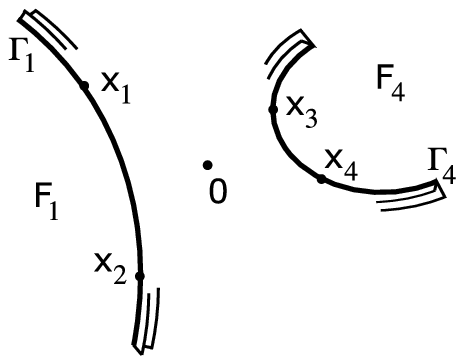}
\hskip 1.2cm
\epsffile{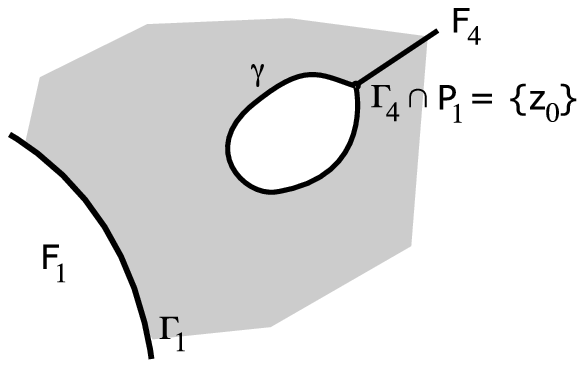}
\medskip
\hskip 1.6cm {\bf Figure 19.1 \hskip 2.8cm Figure 19.2.} Intersection 
with $P_1$
\medskip

First, the intersection $F \cap P_1$ is just a vaguely
round closed curve $\gamma$, with a singularity at the point 
of $\Gamma_4 \cap P_1$ (let us denote it by $z_0$), 
which is actually excluded (because it lies in $F_4$)
and where $\gamma$ has two half tangents 
that make $120^\circ$ angles. See Figure 19.2, 
but do not pay attention to the shaded area yet.
The face $F$ crosses $\gamma$ perpendicularly to $P_1$;
it also approaches $\Gamma_1$ from both sides of $P_1$, with
$60^\circ$ angles.

Call $F_+$ the part of $F$ that lies on the same
side of $P_1$ as $x_3$. Out of $B(0,1)$, $F_+$ coincides with the
union of two faces of $T$ (namely the face $T_{1,3}$ bounded by
$L_1$ and $L_3$, and the face $T_{2,3}$). 
Let us look at $F_+ \cap B(0,1)$. Topologically,
it is a disk, and its boundary is composed as follows. We may start
with the curve $\Gamma_1 \cap B(0,1)$, say from $x_1$ to $x_2$.
Then we follow the arc $\partial B(0,1) \cap T_{2,3}$, up to $x_3$.
Then we follow $\Gamma_4$ from $x_3$ to $z_0$, then we turn around 
$\gamma$, and return along $\Gamma_4$ from $z_0$ to $x_3$. After this
we run along the arc $\partial B(0,1)\cap T_{1,3}$ up to $x_1$,
and we are back where we started. Notice that if we include the
part of $\Gamma_4$ from $x_3$ to $z_0$, $\overline F_+ \cap B(0,1)$ is
a topological tube (because of the gluing along $\Gamma_4$).
Similarly, the whole $\overline F_+$ is a half-infinite tube, 
with a gluing along half of $\Gamma_4$ and only one ``circular" boundary 
along $\gamma$. The reader may try to get a mental picture of 
$\overline F_+$ by looking at the shaded part of Figure 19.2 that lives 
above $\Gamma_1$ and outside of $\gamma$, and imagining that it is
lifted in the third direction. In particular, imagine that
$F_4 \cap P_1$ is lifted until it becomes the upper half
of $\Gamma_4$. Also recall that $\overline F_+$ leaves from
$\Gamma_1$ and both sides of $\Gamma_4$ with $60^\circ$ angles,
and vertically from $\gamma$.

Let $F_-$ denote the part of $F$ that lies on the other side
of $P_1$. Since $\overline F_-$ is obtained from $\overline F_+$
by symmetry with respect to $P_1$, it is also a tube, bounded by $\gamma$
and glued along the other half of $\Gamma_4$. Finally 
$\overline F$ is an infinite tube, obtained from the two pieces
by gluing along $\gamma$. The fact that $F_+$ and $F_-$ live on 
different sides of $P_1$ makes it easy to believe that nothing
weird (like self-intersections) happens when we do the gluing.
It is easy to see that $\overline F$ can be oriented too.

Hopefully the properties above give a precise idea of the set
$F$, and of $E= \overline F_1 \cup \overline F_4 \cup F$. 
The fact that $E$ is locally
equivalent to a plane or a $Y$ is easy to check, and $E$ coincides
with $T$ out of $B(0,1)$.

Let us say a few words about the connected components of 
$\R^3 \setminus E$. It should be expected from the 
proof of Theorem 1.9 that $E$ does not separate $\R^3$
into three or four components, and indeed $\R^3 \setminus E$
has only two components. More precisely, call $y_j=-x_j$ the point of
$\partial B(0,1)$ which is opposite to $L_{j}$, and denote by $U_{j}$ 
the connected component of $\R^3 \setminus E$ that contains $y_{j}$. 
It is easy to see that $U_3 = U_{4}$, because there is a path from
$y_3$ to $y_4$ that does not meet $E$ and crosses $P_1$ inside
$\gamma$. Similarly, $U_1 = U_{2} \,$; rather than trying to use the
fact that $P_1$ and $P_4$ can be interchanged in the construction,
the reader should just look at Figure 19.2 and imagine a curve that runs 
in $P_1$ along $\Gamma_1$, at a short distance from but just opposite to $F_1$. 
This curve does not meet $E$, and it goes from $U_3$ to $U_4$.
Of course the two ways to see that $U_1 = U_{2}$ are similar, since
the analogue of $\gamma$ here is a loop in $P_4 \setminus F_1$
that ends at the point of $P_4 \cap \Gamma_1$ and passes through the
closest point of $\gamma$ on the other side.

It is also true that $U_1 \neq U_3$; this is not easy to see without
a paper model, and not so important, so we skip the proof. Notice
that once we know that $U_1 = U_{2}$ and $U_3 = U_{4}$, there is
no contradiction with the fact that the list ${\cal A}(y)$ of components 
that touch $y$ is locally constant along the curves $\Gamma_i \,$;
in both cases, ${\cal A}(y)= \{ U_1,U_3 \}$ (with one of the 
components counted twice).

This completes our description of $E$. Unless we are extremely lucky,
our example $E$ will not be minimal. We can try to make it minimal
inside $B(0,1)$ under the constraint that $E$ coincides with $T$
out of $B(0,1)$; the author does not know whether this is really possible.
If we can make $H^2(E\cap B(0,1))$ smaller than $H^2(T\cap B(0,1))$,
then we can probably minimize $H^2(E\cap B(0,1))$, but then $E$ is not 
globally minimal in $\R^3$, because a dilation would do better. If 
$H^2(E\cap B(0,1)) > H^2(T\cap B(0,1))$ for $E$ as above, the
only minimum is $T$ (obtained as limit of contractions of sets $E$).
But there is still the unlikely possibility that the minimal value
of $H^2(E\cap B(0,1))$ is precisely $H^2(T\cap B(0,1))$, and then we
could have an exotic minimal set $E$. The case when $E$ is obtained
with curves $\Gamma_i$ that are linked in a complicated way could
even be more fun.

\bigskip
\noindent {\bf 20. Simple computations of density}
\medskip
For the convenience of the reader, we do the simple computations of density
for the minimal cones of types $\Bbb Y$ and $\Bbb T$.

Let us first remind the reader  that $H^2(Y \cap B(0,r)) = 3\pi r^2/2$ when 
$Y$ is a set of type $\Bbb Y$ centered at the origin, simply because 
$Y \cap B(0,r)$ is composed of three half disks of radius $r$.

Next let $T$ be a set of type $\Bbb T$ centered at the origin. 
We want to check that
$$
H^2(T \cap B(0,r)) = d_{+} r^2,
\hbox{ with } d_+ = 3 {\rm Argcos}(-1/3).
\leqno (20.1)
$$
Now $T$ is composed of six equal angular sectors, so if
$\alpha \in (\pi/2,\pi)$ denotes the angle of these sectors, 
$H^2(T \cap B(0,r)) = 3 \alpha r^2$. We just need to show that
$\cos \alpha = -1/3$. 

Recall that $T$ is the cone over the edges of a regular tetrahedron 
centered at the origin; we may obtain such a tetrahedron by starting from 
the cube $[-1,1]^3$, and selecting four vertices (out of the eight) 
that are not mutually contiguous (see Figure 20.1). 
The distance between two vertices is $2\sqrt 2$, and the distance
form a vertex to the origin is $\sqrt 3$. Thus (see Figure 20.2)
$\sin(\alpha/2) = {\sqrt 2 \over \sqrt 3}$ and so
$\cos \alpha = 1 - 2 \sin^2(\alpha/2) = 1 - 4/3 = -1/3$, as needed.
This proves (20.1). 

\medskip 
\hskip 2.1cm
\epsffile{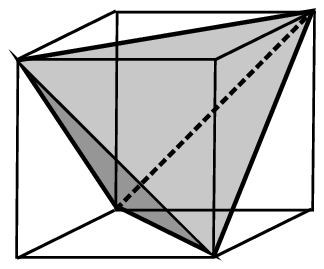}
\hskip 2.7cm
\epsffile{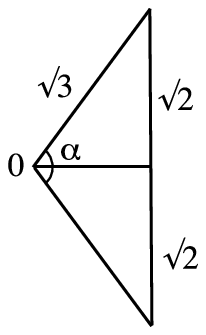}
\medskip
\hskip 1.2cm {\bf Figure 20.1} A tetrahedron in its box
\hskip 1.1cm {\bf Figure 20.2.} A face
\medskip

\bigskip
REFERENCES

\smallskip 

\item {[Al]} F. J. Almgren, Existence and regularity almost everywhere 
of solutions to elliptic variational problems with constraints, 
Memoirs of the Amer. Math. Soc. 165, volume 4 (1976), i-199.
\smallskip

\item {[DMS]} G. Dal Maso, J.-M. Morel, and S. Solimini, 
A variational method in image segmentation: Existence and
approximation results, Acta Math. 168 (1992), no. 1-2, 89--151.

\item {[D1]} G. David, Limits of Almgren-quasiminimal sets, 
Proceedings of the conference on Harmonic Analysis, 
Mount Holyoke, A.M.S. Contemporary Mathematics series, Vol. 320 
(2003), 119-145.
\smallskip
\item {[D2]} G. David, \underbar{Singular sets of minimizers for 
the Mumford-Shah functional},
Progress in Mathematics 233 (581p.), Birkh\"auser 2005.
\smallskip
\item {[D3]} G. David, Epiperimetry and $C^{1+\alpha}$-regularity
for almost-minimal sets in $\Bbb R^3$, 
preprint, Universit\'{e} de Paris-Sud 2006.
\smallskip 
\item {[DDT]} G. David, T. De Pauw, and T. Toro,
A generalization of Reifenberg's theorem in $\Bbb R^3$, preprint.
\smallskip
\item {[DS]} G. David and S. Semmes, Uniform rectifiability and 
quasiminimizing sets of arbitrary codimension, 
Memoirs of the A.M.S. Number 687, volume 144,  2000.
\smallskip
\item {[Do]} A. Dold, \underbar {Lectures on algebraic topology}, Second edition,
Grundlehren der Mathematishen Wissenschaften 200, Springer Verlag 
1980. 
\smallskip 
\item {[Du]} J. Dugundji, \underbar {topology}, Allyn and Bacon, 
Boston, 1966.
\smallskip
\item {[He]} A Heppes, Isogonal sph\"arischen netze, Ann. Univ. Sci. Budapest
E\"otv\"os Sect. Math. 7 (1964), 41-48.
\smallskip
\item {[Fe]} H. Federer, \underbar{Geometric measure theory}, 
Grundlehren der Mathematishen Wissenschaf-ten 
153, Springer Verlag 1969.
\smallskip
\item {[La]} E. Lamarle, Sur la stabilit\'{e} des syst\`{e}mes liquides en 
lames minces, M\'{e}m. Acad. R. Belg. 35 (1864), 3-104.
\smallskip
\item {[LM]}  Gary Lawlor and Frank Morgan,
Paired calibrations applied to soap films, immiscible fluids, and 
surfaces or networks minimizing other norms, 
Pacific J. Math.  166  (1994),  no. 1, 55--83.
\smallskip
\item {[Ma]}  P. Mattila, \underbar{Geometry of sets and 
measures in Euclidean space}, Cambridge Studies in
Advanced Mathematics 44, Cambridge University Press l995.
\smallskip
\item {[Mo1]} F. Morgan, Size-minimizing rectifiable currents, 
Invent. Math. 96 (1989), no. 2, 333-348.
\smallskip
\item {[R]} E. R. Reifenberg, Solution of the Plateau Problem for 
$m$-dimensional surfaces of varying topological type,
Acta Math. 104, 1960, 1--92.
\smallskip
\item {[St]}	E. M. Stein, \underbar{Singular integrals and 
differentiability properties of functions},
Princeton university press 1970 .
\smallskip
\item {[Ta]} J. Taylor, The structure of singularities in 
soap-bubble-like 
and soap-film-like minimal surfaces, 
Ann. of Math. (2) 103 (1976), no. 3, 489--539.
\smallskip

 \bigskip
\vfill \vfill \vfill\vfill
\noindent Guy David,  
\smallskip\noindent 
Math\'{e}matiques, B\^atiment 425,
\smallskip\noindent 
Universit\'{e} de Paris-Sud, 
\smallskip\noindent 
91405 Orsay Cedex, France
\smallskip\noindent 
guy.david@math.u-psud.fr

\bye